\documentclass{article}
\usepackage{amsmath,amssymb,bm,yhmath,mathtools,extarrows}
\usepackage[all]{xy}
\usepackage{caption,graphicx}
\font\ibf=cmbxti10
\font\smallibf=cmbxti8

\usepackage{multicol}

\title{All smooth four-dimensional Schoenflies balls are geometrically simply connected}
\author{by Valentin {\sc Po\'enaru}\footnote{Professor Emeritus at the Universit\'e Paris Sud-Orsay, Math\'ematiques 425, 91405 Orsay Cedex, France. e-mail: valpoe@hotmail.com}}
\date{\sl (June 2016)}

\begin{document}

\maketitle

\vglue 1cm

\section{Introduction}\label{sec1}

This paper will provide a complete proof for the result stated in the title above, namely the following theorem.

\bigskip

\noindent {\bf Theorem 1.} {\it Any smooth $4$-dimensional Schoenflies ball, which we will denote by $\Delta_{\rm Schoenflies}^4$ or simply $\Delta^4$, is geometrically simply-connected (G.S.C.).}

\bigskip

Our proof relies heavily on the celebrated work by Barry Mazur [Ma] from (1958), of which I remind here what is necessary for us right now.

\smallskip
What Barry showed, in the DIFF context is that, for any $n$, if $\Delta^n \equiv \Delta_{\rm Schoenflies}^n$, then
$$
\Delta^n - \{\mbox{a boundary point}\} \underset{\rm DIFF}{=} B^n - \{\mbox{a boundary point}\}.
$$
Soon after 1958, S. Smale, M. Kervaire and J. Milnor [Ke-M], [S] improved this, for all $n \ne 4$ showing that in that case $\Delta^n \underset{\rm DIFF}{=} B^n$.

\smallskip

Of course, for $n \leq 3$ this last result was already known, and less high-power technology was necessary. For $n=3$, this was known through the work of J.W. Alexander in the nineteen twenties, while for $n=2$ any good, strong enough proof of the fundamental conformal representation theorem is complex analysis (sometimes referred to as the ``Riemann mapping theorem''), tells us what we need. But then, in four dimensions, the only known thing, as of today, is just Barry's old result i.e. that
$$
\Delta_{\rm Schoenflies}^4 - \{\mbox{a boundary point}\} \underset{\rm DIFF}{=} B^4 - \{\mbox{a boundary point}\}, \leqno (0)
$$
where of course $B^4$ means the standard 4-ball. One should notice that the (0) is actually equivalent to
$$
{\rm int} \, \Delta^4 \underset{\rm DIFF}{=} R_{\rm standard}^4 
\leqno (1) 
$$
and this is what we will actually use in this paper. [The implication (0) $\Rightarrow$ (1) is obvious, and for the implication (1) $\Rightarrow$ (0) it is sufficient to show the following FACT. Up to diffeomorphism, there is a unique way to add a copy of $R^3$ as boundary of $R^4$. Here is, in a nutshell, the proof of the FACT. Let $Y^4 = \{R^4$ with a copy of $R^3$ glued at infinity, as a boundary$\}$. If we take a tubular neighbourhood of $R^3$ in $Y^4$, its generic fiber defines a PROPER embedding, $[0,\infty) \xrightarrow{ \ \varphi \ } R^4$. It is not hard to show that the diffeomorphism type of the pair $(R^4 , \varphi [0,\infty))$ determines the diffeomorphism of the non compact manifold with non-empty boundary $Y^4$.

\smallskip

At this point, we can invoke the well-known fact that in dimensions $n \geq 4$ there are no wild Artin-Fox arcs, and this clinches the proof that  (1) $\Rightarrow$ (0).]

\smallskip

We shall denote now by $\Delta^2$ the 2-skeleton of $\Delta^4$ and here is now a statement which is equivalent to our Theorem~1.

\bigskip

\noindent {\bf Theorem 2.} {\it The regular neighbourhood of $\Delta^2 \subset \Delta^4$, call it $N^4 (\Delta^2)$, is G.S.C.}

\bigskip

Our arguments, in this paper, will actually provide a proof for Theorem 2.

\smallskip

The rest of this introductory section deals with a bit of history. In July 2006 I have put on line the paper 

\smallskip

\centerline{``On the 3-dimensional Poincar\'e Conjecture and the 4-dimensional Schoenflies Problem''.} 

\centerline{[ArXiv.org/abs/math.GT/0612554]}

\smallskip

There is a big organic connection between that ArXiv paper, which was itself a fast survey of a never published, actually never typed, much longer manuscript, called ``Po V-B'', and our present paper. The Po~V-B and its ArXiv summary, dealt with a simultaneous proof of the following two items, namely
\begin{enumerate}
\item[i)] The proof that, for any homotopy 3-ball $\Delta^3$, we have $\Delta^3 \times I \in {\rm GSC}$. The proof was starting from a previous result of mine [Po1], [Po2], which is the following: ``The open DIFF 4-manifold $\{ \Delta^3 \times I \, \# \, \infty \, \# \, (S^2 \times D^2)$, with all the boundary removed$\}$, is GSC.'' Actually, the proof required not only the [Po1], [Po2] mentioned above, but also some other previous papers by the same author, which I will no longer mention here. But the interested reader may find a synthetical account of these things in [O-Po-Ta]. Then, with very much the same techniques, the Po V-B did contain the following item.
\item[ii)] The proof that $\Delta_{\rm Schoenflies}^4$ is GSC. This time, exactly as the [Po1], [Po2] above was used in the proof of i), one uses Barry's result (0) (or (1)).
\end{enumerate}

If from the ArXiv 2006 paper, as it stands, one carefully removes all the references to $3^{\rm d}$ Poincar\'e Conjecture, then what one gets is a pretty accurate description of the {\ibf plan} of this present work.

\smallskip

The point is that, although the starting points of i) and ii) above look quite different, the arguments for proving them are, in more than one way, very similar. When one looks into the seams, then some points are easier for i) and harder for ii), while for others quite the converse is true. For instance, while $\Delta^3$ is of codimension one inside ${\rm int} \, (\Delta^4 \times I \, \# \, \infty \, \# \, (S^2 \times D^2)$), the $\Delta^4_{\rm Schoenflies}$ is codimension zero in $\Delta^4 \cup \partial \Delta^4 \times [0,1) \underset{\rm DIFF}{=} R^4$, which makes things a bit harder for $\Delta^4$ than for $\Delta^3 \times I$.

\smallskip

But then, while $R^4 - \Delta^4$ is essentially product, the
$$
{\rm int} \, (\Delta^3 \times I \, \# \, \infty \, \# \, (S^2 \times S^2))
$$
contains those infinitely many $S^2 \times D^2$'s, source of some technical complications. But then, also, when viewed from a higher vantage point these are rather minor issues and, more seriously, a certain specific argument for ii) was painfully lacking for a long time. Then, in the Spring 2003, while I was visiting Princeton University, during some homeric working sessions with Dave Gabai and Frank Quinn, I had sort of an illumination and I saw the right ingredient which I was missing before. This ingredient is described in the 2006 ArXiv paper and it occurs here as formula (21.A) in the text which follows.

\smallskip

Anyway, nine years later, when I looked again into these things, in the end of 2014 and the beginning of 2015, then I realized that, for some reasons which may be were understandable in those early years, when I wrote the PO V-B, I had sort of blinders, at that time, failing to see the big forest which was looming behind those nearest trees. That made that my own re-reading, of that old paper of mine, the old Po V-B, was a very painful and rather heavy affair.

\smallskip

The fact is that, years ago, while I was writing that Po V-B I did not know where I would get, until towards the end, and the written paper keeps the traces of the various unnecessary detours and blind alleys. Now I knew exactly what the correct scheme was, and so life was much easier. But once I understood well the whole idea, I thought things over again, and rather than following blindly the old Po V-B, I started by disentangled the Schoenflies part from my Poincar\'e papers, which will no longer be mentioned here from now on, and thereby I managed to vastly simplify and also shorten the arguments for Schoenflies, which I extracted from the Po V-B. I found a lot of simplifying tricks too.

\smallskip

Of course, there is now also the obvious issue of getting from the Theorem 1 stated above, to the full DIFF $4^{\rm d}$ Schoenflies, and the additional step which is still necessary, in order to achieve that, is to show that any smooth $4^{\rm d}$ Schoenflies ball $\Delta^4_{\rm Schoenflies}$, which is in GSC is also standard, i.e. diffeomorphic to $B^4$. That is the object of current joint work with Dave Gabai and I will not discuss it here, except for the following little comment. The fact that $\partial \Delta^4_{\rm Schoenflies} = S^3$ is NOT used in the present paper. This is quite natural, once it is understood that the same technology serves for proving both that $\Delta^4_{\rm Schoenflies} \in {\rm GSC}$ AND that $\Delta^3 \times I \in {\rm GSC}$. But then, the fact that $\partial \Delta^4 = S^3$ should be used big in the joint paper with Dave.

\smallskip

I want to end this introduction with some philosophical-speculative thoughts. The Schoenflies $4^{\rm d}$ DIFF problem is strongly connected to the following two outstanding questions in $4^{\rm d}$ topology: the DIFF $4^{\rm d}$ Poincar\'e Conjecture and the big open gap or chasm, which exists in $4^{\rm d}$, and only in $4^{\rm d}$, between the categories DIFF and TOP.

\smallskip

There are many ways in which our Schoenflies problem is connected with the $4^{\rm d}$ Poincar\'e Conjecture, and I will only mention them here.

\smallskip

If one assumes the full $4^{\rm d}$ DIFF Schoenflies, then it easily follows that every smooth DIFF homotopy sphere $\sum^4$ which is such that $\sum^4 - \{{\rm pt}\} \underset{\rm DIFF}{=} R_{\rm standard}^4$, also comes with $\sum^4 \underset{\rm DIFF}{=} S^4$. For more on this kind of topics, see also the short paper [Po3].

\smallskip

I will not say more, here and now, concerning that big gap, unbridgeable by algebraic topology, between DIFF and TOP, in dimension four, and only in dimension four.

\bigskip

By now years ago, already, I had numerous discussions with Dave Gabai concerning the issues treated in the present paper, which would have probably never existed without his friendly help. I wish to thank him warmly here. Many thanks are also due to Frank Quinn for the many very helpful discussions we had, and to Louis Funar to whom I have lectured about the matters presented here and who came with very useful questions and comments.

\smallskip

The present paper owes a lot to the long continuous connection of the author with the IHES. And then, last but not least, I wish to thank C\'ecile Gourgues for the typing and Marie-Claude Vergne for the drawings.

\section{Geometric preliminaries and the doubling process}\label{sec2}

From now on, $\Delta^4$ will mean the Schoenflies smooth 4-ball $\Delta^4_{\rm Schoenflies}$ and we will denote by $X^4 \underset{\rm DIFF}{=} R^4$ the interior of $\Delta^4$, or some chosen cell-division of it. We will work with the following collection of telescopically nested spaces
$$
\Delta^4 = \Delta^4_{\rm small} \subset \Delta^4 \cup \partial \Delta^4 \times [0,1) \subset \Delta_1^4 \equiv \Delta^4 \cup \partial \Delta^4 \times [0,1],
\leqno (2)
$$
where the middle $\Delta^4 \cup \partial \Delta^4 \times [0,1)$ is, up to diffeomorphism, our $X^4$ while $\Delta_1^4$ is a larger copy of $\Delta^4$.

\smallskip

Here, of course, $\partial \Delta^4  = S^3$ and from now on, $X^4$ will mean the $\Delta^4 \cup \partial \Delta^4 \times [0,1) = R^4$ in (2), or some specific cell-decomposition of it.

\bigskip

\noindent {\bf Lemma 3.} {\it Let $Y^4$ be a smooth compact bounded $4$-manifold, for which by analogy to $(2)$, we consider the telescopic system
$$
Y^4 = Y_{\rm small}^4 \subset Y_{\rm small}^4 \cup \partial Y^4 \times [0,1] \equiv Y_1^4 .
\leqno (2.1)
$$

We also assume that there exists a system of smoothly embedded discs 
$$
\left( \sum_1^P D_i^2 , \sum_1^P \partial D_i^2 \right) \subset (\partial Y^4 \times [0,1] , \partial Y^4 \times \{0\} = \partial Y_{\rm small}^4),
\leqno (3)
$$
which are in {\ibf cancelling position} with the $1$-handles of $Y^4_{\rm small}$. Then $Y^4$ is geometrically simply-connected, we will write GSC.}

\bigskip

By the ``cancelling position'' above, we mean that we can index the 1-handles of $Y^4$ as $H_1^1 , H_2^1 , \ldots , H_p^1$, so that when we consider the geometric intersection matrix $\partial D_j^2$ (exterior 2-handle) $\cdot$ $H_i^1$ (interior 1-handle), we have
$$
\partial D_j^2 \cdot H_i^1 = \delta_{ji}.
$$

The GSC concept makes sense both for cell-complexes and for smooth manifolds, and our hypothesis in this lemma can be also stated as $Y_{\rm small}^4 + \underset{1}{\overset{P}{\sum}} \, D_i^2$ is GSC OR
$$
N^4 \left( Y_{\rm small}^4 + \sum_1^P D_i^2 \right) \ \mbox{is GSC.}
$$
The argument which follows is elementary differential topology.

\bigskip

\noindent {\bf Proof of Lemma 3.} For the collar $\partial Y^4 \times [0,1]$ we have two basic projections
$$
\xymatrix{
\partial Y^4 \times [0,1] \ar[d]^{\pi} \ar[r]_{ \ \ \ \pi_0} &[0,1]. \\
\partial Y^4 = \partial Y_{\rm small}^4
}
$$
We may assume that the smooth map from $2^{\rm d}$ into $3^{\rm d}$ $\pi \ \Bigl\vert \ \underset{1}{\overset{P}{\sum}} \, D_i^2$ is generic, i.e. that this map is a generic immersion outside of a finite, set of isolated singularities, generically denoted by $s \in {\rm int} \, D^2$, which are {\ibf Whitney's umbrellas}. Notice that, because of the embedding condition in (3), the occurrence of clasps as double lines is excluded. For such a clasp there would be a contradiction between the $\pi_0$-values at the two boundary points. We will perform now the following steps which, provisionally, will add 1-handles to $Y_{\rm small}^4 \cup \underset{1}{\overset{P}{\sum}} \, D_i^2$, to be killed later.

\medskip

1) We eliminate the triple points, as follows. At each triple point (of $\pi \ \Bigl \vert \ \underset{1}{\overset{P}{\sum}} \, D_i^2$) we have three branches which, from the viewpoint of the $\pi_0$-values at the triple point, are $B^2 ({\rm lowest})$, $B^2 ({\rm median})$, $B^2 ({\rm highest})$. We create a hole $t$ inside $B^2 ({\rm lowest}) \subset \underset{1}{\overset{P}{\sum}} \, D_i^2$ and, at the same time we add a 1-handle to $N^4 \left( Y_{\rm small}^4 \cup \underset{1}{\overset{P}{\sum}} \, D_i^2 \right)$, by pushing the center of the hole $t$ down to $\partial Y_{\rm small}^4$, so that $\partial ({\rm hole} \, t) \subset  \partial Y_{\rm small}^4$.

\medskip

2) We can also eliminate the closed double curve, proceeding as follows: Start by noticing that each such curve $C$ considered as a subset of $M^2(\pi)$ consists either of two components $C({\rm up})$, $C({\rm down})$, coming with $s(C({\rm up})) = C({\rm down})$, where $s$ is the involution $\underset{1}{\overset{P}{\sum}} \, D_i^2 \times \underset{1}{\overset{P}{\sum}} \, D_i^2 \supset M^2 (\pi) \xrightarrow{ \ s \ } M^2(\pi)$ gotten by exchanging the two factors, or of a single component endowed with the fixed point free involution $C \xrightarrow{ \ s \ } C$, such that $\pi (x) = \pi s(x)$. In this last case we can find two disjoined arcs $A({\rm up})$, $A({\rm down})$, exhausting $C$, with $s \, A({\rm up}) = A({\rm down})$. We can then push a $p \in C({\rm down})$, respectively a $p \in A({\rm down})$, all the way down to $\partial Y_{\rm small}^4$, like we just did it for the $B^2 ({\rm lowest})$ above. With this step, the double curves of $M_2 \left( \pi \ \Bigl\vert \ \underset{1}{\overset{P}{\sum}} \, D_i^2 \right)$ are out of the picture, at the price of our $D_i^2$'s acquiring the holes $t$.

\medskip

So, at the price of changing our source into a bunch of discs with holes $t$, we have gotten rid both of triple points of $\pi$ and of the closed curves of double points, which leaves us, as far as the immersion double points are concerned, only with RIBBONS and with the contribution of the double lines shooting out of the Whitney umbrellas, the two disjoined from each other (since now $M_3 = \emptyset$).

\smallskip

The double lines connecting two Whitney umbrellas in head-on collision can be easily eliminated, without changing the topology of $N^4 \left(Y_{\rm small}^4 \cup \underset{1}{\overset{P}{\sum}} \, D_i^2 \right)$. So, we may assume from now on that from each Whitney umbrella point, one of the two outgoing double lines hits the boundary. From each Whitney singularity $s$ starts now a double line of $\pi \ \Bigl\vert \ \underset{1}{\overset{P}{\sum}} \, D_i^2$ hitting the boundary, like in the Figure~1.

\smallskip

This figure is at the target of the map $\pi$, in $\partial Y_{\rm small}^4$.

$$
\includegraphics[width=120mm]{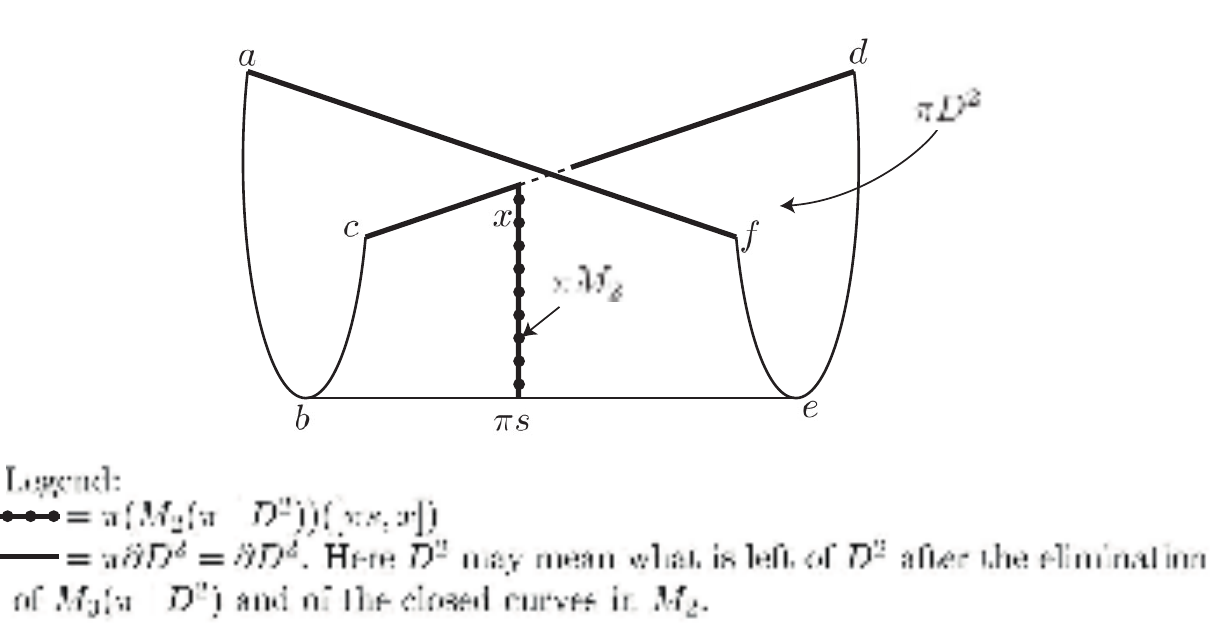}
$$
\label{fig1}
\centerline {\bf Figure 1.} 
\begin{quote}
This figure lives in $\partial Y_{\rm small}^4$. So the line $[\pi s , x]$ joins the Whitney umbrella with $\partial D^2 \subset \partial Y_{\rm small}^4$.
\end{quote}

\bigskip

We will get rid of these last surviving $s$'s by the following steps:
\begin{enumerate}
\item[i)] At the target, inside $\partial Y_{\rm small}^4$ we crush the arc $[\pi s ,x]$ (see Figure~1) to its endpoint $x \in \partial D^2$.
\item[ii)] At the same time as this, at the level of the source $D^2$, we crush to a point the inverse image of $[\pi s , x]$, meaning exactly the set
$$
\pi^{-1} ([\pi s , x]) = [x,s] \cup [s,x'] \, , \quad \mbox{with} \ x \in \partial D^2 \, , \ x' \in {\rm int} \, D^2 .
$$
\end{enumerate}

The topology of $N^4 \left(Y_{\rm small}^4 \cup \underset{1}{\overset{P}{\sum}} \, D_i^2 \right)$ stays unchanged, but obviously not the map $\pi$, which gets simplified.

\smallskip

Notice that the conjunction of these steps i) and ii) changes the topology of $\pi D^2$ but not the topology of $D^2$ and, moreover, they leave intact the closed loop $[a,b,c,d,e,f,a]$ in the Figure~1.

\smallskip

With all these operations, in which we have gotten rid of the triple points, of the closed curves of double points and of the Whitney umbrellas $s$, we have created a new situation, namely

\bigskip

\noindent (3.1) \quad The $\underset{1}{\overset{P}{\sum}} \, D_i^2$ occurring in (3) has been changed into an $\underset{1}{\overset{P}{\sum}} \, E_i^2$ where each $E_i^2$ is a disc with little holes $t$. Instead of (3), we have now 
$$
\left( \underset{1}{\overset{P}{\sum}} \, E_i^2 , \underset{1}{\overset{P}{\sum}} \, \partial E_i^2 \right) \subset (\partial Y_{\rm small}^4 \times [0,1] , \partial Y_{\rm small}^4 \times \{0\} ).
$$

The $\pi \ \Bigl\vert \ \underset{1}{\overset{P}{\sum}} \, E_i^2$ is a {\ibf generic immersion}, the only double points of which are now RIBBONS. Moreover $N^4 \left( Y^4 \cup \underset{1}{\overset{P}{\sum}} \, E_i^2 \right) = N^4 \left( Y_{\rm small}^4 \cup \underset{1}{\overset{P}{\sum}} \, D_i^2 \right) + \{$a collection of  vertical $4^{\rm d}$ 1-{\ibf handles}, which are in canonical bijection with the holes $t$, each 1-handle being a $N^4 (\mbox{vertical arc} \ [t,\pi t]) \subset \partial Y_{\rm small}^4 \times [0,1]\}$.

\smallskip

\noindent End of (3.1)
$$
\includegraphics[width=45mm]{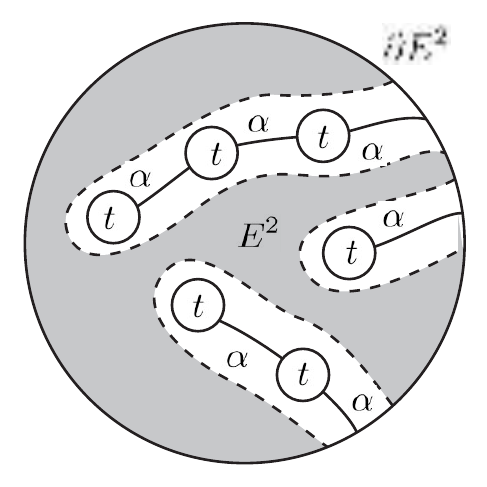}
$$
\label{fig2}
\centerline {\bf Figure 2.} 
\begin{quote}
The disc with holes $E^2$, at the source of $\pi$, with holes $t$ and with arcs $\alpha$ joining them acyclically to $\partial E^2$. This can be achieved so that
$$
\{{\rm arcs} \ \alpha\} \cap \{{\rm RIBBONS}\} = \emptyset .
$$
Moreover, the various $\pi \mid \alpha$ inject and are disjoined from each other. [In order to explain this, notice that any intersection $\pi \alpha_1 \cap \pi \alpha_2$ would be a double point of $\pi \mid \sum E_i^2$. But these are all concentrated in the RIBBONS from which our arcs $\alpha$ are disjoined, by construction.]

We have not tried to draw here the RIBBONS too, so as not to overcharge the figure. But the important point is that the RIBBONS in question do not separate the $t$'s from $\partial E^2$, allowing us to put in the $\alpha$'s.
\end{quote}

Now, we join the holes $t$ by arcs $\alpha$, to $\partial E^2$, {\ibf in an acyclic manner}, as it is suggested  in the Figure~2; this is done without touching the RIBBONS of the immersion $\underset{1}{\overset{P}{\sum}} \, E_i^2 \xrightarrow{ \ \pi \ } \partial Y_{\rm small}^4$.

\smallskip

Next, we enlarge our manifold $Y_{\rm small}^4$ with small dilatations $D^2 (\alpha \cup \pi (\alpha))$, each contained in a vertical plane, and which at the level of $N^4 \left( Y^4 \cup \underset{1}{\overset{P}{\sum}} \, E_i^2 \right)$ (see (3.1)), means addition of 2-handles. These 2-handles, going along the linear chain above, {\ibf exactly cancel} the 1-handles from (3.1). So
$$
N^4 \left(Y_{\rm small}^4 ({\rm enlarged}) \cup \underset{1}{\overset{P}{\sum}} \, E_i^2 \right) \underset{\rm DIFF}{=} N^4 \left( Y^4 \cup \underset{1}{\overset{P}{\sum}} \, D_i^2 \right),
$$
meaning that the result of this construction, with $Y^4({\rm enlarged}) = N^4 \left( Y^4 \cup \underset{\alpha}{\sum} \, D^2 (\alpha \cup \pi(\alpha)\right)$ is GSC.

Moreover, we have now a family of discs $D^2 = D^2 ({\rm small}) \subsetneqq \sum^2$ shaded in Figure~2, replacing the $\underset{1}{\overset{P}{\sum}} \, E_i^2$ in (3.1) and such that
$$
M_2 \left( \pi \ \biggl\vert \ \underset{1}{\overset{P}{\sum}} \, D_i^2 \, ({\rm small}) \right) = M_2 \left( \pi \ \biggl\vert \ \underset{1}{\overset{P}{\sum}} \, E_i^2 \right) .
\leqno (3.2)
$$
It is this family
$$
\left( \underset{1}{\overset{P}{\sum}} \, D_i^2 \, ({\rm small}) , \underset{1}{\overset{P}{\sum}} \, \partial D_i^2 \right) \subset \left(\partial Y_{\rm small}^4 ({\rm enlarged}) \times [0,1], \partial Y_{\rm small}^4\right)
\leqno (3.3)
$$
with which we will work from now on. That is what $Y^4 , D^2$, will mean now.

\smallskip

Now, imagine for a split-second, that we would also know, from the very beginning that actually
$$
M_2 \left( \pi \ \biggl\vert \ \underset{1}{\overset{P}{\sum}} \, D_i^2 \right) = \emptyset .
$$
(This is, of course, a totally irrealistic assumption, at the present stage in the game, which we make here just for the sake of the argument, but let us still imagine~$\ldots$)

\smallskip

Then, by a process of {\ibf raising the see level} inside the collar $\partial Y_{\rm small}^4 ({\rm enlarged}) \times [0,1]$, we could realize a diffeomorphism
$$
N^4 \left( Y_{\rm small}^4 \cup \underset{1}{\overset{P}{\sum}} \, D_i^2 \right) + \{\mbox{$3$-handles}\} \underset{\rm DIFF}{=} Y_{\rm small}^4,
$$
[each $3$-handle gotten by joining $D_i^2$ to $\pi D_i^2$ along the lines $[y,\pi y]$, $y \in D_i^2$].

\smallskip

This diffeomorphism would imply our desired result, but only under the outrageous assumption we just made. Going now back to real life and to (3.3), our RIBBONS are pairs of curves which inside $M_2 \left(\pi \ \Bigl\vert \ \underset{1}{\overset{P}{\sum}} \, D_i^2 \right) \subset \underset{1}{\overset{P}{\sum}} \, D_i^2$ (meaning from now the $D_i^2 ({\rm small})$) look like in Figure~3. What we see there~is~a
$$
\{\mbox{vertical plane}\} \subset \partial Y_{\rm small}^4 \times [0,1],
$$
such that the $\pi^{-1} \pi ({\rm RIBBON})$, lives (for every individual RIBBON) inside such a plane. Moreover, we have exactly
$$
\pi^{-1} \pi ({\rm RIBBON}) \subset \{{\rm RECTANGLE}\} \subset \{\mbox{our vertical plane}\}
$$
$$
\includegraphics[width=110mm]{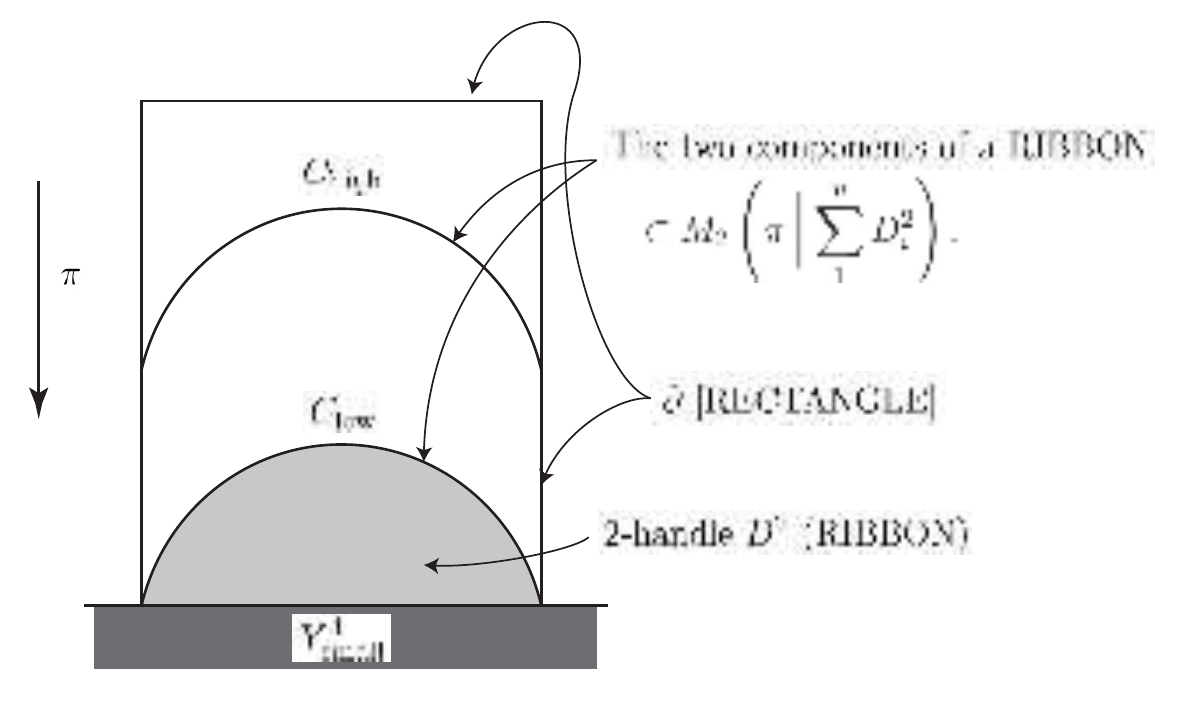}
$$
\label{fig3}
\centerline {\bf Figure 3.} 
\begin{quote}
The RIBBON, at the source, inside the collar $\partial Y_{\rm small}^4 \times [0,1]$, and its projection $\pi$. Here the ${\rm RIBBON} \equiv (C_{\rm high} , C_{\rm low})$, with $C_{\rm high} \subset {\rm int} \, D^2 ({\rm small})$ (3.3) and with a proper embedding $(C_{\rm low} , \partial \, C_{\rm low}) \subset (D^2 , \partial D^2)$.
\end{quote}

Notice now the following facts:
\begin{enumerate}
\item[a)] Since the 2-handles $D^2 (\alpha \cup \pi (\alpha))$ have exactly cancelled the 1-handles in (3.1), we find that

$Z^4 \equiv N^4 \left( Y_{\rm small}^4 \cup \underset{1}{\overset{P}{\sum}} \, E_i^2 (3.1) \right) \cup \underset{\alpha}{\sum} \, D^2 (\alpha \cup \pi (\alpha)) \underset{\rm DIFF}{=}$  (a diffeomorphism already mentioned) $\underset{\rm DIFF}{=} N^4 \left(Y_{\rm small}^4 ({\rm enlarged}) \cup \underset{1}{\overset{P}{\sum}} \, D_i^2 ({\rm small})(3.3) \right) \in {\rm GSC}$.

(This really is a reminder.)
\item[b)] To $Z^4$ we may add the vertical 2-handles $D^2 ({\rm RIBBON})$ which are shaded in the Figure~3 and then we get, as a consequence of a), $V^4 \equiv Z^4 + \underset{\rm RIBBONS}{\sum} \, D^2 ({\rm RIBBONS}) \in {\rm GSC}$.
\item[c)] Consider now
$$
W^4 \equiv N^4 (Y_{\rm small}^4 ({\rm original})) \, \cup \sum \{\mbox{all the vertical $2^{\rm d}$ smooth {\ibf dilatation}} \ D^2 (\alpha \, \cup \, \pi (\infty)) \ {\rm and} \ D^2({\rm RIBBON})\}.
$$
We clearly have $W^4  \underset{\rm DIFF}{=} Y_{\rm small}^4 ({\rm original})$.
\item[d)] The $C_{\rm low}$'s in Figure~3, breaks the system $\underset{1}{\overset{P}{\sum}} \, D_i^2$ (3.3) into a larger system of even smaller discs, call it $\underset{1}{\overset{P}{\sum}} \, D_i^2$ with $Q > P$ (but we will not change the notation for the $D_i^2$'s). This comes with an embedding
$$
\xymatrix{
\left(\underset{1}{\overset{Q}{\sum}} \, D_i^2 , \underset{1}{\overset{Q}{\sum}} \, \partial D_i^2 \right) \  \ar@{^{(}->}[rr] &&(\partial W^4 \times [0,1] , \partial W^4 \times \{0\}) \ar[d]^{\pi} \\
&&\partial W^4
}
$$
which is now such that $M_2 \left( \pi \ \biggl\vert \ \underset{1}{\overset{Q}{\sum}} D_i^2 \right) = \emptyset$. Moreover, clearly also
$$
W^4 \cup \sum_1^Q D_i^2 \underset{\rm DIFF}{=} Z^4 + \sum_{\rm RIBBON} \{\mbox{the 2-handles} \ D^2 ({\rm RIBBON})\} \in {\rm GSC}.
$$
\end{enumerate}

Applying once more the argument of the raising of sea-level, we find that $W^4$ itself is GSC and this proves our lemma. It should be noticed that it is essential for our whole argument that the projections $\partial Y_{\rm small}^4 \times [0,1] \xrightarrow{ \ \pi \ } \partial Y_{\rm small}^4$ and $\partial W^4 \times [0,1]  \xrightarrow{ \ \pi \ } \partial W^4$ are compatible; there is only one notion of VERTICALITY which is involved in both. $\Box$

\bigskip

Once Lemma 3 has been proved, we move back to our main theme. The important fact, at this point, is that in the context of (2) we are presented with two infinite collapses, or collapsing flows.

\bigskip

\noindent (4.A) \quad As a simple consequence of the way in which $\Delta_1^4$ is defined, we have a RED collapse ${\rm int} \, \Delta_1^4 \searrow \Delta_{\rm small}^4$ (actually a compact $\Delta_1^4 \searrow \Delta_{\rm small}^4$ too, but that one we will never use).

\bigskip

\noindent (4.B) \quad From Barry's result $X^4 = {\rm int} \, \Delta_1^4 \underset{\rm DIFF}{=} R^4$ (see (1)), follows the existence of a second, BLUE collapse $X^4 \searrow {\rm pt}$.

$$
\includegraphics[width=11cm]{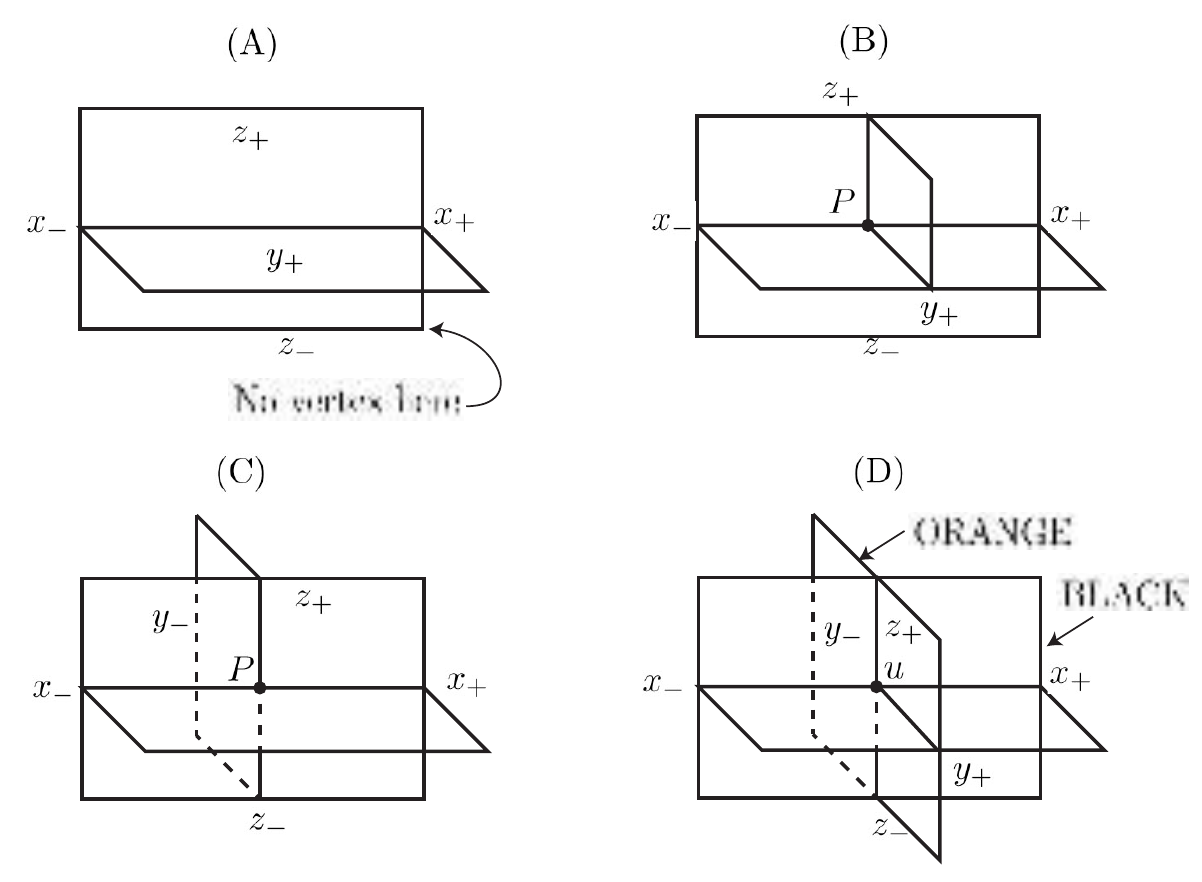}
$$
\label{fig4}
\centerline {\bf Figure 4.} 
\begin{quote}
The local models for $X^3 ({\rm GPS})$. Any permutation of $(x,y,z)$ or of signs, is of course OK too. The (B), (C) are, of course, homeomorphic, but NOT affinely so. So, for us they count as distinct. Of course also, the BLACK (A), (B), (C) may be orange, while in (D) we can also permute BLACK $\Leftrightarrow$ ORANGE. 

The only 1-skeleton present in (D) is $[x_- , x_+] \subset X^1 ({\rm BLACK})$, and this is like in the Figure 4-(A). I mean here the only 1-skeleton of one of the $X^3 ({\rm COLOURS})$'s. But then, when we move to $X^4$ then all the lines in our figure will be in the $X^1 \equiv \{\mbox{1-skeleton of $X^4$}\}$.

Remember that BLACK and ORANGE play a symmetrical role in all our story.
\end{quote}

\bigskip

Both of these collapsing flows will be used and one of the difficulties which we will have to overcome is that, a priori they cut transversally through each other in an uncontrolable fashion, and with a messy situation
$$
(\mbox{RED collapsing flow}) \pitchfork (\mbox{BLUE collapsing flow}),
$$
with which we could not live.

\smallskip 

Our $\Delta^4$ in (2) will be a subcomplex of a specially chosen cell-decomposition of $R^4 = X^4$, something which we will introduce very explicitly, next. We will call it a {\ibf ``general parallelipipedic structure''} of $X^4 = R^4 = R^3 \times R = X^3 \times R$, which I will currently call G.P.S. structure (sorry for the possible confusion!).

\smallskip

In what follows next, we will work with $R^4$ coming with its affine structure and with the coordinate system $(x,y,z,t)$; occasionally, the euclidean metric will also be invoked.

\smallskip

The time axis $R$ will be equipped with the lattice $Z \subset R$, the points of which will be labelled \break $\ldots < t_{-1} < t_0 < t_1 < t_2 < \ldots$ and we will consider the spatial slices $X^3 \times t_i = R^3$. Here comes now the

\bigskip

\noindent {\bf Definition 5.} We will define the {\ibf GPS} structures, first for $R^3 = X^2 \times t_i$ and then for $R^4 = R^3 \times R$. Here is the list of the GPS features.
\begin{enumerate}
\item[a)] Each $X^3 \times t_i$ comes with two cell-decompositions $X^3 ({\rm BLACK})$ and $X^3 ({\rm ORANGE})$, which are $t_i$-independent. The $X^{\varepsilon} ({\rm COLOUR})$ will be the $\varepsilon$-skeleton. We will ask that $X^2 ({\rm BLACK})$ and $X^2 ({\rm ORANGE})$ should cut transversally through each other, so that
$$
X^1 ({\rm BLACK}) \cap X^1 ({\rm ORANGE}) = \emptyset \subset X^3 \times t_i .
$$
\item[b)] The local models for the $X^2 ({\rm COLOUR})$ should be {\ibf generic}, like in Figure~4, where the $X^2 ({\rm COLOUR}) \subset X^3 \times t_i$ are, locally, displayed. What we see in the Figure~4-(D) is the local model for the intersection $X^2 ({\rm COLOUR}) \pitchfork X^2 (\mbox{DUAL COLOUR}) \subset X^3 \times t_i$. Then, the point $u$ is NOT a vertex of any of the two individual $X^2 ({\rm COLOUR})$, but certainly a vertex of the GPS structure of $X^3 \times t_i$, for which the cell-decomposition is canonically defined, starting from $X^2 ({\rm BLACK}) \cup X^2 ({\rm ORANGE})$.
\item[c)] We move now to $R^4$ which, when endowed with its GPS-structure, will be denoted $X^4$, of $2$-skeleton $X^2 ({\rm GPS})$. We want that, after the interiors of some 2-cells, generically called $D^2 (\gamma^1) \subset X^2 ({\rm GPS})$, are deleted, we should have a collapse
$$
X^2 ({\rm GPS}) - \bigcup {\rm int} \, D^2 (\gamma^1) \searrow {\rm pt}.
$$
[Morally speaking the $D^2 (\gamma^1)$'s are to be killed by a $3^{\rm d}$ collapse, but this, so called ``BLUE'' $3^{\rm d}$ collapse, will never be used. It ain't there.]
\item[d)] Inside each $X^3 \times t_i$ the $2^{\rm d}$ collapse from c) above should have {\ibf purely linear trajectories, without bifurcation}, like the blue arrow in Figure~5 may suggest.
\item[e)] When we move from $X^3 \times t_i$ to $X^4$, then each vertex $P \in X^3 \times t_i$ carries at most one of the edges $P \times [t_i , t_{i+1}]$ OR $P \times [t_i , t_{i-1}]$ and NEVER BOTH. But also, sometimes there are no temporal edges.
\item[f)] The $X^4$ has no vertices outside of the $X^3 \times t_i$'s.
\end{enumerate}
[Implicit in this definition, where $P$ is a vertex of $X^4$ the edges coming out of $P$ can only go in one of the directions $\pm \, x$, $\pm \, y$, $\pm \, z$ or $\pm \, t$ and not all of these choices are realized, for any given vertex $P$.] We will denote by $X^4$ the generic GPS subdivision of $R^4$ and by $X^{\varepsilon}$ its $\varepsilon$-skeleton.

\smallskip

We think from now on in terms of $X^4 = (R^4$, with a GPS-structure). Of course there is here a potential source of confusion since $X^3 ({\rm GPS})$ has already been used for $(R^3$, with a GPS-structure). When this confusion is not resolved by the context, we may write $X^{(3)}$ for the 3-skeleton of $X^4$. End of (5).

$$
\includegraphics[width=135mm]{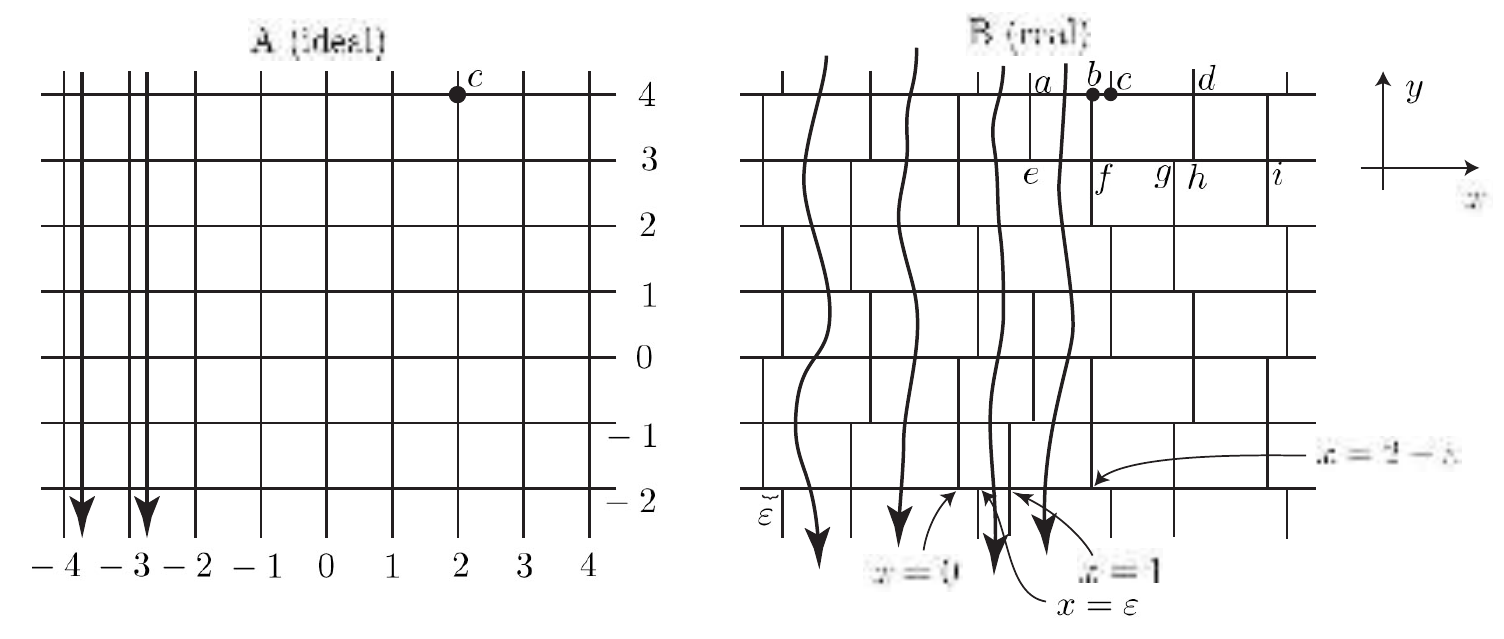}
$$
\label{fig5}
\centerline {\bf Figure 5.} 
\begin{quote}
The $X^1 ({\rm COLOUR})$ in $R^2 (x,y)$, in ideal and in real life. The BLUE arrows stand for the collapsing flow of $R^2$. In A (ideal) vertical lines are black and horizontal ones orange. This goes over to the deformed version B (real).
\end{quote}

\bigskip

Of course, the GPS structure is not uniquely determined by the definition above, and many choices are possible.

\smallskip

So I will describe now the STANDARD GPS STRUCTURE.

\smallskip

We will start with $R^3$ and, ideally, the $X^2 ({\rm COLOUR})$'s should be two perfectly cubical structures cutting transversally through each other, like Figure~5-(A) may suggest, with one dimension less. But then, this does not satisfy the genericity condition in (5)-b). Let us say, more precisely, that in the ideal case we chose, for $(m,n,p) \in Z$, $X^2 ({\rm BLACK}) ({\rm ideal}) = (x = 2m , y = 2n , z = 2p)$, $X^2 ({\rm ORANGE}) = (x = 2 m + 1 , y = 2 n + 1, z = 2p+1)$. Then we go to the real life situation, in two steps.

\medskip

I) We consider first, in $R^2 (x,y)$, the $X^1 ({\rm BLACK})({\rm ideal})$, $X^1 ({\rm ORANGE})({\rm ideal})$, like in the Figure~5-(A). Next, we perturb this, generically, as follows:

\smallskip

The lines $y = n \in Z$, independently of their colour, are left put and the black line $x=2m$ is broken into alternating successive pieces $x=2m$, and $x = 2m + \varepsilon$, with the pattern suggested by Figure~5-(B); something similar if then done for the orange lines $x = 2m + 1$, and this is also represented in Figure~5-(B). After this we have an inclusion

\bigskip

\noindent $(*)$ \quad $\{$the bicoloured structure in Figure~5-(B)$\} \times \{ -\infty < z < +\infty\} \subset \{$2-skeleton of the $X^3 \times t_i\}$. But the 2-skeleton in question contains some horizontal, $z = {\rm const}$, pieces too.

\bigskip

With the structure from Figure~5-(B), we consider all the little edges from the figure in question, like for instance $[a,b], [b,c], [c,d], [a,e], [b,f], [e,f], [d,h]$, and all the thin vertical bands $\{[\mbox{little edge}] \times (-\infty < z < +\infty)$, appropriately divided into little squares by the next step II$\}$ are the $2^{\rm d}$ BLUE collapsing flow lines. These occur in $(*)$, and see here also Figure~35.2. The BLUE $2^{\rm d}$ collapsing flow goes in the direction $\pm z$ (see Figure~35.2).

\bigskip

\noindent {\bf Remarks.} A) The BLUE flow lines in Figure~5, living horizontally in $(x,y)$ are {\ibf not} our real-life BLUE $2^{\rm d}$ collapsing flow. They just explain what we mean by linear trajectories in (5)-d. The real life $2^{\rm d}$ BLUE flow-lines are even straighter than in the Figure~5-(B).

\bigskip

II) If we take the whole $z$-flow of the ORANGE (respectively BLACK) structure in Figure~5-(B) and $(*)$, and then cut it with ORANGE (respectively BLACK) horizontal planes, namely ORANGE planes $z = 2p+1$ and BLACK planes $z = 2p$, the condition b) is still not verified at the intersections of these horizontal planes with the vertical lines going through the vertices from Figure~5-(B). We will name $X^1 ({\rm BLACK} / {\rm ORANGE})$ the $1^{\rm d}$ structure from Figure~5-(B). It has three kind of vertices: pure BLACK vertices, like $b,c$, pure ORANGE vertices, like $g,h$ and then mixed vertices like $a,f$ or $d$. Consider now the BLACK vertex $c$ from Figure~5-(B) and the four $1/4$ planes of BLACK colour, at $z = 2p$, neighbouring it in the ideal case. In the real life case, at $((b), c)$ we perturb the four $1/4$-planes, as follows: $(x > 0 , y > 0) \Rightarrow z=2p$, $(x > 0 , y < 0) \Rightarrow z=2p+\varepsilon$, $(x < 0 , y < 0) \Rightarrow z=2p+2\varepsilon$, $(x < 0 , y > 0) \Rightarrow z=2p+3\varepsilon$. We do a similar change at the ideal ORANGE vertices (with $z=2p$ replaced by $z=2p+1$), i.e. this time we handle $(g,h)$. The BLACK modification renders generic $a$ and also $d$, and the ORANGE modification takes care of $i$ and $f$. Here one uses the fact that $a,b$ are placed between two consecutive black vertices on the same horizontal and then, similarly $f,i$ between two consecutive orange vertices, again in the same horizontal.

\smallskip

So, with all these things we have by now, at level $X^3 \times t_j$ the correct structures $X^2 ({\rm BLACK})$ and $X^2 ({\rm ORANGE})$. Their union, i.e. the $(*)$ above, to which the horizontal BLACK and ORANGE little squares are added {\ibf is} the $X^2 ({\rm GPS}) \times t_j$. But before we can move to $4^{\rm d}$, we will need a refinement of the $X^3 ({\rm GPS})$ we have just introduced.

\smallskip

The $X^1 ({\rm BLACK}/{\rm ORANGE})$ in Figure~5-(B) is clearly a $(Z+Z)$-tesselation of $R^2$ and so are the  $X^1 ({\rm BLACK})$ , $X^1 ({\rm ORANGE})$, independently. Similarly, $X^3 ({\rm GPS}) \approx X^2 ({\rm BLACK}/{\rm ORANGE})$, $X^2 ({\rm BLACK})$, $X^2 ({\rm ORANGE})$ are each of them a $(Z+Z+Z)$-tesselation of $R^3$. The $X^2 ({\rm BLACK} / {\rm ORANGE})$ is a symmetric $\varepsilon$-generic deformation of the standard integral cubical subdivision of $R^3$. We need to identify now, inside $X^2 ({\rm BLACK} / {\rm ORANGE})$ two coarser tesselation (coarser meaning with larger fundamental domain) call them $2X^2 ({\rm COLOUR})$, with the feature that, $2X^2 ({\rm ORANGE})$ and $2X^2 ({\rm BLACK})$, i.e. $2X^2 ({\rm ORANGE}) \subset X^2 ({\rm ORANGE}) \subset X^2 ({\rm BLACK} / {\rm ORANGE}) \supset X^2 ({\rm BLACK}) \supset 2X^2 ({\rm BLACK})$, are transversal to each other. At the $R^2$-level of the ideal Figure~5-(A) we may take $2X^1 ({\rm BLACK})({\rm ideal}) = \{ x = 6m , y = 6n \}$, \break $2X^1 ({\rm ORANGE})({\rm ideal}) = \{ x = 6m + 3 , y = 6n + 3\}$ and from there one continue as before. Importantly, we have $2X^2 ({\rm BLACK}) \pitchfork 2X^2 ({\rm ORANGE})$ (i.e. transversal contact, with $2X^1 ({\rm BLACK}) \cap 2X^1 ({\rm ORANGE}) = \emptyset$, in $X^3 \times t_i$.)

\smallskip

We move now to the $X^4$ and in order to define its standard $4^{\rm d}$ GPS structure, it is sufficient to specify the $X^2$, which afterwards can be filled in canonically (i.e. linearly) with $3^{\rm d}$ and $4^{\rm d}$ cells.

\smallskip

All the features in definition (5) will be satisfied, if we take
$$
X^2 \equiv \sum_i \underbrace{X^2 ({\rm GPS}) \times t_i}_{\rm bicollared} + \sum_i 2X^1 ({\rm BLACK}) \times [t_{2i-1} , t_{2i}] + \sum 2X^1 ({\rm ORANGE}) \times [t_{2i} , t_{2i+1}]. \leqno (6)
$$

This ends the definition of our ``standard'' GPS structure which will be our normal way to put flesh and bones on the DEFINITION (5).

\bigskip

In our context, $\Delta^4 = \Delta^4_{\rm small} \subset X^4$ is a subcomplex and so its $2$-skeleton $\Delta^2$ is a subcomplex of $X^2$. We will use the notation
$$
\Gamma (1) \equiv \{\mbox{1-skeleton of $\Delta^4$}\} \subset \Gamma(\infty) = \{\mbox{the $X^1$, 1-skeleton of $X^2$}\}.
$$

We will express now the basic RED and BLUE collapses from (4-A and B) by the following lemma, which is a detailed version of the (4).

\bigskip

\noindent {\bf The P.L. Lemma 3.1.} {\it For $R^4$ (with its DIFF standard structure), there is a GPS structure $X^4$ such that}
\begin{enumerate}
\item[1)] {\it $\Delta^4_{\rm Schoenflies} \subset \Delta^4 \cup \partial \Delta^4 \times [0,1) = R^4$ is a subcomplex of it. And $X^4$ will have all our desirable BLUE collapsibility properties. Here $\Delta^4_{\rm Schoenflies} = \Delta^4_{\rm small}$, of course.}
\item[2)] {\it The infinite complex $X^4 - {\rm int} \, \Delta^4$ has a RED collapse}
$$
X^4 - {\rm int} \, \Delta^4 \searrow S^3 = \partial \Delta^4 .
$$
\end{enumerate}

\bigskip

This last feature 2) will be called the PROPERTY (P). So the PL-lemma says that there is a GPS structure with property (P) or, said differently, we can realize both features (4-A), (4-B) with the same GPS structure. [Only the $2^{\rm d}$ part of the BLUE collapse will ever be needed.]

\smallskip

Now, what Barry's classical result (0) tells us is that, with the standard affine (and hence DIFF too) structure of $R^4 = \Delta^4 \cup \partial \Delta^4 \times [0,1)$ there is a strictly cubical subdivision, s.t. $\Delta^4_{\rm Schoenflies}$ is a subcomplex. I am stressing here this `` affine'', in view of what is coming next, but the euclidean metric  of $R^4$ will have to be, occasionally, invoked too.

\bigskip

\noindent {\bf Sublemma 3.1-A.} {\it There is a cubical subdivision of $R^4$, like above, s.t. the $R^4 - {\rm int} \, \Delta^4_{\rm Schoenflies}$ has the property} (P).

\bigskip

\noindent {\bf Proof.} Here are the successive steps of the argument.

\medskip

\noindent I) We reformulate explicitly the fact that $R^4 - {\rm int} \, \Delta^4$ has property (P). For the locally affine manifold $R^4 - {\rm int} \, \Delta^4$ there exists a triangulation $T$ by convex affine symplexes, s.t.

\bigskip

\noindent (6.1) $T$ has the property (P).

\bigskip

\noindent (6.2) There is a positive lower bound $\varepsilon > 0$ s.t.
\begin{enumerate}
\item[i)] For each simplex $\sigma$ of $T$, $\Vert \sigma \Vert > \varepsilon$,
\item[ii)] $\{$the angles of the faces of $\sigma$ with the coordinate $1^{\rm d}$, $2^{\rm d}$ or $3^{\rm d}$ planes defined by $(x,y,z,t)\} > \varepsilon$.
\end{enumerate}

\bigskip

Figure 5.1 should suggest the proof of (6.2). Of course, the full detailed argument here has to invoke the DIFF HAUPTVERMUTUNG of J.H.C. Whitehead, but this should be quite standard. There is no cubical structure at all involved in our step II, that will be coming next. Also what we mean by ``cubical'' is not just that the 3-cells are cubes, but that these cubes are nicely aligned with $(x,y,z)$.
$$
\includegraphics[width=80mm]{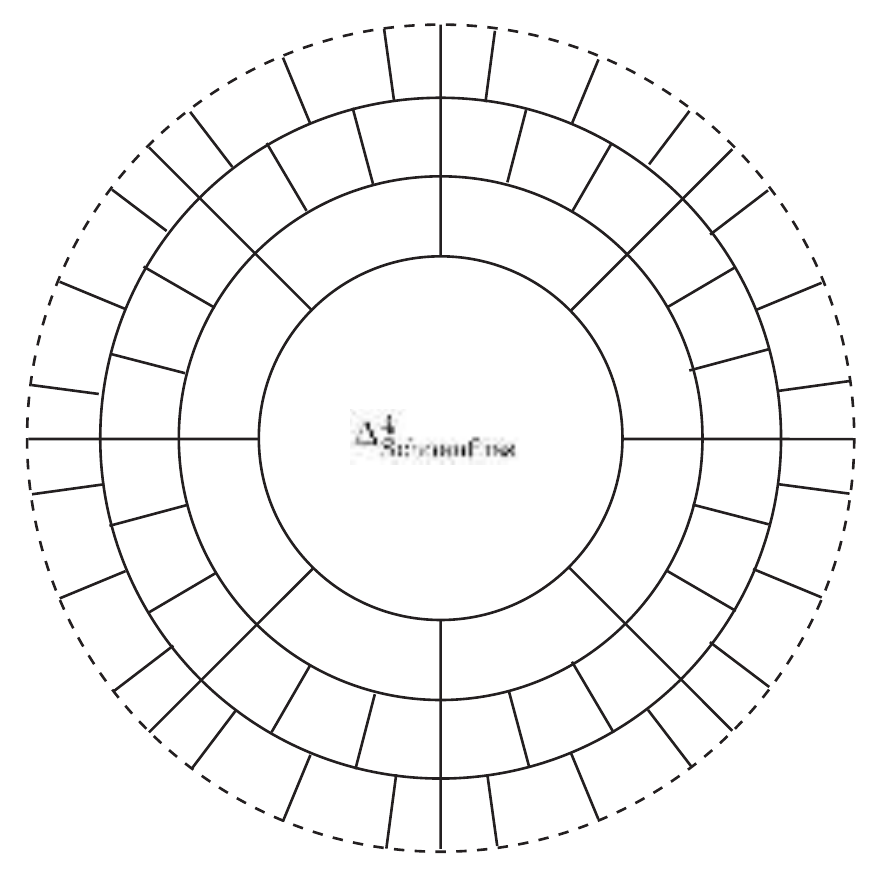}
$$
\label{fig5.1}
\centerline {\bf Figure 5.1.}
\begin{quote}
Idea for the proof of (6.2). The dotted line is the infinity of $R^4$ and the plain circles are equidistant.
\end{quote}

\medskip

\noindent II) Let us denote by $X$ the strictly cubical structure of $R^4$ and by $X'$ the cell-decomposition $X \, \cap \, T$, subdivision of $T$.

\medskip

\noindent {\bf Sub-Sublemma 3.1-B.} {\it When restricted to $R^4 - {\rm int} \, \Delta^4$, the $X \cap T$ continues to have Property {\rm (P)}.}

\medskip

\noindent {\bf Proof.} One gets $X \cap T$ by {\ibf bissecting} the affine-convex symplexes of $T$ by affine planes of dimension $1,2,3$. There are the planes defined by $X$. Our desired conclusion follows easily from this fact. [Property (P) is invariant under barycentric subdivisions, stellar subdivisions and, more up to the point, under {\ibf Siebenmann's bissections}.]

\smallskip

The `` bissections'' mentioned above are to be introduced formally in the next point III).

\medskip

\noindent III) I remind the reader the basic facts of Larry Siebenmann's theory of bissections.

\smallskip

Siebenmann starts by introducing {\it cellulations}, which are an extension of simplicial complexes: instead of using simplexes we use now compact cells $D$ with a linear-convex structure. The notion of (linear) subdivision extends in an obvious way to cellulations and, also, instead of subcomplexes we can introduce now sub-cellulations. What we have gained with this approach is, among other things, the following useful fact: if $Y \subset Z$ is a sub-cellulation, then any subdivision $Y'$ extends canonically to a subdivision of $Z$, not affecting the open cells in $Z-Y$. An important class of subdivisions are the BISSECTIONS. These are localized at the level of an $i$-cell $D^i$ and are obtained by cutting $D^i$ with a hyperplane $H^{i-1} \subset D^i$ and splitting $D^i$ itself and any sub-cell of $D^i$ met by $H^{i-1}$, in the obvious way. Our ``useful fact'' above extends to bissections.

\smallskip

One should notice that no genericity conditions are required here for the position of the hyperplane $H^{i-1} \subset D^i$. As mentioned above, bissections, barycentric subdivisions or Alexander's stellar subdivisions clearly preserve Property (P), but this kind of thing is a priori nor clear for general linear subdivision. But then here comes Siebenmann's very useful version of Alexander's old classical lemma.

\medskip

\noindent {\bf The Siebenmann Lemma.} {\it Let $X$ be a cellulation and $X \to X'$ a linear subdivision of $X$. There exists then a cellulation $X_1$ such that one can go both from $X$ and from $X'$ to $X_1$, via bissections}
$$
\xymatrix{
X \ar[rr]^{\rm linear} \ar[dr]_{\rm bissection} &&X' \ar[dl]^{\rm bissection.} \\
&X_1
}
$$

Now, $X \to X \cap T$ {\ibf is} a linear subdivision of $X$ too, and so one can apply Siebenmann's lemma and get a common bissection $X_1$. Since $X \cap T \in (P)$ and $X \cap T \to X_1$ is a bissection, we also have $X_1 \in (P)$. Diagrammatically we have here
$$
\xymatrix{
X \ar[rr]^{\rm linear} \ar[dr]^{\rm \qquad bissections} &&X \cap T \in (P) \ar[dl] \\
&X_1 \in (P)
}
$$

\medskip

\noindent IV) We consider the full $X$ for $R^4$ with its $X \cap T \mid (X^4 - {\rm int} \,  \Delta^4) \in (P)$ and the bissection $X \to X_1$, considered for the full $X$. (Unlike what happens for, let us say, barycentric subdivisions, a bissection of a subcomplex is, automatically, a bissection of the whole complex.)

\smallskip

From (6.2) it follows that there is a uniform, very fine cubical subdivision of our cubical $X$, call it
$$
X \underset{\mbox{\footnotesize very fine cubical subdivision}}{\xrightarrow{ \qquad \qquad \qquad \qquad \qquad \qquad}} Y,
$$
s.t. inside $Y$ we can MIMICK the bissection $X \to X_1$ by subcomplexes of $Y$ which are very close approximations of convex hyperplanes. The idea of this, is suggested in the Figure~5.2 below.

\smallskip

When each cell $\sigma$ of $X_1$ is replaced by $\sigma' \subset Y$ we get a cell-decomposition $X''$ of $R^4$, which is very close to $\sigma$ and such that:
\begin{enumerate}
\item[a)] The cells $\sigma'$ of $X''$ are nearly-convex.
\item[b)] $X'' \in (P)$.
\end{enumerate}
$$
\includegraphics[width=11cm]{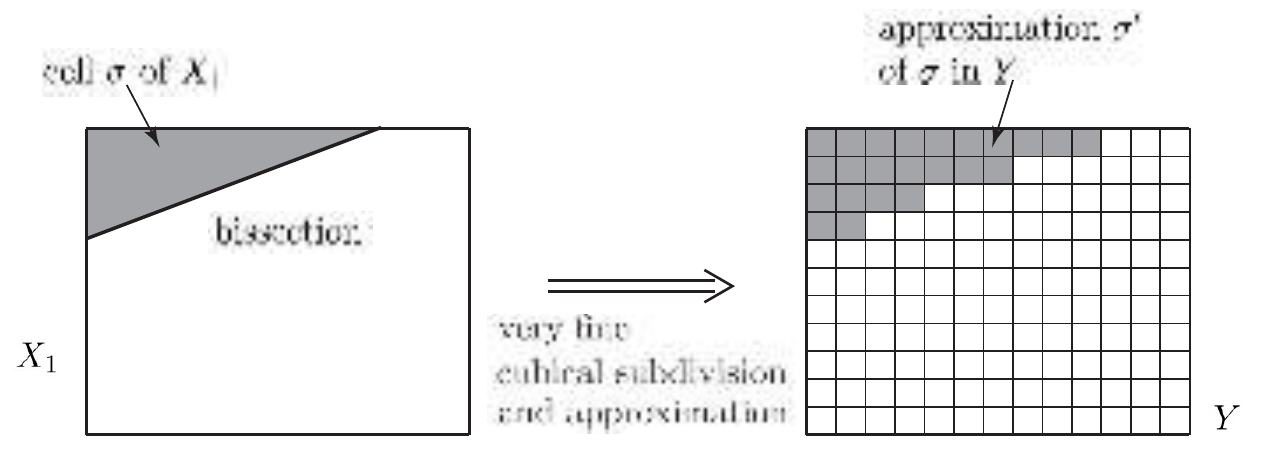}
$$
\label{fig5.2}
\centerline {\bf Figure 5.2.}

\medskip

\centerline{ \ MIMICKING}

\bigskip

Consider now the subdivision $X'' \to Y$ which is gotten by cutting the $\sigma'$'s with the hyperplanes suggested in the RHS of Figure~5.2. Because of a) above, for the same reasons as in the proof of the Sublemma 3.1-(B) we have the implication
$$
X'' \in (P) \Longrightarrow Y \in (P).
$$
This proves our Sublemma 3.1-A. $\Box$

\bigskip

In order to clinch the proof of the PL Lemma 3.1 we have to show how to go from cubical subdivisions to GPS structures without loosing neither the BLUE $2^{\rm d}$ collapsibility nor the RED property (P).

\smallskip

We start with the cubical subdivision $Y$ which has property (P) and which hence satisfied the Sublemma 3.1-A.

\bigskip

\noindent {\bf Sublemma 3.1-C.} {\it Via appropriate slidings and appropriate subdivisions, none of which violate, neither the BLUE $2^{\rm d}$ collapsibility nor the RED property (P), we can change $Y$ into the standard GPS structure:}

$$
\includegraphics[width=125mm]{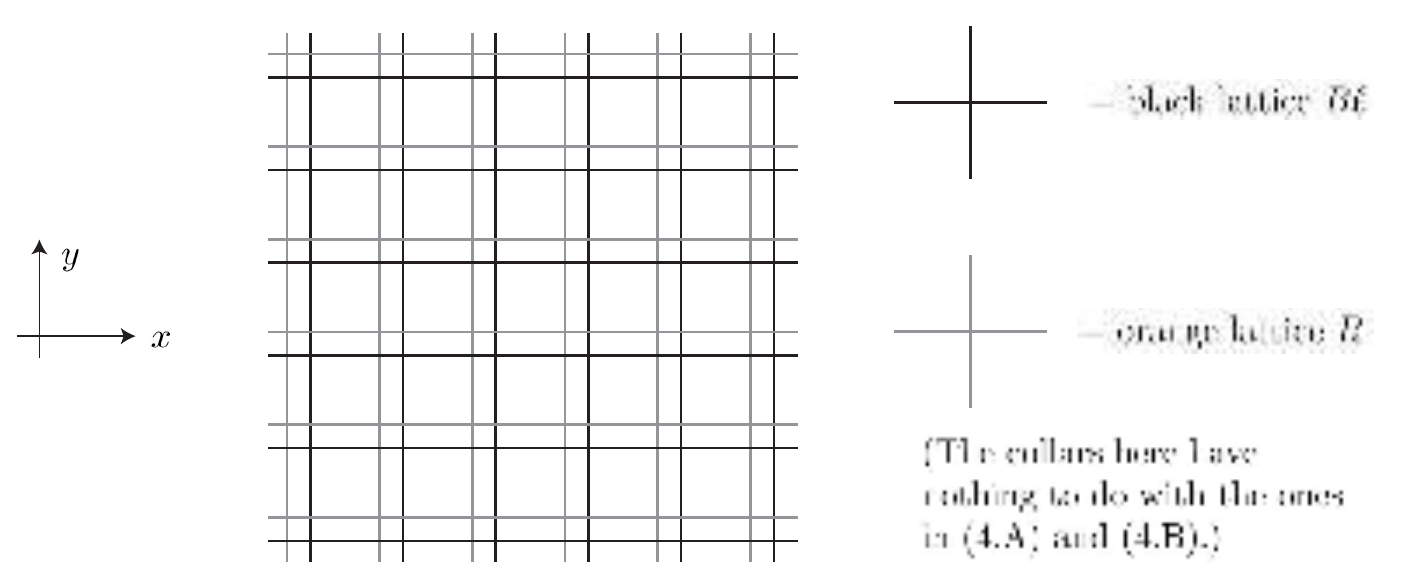}
$$
\label{fig5.3}
\centerline {\bf Figure 5.3.}

\medskip

\centerline{This figure goes with the explanations for DILATATING SLIDINGS.}

\bigskip

\noindent {\bf Proof.} All our manipulations now will be $2^{\rm d}$ and the extension to $3^{\rm d}$ and $4^{\rm d}$ will be canonically automatic. It is well-known that, generically speaking, $2^{\rm d}$ sliding moves of J.H.C. Whitehead violates collapsibility. But there is a category of sliding moves, which I will call {\ibf dilatating slidings} which do {\ibf not} violates collapsibility.

\bigskip

For these special slidings there is no disconnecting of smooth strata of maximal dimension. Here is our typical paradigmatic example.

\smallskip

Consider, to begin with the following $2^{\rm d}$ infinite complex $K$.

\bigskip

$K \equiv \{$the plane $z=0$ from Figure 5.3$\} \cup \{$for each line $L$ of the Black lattice $B\ell$ we consider the plane $L \times (-\infty < z < \infty)$, which we add to $(z=0)\}$.

\bigskip

\noindent This obviously has an infinite collapse $K \searrow {\rm pt}$. Move now from $K$ to the following

\bigskip

$K' \equiv \{$the same $z=0$ as above$\} \cup \{$for each line $L \in B\ell$ we add the $\frac12$-plane $L \times [-\infty < z \leq 0]\} \cup \{$for each line $\Lambda \in R$ we add the $\frac12$-plane $\Lambda \times [0 \leq z < \infty]\}$.

\bigskip

\noindent The $K \Rightarrow K'$ is our typical dilatating slide and clearly $K' \searrow {\rm pt}$ too. There are of course much simpler and more simple-minded examples.

\smallskip

We look now at the $2^{\rm d}$ skeleton of $Y$, call it $Y^2$, and which has the following obvious structure. We have a very fine cubical subdivision of the time axis $t \in R$, call it
$$
\ldots < t_{-1} < t_0 < t_1 < t_2 < \ldots
$$
Then, there is a time independent cubical subdivision $Z^2$ of $R^3 = (x,y,z)$ with 1-skeleton $Z^1$. With this, we have the following structure for our $Y^2$
$$
Y^2 = \left( \sum_i Z^2 \times t_i \right) \cup \left( Z^1 \times (-\infty < t < \infty) \right).
$$

$$
\includegraphics[width=10cm]{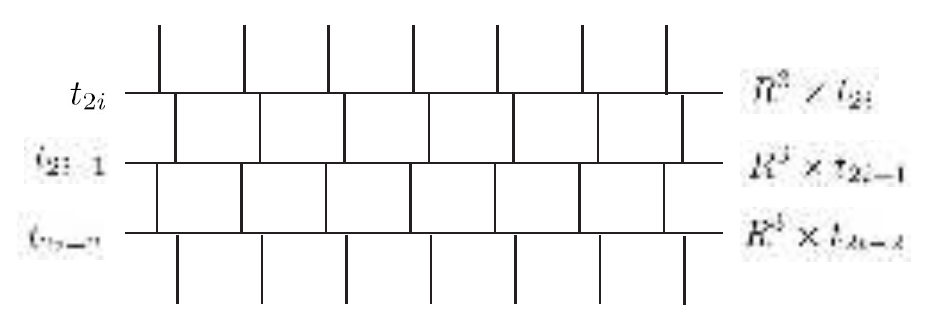}
$$
\label{fig5.4}
\centerline {\bf Figure 5.4.}

\centerline{
The vertical arcs of the form $[t_{2j-1} , t_{2j}]$ are here ORANGE.
}

\bigskip

We proceed now with the dilatating sliding suggested in Figure 5.4. This implements the e) from (5), and also f) from the same (5). On each individual $R^3 \times t_j$ we have a strictly cubical subdivision. From this structure one can then go to the standard GPS structure by more cubical subdivisions on the $R^3 \times t_j$'s, followed by small slidings of the dilatating type.

\smallskip

This ends the proof of our PL Lemma 3.1. $\Box$

\bigskip

We work now with the structure produced by the PL Lemma 3.1.

\bigskip

\noindent {\bf Lemma 4.} 1) {\it Inside $\Gamma (\infty)$ we will consider ``{\ibf spots}'', small open arcs far from any vertex. There will be two sets of them
$$
R(\mbox{for RED}) \subset \Gamma (\infty) \supset B (\mbox{for BLUE}), \leqno (7)
$$
and we should also think of them as $1$-handles. A given edge $e \in \Gamma (\infty)$ will carry at most one of the $r_i \in R$ or $b_j \in B$. Also, if by any chance we find $r_i \in e \ni b_j$, same $e$, then $r_i = b_j$, and this is how $R \cap B$ is generated.

\smallskip

We will introduce the more appropriate notations
$$
R \cap \Gamma (1) = \{ R_1 , R_2 , \ldots , R_n \} , \quad R - \{ R_1 , R_2 , \ldots , R_n \} = \{ h_1 , h_2 , \ldots \}. \leqno (8)
$$
The graphs $\Gamma (1) - R \cap \Gamma (1)$, $\Gamma (\infty) - R$, $\Gamma (\infty) - B$ are trees but, generically, the $\Gamma (1) - B$ is a disconnected union of {\ibf several} trees.}

\medskip

2) {\it The $X^2$ will be gotten by adding $2$-cells along a $\{\mbox{link}\} \subset \Gamma (\infty)$. We will have two disjoined partitions for this link, namely
\bigskip

\noindent $(9)$
\vglue -1cm
\begin{eqnarray}
\{\mbox{link}\} &= &\sum_1^{\bar n} \Gamma_i + \sum_1^{\infty} C_j + \sum_1^{\infty} \gamma_k^0 \ \mbox{(RED partition)} \nonumber \\
&= &\sum_1^{\infty} \eta_{\ell} + \sum_1^{\infty} \gamma_m^1 \ \mbox{(BLUE partition)} \nonumber
\end{eqnarray}
where, to begin with, $\underset{1}{\overset{\bar n}{\sum}} \ \Gamma_i \subset \Gamma(1)$ and the $2$-skeleton of $\Delta^4_{\rm Schoenflies}$ is here
$$
\Delta^2 = \Gamma(1) + \sum_1^{\bar n} D^2 (\Gamma_i);
$$
next,
\begin{eqnarray}
X^2 &= &\Gamma (\infty) \cup \sum_1^{\bar n} D^2 (\Gamma_i) \cup \sum_1^{\infty} D^2 (C_j) \cup \sum_1^{\infty} D^2 (\gamma_k^0) = \nonumber \\
&= &\Gamma (\infty) \cup \sum_1^{\infty} D^2 (\eta_e) \cup \sum_1^{\infty} D^2 (\gamma_m^1). \nonumber 
\end{eqnarray}
Here $D^2 (\mbox{curve})$ is the $2$-cell attached along the respective curve. Moreover, clearly $\bar n = \# \, D^2 (\Gamma) \geq n = \# \, R \cap \Gamma(1)$, $\# \, B \cap \Gamma (1) \geq n$.

\smallskip

The $D^2 (\gamma^0)$ (respectively $D^2 (\gamma^1)$) are exactly the $2$-cells which are killed by the $3^{\rm d}$ part of the RED (respectively BLUE collapse), in $(4)$. For the time being, at least, that is all we will say concerning the $3^{\rm d}$ collapses.

\smallskip

We have
$$
X^2 = \Gamma (\infty) \cup \sum_1^{\bar n} D^2 (\Gamma_j) + \sum_1^{\infty} D^2 (C_i) + \sum_1^{\infty} D^2 (\gamma_k^0) \supset X_0^2 \equiv \Gamma (\infty) \cup \sum_1^{\bar n} D^2 (\Gamma_j) + \sum_1^{\infty} D^2 (C_i).
$$
Actually, we will leave alive in $X_0^2$ a thin boundary collar, for each of the deleted $D^2 (\gamma_k^0)$'s.}

\medskip

3) {\it The $2^{\rm d}$ part of our RED and BLUE collapsing flows are expressed by the following two geometric intersection matrices

\medskip

\noindent $(10)$
\vglue -1cm
\begin{eqnarray}
C_j \cdot h_p &= &\delta_{jp} + \xi_{jp}^0 , \ \mbox{when $\xi_{jp}^0 \ne 0$ implies $j > p$}, \nonumber \\
\eta_{\ell} \cdot b_q &=& \delta_{\ell q} + \zeta_{\ell q}^0 , \ \mbox{when $\zeta_{\ell q}^0 \ne 0$ implies $\ell > q$}. \nonumber
\end{eqnarray}
} End of Lemma.

\bigskip

\noindent [{\bf An important digression.} The kind of infinite matrices which occur in (10) will be called of the {\ibf ``easy id + nilpotent''} type. The easy id $+$ nil implies GSC, as we will explain more in detail soon. But then, there is also the ``id $+$ nilpotent of the difficult type'', with the final inequalities occurring in (10) reversed. This occurs, for instance, for the classical Whitehead manifold ${\rm Wh}^3$, which certainly is {\ibf not} GSC. [It is not GSC for many reasons, the first one coming to my mind being that for open 3-manifolds, GSC implies $\pi_1^{\infty} = 0$.] Let us be a bit more precise concerning our condition of easy id $+$ nilpotent. To simplify the exposition, I will concentrate on $\eta_{\ell} \cdot b_q$, when $\Delta^2$ does not play any special role, but, with appropriate changes similar things can be said for $C_j \cdot h_p$ too.

\smallskip

To begin with, our indices belong to an appropriate ordered set, NOT necessarily totally ordered, and with $\zeta_{\ell q}^0 \ne$ in (10) implying $\ell > q$, satisfying also the following additional condition: From any element in our set of indices starts at least one back-going trajectory which is INFINITE. Trajectories mean things like $\ell \to q$, in the context above. Anyway, for the open manifolds concerned right now, this condition is natural.

\smallskip

Now, we are in a countable context for the $\eta_{\ell} \cdot b_q$ and so , we can re-index things compatibly with $\eta_{\ell} \cdot b_q =$ easy id $+$ nil, and make our set of indices be $Z_+$. [We have a monomorphism
$$
\{\mbox{our not necessarily totally ordered, but still ordered set}\} \to Z_+ .]
$$
Our geometric intersection matrix concerns the 2-skeleton $X^2$ of our $X^4 = R^4$ and, with $1$-skeleton $X^1 \subset X^2$, we have a tree $T$ to which the $h_q$'s and $D^2 (\eta_{\ell})$'s get attached, so as to get an $X^2 - \Sigma \, {\rm int} \, D^2 (\gamma^1) \subset X^2$. With this, what $C \cdot h =$ easy id $+$ nilpotent means, is that one can go from $T$ to the $X^2 - \underset{1}{\overset{\infty}{\sum}} \ \overset{\!\!\circ}{D^2}  (\gamma^1)$, by an infinite sequence of dilatation. Alternatively, this can be also expressed in one of the following ways:
\begin{enumerate}
\item[$\bullet$] There is an `` infinite collapse'' $X^2 - \underset{1}{\overset{\infty}{\sum}} \, \overset{\!\!\circ}{D^2}  (\gamma^1) \searrow T$, OR
\item[$\bullet\bullet$] (The GSC condition) The 1-handles cancell with the 2-handles in the manner suggested by the drawing below.
\end{enumerate}
$$
\includegraphics[width=10cm]{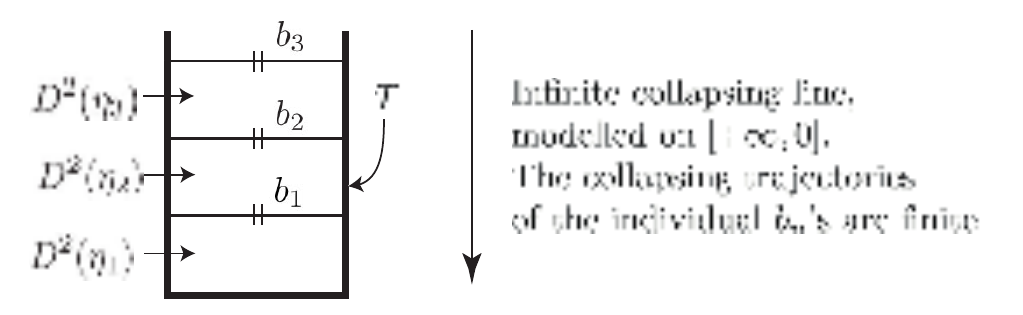}
$$
\label{fig5.5}

End of digression.]

\bigskip

In the digression above, we have worked with $X^2$ ($-$ the $D^2 (\gamma^1)$'s), but we will also need
$$
X_0^2 \equiv \Gamma (\infty) \cup \sum_1^{\bar n} D^2 (\Gamma_j) \cup \sum_1^{\infty} D^2 (C_j) \subset X^2 \ (X_0^2 = \{ X^2 \ \mbox{with the $D^2 (\gamma^0)$'s deleted}\}), 
\leqno (11)
$$
which comes with a RED $2^{\rm d}$-collapse $X_0^2 \searrow \Delta^2$. But then $X_0^2$ is clearly limping BLUE-wise; there the full $X^2$ is needed, and then the $\gamma^1$'s get deleted,~$\ldots$

\smallskip

We will need to extend this little theory, from $X^2$ to two successive, more elaborate contexts, the $X^2 [{\rm new}]$ and the $2X^2$. So, we will call, from now on, ``old'', the context of Lemma 4, the role of which was to be a pedagogical introduction. But then, our real-life contexts will be, successively, the NEW context, and then the ``DOUBLED'' one. Here is what these will do for us.

\smallskip

On the one hand, they will put our $\Delta^2$ in a protected position with respect to the BLUE flow. Then, they will eliminate the horrible complications stemming from the intersection of the two collapsing flows
$$
(\mbox{RED $2^{\rm d}$ collapsing flow}) \pitchfork (\mbox{BLUE $2^{\rm d}$ collapsing flow}) = \ ?
$$
Each of the individual flows is devoid of closed oriented orbits, but these may happily occur when we put them together. And, it is so that we could not live with such a thing.

\bigskip

\noindent {\bf A first change of viewpoint, old $\Rightarrow$ new.} The $X^2$ from Lemma 4 is called $X^2 ({\rm old})$ from now on.

\smallskip

And then, for the time being purely abstractly, we will introduce a fifth coordinate called $\xi_0$, in addition to $(x,y,z,t)$. We define
$$
X_0^2 [{\rm new}] \equiv \left( X_0^2 ({\rm old}) - \sum_1^{\bar n} \, \overset{\!\!\circ}{D^2} (\Gamma_j) \times (\xi_0 = 0) \right) \underset{\overbrace{\mbox{\footnotesize$\Gamma (1) \times (\xi_0 = 0)$}}}{\cup} [\Gamma (1) \times [0 \geq \xi_0 \geq -1]] 
$$
$$
\underset{\overbrace{\mbox{\footnotesize$\Gamma (1) = \Gamma (1) \times (\xi_0 = -1)$}}}{\cup} \Delta^2 \times (\xi_0 = -1).
\leqno (12)
$$

Forgetting for the time being about the BLUE collapse and keeping only the RED $2^{\rm d}$ collapsing flow alive, we have here a (RED) $2^{\rm d}$ collapse
$$
X_0^2 [{\rm new}] \searrow \Delta^2 \times (\xi_0 = -1) .
\leqno (12.0)
$$
This should be obvious, but we will write explicitly below the relevant geometric intersection matrices. Here is now OUR CHANGE OF VIEWPOINT: We decree, from now on, that the $\Delta^2 \times (\xi_0 = -1)$ {\ibf is} our $\Delta^2$ of interest, its 2-cells {\ibf are} the $D^2 (\Gamma_i)$; the 2-cells of $\Delta^2 \times (\xi_0 = 0)$ (which are anyway absent in $X_0^2 [{\rm new}]$ and which will afterwards reappear in $X^2 [{\rm new}] \supset X_0^2 [{\rm new}]$, see (12.1) below), become new $D^2 (\gamma^0)$'s and the new 2-cells involving an edge $[0 \geq \xi_0 \geq -1]$ become new $D^2 (C_j)$'s to be added to the previous $D^2 (C_j) \subset X_0^2 ({\rm old})$. When we talk about the $X_0^2 [{\rm new}]$, keep in mind that both the $\overset{\!\!\circ}{D^2}  (\gamma ^0) \subset X^2 ({\rm old})$ AND the  $\overset{\!\!\circ}{D^2}  (\Gamma_i) \times (\xi_0 = 0) \subset \Delta^2 \times (\xi_0 = 0)$, are absent, replaced by boundary collars.

\smallskip

Together, they will be the NEW family $D^2 (\gamma^0)$ of $X_0^2 [{\rm new}]$. Explicitly,
$$
D^2 (\gamma^0)[{\rm new}] = \{{\rm The} \ D^2 (\Gamma_i) \times (\xi_0 = 0)\} + \{{\rm The} \ D^2 (\gamma^0) ({\rm old})\}.
$$

The sites $R_i \times (\xi_0 = -1)$ for $i = 1,2,\ldots , n$ are promoted as NEW $R_i \subset \Delta^2 = \Delta^2 \times (\xi_0 = -1)$ and they live in $\Gamma(1) \times (\xi_0 = -1)$, of course. The old $h$'s stay put and, additionally to them, {\ibf every} edge $e \subset \Gamma (1) \times (\xi_0 = 0)$ acquires an $h$, out of the blue. [The old $R^i \times (\xi_0 = 0)$ are thereby all retrogrades as $h$'s.] BUT there is no $h$ on the new edges $P \times [0 \geq \xi_0 \geq -1]$, $P \in \Gamma_0 (1) \equiv 0$-skeleton of $\Delta^2$. The RED collapse of $X_0^2 [{\rm new}]$ proceeds as follows: We first collapse, normally $X_0^2 [{\rm new}] \supset X_0^2 ({\rm old}) \searrow \{ \Delta^2 \times (\xi_0 = 0)$ with the interiors of its 2-cells all deleted, they are now $D^2 (\gamma^0)$'s$\}$, then starting from $\Gamma(1) \times (\xi_0 = 0)$ we erase $\Gamma(1) \times [ 0 \geq \xi_0 > -1]$ and we are left with our
$$
\Delta^2 \equiv \Delta^2 \times (\xi_0 = -1), \ \mbox{as we should.}
$$

The RED story for $X_0^2 [{\rm new}]$, as told above, is OK, but since things are limping BLUE-wise, we introduce the following space, on the lines of (12) above
$$
X^2 [{\rm new}] = X^2 ({\rm old}) \cup [\Gamma(1) \times [0 \geq \xi_0 \geq -1]] \cup \Delta^2 \times (\xi_0 = -1). 
\leqno (12.1)
$$

Notice that each $e_i \subset \Gamma (1) \times (\xi_0 = 0)$ carries now an $h_i = h(e_i)$ its RED dual 2-cell is $D^2 (C_i) = e_i \times [0 \geq \xi_0 \geq -1]$. Similarly, the BLUE context of $X^2 [{\rm new}]$, each $e_i \subset \Gamma(1) \times (\xi_0 = -1)$ carries a $b(e_i)$ with dual the
$$
D^2 (\eta (e_i)) = e_i \times [-1 \leq \xi_0 \leq 0].
$$
These things generate the little schematical figure below.

$$
\includegraphics[width=9cm]{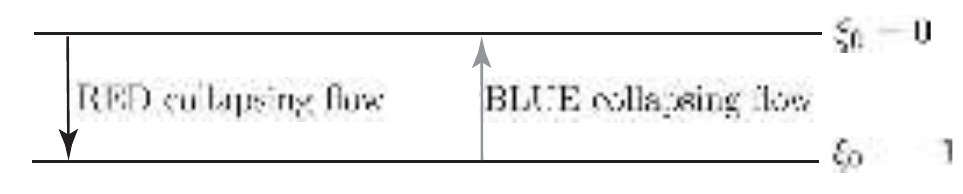}
$$
\label{fig5.bis}
\centerline {\bf Figure 5.bis}
\begin{quote}
The BLUE and RED collapsing flows inside 
$$
\Gamma (1) \times (0 \geq \xi_0 \geq -1) \subset X_0^2 [{\rm new}] \subset X^2 [{\rm new}].
$$
\end{quote}

Something should be stressed at this point. In the context of the collapse from (4-R, B) we had 2-cells $D^2 (\gamma^0)$, respectively $D^2 (\gamma^1)$ killed by the corresponding $3^{\rm d}$ collapses. In our present, purely $2^{\rm d}$ abstract context, before the RED (respectively the BLUE) $2^{\rm d}$ collapses, the $D^2 (\gamma^0)$'s (respectively the $D^2 (\gamma^1)$'s) have to be deleted, {\ibf by decree} (not by $3^{\rm d}$ collapse, there ain't any, at least not in the BLUE context, when our game will always be purely 2-dimensional). Of course all our $2^{\rm d}$ objects will eventually be thickened into $4^{\rm d}$ manifolds, but we will not talk about that now.

\smallskip

For expository purposes, a RED $3^{\rm d}$ collapse will be introduced below, and very transiently only, in the NEW context. Afterwards, a ghostly memory of it will survive in our really serious, DOUBLED context. But never from now on will then be any trace of any $3^{\rm d}$ BLUE collapse. Now, we define for the $X^2 [{\rm new}]$ in (12.1), the family
$$
D^2 (\gamma^1) [{\rm NEW}] = \{{\rm The} \ D^2 (\Gamma_i) \times (\xi_0 = -1)\} + \{{\rm The} \ D^2 (\gamma^1)({\rm old})\}.
\leqno (12.2)
$$
With this, we also introduce, for the same $X_0^2 [{\rm new}]$, the family
$$
B[{\rm new}] \equiv B({\rm old}) + \{\mbox{a $b_i$ on {\ibf every} edge} \ e_i \subset \Gamma(1) \times (\xi_0 = -1)\}.
$$
For these $b \in B[{\rm new}]$ the dual $\eta$'s are: For $b_j \in B({\rm old})$, the same $\eta_j$ as before, in $X^2 ({\rm old})$ and, for $b_i \subset e_i \subset \Gamma (1) \times[\xi_0 = -1]$ the dual object is $D^2 (\eta_i) = e_i \times [0 \geq \xi_0 \geq -1]$ with, of course $\eta_i = \partial D^2 (\eta_i)$, as already said.

\smallskip

The BLUE collapse $X^2 [{\rm new}] \searrow {\rm pt}$ proceeds as follows.
\begin{enumerate}
\item[i)] We delete the $D^2 (\gamma^1) [{\rm NEW}]$'s introduced above.
\item[ii)] Then, starting from $\Gamma(1) \times (\xi_0 = -1)$ we collapse away $\Gamma(1) \times [0 > \xi_0 \geq 1]$ which leaves us with $X^2 ({\rm old})$ which we collapse, then, normally.
\end{enumerate}

\bigskip

\noindent {\bf Addendum.} It also makes sense to introduce the following NEW RED $3^{\rm d}$ object
$$
X^3 [{\rm new}] \equiv \{ X^3 ({\rm old}) \underset{\overbrace{\mbox{\footnotesize$\Delta^2 \times (\xi_0 = 0)$}}}{\cup} \Delta^2 \times [0 \geq \xi_0 \geq -1] - \{\mbox{the interiors of the 3-cells of} \ \Delta^3 \subset X^3 ({\rm old})\}\},
$$
which comes endowed with a RED $3^{\rm d}$ collapse proceeding as follows.
\begin{enumerate}
\item[i)] Start with the normal collapse $X^3 ({\rm old}) \searrow \Delta^2 \times (\xi_0 = 0) \subset X^3 ({\rm old})$.
\item[ii)] Then, continue with $\Delta^2 \times [0 \geq \xi_0 \geq -1] \searrow \Delta^2 \times (\xi_0 = -1)$.
\end{enumerate}
The 2-cells which are demolished by this RED $3^{\rm d}$ collapse are exactly the
$$
\{{\rm old} \ D^2 (\gamma^0) \subset X^3 ({\rm old})\} + \{\mbox{the $D^2 (\Gamma) \times (\xi_0 = 0)$, retrograded as $D^2 (\gamma^0) [{\rm NEW}]$'s}\}.
$$
End of the description of the transformation OLD $\Rightarrow$ NEW. $\Box$

\bigskip

We will now move to the DOUBLED CONTEXT, where the objects introduced will be 2-dimensional only, for the time being. We continue to talk about sites $B,R$ and of 2-cells $D^2 (\eta)$, $D^2 (\gamma^1)$, $D^2(\gamma^0), \ldots$. There are no $3^{\rm d}$ collapses alive but, before the RED (respectively BLUE) $2^{\rm d}$ collapse can start, the $D^2 (\gamma^0)$'s (respectively the $D^2 (\gamma^1)$'s) have to be deleted, by {\ibf decree}.

\smallskip

[We do not throw away, yet, the RED $3^{\rm d}$ collapse, its ghost will be still of great use to us, although the collapse, as such will not be put into effect.]

\smallskip

Next comes the really important change of point of view, namely THE DOUBLING PROCESS. Still another abstract sixth axis is to be introduced, call it $-\infty < \zeta < +\infty$ and on it we fix three points labelled $r < \beta \ll b$ with $\vert \beta - r \vert \ll \vert b-r \vert$. Two copies of $X^2 = X^2 [{\rm new}]$ are to be considered now, call them respectively $X^2 \times r$ and $X^2 \times b$. The $X^2$ is here our friend $X^2 ({\rm new})$. What follows next will be, for the time being, a purely abstract $2^{\rm d}$ construction inside the space
$$
(X^2 \times r) \cup (\Gamma (\infty) \times [r,b]) \cup (X^2 \times b),
\leqno (13)
$$
where $[r,b] \subset (-\infty < \zeta < \infty)$ and where the first ``$\cup$'' is along $\Gamma(\infty) = \Gamma(\infty) \times r$ and the second along $\Gamma(\infty) = \Gamma(\infty) \times b$.

\smallskip

Notice that there are two kinds of edges $e \subset \Gamma(\infty)$, the edges $e(B)$ which contain a $b_i \in e$ (call them, specifically $e(b_i)$) and all the others, which I will chose to call, generically, $e(r)$. Among the $e(r)$'s we will find all the edges
$$
e = P \times [0 \geq \xi_0 \geq -1] , \quad P \in \{\mbox{vertices of} \ \Gamma(1)\} \subset \Delta^2 .
$$

In the Figure 7 below, the last edges above are in the category (IV).

\smallskip

In the context of (13), to each edge $e \subset \Gamma(\infty)$ corresponds a 2-cell $e \times[r,b] \subset \Gamma(\infty) \times [r,b]$.

\smallskip

We will use the notations
$$
(D^2 (c(b_i)) , c(b_i)) \equiv (e(b_i) \times [r,b] , \partial (e(b_i) \times [r,b]) \leqno (14)
$$
$$
(D^2 (c(r)) , c(r)) \equiv (e(r) \times [r,b] , \partial (e(r) \times [r,b]) , \ {\rm i.e.} \ c = \partial D^2 (c). 
$$

Remember that, by now the old context has been replaced by the new one and, in turn, this will be replaced itself by the DOUBLED context, i.e. we perform successively the changes
$$
{\rm old} \Longrightarrow {\rm new} \Longrightarrow {\rm DOUBLE}.
\leqno (15)
$$
But while the old context was only a pedagogical gimick, to be forgotten, both the NEW and the DOUBLED context (which does not quite superside the NEW one, for instance the BLUE flow on $X^2 [{\rm NEW}]$ is not the restriction of the BLUE flow on $2X^2$), will have to be used, both of them.

\smallskip

When we go to the DOUBLED context, then the higher analogue of $X_0^2$, the space of the $2^{\rm d}$ RED collapsing flow, will be the following $2^{\rm d}$ object where, of course $X_0^2 \times r = X_0^2 [{\rm new}] \times r$,

\bigskip

\noindent (16) \quad $2X_0^2 \equiv (X_0^2 \times r) \cup \{\Gamma (\infty) \times [r,b] \ \mbox{with any 2-cell} \ D^2 (C(b_i))$ deleted and replaced by a very thin tubular neighbourhood of its boundary. The newly created boundary component is re-baptized
$c(b_i)\} \cup \Bigl\{ \left( \underset{1}{\overset{\infty}{\bigcup}} \, D^2 (\eta_{\ell}) \right) \times b$, a space which I will call $X_b^2 \Bigl\}$.

\bigskip

Remember that at the level of our $X_0^2 \times r \approx X_0^2 [{\rm new}]$ all the $D^2 (\gamma^0)[{\rm NEW}] = D^2 (\Gamma_i) \times (\xi_0 = 0)$ are {\ibf deleted} and actually replaced by thin tubular neighbourhoods of their boundaries. This allows us to write generically
$$
\partial (2X_0^2) = \sum_k \gamma_k^0 [{\rm NEW}] + \sum_{b_i \in B[{\rm NEW}]} c(b_i) \supset \sum_j \Gamma_j \times (\xi_0 = 0) \ \mbox{(on the $r$-side)},
\leqno (16.1)
$$
with an obvious small twist of notation. [Any 2-cell which gets deleted is replaced by a boundary collar, the exterior frontier of which occurs now in (16.1).]

\smallskip

The 1-skeleton of $2X_0^2$ is the following object
$$
2\Gamma (\infty) \equiv (\Gamma (\infty) \times r) \cup (\Gamma_0 (\infty) \times [r,b]) \cup (\Gamma (\infty) \times b),
$$
when $\Gamma_0 (\infty) \subset \Gamma (\infty)$ is the 0-skeleton of $\Gamma(\infty)$.

\smallskip

BLUE-wise, our $2X_0^2$ is limping, reason for introducing a higher analogue of $X^2$ too, namely the
$$
2X^2 \equiv (X_0^2 \times r) \cup (\Gamma (\infty) \times [r,b]) \cup X_b^2 \supset 2X_0^2 .
\leqno (17)
$$
The reason why, inside $2X^2$ we use $X_0^2 \times r \approx X_0^2 [{\rm new}]$, just like for $2X_0^2$, and {\ibf not} the seemingly more appropriate $X^2 [{\rm new}] \times r$, will soon be crystal clear.

\smallskip

In the context of $X^2 [{\rm new}]$, the $2^{\rm d}$ RED and BLUE collapses were expressed in terms of sites $R,B$ which, for $X^2 [{\rm new}]$ I will denote them now by $R_0 , B_0$ and the corresponding curves $\Gamma , C , \gamma^0$, respectively $\eta , \gamma^1$. [We speak now about $X^2 [{\rm new}]$ while the initial $X^2$ call it $X^2 ({\rm old})$ will never play any role any longer, and will be forgotten, whenever the contrary is not explicitly said.]

\bigskip

\noindent {\bf Lemma 5.} 1) {\it We can endow $2X^2 \supset 2X_0^2$ with RED and BLUE sites $R_1 \supset R_0$, $B_1 \supset B_0$ and with extended system of curves so that $2X_0^2$ should carry a $2^{\rm d}$ RED collapsing flow $2X_0^2 \searrow \Delta^2 = \Delta^2 \times (\xi_0 = -1) \subset X_0^2 \times r$ and a $2^{\rm d}$ BLUE collapse $2X^2 - \sum \, {\rm int} \, D^2 (\gamma^1)$ (to be defined in the doubled case) $\searrow$ {\rm pt}. It will turn out that the $2X_0^2$, as defined, is already free of $D^2 (\gamma^0)$'s so we did not have to take, for defining its RED collapse, the same precautions as for the BLUE collapse of $2X^2$.

\smallskip

No $3^{\rm d}$ collapse will ever be considered at the doubled level but, at the level of $X_0^2 [\mbox{new}]$, but for technical reason, we will still need to invoke later on, the RED $3^{\rm d}$ collapse and a ghostly memory of it will linger at the level of our $2X_0^2$.}

\medskip

2) {\it Here are the explicit spots ($1$-handles) and curves (attaching zones of $2$-handles or simply boundary components), contained in $\Gamma (2\infty)$, at the double level. The $1$-handles (or spots) are 

\bigskip

\noindent {\rm (18)} \quad $R_1 \equiv \{\mbox{The family} \ R_0 \subset X_0^2 \times r \} + \{ b_i \times b \} (= \{e(b_i) \times b \}) + \{e(r) \times b \} = \left[\underset{1}{\overset{n}{\sum}} \, R_i \times (\xi_0 = -1) \right] + \underset{1}{\overset{\infty}{\sum}} \, h_n$, with $b_i \in B_0$ and with $e(r)$ like in the Figure {\rm 7-(I, IV)}, $e(b)$ like in Figure {\rm 7-(II, III)}, and, of course, with $\{b_i \times b \} + \{ e(r) \times b \} \subset X^2 \times b$; then
$$
B_1 = \underbrace{\{ b_i \times r \}}_{\mbox{\footnotesize The $B_0$, living in $X_0^2 \times r$}} + \, \{ b_i \times b \} + \{ e(r) \times b \},
$$
the last two items being common for $R_1$ and $B_1$, and see here the Figures {\rm 7} and {\rm 7-bis}. We have $B_0 = B[\mbox{new}]$ (introduced in the context $X^2 [\mbox{new}]$ {\rm (12.1)}. End of {\rm (18)}.

\bigskip

Next we have

\bigskip

\noindent {\rm (18.1)} \quad $\{$extended set of $C$'s$\} = \{ C \subset X_0^2 \times r\} + \{ c(r)$ (see {\rm (14)}) $+ \{\eta \times b \subset X_b^2 \}$; $\{$extended set of $\gamma^0$'s$\} = \{$the $\gamma^0 \subset X_0^2 \times r$ with the $\partial \Gamma_i \times (\xi_0 = 0)$ included$\} + \{c(b)\}$ (see {\rm (14)}); $\{$extended set of $\eta$'s$\} = \{\eta \times b \subset X_b^2 \} + \underset{\mbox{in} \ \Gamma (\infty) \times [r,b]}{\underbrace{\{ c(r) \} + \{c(b)\}}}$; $\{$extended set of $\gamma^1$'s$\} = \{\Gamma_i \subset (\xi_0 = 0$ and $-1)\} + \{ C_j \}$.

\bigskip

So, {\ibf all} $2$-cells of $X_0^2 \times r$ are, BLUE-wise speaking, $D^2 (\gamma^1)$'s, when we go to $2X^2$, and this is the reason to define the $2X^2$ in {\rm (17)} as we did. The $2^{\rm d}$ BLUE flow of $2X^2$ is now mute on the piece $X_0^2 \times r$.

\smallskip

Here are some EXPLANATIONS concerning the lines $2$ (the $\gamma^0$'s) and $4$ (the $\gamma^1$'s) in {\rm (18.1)}. None of the $\gamma^0 \in \gamma^0$ (of $X_0^2 \times r = X_0^2 [\mbox{new}]$) are physically there, neither in $2X_0^2$ nor in $2X^2$. But the ghost of the RED $3^{\rm d}$ collapse (where the $\gamma^0$'s do appear) will be used, reason for including those $\gamma^0 [\mbox{NEW}] \subset X_0^2 [\mbox{NEW}]$ in our formulae. The $\{ D^2 (\gamma^0)\} \ni D^2 (c(b)) \subset 2X^2$ and they have to be deleted in $2X_0^2$ before the RED $2$-collapse can start. [The $D^2 (c(b))$'s are $\{$extended $D^2 (\gamma^0)$'s$\}$.] Similarly, the $D^2 (\mbox{extended} \ (\gamma^1))$ have to be deleted so as to proceed to the BLUE $2^{\rm d}$ collapse of $2X^2 \searrow {\rm pt}$. With all these things, here are the two basic RED and BLUE geometric intersection matrices, at levels $2X_0^2 \subset 2X^2$. It is these matrices which define our two $2^{\rm d}$ collapsing flows

\bigskip

\noindent {\rm (19)} \quad $\{$extended set of $C$'s$\} \cdot \{$extended set of $h$'s$\} =$ easy id $+$ nilpotent, where $\sum \, h$ (extended) $\equiv R_1 - \Gamma(1) \times (\xi_0 = -1)$, and remember that $\underset{1}{\overset{n}{\sum}} R_i \times (\xi_0 = -1)$ are not to be mixed up with the bigger set $\sum \, h_n \subset R_1$; then $\{$extended set of $\eta$'s$\} \cdot \{$extended set of $B_1$'s$\} =$ easy id $+$ nil. End of {\rm (19).}

\bigskip

Here are the explicit dualities establised by the diagonal $\delta_{ij}$'s of our matrices, between sets of $1$-handles and sets of attaching curves of $2$-handles.

\bigskip

\noindent {\rm (20)} \quad (RED duality) $\{ C_i \subset X_0^2 \times r\} \underset{\approx}{\longleftrightarrow} \left( \underset{1}{\overset{\infty}{\sum}} \, h_n \right) \cap X_0^2 \times r$, occurring in the $R_0 - B_0$, or $R_0 \cap B_0$, and see here Figure {\rm 7-(I)}, respectively {\rm 7-(II)},
$$
\{\eta_i \times b \} \underset{\approx}{\longleftrightarrow} \{ b_i \times b \} \quad \mbox{and} \quad \{ c(r) \} \underset{\approx}{\longleftrightarrow} \{ e(r) \times b \}.
$$

Here the $\{ b_i \times b \}$, $\{ e(r) \times b \}$ are the $\left( \underset{1}{\overset{\infty}{\sum}} \, h_n \right) \cap X_b^2$ and they occur, respectively, the Figures {\rm (7-II, III) $+$ (7.bis)} and {\rm 7-(I, IV)}. All of them are $R_1 \cap B_1$'s.

\bigskip

\noindent {\rm (21)} \quad (BLUE duality) $\{ \eta_i \times b \} \underset{\approx}{\longleftrightarrow} \{ b_i \times b \}$ (and this is the BLUE duality of $X^2 [\mbox{new}]$, {\ibf transported} from $X^2 \times r$ to $X^2 \times b$), then
$$
\{ c(r)\} \underset{\approx}{\longleftrightarrow} \{ e(r) \times b \} , \quad \{ c(b_i)\} \underset{\approx}{\longleftrightarrow} \{ b_i \times r \in e(b) \times r \}.
$$

The first of these occurs in Figure {\rm 7-(I, IV)} and the second in {\rm 7-(II, III) $+$ 7.bis}. With all these things the two $2^{\rm d}$ RED and BLUE flows, on $2X_0^2$, respectively on $2X^2$, are completely defined. But we will explicit them even more.

\bigskip

\noindent {\rm (21.1)} \quad The BLUE $2^{\rm d}$ flow goes like follows, for $2X^2$
$$
\underbrace{\{ b_i \times r \}}_{{\mbox{\footnotesize no internal blue} \atop \mbox{\footnotesize flow lines among}} \atop \mbox{\footnotesize these guys}} \longrightarrow \underbrace{\{ b_i \times b \}}_{{\mbox{\footnotesize All the blue flow lines of} \atop \mbox{\footnotesize $X^2 [\mbox{new}]$, given by $\eta_j \cdot b_j$}} \atop \mbox{\footnotesize concerns {\smallibf these} guys}} \longrightarrow \underbrace{\{e(r) \times b \}}_{{\mbox{\footnotesize no internal blue} \atop \mbox{\footnotesize flow lines among}} \atop \mbox{\footnotesize these guys}}.
$$

To make things completely explicit here are also the off-diagonal terms in our geometric intersection matrices.

\bigskip

\noindent {\rm (21.2)} \quad (RED case) On $X_0^2 [\mbox{new}] = X_0^2 \times r$ we have the $C \cdot h$ of $X_0^2 [\mbox{new}]$, with the corresponding off-diagonal terms. On $X_b^2$ we have the $\eta \cdot b$ of $X^2 [\mbox{new}]$, with the corresponding off-diagonal terms. For the part of the matrix concerning $c(r)$ the situation is completely readable in Figure {\rm 7-(I, IV)}. In the cases of {\rm 7-(I, IV)} we have off-diagonal terms $e(r) \times b \xrightarrow{ \ C \cdot h \ } e(r) \times r$ and, while in {\rm (I)} things continue with $C \cdot h \mid X_r^2$, at $e(r) \times r$ in {\rm (IV)} the RED flow stops (Dead End).

\bigskip

\noindent {\rm (21.3)} \quad (BLUE case) On $X_b^2$ we have the $\eta \cdot b$ of $X^2 [\mbox{new}]$ with the same off-diagonal terms as in the RED case. For $c(r)$, $c(b)$ the situation is readable respectively as Figure {\rm 7-(I, IV)} OR {\rm ((7-(II, III)) $+$ (7.bis))}. In the cases {\rm 7-(I)} and {\rm 7-(IV)} things are exactly as in the RED case. In the situation {\rm ((7-(II, III)) $+$ (7.bis))}, we have $e(b) \times r \underset{\eta \cdot B}{-\!\!\!-\!\!\!-\!\!\!\longrightarrow} e(b) \times b$, as off-diagonal term.

\bigskip

So, keep in mind that, RED-wise, at level $2X_0^2$, on the $X_0^2 \times r$ side we just leave the $C \cdot h$ of $X_0^2 [\mbox{new}]$, while the $\eta \cdot B$ of $X^2 [\mbox{new}]$ is transported on $X_b^2$, RED-wise, not only BLUE-wise.

\smallskip

Formulae {\rm (20), (21)} and Figure {\rm 7} should help making explicit our two collapses. Here is, in detail, the RED $2^{\rm d}$-collapse of $2X_0^2$ (at the level of which the $D^2 (\gamma^0)$'s are already deleted).}

\begin{enumerate}
\item[$\bullet$)] {\it We collapse away $X_b^2$ using the sites $e(b) \times b = b_i \times b$ in Figure~{\rm 7-(II, III), 7.bis}, which are free on their left side. This leaves us free to unleash the flow $\eta \cdot B$ [of $X^2 [\mbox{NEW}]$] on the $X_b^2$, where it has been transported and demolish the $X_b^2$. There is, of course, no RED action on $\Delta^2 \times (\xi_0 = -1)$ but, for the $b \in B_0 \cap \Delta^2 \times (\xi_0 = -1)$, the $b \times b \in \Gamma (\infty) \times b$ partake into the RED flow of $X_b^2$. [See here the $[BCGF]$ in Figure~{\rm 19.1}.] Also, we talk here in terms of a mythical ``infinite collapse''; but what one should read, in real life is actually the following story: when we consider $\Gamma (\infty) \times b \supset \{ b_i \times b \}$, for all the $b_i \in B_0$, and to which the $D^2 (\eta) \times b$ are attached, this is the same thing, up to isomorphism as $X^2 ({\rm old}) - \sum \, \overset{\!\!\circ}{D^2}  (\gamma^1)$ with $B_i \subset \Gamma (\infty)$. So all the $\eta \cdot B$ game can be played on the $X_b^2$ side now, and leave us, on the $b$-side with only $\Gamma (\infty) \times b - \{ b_i \times b \}$ alive. From here on we can move to the next $\bullet\bullet)$.}
\item[$\bullet$$\bullet$)] {\it Next we collapse the pairs $(e(r) \times b , e(r) \times [r,b])$, occurring in the Figures~{\rm 7-(I and IV)}. They are now free on their right side, because of $\bullet)$. By now, also, after the collapse in $\bullet\bullet)$, only $(\Gamma (\infty) \times r) \, \cup \, \Gamma_0 (\infty) \times [r,b]) \subset 2\Gamma (\infty)$, is alive, as far as the $1$-skeleton is concerned, plus the $X_0^2 \times r$, of course.}
\item[$\bullet$$\bullet$$\bullet$)] {\it Finally, we collpase normally
$$
\{ X_0^2 \times r = X_0^2 [\mbox{new}]\} - \sum_1^{\bar n} D^2 (\Gamma_i) \times (\xi_0 = 0) \searrow \Delta^2 = \Delta^2 \times (\xi_0 = -1).
$$
End of the RED collapse.}
\end{enumerate}

{\it Here is now, again in detail, the BLUE $2^{\rm d}$ collapse of $2X^2 \supsetneqq 2X_0^2$:}

\begin{enumerate}
\item[$\bullet$)] {\it We start by eliminating all the ${\rm int} \, D^2 (\mbox{extended} \ \gamma^1)$, which includes now all the $D^2$'s of $X^2_0 [\mbox{NEW}] \approx X_0^2 \times r$, including the $D^2 (\Gamma_j) \times (\xi_0 = -1)$, extending thus the {\rm (12.2)}. This is an important feature of the DOUBLING.}
\item[$\bullet$$\bullet$)] {\it Then we collapse away the pairs $(e(b) \times r , e(b) \times [r,b])$ which are now free on their left side, Figures}~7-(II, III), 7-bis.
\item[$\bullet$$\bullet$$\bullet$)] {\it Afterwards, we collapse normally $X_b^2$ starting at the site $e(b) \times b$, which are by now free. In other terms, we perform the BLUE collapse of $X^2 [\mbox{new}]$, transported to the $X_b^2$-side.}
\item[$\bullet$$\bullet$$\bullet$$\bullet$)] {\it Finally, we collapse away the $(e(r) \times b, e(r) \times [b,r])$, in Figure~{\rm 7-(I, IV)}, and this leaves only with a tree alive. End of the collapsing story.}
\end{enumerate}

\medskip

3) {\it ({\ibf Punch line}) Here are the two goals which the DOUBLING has achieved, namely}

\smallskip

A) {\it We have now, in agreement with {\rm (12.2)}
$$
\{ \Gamma \times (\xi_0 = -1)\} \subset \{\mbox{extended $\gamma^1$'s}\}. \leqno (21.{\rm A})
$$
[Actually already at the level $X^2 [\mbox{new}]$, $\{ \Gamma \times (\xi_0 = -1)\} \subset \{\gamma^1\}$, and this is not violated  by the passage NEW $\Longrightarrow$ DOUBLE.]}

\bigskip

\noindent (21.B) \quad {\it At the level of $2X^2$ the RED and BLUE never cut transversally through each other. When a $2$-cell $\sigma^2 \subset 2X^2$ contains both a RED and a BLUE collapsing arrow, they coincide. When before the DOUBLING they might have happily cut through each other, the RED and BLUE arrow become now parallel, or disjoined, if you wish.}

\medskip

4) {\it ({\ibf A strategical decision}) All the edges $e \subset \Gamma (1) \times (\xi_0 = -1)$ contain a $b = b(e) \in B_1$, so they are in the same boat as the edges $e(b)$ from the Figures~{\rm 7-(II, III)}. Moreover, there are no other edges at $\xi_0 = -1$.

\smallskip

Now, our {\ibf strategic decision} is that for exactly all the edges at $\xi_0 = -1$ we completely erase the shaded collar which is visible in the Figures~{\rm 7-(II, III)}. What this achieves is the following (and see here Figure~{\rm 7.bis}, replacing {\rm 7-(II, III)} at $\xi_0 = -1$).}

\bigskip

\noindent (21.C) \quad {\it The only $2^{\rm d}$ pieces not in $\Delta^2 \times (\xi_0 = -1)$ but adjacent to it, come from $\Gamma (1) \times [0 \geq \xi_0 \geq -1$]. In particular, there is no $2^{\rm d}$ piece $A^2 \subset \Gamma (\infty) \times [r,b]$ adjacent to $\Delta^2 \times (\xi_0 = -1)$.}

\smallskip

\noindent End of Lemma 5.

\bigskip

\noindent {\bf Very important remark.} We may sometimes have to use subdivisions which normally would have to preserve the GPS system. But, when we apply such a subdivision, let us say to $X^2 [{\rm new}]$, then $D^2 (\Gamma) , D^2 (C) , D^2(\eta)$ break into several 2-cells with the same labels, while for a normal subdivision we find
$$
D^2 (\gamma^0) \Longrightarrow \{\mbox{a unique smaller $D^2 (\gamma^0)$ and many small $D^2 (C)$'s}\},
$$
$$
D^2 (\gamma^1) \Longrightarrow \{\mbox{a unique smaller $D^2 (\gamma^1)$ and many small $D^2 (\eta)$'s}\}.
$$
The last line above comes with the big potential danger of destroying the feature (21.A) which we will very much need. But we will {\ibf need} to subdivide $\Delta^2 \times (\xi_0 = -1)$ too, in section V below (CONFINEMENT). And then, so as to preserve (21.A) we will similarly subdivide $\Delta^2 \times (\xi_0 = 0)$ too, and for {\ibf any} vertex $P \in \Delta^2$ we will add a line $P \times [0 \geq \xi_0 \geq -1]$, see here how we proceed in Figure 30 below. The whole of $\Delta^2 \times [0 \geq \xi_0 \geq -1]$ is subdivided uniformly there. The same important remark is valid after DOUBLING, of course. $\Box$

\bigskip

$$
\includegraphics[width=144mm]{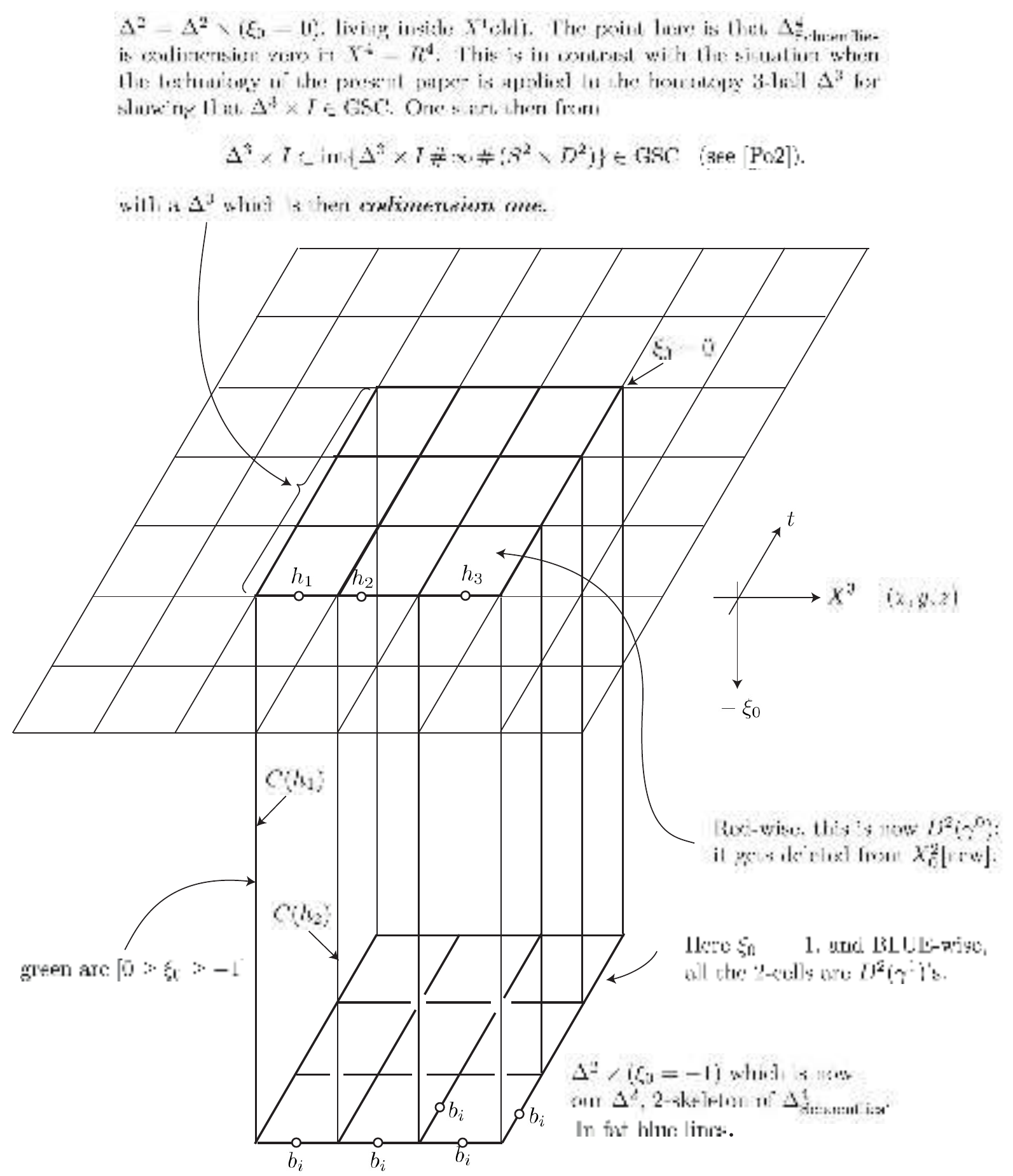}
$$
\label{fig6}
\centerline {\bf Figure 6.}
\begin{quote}
Schematical representation of $X^2 [{\rm new}] \supset X_0^2 [{\rm new}]$. Our drawing  is, oversimplified, at least in two respects: we have failed to represent the GPS cellular structure and, of course also, the dimensions are highly irrealistic. But it should be stressed that {\ibf each} vertex of $\Delta^2 \times (\xi_0 = -1)$ is joined via a green arc to correspondant vertex in $\Delta^2 \times (\xi_0 = -1)$. In this figure we see the $\Gamma (1) \times (\xi_0 = 0)$ represented in fat black lines ($={\bm -\!\!\!\!\!\bm -\!\!\!\!\!\bm -\!\!\!\!\!\bm -\!\!\!\!\!\bm -}$), the $[0 \geq \xi_0 \geq -1] \times P$'s in green lines ($= -\!\!\!-\!\!\!-)$ and the $\Gamma (1) \times (\xi_0 = -1)$ in fat blue lines ($={\bm -\!\!\!\!\!\bm -\!\!\!\!\!\bm -\!\!\!\!\!\bm -\!\!\!\!\!\bm -}$). Every edge $e_i \times (\xi_0 = -1)$ carries a $b_i \in B_0$, every $e_j \times (\xi_0 = 0)$ carries a $h_j$ with the dual cell being $D^2 (C_j) = e_j \times [0 \geq \xi_0 \geq -1]$ and the edges $P \times [0 \geq \xi_0 \geq -1]$ are free of $R \cup B$. The edges $e \times (\xi_0 = 0)$ may carry $b \in B_0$'s too.
\end{quote}

\bigskip

Figure 7 illustrates what  goes on along $\Gamma (\infty) \times [r,b] \subset 2X^2 \supset 2X_0^2$. We see a 2-cell $e \times [r,b]$, when $e \subset \Gamma (\infty)$ is an edge, with $e \subset X^2 [{\rm new}]$ and there are four cases to consider, namely the following
$$
\mbox{(Case I) $R_0 - B_0 \in e$, \ (case II) $R_0 \cap B_0 \in e$, \ (case III) $B_0 - R_0 \in e$, \ (case IV) $R_0 \cup B_0 \notin e$.}
$$
$$
\includegraphics[width=15cm]{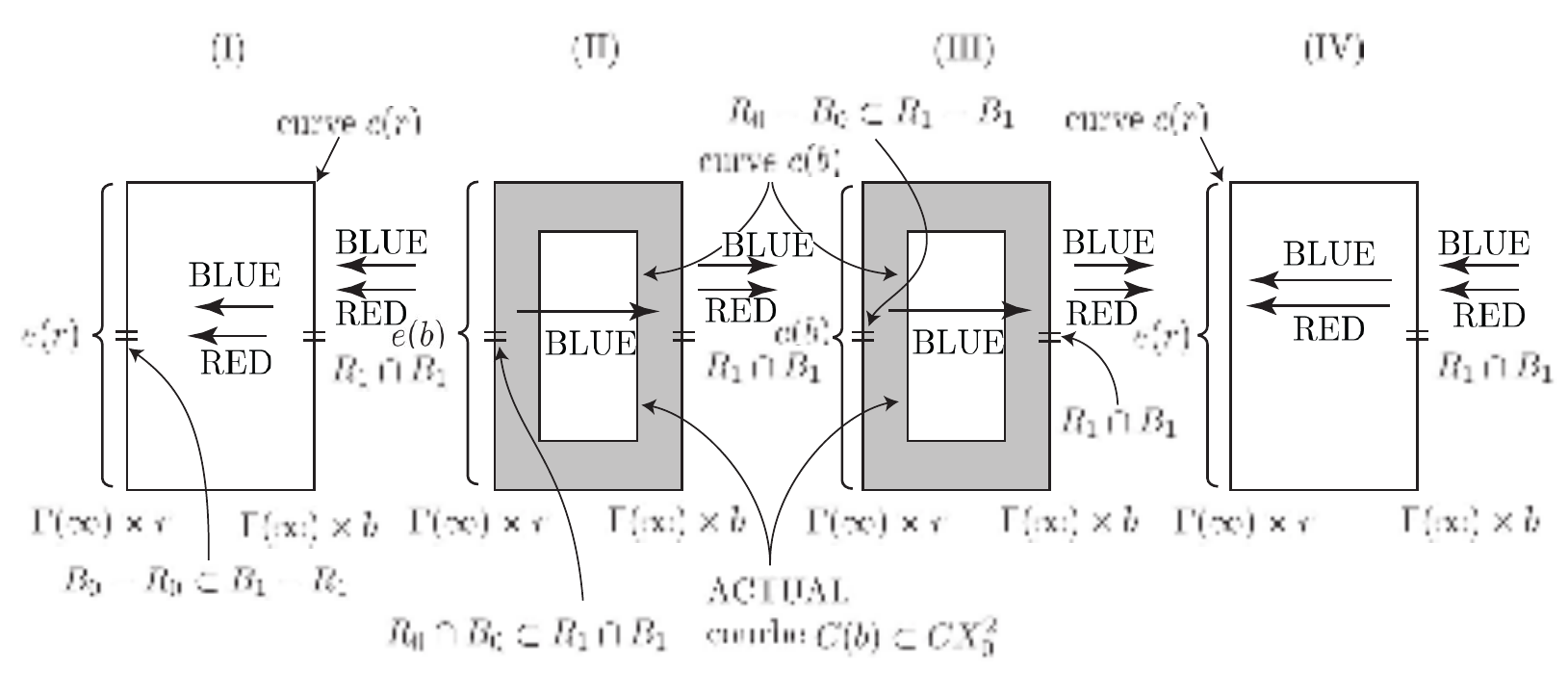}
$$
\label{fig7}
\centerline {\bf Figure 7.}
\begin{quote}
The discs $e \times [r,b] \subset 2X^2$, when far from $\xi_0 = -1$, where Figure 7.bis replaces the present Figures 7-(II, III). At the level $2X_0^2$, only the shaded part of II, III survives.

LEGEND: $\Relbar\!\!\!\!\!\vert \ =$ spot, in $R_1 \cup B_1$, $\xleftarrow{ \ {\rm BLUE} \ } \, =$ BLUE $2^{\rm d}$ collapsing flow (for $2X^2$), $\xleftarrow{ \ {\rm RED} \ } \, =$ RED $2^{\rm d}$ collapsing flow (for $2X_0^2$). These flows are represented here only to the extent they affect the $\Gamma (\infty) \times [r,b] \subset 2X^2$. More RED arrows may actually also leave from $e(r)$ in (I), (IV) to the left, in $X_0^2 \times r$, but we have failed to represent this here. The $X_b^2$ lives to the right of these figures. In (I), (IV), the outer curve of the rectangle is of the type $c(r)$ while in (II), (III) it is of type $c(b)$; but then when going from $2X^2$ to the smaller $2X_0^2$, the name is inherited by the corresponding smaller curve in $\partial (2X_0^2)$. The RED arrows inside (I), (IV) correspond, in terms of the dualities from (20), to $\{ c(r)\} \underset{\approx}{\longleftrightarrow} \{ e(r) \times b \}$. The BLUE arrows inside (I), (IV) corresponds to the same $\{ c(r)\} \underset{\approx}{\longleftrightarrow} \{ e(r) \times b \}$ (which occurs in both geometric intersection matrices). The outgoing BLUE and RED arrows of (II), (III) correspond to the $\{ \eta_i \times b \} \underset{\approx}{\longleftrightarrow} \{ b_i \times b \}$ in our matrices. In the context of (I), (IV) the BLUE and RED incoming arrows in $e(r) \times b \in R_1 \cap B_1$, just means that the corresponding edges $e(r) \times b$ may receive such arrows from $X_b^2$, without sending themselves any, back into $X_b^2$. To the left of the $e \times r$ in (I) or (II) we just have the normal RED flow inside $X^2 [{\rm new}]$, which we did not represent here by arrows, as already said.

Any edges $e \subset \Delta^{(2)} \times (\xi_0 = -1)$ would find itself normally as an $e(b)$ in (II), (III). But since now the whole collar is deleted, in agreement with 4) in Lemma 5, we have redrawn the situation at $(\xi_0 = -1)$ in Figure 7.bis. Finally, for any vertex $P \in \Gamma (1)$, the $e = P \times [0 \geq \xi_0 \geq -1]$ is among the $e(r)$'s in (IV).
\end{quote}
$$
\includegraphics[width=9cm]{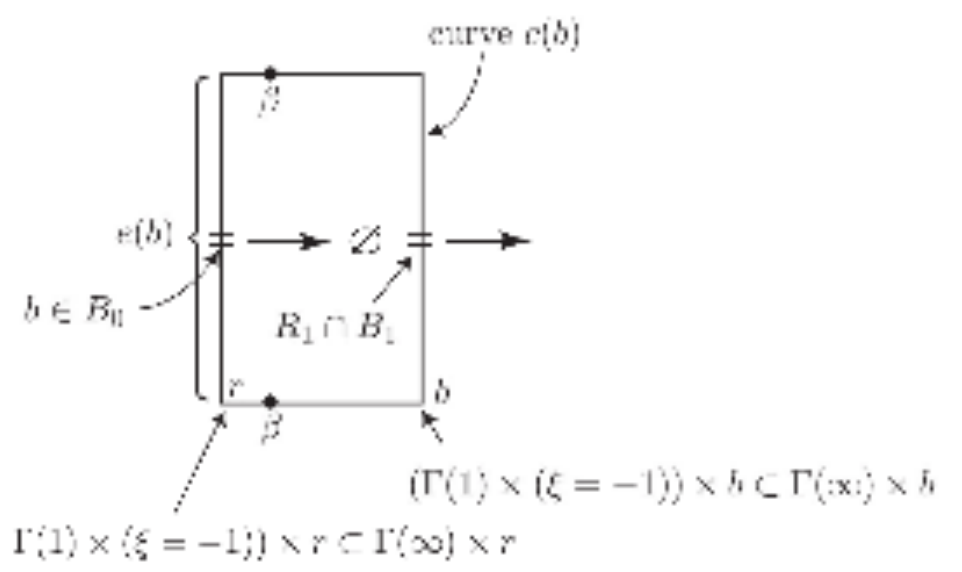}
$$
\label{fig7.bis}
\centerline {\bf Figure 7.bis}
\begin{quote}
Here is what the Figures 7-(II, III) become when the corresponding $e(b) \subset \Delta^2 \times (\xi_0 = -1)$. {\ibf All} the edges of $\Delta^2 \times (\xi_0 = -1)$ are like this $e(b)$.

LEGEND: ${\bm -\!\!\!\!\!\bm -\!\!\!\!\!\bm -\!\!\!\!\!\bm -\!\!\!\!\!\bm -} \, =$ the curve $c(b)$ for $b \in \Delta^2 \times (\xi_0 = -1)$. $\longrightarrow \, =$ BLUE flow at doubled level. No collar is surviving here, unlike in Figures 7-(II, III). All arrows are here BLUE.
\end{quote}


\section{Four dimensional thickenings and compactifications}\label{sec3}

We start now from
$$
2X_0^2 \supset X_0^2 \times r \approx X_0^2 [{\rm NEW}] \supset \Delta^2 \times (\xi_0 = -1) = \{ \Delta^2 , \ \mbox{2-skeleton of} \ \Delta_{\rm Schoenflies}^4\},
$$
and we want to define $4^{\rm d}$ regular neighbourhoods for $2X_0^2$. Now, there is no a priori given DIFF 4-manifold into which $2X_0^2$ is embedded, or even immersed, and we will start by defining, {\ibf by decree}, local $4^{\rm d}$ thickenings. Then, going from local to global, we glue together these pieces, making use of appropriate framings, when 2-handles are to be added. This construction willl be restricted by two conditions which it will {\ibf have to satisfy}

\medskip

\noindent (22.A) \quad (Compatibility with Schoenflies) So as not to loose the connection with the Schoenflies issue, we will insist that $\{ N^4 (2X_0^2)$ defined by glueing the chosen local piece above$\} \mid \Delta^2 \underset{\rm DIFF}{=} \{$the normal $N^4(\Delta^2)$ from Theorem 2, i.e. the normal $4^{\rm d}$ regular neighbourhood of the $\{2^{\rm d}$ skeleton $\Delta^2$ of $\Delta_{\rm Schoenflies}^4\} \subset \Delta_{\rm Schoenflies}^4\}$.

\medskip

\noindent (22.B) \quad Compatibility with the $2^{\rm d}$ RED collapsing: We do insist that when the $N^4 (2X_0^2)$ is expressed as
$$
\underset{\overbrace{{\mbox{\footnotesize $N^4$ (infinite tree) +} \atop \mbox{\footnotesize $\Bigl\{$The 1-handles $\underset{1}{\overset{n}{\sum}} R_j + \underset{1}{\overset{\infty}{\sum}} h_i\Bigl\}$}} \atop \mbox{\footnotesize (see 19)}}}{N^4 (2X_0^2) = N^4 (2\Gamma (\infty))} + \sum \{\mbox{2-handles $D^2 (c)$ corresponding to the}
$$
\vglue -18mm
$$
\qquad \qquad \qquad \qquad \qquad \qquad \qquad \qquad \quad \mbox{ $\{$extended set of $D^2 (c)$'s$\}$ (18.1)}\} + \sum_1^{\bar n} D^2(\Gamma_i),
$$

\vglue 5mm

\noindent then we should find exactly the same geometric intersection matrix
$$
C \cdot h = {\rm id} + {\rm nilpotent} \ \mbox{(of the easy type)},
$$
as in the first line of (19).

\bigskip

\noindent {\bf Remark.} There is, a priori, another possibility to proceed, in order to get our $N^4 (2X_0^2)$, namely to start by choosing some appropriate immersion of $2X_0^2$ into some ambient 4-manifold ($\Delta_1^4$ in (2), for instance) and then to take its $4^{\rm d}$ regular neighbourhood.

\smallskip

But I found the present manner of proceeding much more convenient, in particular more efficient when discontinuous {\ibf changes on local topology} always compatible with (22.A) $+$ (22.B), hence NOT contradicting the {\ibf global} topology, will be necessary. End of Remark.

$$
\includegraphics[width=13cm]{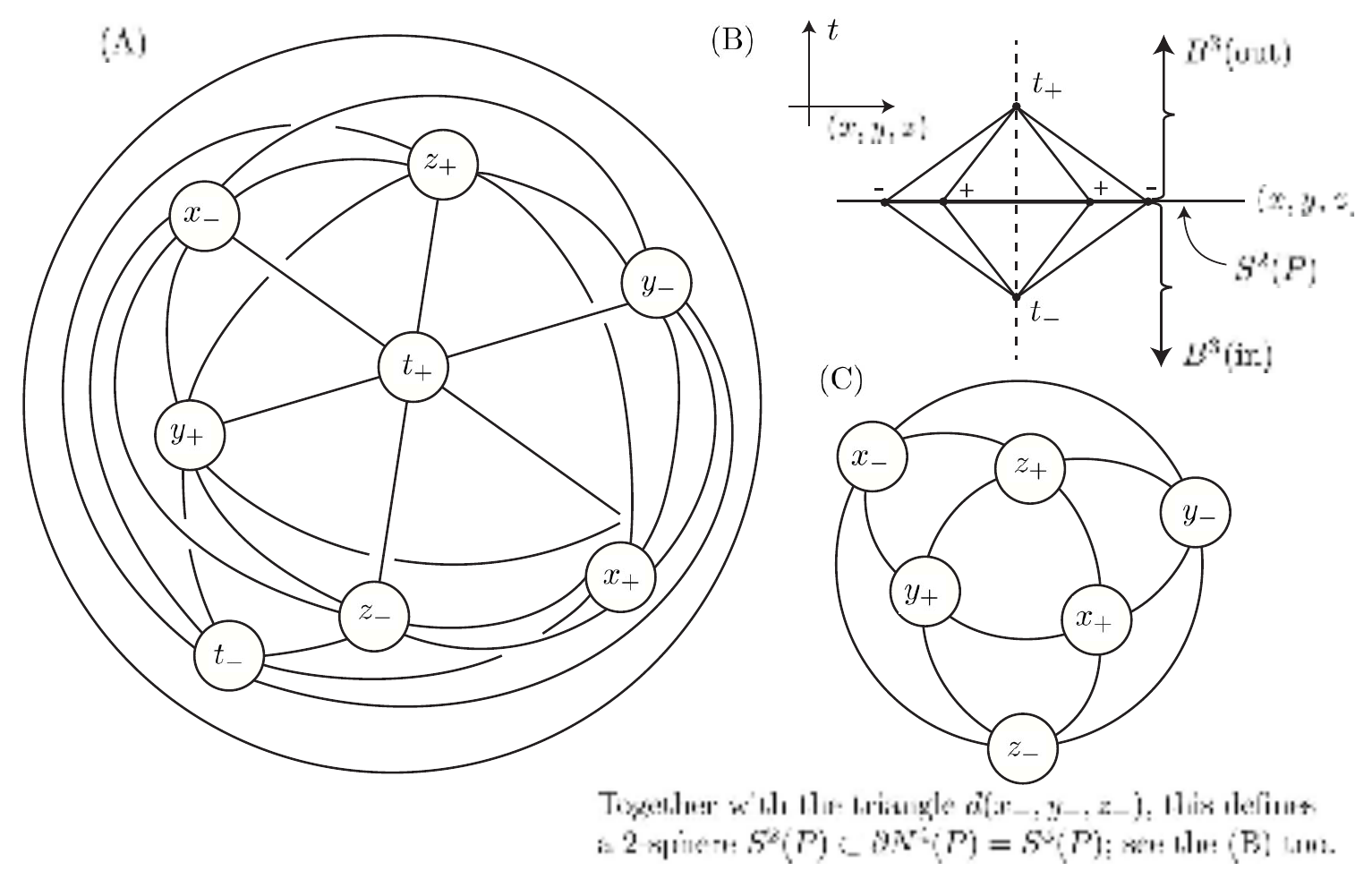}
$$
\label{fig8}
\centerline {\bf Figure 8.}
\begin{quote}
We see here, in (A), the trace $(X_{\rm cubical}^2 \mid P) \cap S^3 (P)$.

EXPLANATIONS: The (A) is gotten by suspending the $S^2(P)$ from (C), from the 
$\xy *[o]=<15pt>\hbox{$+t$}="o"* \frm{o}\endxy$ and $\xy *[o]=<15pt>\hbox{$-t$}="o"* \frm{o}\endxy$, like it is suggested in (B). Actually, what (B) suggests is a SPLITTING
$$
S^3 (P) = \partial N^4 (P) = B^3 ({\rm out}) \underset{\overbrace{\mbox{\footnotesize$S^2(P)$}}}{\cup} B^3 ({\rm in}), \ {\rm with} \ \xy *[o]=<15pt>\hbox{$t_-$}="o"* \frm{o}\endxy \subset B^3 ({\rm in}), \ \xy *[o]=<15pt>\hbox{$t_+$}="o"* \frm{o}\endxy \subset B^3 ({\rm out}),
$$
while in (C) we have redrawn a detail of the main (A), which corresponds exactly to
$$
S^2 (P) \equiv \bigcup \{\mbox{the eight triangles $d^2$ (space, space, space)}\}.
$$

But there is also a second SPLITTING, by $S_{\infty}^2 (P) \equiv S^3 (P) \cap \sum_{\infty}^2$, with $\sum_{\infty}^2$ like in (27), to be discussed later,
$$
S^3 (P) = B^3 (-) \underset{\Sigma_{\infty}^2}{\cup} B^3 (+) ;
$$
what we see in (A) is actually $B^3 ({\rm in}) \cup B^3 (-)$, with the interaction between these two splittings of $S^3 (P)$, by $S^2 (P)$ and by $\sum_{\infty}^2$ (or rather its trace $S_{\infty}^2$ on $S^3 (P)$), is suggested in the Figure~9.bis.
\end{quote}

\bigskip

As a preliminary for the $N^4 (2X^2_0)$ which will have to be constructed from scratch for our $2X_0^2$, hence not coming with any a priori God-given embedding into some given 4-manifold, I will consider now
$$
X_{\rm cubical}^2 \equiv \{\mbox{the 2-skeleton of {\ibf the standard cubical} cell-decomposition of $R^4 = R^3 \times R$} , 
$$
$$
{\rm with} \ R^4 = (x,y,z,t)\},
$$
which certainly comes with the obvious embedding $X^2_{\rm cubical} \subset R^3 \times R$, from which the $N^4 (X^2_{\rm cubical})$ is defined. If $P \in X^2_{\rm cubical}$ is a vertex, we have
$$
N^4 (P) \equiv N^4 (X^2_{\rm cubical}) \mid P = \{{\rm the} \ N^4 ({\rm germ} \ X^2_{\rm cubical} \mid P)\}.
$$
This comes with $N^4(P) = B^4$, $S^3(P) \equiv \partial N^4(P)$.

\smallskip

The $N^3 ((X^2_{\rm cubical} \mid P) \cap S^3 (P)) \subset S^3(P)$ is a collection of small 3-balls $b^3 (x_{\pm}) , \ldots , b^3 (t_{\pm})$, occuring as $\xy *[o]=<15pt>\hbox{$x_{\pm}$}="o"* \frm{o}\endxy , \ldots , \xy *[o]=<15pt>\hbox{$t_{\pm}$}="o"* \frm{o}\endxy$ in Figure~8. Here
$$
S^2(P) \equiv S^3 (P) \cap \{\mbox{the $(x,y,z)$-coordinate hyperplane}\}
$$
induces the splitting $S^3 (P) = B^3 ({\rm out}) \underset{\overbrace{\mbox{\footnotesize$S^2(P)$}}}{\cup} B^3 ({\rm in})$, with $b^3 (t_+) \subset B^3 ({\rm out})$, $b^3 (t_-) \subset B^3 ({\rm in})$. Figure 8-(A) is a rigorously {\ibf correct} representation of $(S^3 (P) , (X^2_{\rm cubical} \mid P) \pitchfork S^3 (P))$; one can think of the 3-ball inside which the drawing lives as being $B^3 ({\rm in}) \cup B^3 (-)$, see here Figure 9.bis too.

\smallskip

Of course, also, nice small isotopies are allowed. For instance, the $b^3 (t_-)$ which, for reasons of graphical commodity has been pulled a bit to the side, could come just under $b^3 (t_+)$ so that these should be a line: $\{$observer's eye$\} -\!\!\!-\!\!\!-\!\!\!- \ b^3 (t_+)  -\!\!\!-\!\!\!-\!\!\!- \ b^3 (t_-)$ (see (B) too). It is this move which is used when in Figure 9 we change $(B_+)$ into $(B_-)$. The $(X^2_{\rm cubical} \mid P)$ occuring above, is the germ at the origin $0 \in R^4$ of the six coordinate planes. Each triangle $d^2$ (space-time, space-time, space-time) stands for a corner of $X^2_{\rm cubical} \mid P$, and the $d^2$'s are disjoined except for their common edges or vertices. In our Figure 8-(A), the {\ibf arcs} joining two $b^3$'s are pieces of the curves along which 2-handles get attached.

\bigskip

\noindent {\bf The reconstruction of} $N^4 (X^2_{\rm cubical})$, out of the correct Figures 8-(A), for the various vertices $P$. This is done in two steps.

\medskip

I) If $P_1 , P_2$ are two adjacent vertices, then in the 1-skeleton $X^2_{\rm cubical}$ some $b^3 (u_{\pm}) \subset S^3 (P_1)$ communicates with $b^3 (u_{\mp}) \subset S^3 (P_2)$, and we will join the two via a 1-handle. This way we have reconstructed the $N^4 (X^1_{\rm cubical})$, which should be, of course orientable.

\medskip

II) From the arcs $\subset \underset{P}{\bigcup} S^3 (P)$, joined along the lateral surfaces of the 1-handles above, in a canonical manner, we get a
$$
\{{\rm link}\} \subset \partial N^4 (X^1_{\rm cubical})
$$
to which, with appropriate framings we add 2-handles. This {\ibf is} the reconstruction of $N^4 (X^2_{\rm cubical}) \subset R^3 \times R = R^4$. But the $X^2_{\rm cubical}$ is only an intermediary tool, since what really interests us is the real-life GPS structure, and that will be used for $2X_0^2$. But then, the GPS structure is compatible with the cubical one: One moves from one to the other (cubical $\Rightarrow$ GPS) via some deletions and simple isotopies. And so, as long as $N^4 (\Delta^2) = N^4 (2X_0^2) \mid \Delta^2$ stays compatible with $N^4 (X^2_{\rm cubical})$, we are sure that the (22.A) is satisfied. Otherwise, there are no restriction, outside of $\Delta^2 = \Delta^2 \times (\xi_0 = -1)$, as far as the framings are concerned.

\bigskip

$$
\includegraphics[width=125mm]{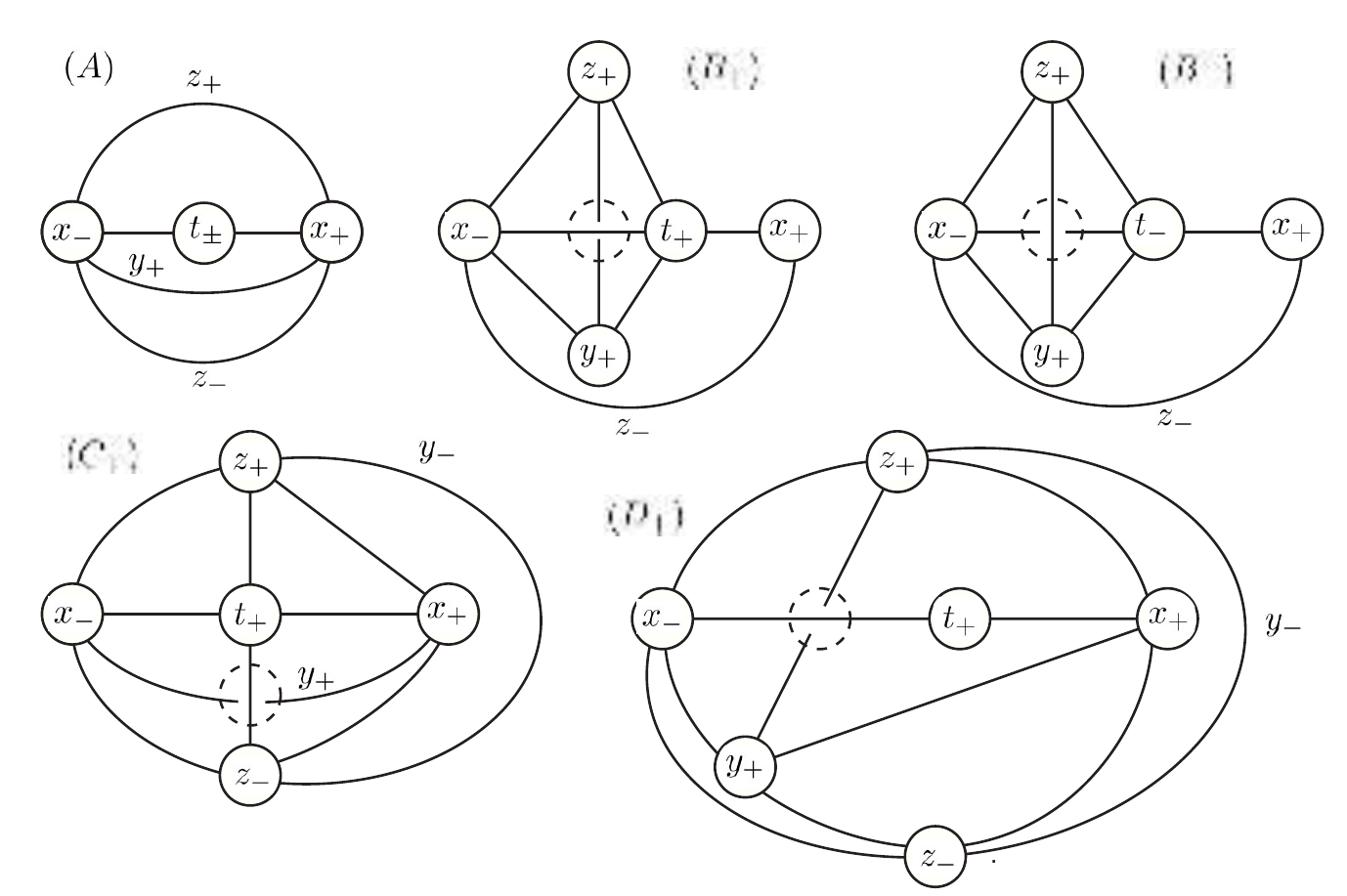}
$$
\label{fig9}
\centerline {\bf Figure 9.}
\begin{quote}
The local models for the real-life GPS $N^4 (P) \equiv N^4 (2X_0^2 ({\rm GPS})) \mid P$, with a $N^4 (2X_0^2 ({\rm GPS}))$ which is still to be defined but for which the present figures settles the local models. By appropriate permutation of space-time letters, all the local GPS models can be deduced from the ones drawn here. The present figure should be compared to Figure 8-(A), but for typographical commodity we have omitted the outer surrounding circle. The present drawings only refer to the $(x,y,z,t)$ part of the $N^4 (2X_0^2)({\rm GPS})$, and they are certainly compatible with Figure 8-(A). This figure will have to be embellished with the contributions coming from the axes $\xi_0$ and $\zeta$. The typical embellished figure is 24.

The present $(A)$, $(B_{\pm})$, $(C_+)$, $(D)$ correspond, respectively, to the (A), (B), (C), (D) in Figure~4, with the $t_{\pm}$ added. The lines $[z_+, y_+]$, $[y_+ , z_-]$, $[z_- , z_+]$ are ORANGE, the others are BLACK.
\end{quote}

\bigskip

\noindent {\bf Additional explanations for the Figure 9.} Each of the five little drawings in the figure which is displayed here, once filled in with the appropriate 2-cells, with disjoined interiors, bounded by the closed curves in the figure, can be recognized as being a detail which appears in the Figure 8-(A) too. So, Figure 9 is compatible with 8-(A).

\smallskip

All the configurations drawn here embed, isometrically in Figure 8-(A) respecting the presently drawn $b^3$ (space-time)'s. When the corresponding line in Figure 8-(A) has a $b^3 \, \mbox{(space-time)}$ not part of the GPS local model, then one wrote the corresponding letter next to the corresponding line.

\smallskip

We should also have Figures 9-$(C_-)$, $(D_-)$ gotten from the drawn $(C_+)$, $(D_+)$, by letting the lines with a dotted loop over them cross, in the following pattern

$$
\includegraphics[width=8cm]{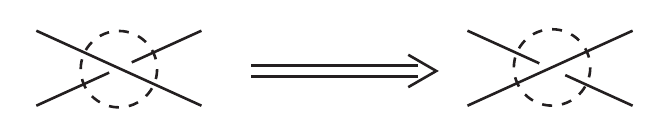}
$$

Figures 9 are used for reconstructing the GPS $N^4 (X_0^2 \times r)$, or at least the $N^4 (\Delta^2 \times (\xi_0 = -1))$, on the same lines as the reconstruction of the $N^4 (X_{\rm cubical}^2)$. And, as already said, a long as our local models are compatible with Figure 8-(A), {\ibf and} we respect the framings of $N^4 (\Delta^2)$, (see the context of the law II for reconstructing $N^4 (X_{\rm cubical}^2)$), our condition (22.A) is automatically satisfied.

\smallskip

In the context of Figure 9 there is also a $(C_-)$, which we have not drawn, and this is gotten from $(C_+)$ by changing $t_+ \Rightarrow t_-$ and then using the obvious crossing operation at the spot surrounded by the dotted circle. Moving now to $(D_+)$, which is related to Figure~4-(D), I stress that we are here in the context of formula (6), with the time directions past/future obeying the scheme presented in Figure~9.1. In this context, $(D_+)$ corresponds to $t = t_{2i-1}$ and if we move to $(D_+^*)$, we get to $t = t_{2i}$, with a change BLACK $\Leftrightarrow$ ORANGE. 

\bigskip

$$
\includegraphics[width=13cm]{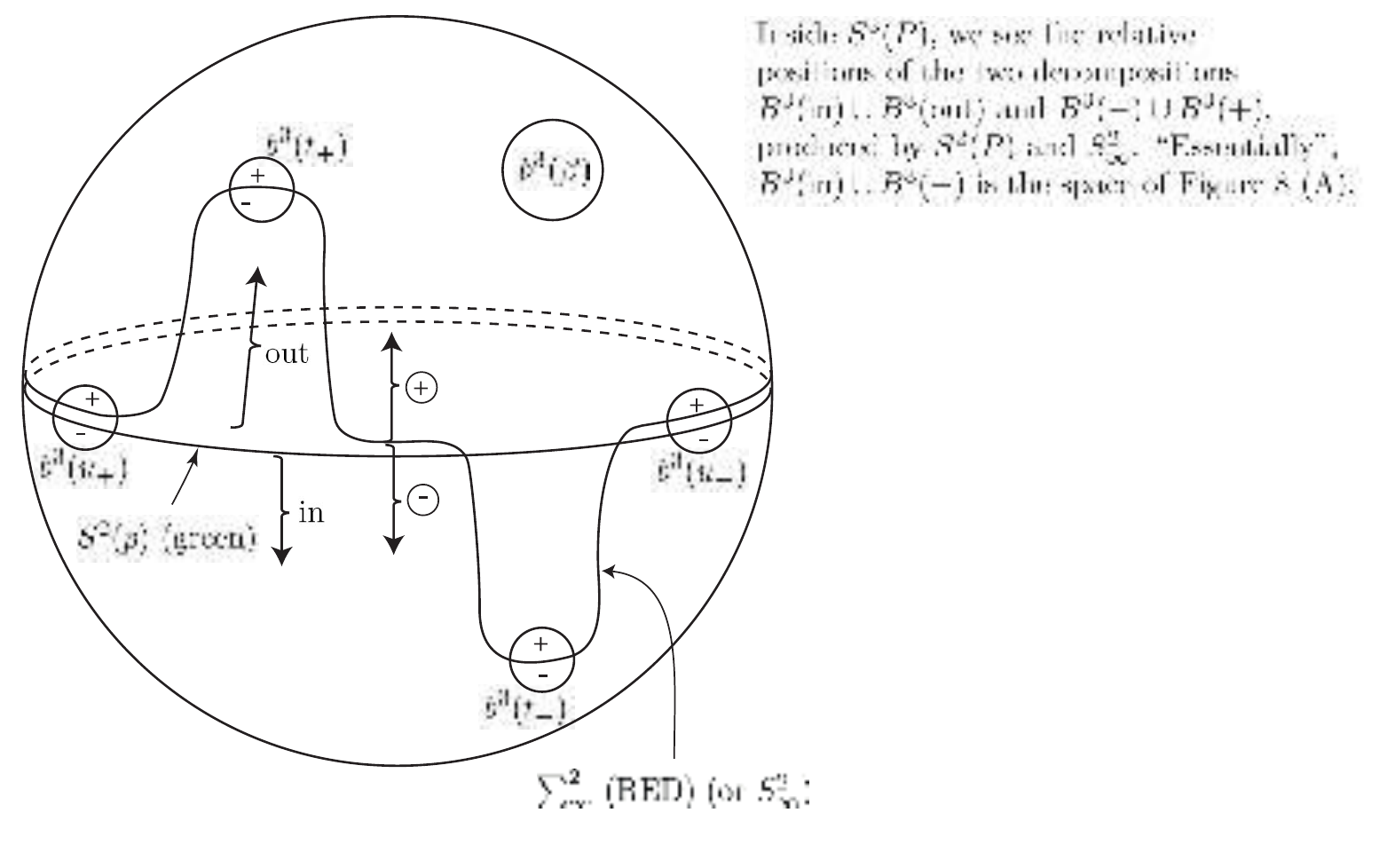}
$$
\label{fig9bis}
\centerline {\bf Figure 9.bis.}
\begin{quote}
With one dimensionless we suggest here the two splittings $S^3(P) = B^3 ({\rm in}) \underset{\overbrace{\mbox{\footnotesize$S^2(p)$}}}{\cup} B^3 ({\rm up})$, $S^3 (P) = B^3 (-) \underset{\sum_{\infty}^2}{\cup} B^3 (+)$, together.
\end{quote}

\bigskip

Then there are also $(D_{\pm}^{(*)})$ figures, gotten like $(C_+) \Rightarrow (C_-)$. The present 9-$(D_+)$ corresponds to 4-(D), with the same colours and with a $t_+$ added. So we go now from $3^{\rm d}$ to $4^{\rm d}$. With this, also, the ``$u$'' in Figure 4-(D) becomes now a vertex, bringing an ORANGE CONTRIBUTION to the 1-skeleton. End of explanations.

\bigskip

We move now to the still to be defined $N^4 (2X_0^2)_{\rm GPS}$, actually an object to be defined ex-nihilo, with the only restriction that the conditions (22.A), (22.B) should be satisfied. We start by looking at
$$
X_0^2 [{\rm new}] = X_0^2 \times r \supset X_0^2 [{\rm new}] ({\rm truncated}) \equiv X_0^2 [{\rm new}] - ((\Gamma (1) \times [0 > \xi_0 \geq -1]) \cup \Delta^2 \times (\xi_0 = - 1)).
$$ 
As long as we forget the coordinate $t$, the local models for $X_0^2 [{\rm new}] ({\rm truncated})$ are the (A), (B), (C), (D) from Figure 4. Realistically, a $t_+$ OR a $t_-$ will have to be added, never both, according to the pattern from Figure~9.1. With this addition, the Figures 4 turn into 9.

\smallskip

The $X_{\rm cubical}^2 \mid P$ is the standard form
$$
\{ [x,y] \cup [y,z] \cup [x,z] \cup [x,t] \cup [y,t] \cup [z,t]\} \mid P.
$$

\smallskip

$$
\includegraphics[width=12cm]{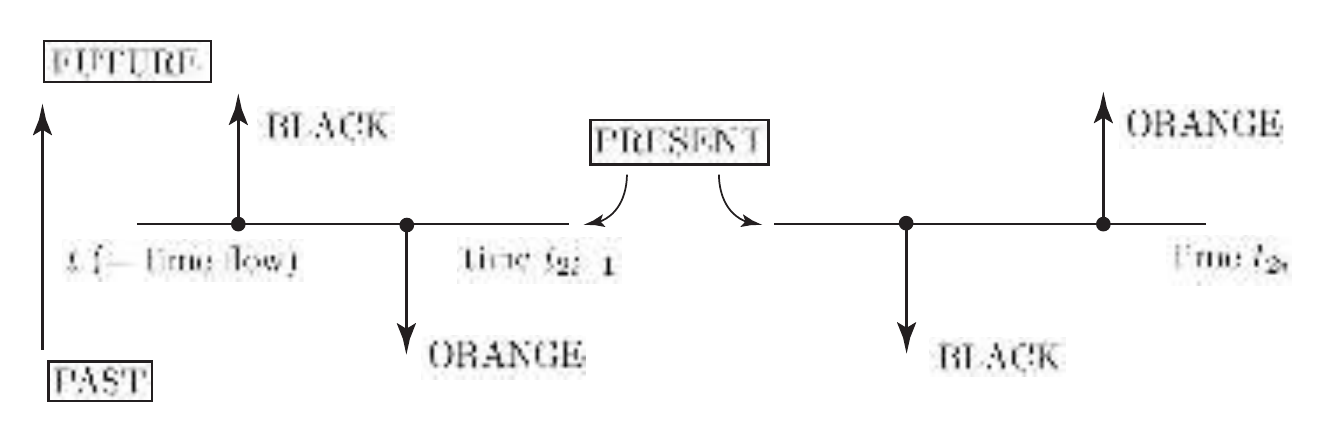}
$$
\label{fig9.1}
\centerline {\bf Figure 9.1.}

\smallskip

\centerline{\ \ Illustration for (6).}

\bigskip

Here is a more or less obvious fact: there is a simple transformation
$$
X_{\rm cubical}^2 \mid P \Longrightarrow \{ X_0^2 [{\rm new}] ({\rm truncated}) \mid P , \ \mbox{with GPS structure}\}, \leqno (23)
$$
proceeding by obvious minor deletions and embedded $2^{\rm d}$ slidings of the dilatating type.

\smallskip

The local models for $X_0^2 [{\rm new}] ({\rm truncated})$ are our (A), (B), (C), (D) above (with $t_{\pm}$ added). What (23) implies is that, of we take for $P \in X_0^2 [{\rm new}] ({\rm truncated})$, {\ibf by decree}, the local models
$$
(N^4 (P)) = B^4 , \partial N^4 (P) \cap (X_0^2 [{\rm new}] ({\rm truncated})  \mid P)) \leqno (24)
$$
from Figure 9, we are {\ibf compatible} with the local model from Figure 8. This is the main fact behind the (22.A).

\smallskip

When moving from $X_0^2 [{\rm new}] ({\rm truncated})$ to the full $2X_0^2$ our local models are incomplete since they lack the following two basic ingredients:

\bigskip

\noindent (25) \quad i) The contribution $[0 \geq \xi_0 \geq -1]$ joining two vertices $P \times (\xi_0 = 0)$ and $P \times (\xi_0 = -1)$, when $P \in \Delta^2$; ii) Then, for {\ibf any} $P \in X_0^2 \times r \approx X_0^2 [{\rm new}]$ the $[r \leq \zeta \leq b]$ joining $P \times r$ to $P \times b$ is also to be added.

\bigskip

At this point an important prentice has to be opened, and although our $N^4 (2X_0^2)$ has not yet here fully constructed, we can still talk about it. Our $\beta \in [r,b]$ will be taken sufficiently close to $r$, so that for each edge $e \subset \Gamma (\infty) = \Gamma (\infty) \times r$ we should have $e \times [r,\beta] \subset 2X_0^2$, and

\bigskip

\noindent (26) \quad $\Gamma (\infty) \times \frac{\beta}{2} \subset 2X_0^2$ splits $2X_0^2$ into a ``$-$'' part, towards $r$ and a ``$+$''  part towards $b$. For $e = e(b)$ (see Figure 7), the situation of $e \times \frac{\beta}2$ and $e \times \beta$ is displayed in Figure 10.

\bigskip

$$
\includegraphics[width=6cm]{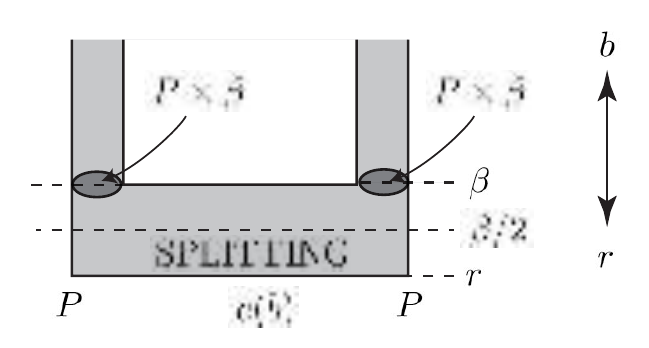}
$$
\label{fig10}
\centerline {\bf Figure 10.}

\smallskip

\centerline{
We have shaded here the collar which is alive in the Figures 7-(II, III).
}

\bigskip

We introduce now 
$$
\Gamma_1 (\infty) \equiv \Gamma (\infty) \times r \cup \Gamma_0 (\infty) \times [r,\beta] \leqno (26.1)
$$
where $\Gamma_0 (\infty) \subset \Gamma_1 (\infty)$ is the set of vertices), and $N^4 (\Gamma_1 (\infty)) \equiv \{$the (abstract) $4^{\rm d}$ regular neighbourhood of $\Gamma_1 (\infty) \cong \Gamma_1 (\infty) \times r \subset 2X_0^2 \}$. I will give now a useful 

\bigskip

\noindent {\bf Model for $N^4 (2\Gamma (\infty))$.} We will concentrate now on $\Gamma_1 (\infty)$ (27), the extension to the full $2\Gamma (\infty)$, afterwards, should be transparent. Our aim is to present now, in a manner which is consistent with everything said so far, the $N^4 (\Gamma_1 (\infty))$ as a product
$$
N^4 (\Gamma_1 (\infty)) = N^3 (\Gamma_1 (\infty)) \times [0,1], \leqno (27)
$$
and we shall think of the factor $[0,1]$ above as being $[r,\beta] \subset [r,b]$. Here is a very easy general fact. Consider any graph $\Gamma$ (like $\Gamma = \Gamma_1 (\infty)$). Then we have the following

\bigskip

\noindent (27.0) \quad For every $N^n (\Gamma)$, $n \geq 4$, there is a regular neighbourhood $N^3 (\Gamma)$, such that $N^n (\Gamma) = N^3 (\Gamma) \times B^{n-3}$.

\bigskip

Hence, there is no problem with the bare (27), but we have to chose the correct $N^3 (\Gamma_1 (\infty))$ and check the compatibilities, in particular with the basic (22.A).

\bigskip

\noindent {\bf Construction of the correct $N^3 (\Gamma_1 (\infty))$.} With the not yet defined $N^4 (\Gamma_1 (\infty))$ the {\ibf embellished Figure~9} (with the contributions of $\xi_0 , \zeta$ added), i.e. Figure 24, we have local models of $(X^2 [{\rm NEW}] \pitchfork \partial N^4 (P))$ (embellished), for the vertices $P \in \Gamma (\infty) \subset \Gamma_1 (\infty)$. The embellishment means that we are actually talking about $(2X_0^2 \mid P) \pitchfork \partial N^4 (P)$, for $P \in \Gamma (\infty) \times r \subset 2\Gamma (\infty)$. We can think of the plane of the figures of type 9 as being a $S_0^2(P)$ on which all the $b^3 (u_{\pm})$ others than $b^3 (t_{\pm})$ ride already, as in Figure 9.bis and on which all the rest of Figure 9 is projected like in a link diagram, as an immersion
$$
(X^2 [{\rm NEW}] \pitchfork N^4 (P))_{\rm embellished} \xrightarrow{ \ \psi \ } S_0^2 (P). \leqno (27.1)
$$
This comes with sites $u \in \{\pm \, x , \pm \, y , \pm \, z , \pm \, t , \pm \, \xi_0 , \beta \in [r,b]\}$ and little discs $b^2 (u)$. Our {\ibf link diagram} contains curves joining the $b^2 (u)$'s, and at the double points of $\psi$ ($=$ the CROSSINGS), UP/DOWN's are specified.

\smallskip

We think of $S_0^2 (P)$ which, until further notice, is simply the plane of the figures of type 9, as being the $S_0^2 (P) = \partial B_0^3 (P)$, with the $B_0^3 (P)$ living on the other side of the Figures 9, 24, with respect to the observer.

\smallskip

In order to construct our $N^3 (\Gamma_1 (\infty))$, we start with
$$
\sum_{\overbrace{\mbox{\footnotesize$P \in \{$vertices of $\Gamma (\infty)\}$}}} B_0^3 (P) , \ \mbox{coming with $b^2 (u)$'s $\subset$ $S_0^2 (P)$,}
$$
and notice that, if for two adjacent vertices $P_1 , P_2$ we have $b^2 (u_{\pm}) \subset S_0^2 (P_1)$, then we certainly find $b^2 (u_{\mp}) \subset S_0^2 (P_2)$. If we join each $B_0^3 (P_1)$ to $B_0^3 (P_2)$, in an orientable manner, along a $b^2 \times [0,1]$ with $b^2 \times \{ 0 \} = b^2 (u_{\pm}) \subset S_0^2 (P_1)$, $b^2 \times \{1\} = b^2 (u_{\mp}) \subset S_0^2 (P_2)$, then we get a model for $N^3 (\Gamma (\infty))$, endowed with a $b^2 (\beta) \subset \partial N^3 (\Gamma (\infty))$ for every $B_0^3 (\beta) \subset N^3 (\Gamma (\infty))$. This structure, where we will think of the $b^2 (\beta)$'s as sticking out, a bit, like in the Figure 11 below, is our $N^3 (\Gamma_1 (\infty))$, the construction of which is by now finished, and, by construction, it is compatible with the figures of type 9. 

\smallskip

We come now with our additional factor $[r,\beta] \subset [r,b]$ and it is the $[r,\beta]$ which will be our factor $[0,1]$ in the context of (27).

\smallskip

As already said, when it  comes to the compatibility with (22.A), the bare (27), which is by now in place, is irrelevant. What counts is the link diagrams (27.1) for various $P \in \Gamma (\infty) \times r$, and their lifts to $4^{\rm d}$.

\smallskip

It should be understood here that $\beta \in (r,b)$ is very close to $r$, so that all the thin band $\Gamma (\infty) \times [r,\beta] \subset 2X_0^2$, is not affected by the deletion of the non-shaded areas in Figures 7-(II, III). A more detailed version of Figure~7 is presented in Figure 45, where the splitting line occurs as a dotted line. Figure 7.bis is not concerned here; it is split by the two sites marked $\beta$ into $\pm$ halves, $-$ on the $r$-side and $+$ on the $b$-side. Let us also notice here the obvious embedding
$$
\Gamma_1 (\infty) \subset \Gamma (\infty) \times [r,\beta] \leqno (27.2)
$$
and one should notice that the pairs (germ of $2X_0^2$ at $\Gamma (\infty) , \Gamma (\infty)$) and (germ of $\Gamma (\infty) \times [r,\beta]$ at $\Gamma (\infty) \times \frac\beta2$, $\Gamma (\infty) \times \frac\beta2 \approx \Gamma (\infty)$) are completely different from each other. In terms of Figure 24, the first germ contains the complete information, while the second one only knows about those things which are adjacent to $\beta$. Also, while the germ of $2X_0^2$ lives in $N^4 (\Gamma (\infty))$ and cannot be pushed into $\partial N^4 (\Gamma (\infty))$, the germ of $\Gamma (\infty) \times [r,\beta]$ mentioned above, can. This last fact will be used in this paper.

\smallskip

We finally can introduce the
$$
N^4 (\Gamma_1 (\infty)) = \{ N^3 (\Gamma (\infty)) \times [r,\beta] , \ \mbox{with a} \ b^2 (\beta) \times [r,\beta] \leqno (27.3)
$$
$$
\mbox{added to $\partial N^4 (\Gamma (\infty))$ at each $P\}$, see Figure 11.} 
$$
We have here the commutative diagram
$$
\xymatrix{
N^4 (\Gamma_1 (\infty)) \ar[d] &\underset{\rm DIFF}{=} &N^3(\Gamma (\infty)) \times [r,\beta] \ar[d]^p \\
\Gamma_1 (\infty) &\subset &\Gamma(\infty) \times [r,\beta].
} \leqno (27.4)
$$
At the level of $N^4 (\Gamma_1 (\infty))$, the $b^2 (u)$'s of $\partial B_0^3 (P)$ become $B^3 (u) \subset \partial N^4 (P)$. Concerning the $S_0^2 (P)$ from (27.1), when we go to $4^{\rm d}$, it gives rise to $S_0^2 (P) \times \frac\beta2 = \{$the $S_{\infty}^2 (P)$ below$\}$.

\smallskip

The compatibility with (22.A) is taken care of by the explicit liftings of the link diagrams, and what we can perceive on our model (27.3) is a PROPER and proper embedding $N^3 (\Gamma(\infty)) \times \frac\beta2 \subset N^4 (\Gamma_1 (\infty))$. This {\ibf defines} a SPLITTING SURFACE
$$
\Sigma_{\infty}^2 \equiv \partial (N^3 (\Gamma_1 (\infty)) \times \frac\beta2 \subset \partial N^4 (\Gamma_1 (\infty)),
\leqno (27.5)
$$
coming with the SPLITTING
$$
\partial N^4 (\Gamma_1 (\infty)) = \partial N^4_- (\Gamma_1 (\infty)) \underset{\Sigma_{\infty}^2}{\cup} \partial N_+^4 (\Gamma_1 (\infty)) . \leqno (27.6)
$$
Notice the connection with (27.4), namely the equation
$$
p^{-1} \left(\Gamma (\infty) \times \frac\beta2 \right) \cap \partial N^4 (\Gamma_1 (\infty)) = \Sigma_{\infty}^2 .
$$

We will also use here the notations
$$
Z^3 (-) \equiv \partial N_-^4 (\Gamma_1 (\infty)) , \quad Z^3 (+) \equiv \partial N_+^4 (\Gamma_1 (\infty)).
$$
The $\pm$ convention will be fixed by setting $\Sigma \, b^3 (\beta) \subset \partial N_+^4 (\Gamma_1 (\infty))$, with the rest of the $b^3 (\pm \, u)$'s (see Figure 24) being each split by $\Sigma_{\infty}^2$ into two halves, like Figure 9.bis may suggest. All the corresponding interconnections, inside the $\partial N^4 (P)$'s, inbetween the $b^3 (\pm \, u)$'s live, for the time being at least, inside the $\partial N_-^4 (\Gamma_1 (\infty))$. Later, some of them will have to be pushed into $\partial N_+^4 (\Gamma_1 (\infty))$, in a very controlled manner. Our SPLITTING
$$
\partial N^4 (\Gamma_1 (\infty)) = Z^3 (-) \underset{\Sigma_{\infty}^2}{\cup} Z^3 (+)
$$
is suggested very schematically in the Figure 11.

\bigskip

$$
\includegraphics[width=7cm]{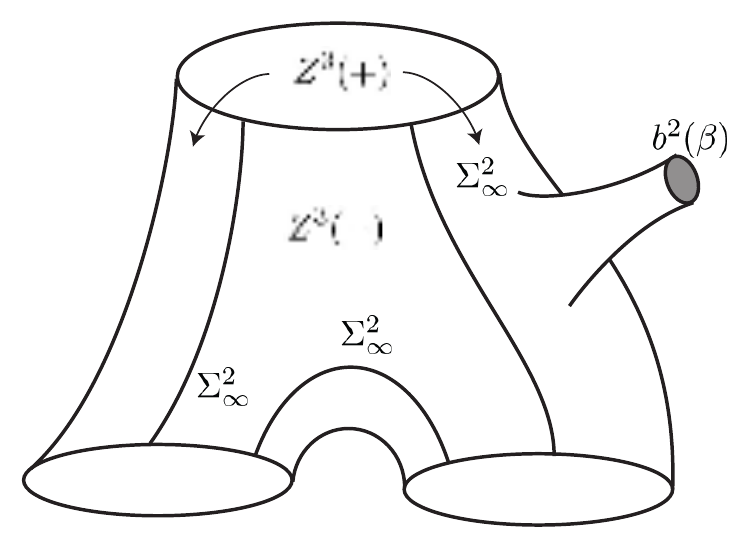}
$$
\label{fig11}
\centerline {\bf Figure 11.}
\begin{quote}
The splitting by $\Sigma_{\infty}^2$, see (27.6) is suggested here, with one dimension less. Via the sites $b^3 (\beta)$ our $N^4 (\Gamma_1 (\infty))$ communicates with the rest of $N^4 (2\Gamma (\infty))$.
\end{quote}

\bigskip

For each vertex $P \in X_0^2 \times r$, the $\Sigma_{\infty}^2$ induces a splitting (see Figure 9.bis)
$$
S^3 (P) = B^3 (-) \underset{S_{\infty}^2}{\cup} B^3 (+) , \leqno (28)
$$
where we use the notation $S_{\infty}^2 \cong \Sigma_{\infty}^2 \cap \partial N^4 (P)$. This splits each $b^3 (\pm \, r)$ in two halves. We mean here all the $b^3 (\pm \, u)$'s from the embellished Figure 24. Next, if $P_1 , P_2$ are two adjacent vertices of $\Gamma_1 (\infty)$ and $b^3 (\pm \, u) \subset \partial N^4 (P_1)$, then $b^2 (\mp \, u) \subset \partial N^4 (P_2)$ and $b^2 (u) \times [-1,+1]$ joins the two $N^4 (P_i^3)$'s. The splitting $\Sigma_{\infty}^2$ continues along the $\partial b^3 (u) \times [-1,+1]$. This is what Figure 11 is supposed to suggest. The splitting (28) and the splitting from the Figure 8, the $\partial N^3 (P) = B^3 ({\rm in}) \underset{S^2(P)}{\cup} B^3 ({\rm out})$ are articulated with each other like in the Figure 9.bis, which shows exactly how the splittings of $S^3 (P) = \partial N^4 (P)$ by $S^2(P)$ and by $S_{\infty}^2$ ($=$ trace of $\Sigma_{\infty}^2$ on $\partial N^4 (P)$) interact with each other. In view of Figure 9.bis, we may also happily assume that the $S_0^2 (P) = \partial B_0^3 (P)$, which occurs in (27.1), is the $S_{\infty}^2$.

\smallskip

When we move from $\Gamma_1 (\infty)$ to the larger $2 \Gamma (\infty)$, we define the following big SPLITTING which naturally extends (27.5)
$$
\partial N^4 (2\Gamma (\infty)) = \partial N_-^4 (\Gamma_1 (\infty)) \underset{\overbrace{\mbox{\footnotesize $\Sigma_{\infty}^2$}}}{\cup} \partial N_+^4 (2\Gamma (\infty)), \leqno (28.1)
$$
where $\partial N_+^4 (2\Gamma (\infty)) \equiv \partial N^4 (2\Gamma (\infty)) - \partial N_-^4 (\Gamma_1 (\infty))$. You may happily use the terminology $\partial N_-^4 (2\Gamma (\infty)) \equiv \partial N_+^4 (2\Gamma_1 (\infty))$.

\smallskip

We will denote by $\{{\rm link}\}[{\rm NEW}]$ the obvious extension of the $\{{\rm link}\}$ from (9) to $X^2 [{\rm NEW}]$ with the $D^2 (\Gamma_i) \times (\xi_0 = 0)$ becoming now, RED-wise, $D^2 (\gamma_k^0)$'s and continuing to stay, BLUE-wise, according to the case $D^2 (\eta)$ or $D^2 (\gamma^1)$. Also, of course, the contribution of $(\Gamma (1) \times [0 \geq \xi_0 \geq -1]) \cup (\Delta^2 \times (\xi_0 =-1))$ is to be included now too. And, at this point, remember that the $\Gamma (1) \times [0 \geq \xi_0 \geq -1]$ is made out of $D^2 (C)$'s RED-wise and of $D^2 (\eta)$'s BLUE-wise, while the $D^2 (\Gamma_i) \times (\xi_0 = -1)$ (our $D^2 (\Gamma_i)$'s, now) are $D^2 (\gamma^1)$'s BLUE-wise. With this, we get, for the time being, the inclusion
$$
\{{\rm link}\}[{\rm new}] \subset \partial N_-^4 (\Gamma_1 (\infty)). \leqno (28.2)
$$

\bigskip

\noindent {\bf Reconstruction of $N^4 (X_0^2 \times r) \supset N^4 (\Gamma_1 (\infty))$.} We will imitate the same steps as in the RECONSTRUCTION of $N^4 (X_{\rm cubical}^2)$ above, but now we are GPS, not cubical.

\medskip

\noindent STEP (0). We need, for each $P \in \{$vertex of $X_0^2 \times r\} \in \Gamma (\infty) \times r$, a pair $(N^4 (P) , \{ \partial N^4 (P) \cap [((X_0^2 \times r) \cup (\Gamma (\infty) \times [r,\beta]))\} \mid P])$. When the axes $\xi_0 , \zeta$ are disregarded, then the corresponding restrictions will be exactly  the ones from Figure 9. If we keep this condition strictly for the $P \in \Delta^2 \times (\xi_0 = -1)$, then (22.A) is {\ibf automatically satisfied}, provided we proceed correctly along the edges $[P_1 , P_2] \subset \Delta^2 \times (\xi_0 = -1)$ too. To the $N^4 (P)$'s we add now, for all $P$'s, the contribution $b^3 (\beta) \subset B^3 (+) \subset \partial N^4 (P)$, and for any $P \in \Delta^2 \times (\xi_0 = 0)$, respectively $P \in \Delta^2 \times (\xi_0 = -1)$, the contribution $b^3 (- \, \xi_0)$, respectively $b^3 (+ \, \xi_0)$, with $b^3 (\pm \, \xi_0) \subset B^3 (-)$. We are now at the level of the extended figures \`a la 9 and, in these, the $b^3 (\pm \, \xi_0)$'s are riding on top of the $\Sigma_{\infty}^2$, like the $b^3 (\pm \, u)$, $u \in \{ x,y,z,t \}$. Of course, in our extended figures we add lines joining $b^3 (\beta)$ to $b^3 (\pm \, \xi_0)$ and $b^3 (\beta \ {\rm or} \ \pm \, \xi_0)$ to all the $b^3 \, \mbox{(space-time)}$. Any crossing of lines $\includegraphics[width=3cm]{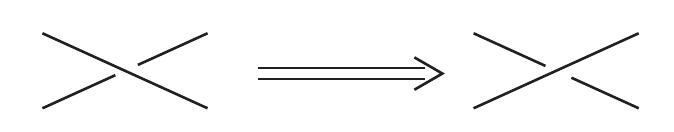}$ inside $\partial N^4 (P)$ is allowed EXCEPT for the crossing of lines $b^3 (\mbox{space-time}) -\!\!\!-\!\!\!-\!\!\!-\!\!\!- b^3 (\mbox{space-time})$, belonging to $\Delta^2 \times (\xi_0 = -1)$. This way we stay compatible with (22.A), (22.B).

\smallskip

[Repeating myself once, more, compatibility with (22.A) also requires paying attention to what we do along the edges, in particular paying attention to the framings. That should be obvious.] Concerning the allowed crossings notice that they do not change the geometric intersection matrices $C \cdot h$, $\eta \cdot B$. They also leave the topology of $N^4 (\Delta^2 \times (\xi_0 = -1))$ intact.

\smallskip

This ends STEP (0) and next come the STEPS I, II, which are exactly like for the reconstruction of $N^4 (X_{\rm cubical}^2)$. This can continue now to a whole reconstruction of $N^4 (2X_0^2)$, without caring any longer about (22.A), on the $X_b^2$-side. Our construction of $N^2(2X_0^2)$ automatically comes with $C \cdot h = {\rm id} + {\rm nil}$, implying (22.B).

\smallskip

Let us consider now $N^4 (\Gamma_1 (\infty)) \subset N^4 (2\Gamma (\infty))$, endowed with the SPLITTING (28.1) and with the $\{$link$\}$ from (18.1). But we will restrict our attention to $N^4 (\Gamma_1 (\infty))$, with the splitting restricted to it and with the ${\rm link} \, [{\rm new}] \subset \partial N_1^4 (\Gamma_1 (\infty))$ from (28.2). For each $P \in X_0^2 \times r$, the $\Sigma_{\infty}^2$ from (27.5) splits $S^3 (P) = \partial N^4 (P)$, via $S_{\infty}^2 = \Sigma_{\infty}^2 \cap S^3 (P)$, into the $B^3 (-)$ and $B^3 (+)$ which we have already met. And, as already mentioned earlier $S^2 (\infty) = S_0^2 (P) \times \frac\beta2$, with the space of the Figures 9, 24 being $B^3 (-) \cup B^3 ({\rm in})$ (see Figure 9.bis).

\smallskip

With this, from now on our {\ibf link diagram} (27.1) is actually the restriction of a larger link diagram
$$
\{{\rm link}\} [{\rm NEW}] \xrightarrow{ \ \ \Lambda \ \ } \Sigma_{\infty}^2 , \leqno (28.3)
$$
a generic immersion with prescribed UP/DOWN labels at each crossing. The $\{{\rm link}\} [{\rm NEW}]$ is generated, then, from this LINK DIAGRAM and, next, the canonical framings come from the neighbourhood of (28.3) in $\Sigma_{\infty}^2$.

\smallskip

We move next to the topic of EXTENDED COCORES. Here what I CLAIM is the following:

\bigskip

\noindent (28.4) \quad Inside $2X_0^2$, for any point $x \in 2X_0^2 - \underset{\Gamma_i \subset \Delta^2}{\sum} {\rm int} (D^2 (\Gamma_i) \times (\xi_0 = -1))$ there is an
$$
\{\mbox{extended cocore $x$}\} \subset 2X_0^2 ,
$$
which is a tree based at $x$, PROPERLY embedded inside $2X_0^2$, cutting transversally the $2\Gamma (\infty) \subset 2X_0^2$ and which has the feature that  at any point $\{$extended cocore $x\} - \{x\}$, the $2X_0^2$ is locally split in two. [This means that if $y \in \{\mbox{extended cocore $x$}\} \pitchfork \Gamma (2\infty)$ ($\not\ni x$), then for each $2^{\rm d}$ branch of $2X_0^2$ at $e \equiv \{$edge of $2\Gamma (\infty)$ containing $y\}$, there is a $1^{\rm d}$ branch of the $\{$extended cocore $x\}$ at $y$. We will also say that such a tree like our extended cocore, which locally splits $2X_0^2$ into two sides, is a {\ibf complete tree}.]

\smallskip

Now, we start by considering any $b_i \in B \cap (\Gamma (1) \times (\xi_0 = -1))$ and show how the all-important $\{$extended cocore $b_i\} \subset 2X_0^2$ is to be constructed. From its root at $b_i = b_i \times (\xi_0 = -1)$, the extended cocore anyway starts with the edge $b_i \times [-1 \leq \xi_0 \leq 0]$. When we get at $b_i \times (\xi_0 = 0) \in X_0^2 [{\rm new}] \mid ({\rm truncated}) \equiv X_0^2 [{\rm new}] - (\Gamma (1) \times (0 > \xi_0 \geq -1)) \cup \Delta^2 \times (\xi_0 = -1)$, we have two possibilities.

\medskip

i) (The trivial case) The $b_i \times (\xi_0 = 0)$ is an ``isolated'' point of $X_0^2 [{\rm new}] \mid ({\rm truncated})$, i.e. it is not touched by the 2-skeleton. [And remember that the $D^2 (\Gamma_i) \times (\xi_0 = 0)$'s are NOT part of the 2-skeleton in question.] In this case
$$
\{\mbox{extended cocore $b_i$}\} \equiv b_i \times [-1 \leq \xi_0 \leq 0] .
$$
It will turn out that most of the theory in the present paper becomes TRIVIAL in this situation and we will not dwell on it any longer. [When, via the rest of this paper we will understand how to deal with the generic case below, then automatically and moreover trivially so, one will also know how to deal with our trivial case i).]

\medskip

ii) The {\ibf generic} case, when $b_i \times (\xi_0 = 0) \subset e \subset \Gamma (1) \times (\xi_0 = 0)$ is incident, in $X^2 [{\rm NEW}]$ to 2-cells NOT in $\Delta^2 \times (\xi_0 = 0)$ and then $\{$extended cocore $b_i\}$ is an infinite tree. 

\smallskip

In this case, our $b_i \times (\xi_0 = 0)$ is touched by the infinite collapsing flow unleashed by $C \cdot h = {\rm id} + {\rm nilpotent}$ inside $X_0^2 [{\rm NEW}]$.

\smallskip

So, in the generic case 
$$
\{\mbox{extended cocore of $b_i \times (\xi_0 = 0)\} =$} 
$$
$$
\mbox{$\bigcup \, \{$all the trajectories of the $2^{\rm d}$ RED collapsing flow of $2X_0^2$, hitting $b_i \times (\xi_0 = 0)\}$}.
$$

\smallskip

In this case, we also have, of course
$$
\{\mbox{extended cocore $b_i \times (\xi_0 = -1)\} = b_i \times [-1 \leq \xi_0 \leq 0] \cup \{$extended cocore of $b_i \times (\xi_0 = 0)\}$.}
$$

This definition easily extends to all points $x \in 2X_0^2$ {\ibf not} in $\bigcup \overset{\!\!\circ}{D^2} (\Gamma_i) \times (\xi_0 = -1)$. Retain that all $x \in R_1 = \{ R_1 , R_2 , \ldots , R_n \subset \Delta^2 \times (\xi_0 = -1)\} + \underset{1}{\overset{\infty}{\sum}} \, h_n$, or $x \in B_1$ or $x \in {\rm int} \, D^2 (C)$, with $C$ like in (18.1) have extended cocores, but certainly {\ibf not} the $x \in {\rm int} \, (D^2 (\Gamma) \times (\xi_0 = -1))$.

\smallskip

Even better, let us go now $4$-dimensional, then we get a PROPERLY embedded copy of $B^3 - \{$a tame Cantor set of $\partial B^3 \}$, denote it with the same notation as above,
$$
\{\mbox{extended cocore of $x$}\} \subset N^4 \left(2\Gamma (\infty) \cup \sum_1^{\infty} D^2 (C_i) \right) \subset N^4 (2X_0^2).
$$

Here $x \in \partial \{\mbox{extended cocore of $x$}\}$ (conceived now as a $3^{\rm d}$ object) and the embedding $\{$exterior cocore $x\} - \{x\} \subset N^4 (2X_0^2)$ is proper. [Reminder on the terminology. PROPER means inverse image of compact is compact, while proper means interior to interior and boundary to boundary.]

\smallskip

The $x \in {\rm int} \, D^2 (\Gamma_i)$'s DO NOT have such extended cocores. $\Box$

\bigskip

\noindent {\bf Compactification.} We have given a general idea how to get our $N^4 (2X_0^2)$, with some details, like the explicit contributions of the vertices $b^3 (\beta)$, $b^3 (\pm \, \xi_0)$ to be made more explicit later.

\smallskip

But what we have just said above should suffice, for right now. In particular, notice the following decomposition

\bigskip

\noindent (29) \quad $2X_0^2 = \Delta^2$ (always meaning $\Delta^2 \times (\xi_0 = -1)) \cup \underset{1}{\overset{K}{\sum}} \, T_i \ (=$ a finite union of trees, each resting with its foot on $\Delta^2) + \underset{1}{\overset{\infty}{\sum}} \, h_j + \underset{1}{\overset{\infty}{\sum}} \, D^2 (C_j)$.

\bigskip

Here, of course, the $h_j$, $D^2 (C_j)$ are our RED 1-handles and 2-handles, coming with $C_j \cdot h_i = {\rm id} + {\rm nil}$ (of the easy type), like in (19). But, for expository purposes, let us be for a while a bit more general, and just assume, in the context of (29), that we are given a canonical isomorphism between the sets of 1-handles and 2-handles, such that
$$
\partial D^2 (C_i) \cdot h_j = \delta_{ij} + \{\mbox{off diagonal} \ \eta_{ij} \in Z_+ \} . \leqno (30)
$$
We will make it clear, explicitly, when we go from the general (30) to the easy id $+$ nil of (19). But we will stay general, for a short while.

\smallskip

Our $h_i$ and $D^2 (C_i)$ are actually $4^{\rm d}$ handles of index $\lambda = 1$ and $\lambda = 2$, respectively; we will set for them

\bigskip

\noindent (31) \quad $h_i = ($arc $J_i =$ core of the $1$-handle) $\times$ ($3^{\rm d}$ ball $B_i^*$, the cocore) $= J_i \times B_i^*$, $D^2 (C_i) = $ (disc $D_i =$ core) $\times$ (another $2^{\rm d}$ disc $D_i^*$, the cocore) $= D_i \times D_i^*$.

\bigskip

With this, the canonical $\delta_{ij}$ in (30) comes with inclusions:
$$
J_i (\mbox{core of $h_i$}) \subset \partial D_i^2 = C_i \ (\mbox{we think here in terms of} \ D^2 (C_i) \approx D_i ({\rm core})) , \partial (B_i^* ({\rm cocore})) \supset D_i^* ({\rm cocore}).
$$

The $4^{\rm d}$ context we talk about here, means the not yet fully explicit $N^4 (2X_0^2)$. The $h_i$, $D^2 (C_i)$ are now $4^{\rm d}$ cells with disjoined interiors, out of which we will put together the following connected sum along the respective boundaries
$$
h_i \cup D^2 (C_i) = h_i \underset{\overbrace{\mbox{\footnotesize (core of $h_i$) $\times$ (cocore of $D^2 (C_i)) = J_i \times D_i^*$}}}{\#} D^2 (C_i) \underset{\rm DIFF}{=} B^4 \equiv \mbox{\ibf the ``state'' $i$}.
$$

\bigskip

Figure 12 should show what we are talking about, with one ambient dimension less.

$$
\includegraphics[width=17cm]{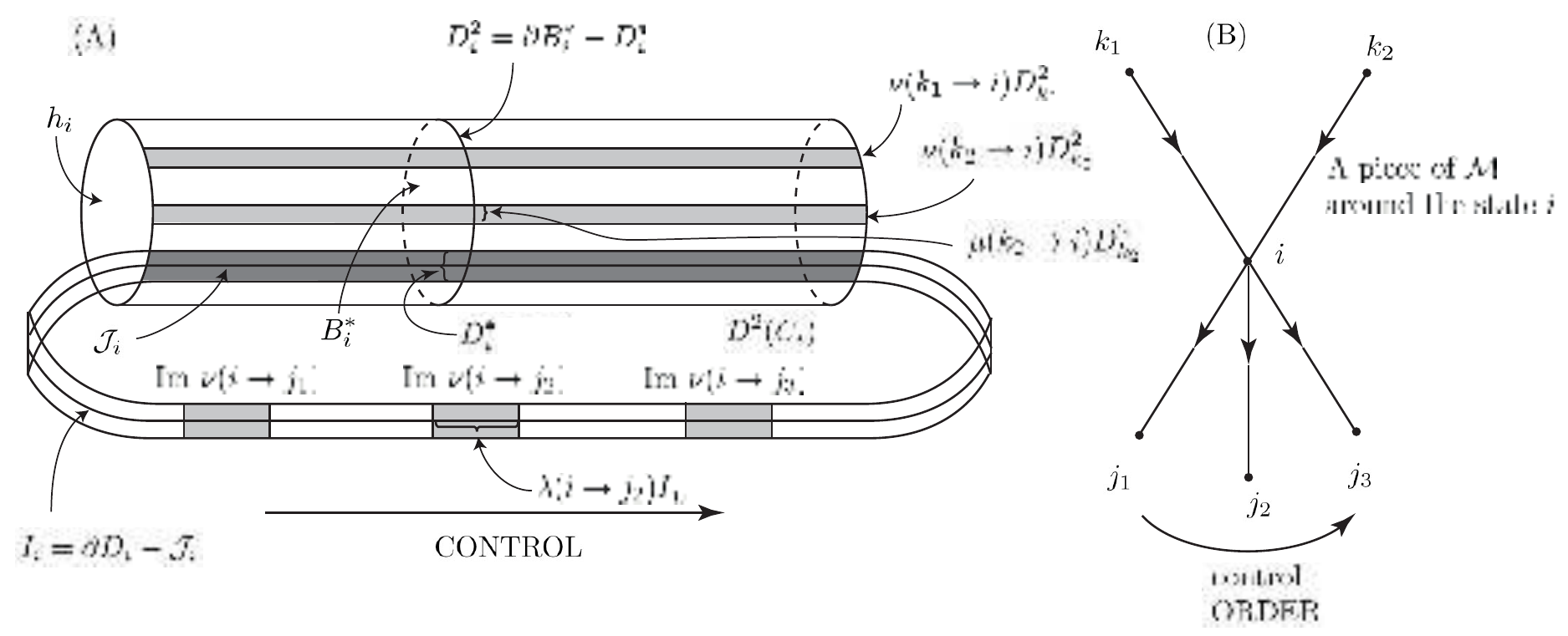}
$$
\label{fig12}
\centerline {\bf Figure 12.}
\begin{quote}
We display here, with one dimensionless, the state $i$.

LEGEND: $\includegraphics[width=1cm]{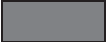} = J_i \times D_i^* \approx {\rm Box}  (i)$; $\includegraphics[width=1cm]{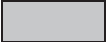} = \nu (i \to j_k) \, {\rm Box} (i) \subset {\rm Box} (j_k)$.

$\longrightarrow \ =$ CONTROL arrow. Here, because $\dim I_i = 1$ there is a natural linear order (well-defined up to a global orientation, which can be chosen arbitrarily, on the set of arrows $i \to j$, which are outgoing from $i$. In real life, $\dim D_i^2 = 2$ and there is NO such order for the arrows incoming into $i$.
\end{quote}

\bigskip

To the matrix (30) with its leading $\delta_{ij}$-term followed by off-diagonal terms, we will associate an oriented graph denoted ${\mathcal M}$ (30). The vertices of this ${\mathcal M}$ are the various states $i$ above, i.e. the $h_i \cup D^2 (C_i)$'s but another incarnation of $i$ will be the ${\rm Box} (i)$ defined in (33) below. The off-diagonal part of (30) defines the oriented edges of ${\mathcal M}$, via the recipe
$$
\# \, \{{\rm edges} \ i \to j \ {\rm in} \ {\mathcal M}\} = \{ \eta_{ij} \ \mbox{in (30)}\}.
$$

Such oriented edges will also be called ``arrows''.

\smallskip

In connection with the formulae (31) we introduce the following objects, displayed in the Figure 12:
$$
D_i^2 ({\rm disc}) = \partial \, {\rm cocore} \, h_i - {\rm cocore} \, D^2 (C_i) = \partial B_i^* - D_i^* \leqno (32)
$$
$$
I_i ({\rm arc}) = \partial \, {\rm core} \, D^2 (C_i) - {\rm core} \, h_i = \partial D_i - J_i .
$$
These formulae suggest obvious identifications, rel $\partial I_i = \partial J_i$, $\partial D_i^2 = \partial D_i^*$,
$$
I_i \approx J_i \, (= 1^{\rm d} \, {\rm core} \, h_i), \quad D_i^2 \approx D_i^* (= 2^{\rm d} \, \mbox{cocore of the $4^{\rm d}$ handle of index $2$,} \, D^2 (C_i)). \leqno (32.1)
$$

More explicitly, $D_i^2$ and $D_i^*$ are the two hemispheres of the $2$-sphere $\partial B_i^*$ and, like wise, $I_i$ and $J_i$ the two halves of the circle $C_i = \partial D_i$.

\smallskip

I will also introduce now the following purely abstract object
$$
{\rm Box} (i) \equiv I_i \times D_i^2 = (\partial D^2 (C_i) - J_i) \times (\partial B_i^* - D_i^*). \leqno (33)
$$
The formula (33) can be visualized on Figure 12 and ${\rm Box} (i)$ is a $3^{\rm d}$ cell. We say it is ``abstract'' because it does not come with any a priori canonical embedding into some ambient space.

\smallskip

Of course, in (33) $I_i , D_i^2$ are like in (32) and there is an abstract identification ${\rm Box} (i) \approx J_i \times D_i^* \subset \partial h_i$, an embedding suggested by the double shading in Figure 12.

\bigskip

\noindent {\bf Lemma 6.} {\it To each of the $\eta_{ij}$ arrows $i \to j$, coming with the matrix ${\mathcal M}$ {\rm (30)}, correspond embeddings, independant from each other,
$$
I_j \xrightarrow{ \ \lambda (i \to j) \ } I_i \approx J_i \quad \mbox{and} \quad D_i^2 \approx D_i^* \xrightarrow{ \ \mu (i \to j) \ } D_j^2 ;
$$
see Figure {\rm 12}.

\smallskip

Out of these two we can extract the following map, which is a smooth embedding when restricted to its domain of definition
$$
{\rm Box} (i) = I_i \times D_i^2 \supset (\lambda (i \to j) I_j) \times D_i^2 \underset{\nu (i \to j) \, \equiv \, {\rm id} \times \mu (i \to j)}{-\!\!\!-\!\!\!-\!\!\!-\!\!\!-\!\!\!-\!\!\!-\!\!\!-\!\!\!-\!\!\!-\!\!\!-\!\!\!-\!\!\!-\!\!\!\longrightarrow} I_j \times D_j^2 = {\rm Box} (j). \leqno (34)
$$
Notice here that the $\lambda (i \to j) I_j \subset I_i$ is an isomorphic copy of $I_j$.}

\bigskip

It is best to think of the $\nu (i \to j)$ as a purely abstract map. Figures 12 and 13 should help visualizing these things. In the formula (34), for $x \in I_j$, the id factor of $\nu (i \to j)$, takes $\lambda (x) = \lambda (i \to j)(x)$ back to $x$. Whatever minor ambiguity there might be in the definition of the map $\nu (i \to j)$, the point here is that the ${\rm Box} (i)$ goes {\ibf cleanly through} ${\rm Box} (j)$, like in a standard {\ibf Markov partition}. This is suggested in the Figure 13.
$$
\includegraphics[width=13cm]{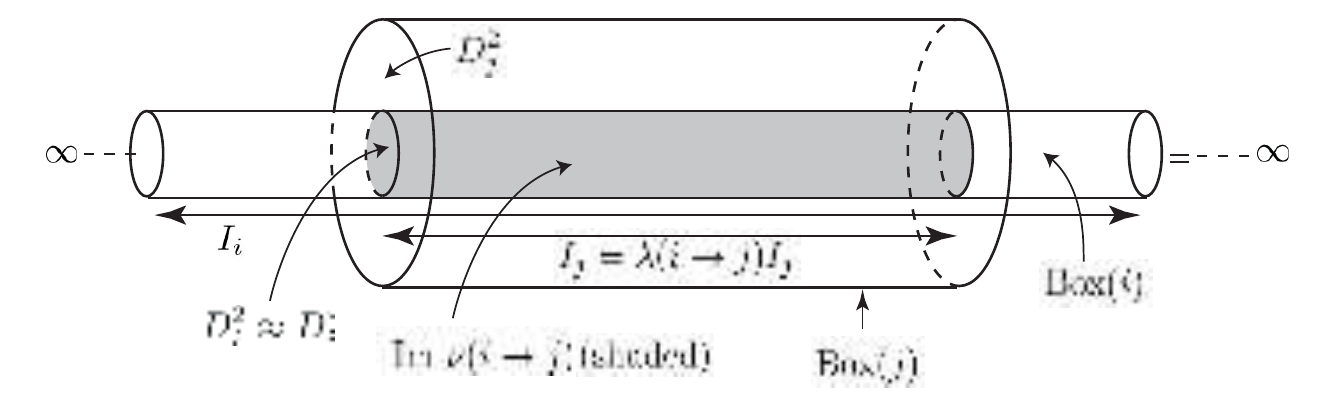}
$$
\label{fig13}
\vglue -7mm \centerline {\bf Figure 13.}

\begin{quote}
This figure accompanies the arrow $i \to j$ of ${\mathcal M}$, and the points $p \in {\rm Box} (i)$ and $q \in {\rm Box} (j)$ which are such that $q = \nu (i \to j) p$, get identified. See here (32.1) too. Our figure shows what the map $\nu (i \to j)$ defined in (34) does.
\end{quote}

\bigskip

Let ${\mathcal M} = \underset{\alpha}{\sum} \, {\mathcal M}_{\alpha}$ be the decomposition into connected components and $A = A(\alpha) \equiv \underset{i \in \alpha}{\sum} {\rm Box} (i)$. On the set $A (\alpha)$ we introduce the equivalence relation ${\mathcal R} = {\mathcal R} (\alpha)$, generated as follows: For each arrow $\{ i \to j \} \in {\mathcal M}_\alpha$ and for each $p \in \{ \lambda (i \to j) I_j \times D_i^2$, the domain of definition of the map $\nu (i \to j)$ (34)$\} \subset {\rm Box} (i)$, with $q = \nu (i \to j)(p) \in {\rm Box} (j)$ one has, to begin with,
$$
(p,q) \in {\mathcal R}.
$$
The whole equivalence relation ${\mathcal R}$ is, afterwards, generated by these pairs, by adding the obvious things like
$$
\{ (p,q) , (q,r) \in {\mathcal R}\} \Longrightarrow (p,r) \in {\mathcal R}, \ \mbox{a.s.o.}
$$
With this we consider the quotient space
$$
X^3 (\alpha) \equiv \left( \sum_{i \in \alpha} {\rm Box} (i) \right) / \, {\mathcal R} (\alpha). \leqno (35)
$$

At the (temporary) level of whole generality, there is a little problem here. In the general case ${\mathcal M}$ might have oriented {\ibf cycles}. Then, there can be orbits of ${\mathcal R}$ which, as subsets of $A({\alpha})$ may fail to be closed sets. Then, continuous functions no longer separate points, and $X^3 (\alpha)$ only makes sense as a {\ibf non-commutative space}, with the algebra of functions which is natural for quotient-spaces, namely matrices of functions, of the form
$$
\begin{pmatrix}
f(p,p), &f(p,q) \\ f(q,p), &f(q,q)
\end{pmatrix},
$$
with the matrix algebra multiplication law, i.e. convolution product
$$
(f * g)(p,q) = \sum_{\overbrace{\mbox{\footnotesize$(p,r) , (r,q) \in R$}}} f(p,r) \, g(r,q) \qquad \mbox{(see [C]).}
$$

Fortunately, in our real-life case, when (30) is of the easy id $+$ nil type these complications vanish into thin air. In our real-life case of the easy id $+$ nil matrix (30) we have the following

\bigskip

\noindent {\bf Lemma 7.} {\it When the graph ${\mathcal M} (30)$ is of the type easy id $+$ nil, then the $X^3 (\alpha)$ from $(35)$ is a mundane connected non-compact $3$-manifold with non-empty boundary, which comes naturally equipped with a {\ibf lamination} ${\mathcal L}_{\alpha}$ by lines.}

\bigskip

The proof of Lemma 7 will be given later. [For the time being I will only give some comments. With any ${\mathcal M} (30)$ the Markov type Figures 13 are always there. But the issue is what goes on when one puts them together. In the easy id $+$ nilpotent case, they pile up nicely, we have things like (55) below and one gets a smooth final object like in Lemma 7. In particular, with ${\mathcal M} (30)$ of the type easy nil $+$ nilpotent, there are neither closed orbits nor dangerous going-back trajectories \`a la ${\rm Wh}^3$, for ${\mathcal R} (\alpha)$ and our $X^3 (\alpha)$ from (35) is then a mundane $3$-manifold, useful for our present purposes. But let us stop here for one minute and look back. In the most general case, when for some manifold we give a handle-decomposition coming with a canonical identification between the following two sets
$$
\{\mbox{handles of index $\lambda + 1$}\} \approx \{\mbox{handles of index $\lambda$}\},
$$
dictated by a corresponding geometric intersection matrix \`a la (30), then strange objects, quotient spaces which may be non-commutative spaces, in no way directly related with the manifold under investigation, are lurking around.

\smallskip

Of course, with our easy id $+$ nilpotent condition, for the geometric intersection matrix, these additional spaces are no longer {\ibf non}-commutative, but mundane manifolds which we will happily use in our constructions. Notice here the following paradigm, the first part of which is a proved result (our Lemma 7) while the second part more like a hunch which deserves, I think, to be investigated:
$$
\mbox{GSC (i.e. easy id $+$ nilpotent)} \Longrightarrow \mbox{an $X^3 (\alpha)$ which is a nice smooth manifold},
$$
and then,
$$
\mbox{general case (non-GSC)} \Longrightarrow \mbox{an $X^3 (\alpha)$ which is a horrible quotient-space possibly}
$$
$$
\mbox{a {\ibf non}-commutative space (\`a la Alain Connes)}.
$$
It should be stressed that all this little story does not function only in the non-compact context which interests us here, but in the compact case too.] From now on we go back to our real life easy id $+$ nil situation.

\smallskip

In the context of Lemma 7, we denote by $\bar\varepsilon (X^3 (\alpha)) \subset X^3 (\alpha)$ the closed set of our lamination ${\mathcal L}_{\alpha}$. Transversally to $\bar\varepsilon (X^3 (\alpha))$, the $X^3 (\alpha)$ is a pair

\bigskip

\noindent ($R^2$, a {\ibf tame} Cantor set $\subset R^2$, i.e. a Cantor set $C$ such that its inclusion into $R^2$ factorizes through a smooth line $C \subset L \subset R^2$).

\bigskip

Let us be more specific. Every $D_i^2$ is a transversal to the lamination ${\mathcal L}_{\alpha}$, the closed set of which is $\bar\varepsilon \, X^3 (\alpha)$ and, we have
$$
(\bar\varepsilon \, X^3 (\alpha)) \cap D_i^2 = \bigcup_{\overbrace{{{ i \, = \, i_0 \leftarrow i_1 \leftarrow i_2 \leftarrow \ldots , \atop \mbox{\footnotesize all the trajectories of ${\mathcal M}$}} \atop \mbox{\footnotesize ending at $i$}}}} \ \bigcap_{j=1}^{\infty} \mu (i_j \to i) \, D_{i_j}^2.
$$

\smallskip

Let's go now to the smooth $4^{\rm d}$ version of (29),

\bigskip

\noindent (36) \quad $N^4 (2X_0^2) = \left(N^4 (2\Gamma (\infty)) - \underset{1}{\overset{\infty}{\sum}} \, h_i \right) + \underset{1}{\overset{\bar n}{\sum}} \, D^2 (\Gamma_j) + \underset{1}{\overset{\infty}{\sum}} \, (h_i \cup D^2 (C_i)) = N^4 (\Delta^2) \cup \underset{\ell = 1}{\overset{K}{\sum}} \, \{$the $4^{\rm d}$ regular neighbourhood of a tree $T_{\ell}$ resting on $\partial N^4 (\Delta^2)$, and which we continue to denote by $T_{\ell}\} + \underset{1}{\overset{\infty}{\sum}} \, h_i + \underset{1}{\overset{\infty}{\sum}} \, D^2 (C_i)$.

\bigskip

Here $N^4 (2\Gamma (\alpha)) - \underset{1}{\overset{\infty}{\sum}} \, h_i \equiv \biggl\{$the $N^4 (2\Gamma (\infty))$ split along $\underset{1}{\overset{\infty}{\sum}} \, {\rm cocore} \, h_i \biggl\} \supset N^4 (\Gamma(1)) \supset \underset{1}{\overset{n}{\sum}} \, R_i$ (the RED $1$-handles if $N^4 (\Delta^2))$.

\smallskip

Let $\varepsilon (T_{\ell})$ denote the set of end-points of the tree $T_{\ell}$. At $4^{\rm d}$ level we have an obvious smooth compatification, albeit a completely trivial one

\bigskip

\noindent (37) \quad $T_{\ell}^{\wedge} \equiv \{$The $4^{\rm d}$ $T_{\ell}\} \cup \varepsilon (T_{\ell}) \equiv B^4 (T_{\ell}) \underset{\rm DIFF}{=} B^4$, with $\varepsilon (T_{\ell}) \subset \{$a tame Cantor set inside $\partial B^4 (T_{\ell}\}$.

\bigskip

Of course, $\varepsilon (T_{\ell})$ lives at the infinity of $T_{\ell}$.

\smallskip

So, we get here another diffeomorphic model for $N^4 (\Delta^2)$
$$
N^4 (\Delta^2)^{\bullet} = N^4 (\Delta^2) \cup \sum_{\ell = 1}^K \widehat T_{\ell} = \left( N^4 (2\Gamma (\infty)) - \sum_1^{\infty} h_n \right)^{\wedge} . \leqno (37.1)
$$
It is essential here that we have $\left( \underset{1}{\overset{\bar n}{\sum}} \, \Gamma_j \right) \cap h_i = \emptyset$, a feature which will be lost when, for the so-called BALANCING purpose, we will have to extend the $\underset{1}{\overset{\bar n}{\sum}} \, \Gamma_j$.

\smallskip

For each $\alpha$ in ${\mathcal M} = \underset{\alpha}{\sum} \, {\mathcal M}_{\alpha}$, the following object is a smooth non-compact 4-manifold with non-empty boundary
$$
{\rm LAVA}_{\alpha} \equiv \sum_{i \in \alpha} (h_i \cup D^2 (C_i)) , \leqno (38)
$$
and with this we also set
$$
{\rm LAVA} \equiv \sum_{\alpha} {\rm LAVA}_{\alpha}.
$$
This LAVA may well touch the $\underset{1}{\overset{n}{\sum}} \, R_i \times (\xi_0 = -1)$ but not to $\underset{1}{\overset{\bar n}{\sum}} \, D^2 (\Gamma_j)$. We also have
$$
\delta \, {\rm LAVA}_{\alpha} \equiv (\partial \, {\rm LAVA}_{\alpha}) \cap \partial \left( N^4 (2X_0^2) - \underset{1}{\overset{\infty}{\sum}} \, h_k \right) , \ \delta \, {\rm LAVA} \equiv \sum_{\alpha} \delta \, {\rm LAVA}_{\alpha} . \leqno (38.1)
$$

Here is how $N^4 (2X_0^2)$ is structured
$$
N^4 (2X_0^2) = \left[ \left( N^4 (2\Gamma (\infty)) - \sum_1^{\infty} h_k \right) \underset{\overbrace{\mbox{\footnotesize$\delta \, {\rm LAVA}$}}}{\cup} {\rm LAVA} \right] + \sum_1^{\bar n} D^2 (\Gamma_j). \leqno (39)
$$

In the RHS of this formula, LAVA is a PROPER codimension zero submanifold. Each component ${\rm LAVA}_{\alpha} \subset {\rm LAVA}$ is actually itself PROPER, and it is split from the rest by the connected codimension one submanifold $\delta \, {\rm LAVA}_{\alpha}$. Very importantly, in the context of (39), the 2-handles $D^2 (\Gamma_j)$ are directly attached to $N^4 (2\Gamma (\infty)) - \underset{1}{\overset{\infty}{\sum}} \, h_k$.

\bigskip

\noindent {\bf Lemma 8.} 1) {\it As far as the ends are concerned, we have
$$
\varepsilon \left( N^4 (2\Gamma (\infty)) - \sum_1^{\infty} h_k \right) = \sum_1^K \varepsilon (T_i),
$$
a formula which would of course no longer be true with the LHS replaced by $\varepsilon (N^4 (2\Gamma (\infty)))$, (in which case we would only find a quotient space projection $\underset{1}{\overset{K}{\sum}} \, \varepsilon (T_i) \twoheadrightarrow \varepsilon (N^4 (2\Gamma (\infty)))$. Also, there is an easy smooth compactification
$$
\left( N^4 (2\Gamma (\infty)) - \sum_1^{\infty} h_k \right)^{\wedge} = \left( N^4 (2\Gamma (\infty)) - \sum_1^{\infty} h_k \right) \cup \varepsilon \left( N^4 (2\Gamma (\infty)) - \sum_1^{\infty} h_k \right) \underset{\rm DIFF}{=} \leqno (40)
$$
$$
\underset{\rm DIFF}{=} n \, \# \, (S^1 \times B^3) \quad \mbox{(with $n = \, \# \, R \cap \Gamma (1)$)}.
$$

Moreover, in the same vein, we have the diffeomorphism
$$
N^4 (\Delta^2) \underset{\rm DIFF}{=} \left( N^4 (2\Gamma (\infty)) - \sum_1^{\infty} h_k \right)^{\wedge} + \sum_1^{\bar n} D^2 (\Gamma_j) , \leqno (40.1)
$$
with the RHS of $(40.1)$ being the same object as the $N^4 (\Delta^2)^{\bullet}$ above (see {\rm (37.1)}).}

\medskip

1-bis) {\it Without any loss of generality, $\left( N^4 (2\Gamma (\infty)) - \underset{1}{\overset{\infty}{\sum}} \, h_k \right)^{\wedge}$ is diffeomorphic to the closure of $ N^4 (2\Gamma (\infty))$ $-$ $\underset{1}{\overset{\infty}{\sum}} \, h_k$ in the ambient space $N^4 (2X_0^2)$; so we have the following diffeomorphism between the abstract compactification occurring in the LHS of the formula below, and the mundane closure inside the ambient space $N^4 (2X_0^2)$ {\rm (39)}
$$
\left( N^4 (2\Gamma (\infty)) - \sum_1^{\infty} h_n \right)^{\wedge} \underset{\rm DIFF}{=} \overline{\left( N^4 (2\Gamma (\infty)) - \sum_1^{\infty} h_n \right)} .
$$
This kind of equalities abund in the sequel of this paper, but we will not always bother to write them down explicitly. The reader should not have any difficulty in guessing them.

\smallskip

The $\varepsilon \left( N^4 (2\Gamma (\infty)) - \underset{1}{\overset{\infty}{\sum}} \, h_k \right) \subset \partial \left( N^4 (2\Gamma (\infty)) - \underset{1}{\overset{\infty}{\sum}} \, h_k \right)^{\wedge}$ is a tame Cantor set and we also have the following equality between spaces of ends $\varepsilon \left( N^4 (2\Gamma (\infty)) - \underset{1}{\overset{\infty}{\sum}} \, h_k \right) = \varepsilon \left( \partial \left( N^4 (2\Gamma (\infty)) - \underset{1}{\overset{\infty}{\sum}} \, h_k \right)\right)$.}

\medskip

2) {\it There is a PROPER {\ibf smooth embedding}
$$
\sum_{\alpha} X^3 (\alpha) \xrightarrow{ \ J \ } \partial \left( N^4 (2\Gamma (\infty)) - \sum_1^{\infty} h_k \right) = \partial \left( N^4 (2\Gamma (\infty)) - \sum_1^{\infty} h_k \right)^{\wedge} - \varepsilon \left(\partial \left( N^4 - \sum_1^{\infty} h_k \right)\right). \leqno (41)
$$
Moreover, in terms of {\rm (38.1)} we have here the following equality of sets}
$$
\delta \, {\rm LAVA}_{\alpha} = J X^3 (\alpha) \approx X^3 (\alpha). \leqno (42)
$$

3) {\it For each $\alpha$ there is a diffeomorphism of pairs
$$
({\rm LAVA}_{\alpha} , \delta \, {\rm LAVA}_{\alpha}) \underset{\rm DIFF}{=} (X^3 (\alpha) \times [0,1] - \left( \underbrace{\bar\varepsilon (X^3 (\alpha))}_{\mbox{\footnotesize the closed set of the lamination ${\mathcal L}_{\alpha}$}} \times \{1\} \right) , X^3 (\alpha) \times \{ 0 \}). \leqno (43)
$$
}

\bigskip

The proof and various embellishments of this lemma, will follow in this section, later on. We introduce now the endpoint compactification
$$
(\delta \, {\rm LAVA}_{\alpha})^{\wedge} \equiv \delta \, {\rm LAVA}_{\alpha} \cup \varepsilon (\delta \, {\rm LAVA}_{\alpha}).
$$
As a consequence of (42), (43), we have an (almost) SMOOTH COMPACTIFICATION OF LAVA
$$
{\rm LAVA}_{\alpha}^{\wedge} = (\delta \, {\rm LAVA}_{\alpha} \times [0,1]) \underset{\overbrace{\mbox{\footnotesize $\delta \, {\rm LAVA}_{\alpha}$}}}{\cup} (\delta \, {\rm LAVA}_{\alpha})^{\wedge} = (\delta \, {\rm LAVA}_{\alpha} \times [0,1]) \cup \varepsilon ({\rm LAVA}_{\alpha}), \leqno (44)
$$
and here the $\delta \, {\rm LAVA}_{\alpha} \times [0,1] \subset {\rm LAVA}_{\alpha}^{\wedge}$ is a bona fide smooth manifold. But it should be stressed here that (44) as well as the
$$
(\delta \, {\rm LAVA}_{\alpha})^{\wedge} \equiv \delta \, {\rm LAVA}_{\alpha} \cup \varepsilon (\delta \, {\rm LAVA}_{\alpha})
$$
are provisional formulae, soon to be superseded by more appropriate and more accurate versions. Also, (43), (44), $\ldots$ are ABSTRACT, meaning that they come, a priori, without any connections with the embedding ${\rm LAVA}_{\alpha} \subset N^4 (2X_0^2)$; all this will be corrected in time.

\smallskip

In the context of (43), the induces diffeomorphism
$$
\delta \, {\rm LAVA}_{\alpha} = X^3 (\alpha) \times \{ 0 \},
$$
is nothing else but the following equality among subsets of $\partial \left( N^4 (2\Gamma (\infty)) - \underset{1}{\overset{\infty}{\sum}} \, h_n \right)$, the (42), which I will re-write here
$$
\delta \, {\rm LAVA}_{\alpha} = JX^3 (\alpha) \approx X^3 (\alpha).
$$
One should not mix up the $\varepsilon (\delta \, {\rm LAVA}_{\alpha} ) \equiv \{$the space of ends in the sense of Freudenthal, Hopf et al$\}$ and the $\varepsilon ({\rm LAVA}_{\alpha})$ which will NOT mean ends, but the set of points at infinity of ${\rm LAVA}_{\alpha}$, to be made explicit below, and which is certainly not totally discontinuous, since the lamination ${\mathcal L}_{\alpha}$ contributes to $\varepsilon ({\rm LAVA}_{\alpha})$. Each leaf of it joins a pair of points at infinity, living in $\varepsilon (\delta \, {\rm LAVA}_{\alpha})$. With this
$$
\varepsilon ({\rm LAVA}_{\alpha}) = \bar\varepsilon (X^3 (\alpha)) \cup \varepsilon (\delta \, {\rm LAVA}_{\alpha}) \ \mbox{(and this is a set of ends).} \leqno (44.1)
$$

The embedding $J$ from (41) induces, automatically, a continuous map, at the level of the ends \`a la Hopf-Freudenthal
$$
\sum_{\alpha} \varepsilon (X^3 (\alpha)) \xrightarrow{ \ E \ } \varepsilon \left( N^4 - \sum_1^{\infty} h_k \right) = \varepsilon \left(\partial \left( N^4 - \sum_1^{\infty} h_k \right)\right). \leqno (45)
$$
For each $\alpha$ we consider the closed subspace $\lambda_{\alpha} \equiv E (\varepsilon (X^3 (\alpha))) \subset \varepsilon \left( N^4 - \underset{1}{\overset{\infty}{\sum}} \, h_k \right)$ and the quotient-space projection
$$
\varepsilon (\delta \, {\rm LAVA}_{\alpha}) = \varepsilon (X^3 (\alpha)) \twoheadrightarrow \lambda_{\alpha} . \leqno (45.1)
$$
\vglue-6mm
$$
\mbox{\hglue -5mm}{\uparrow\mbox{\hglue -2mm}}_{-\!\!-\!\!-\!\!-\!\!-\!\!-\!\!-\!\!-\!\!-\!\!-\!\!-\!\!-\!\!-\!\!-\!\!-\!\!-\!\!-\!\!-}{\mbox{\hglue -2mm}\uparrow}
$$
\vglue-5mm
$$
\mbox{\hglue -4mm\footnotesize $\approx$ (see (42))}
$$
The compact space $\lambda_{\alpha}$ are not necessarily disjoined and we have
$$
\overline{{\rm Im} \, E} = \overline{\bigcup_{\alpha} \lambda_{\alpha}} = \bigcup_{\alpha} \lambda_{\alpha} \cup \lambda_{\infty}, \leqno (46)
$$
where $\lambda_{\infty} \subset \varepsilon \left( N^4 - \underset{1}{\overset{\infty}{\sum}} \, h_k \right)$ is the closed set of points $p_{\alpha}$ such taht we can find a sequence $p_{\alpha} \in \lambda_{\alpha}$ with $\underset{\alpha = \infty}{\lim} \, p_{\infty} = p_{\infty}$.

\smallskip

We will sharpen now the definition of ${\rm LAVA}_{\alpha}^{\wedge}$ given in (44) by setting, from now on the following improved (and/or corrected) form of (44), namely
$$
{\rm LAVA}_{\alpha}^{\wedge} \equiv (\delta \, {\rm LAVA}_{\alpha} \times [0,1]) \cup \lambda_{\alpha} . \leqno (47)
$$
Notice that the RHS of this formula, which supersedes (44) is a quotient-space of the RHS of (44). Similarly, we will set
$$
(\delta \, {\rm LAVA}_{\alpha})^{\wedge} \equiv (\delta \, {\rm LAVA}_{\alpha}) \cup \lambda_{\alpha} . \leqno (47.1)
$$

So, from now on, we use the sharper $\lambda_{\alpha}$ instead of the ABSTRACT $\varepsilon (\delta \, {\rm LAVA}_{\alpha})$.

\bigskip

\noindent {\bf The compactification Lemma 9.} 1) {\it Via a smooth, not necessarily ambient isotopy of the $N^4 (\Delta^2) \cup \left( N^4 (2\Gamma (\infty)) - \underset{1}{\overset{\infty}{\sum}} \, h_k \right)$, inside the ambient space $N^4 (2 X_0^2)$, we can get, for each $\alpha$, the following equality, connecting the LHS CONCRETE pair to the RHS ABSTRACT pair, as conceived from now on}
$$
\left( \overline{{\rm LAVA}_{\alpha}} , \overline{\delta \, {\rm LAVA}_{\alpha}} \right) = \left( \underbrace{(\delta \, {\rm LAVA}_{\alpha} \times [0,1]) \cup \lambda_{\alpha}}_{{\rm LAVA}_{\alpha}^{\wedge}} \, , \underbrace{(\delta \, {\rm LAVA}_{\alpha} \cup \lambda_{\alpha}) \times \{0\}}_{(\delta \, {\rm LAVA}_{\alpha})^{\wedge}} \right) . \leqno (48)
$$

2) {\it The smooth non-compact $4$-manifold, part of $N^4 (2X_0^2)$
$$
\left(N^4 (2\Gamma (\infty)) - \sum_1^{\infty} h_k \right) \underset{\overbrace{\mbox{\footnotesize $\delta \, {\rm LAVA}$}}}{\cup} {\rm LAVA} \qquad \mbox{(see $(39)$)}
$$
has the following canonical smooth compactification

\bigskip

\noindent {\rm (49)} \quad $\left[ \left( N^4 (2\Gamma (\infty)) - \underset{1}{\overset{\infty}{\sum}} \, h_k \right) \cup {\rm LAVA} \right]^{\wedge} \equiv \Biggl\{\Biggl\{ \left(N^4 (2\Gamma (\infty)) - \underset{1}{\overset{\infty}{\sum}} \, h_k \right)^{\wedge} \cup \underset{\alpha}{\bigcup} \ {\rm LAVA}_{\alpha}^{\wedge}$, where each $\lambda_{\alpha} \subset {\rm LAVA}_{\alpha}^{\wedge}$ is identified with its counterpart in $\varepsilon \left(\left( N^4 (2\Gamma (\infty)) - \underset{1}{\overset{\infty}{\sum}} \, h_k \right)^{\wedge} \right) \Biggl\}\Biggl\}^{\wedge} = \left( N^4 (2\Gamma (\infty)) - \underset{1}{\overset{\infty}{\sum}} \, h_k \right)^{\wedge}$ $\cup$ $ \underbrace{\cup \left(\underset{\alpha}{\bigcup} \ {\rm LAVA}_{\alpha}^{\wedge} \cup \lambda_{\infty} \right).}_{\mbox{\footnotesize This is now a closed subset of $\left(\left( N^4 (2\Gamma (\infty)) - \underset{1}{\overset{\infty}{\sum}} \, h_k \right) \cup {\rm LAVA} \right)^{\wedge}$}}$

\bigskip

Here the various $\delta \, {\rm LAVA}_{\alpha} \times [0,1] \subset {\rm LAVA}_{\alpha}^{\wedge}$ are $2$-by-$2$ disjoined, and also ${\rm LAVA}_{\alpha} \subset \left( N^4 (2\Gamma (\infty)) - \underset{1}{\overset{\infty}{\sum}} \, h_k \right)$ $\cup$ ${\rm LAVA}$ is a PROPER submanifold, unlike the $\delta \, {\rm LAVA}_{\alpha} \times [0,1] \subset \left(\left( N^4 (2\Gamma (\infty)) - \underset{1}{\overset{\infty}{\sum}} \, h_k \right) \cup {\rm LAVA} \right)^{\wedge}$, where the LHS lacks the $\lambda_{\alpha}$.

\smallskip

In the same vein, the $\underset{\alpha}{\bigcup} \ {\rm LAVA}_{\alpha}^{\wedge} \subset \left(\left( N^4 (2\Gamma (\infty)) - \underset{1}{\overset{\infty}{\sum}} \, h_k \right) \cup {\rm LAVA} \right)^{\wedge}$ is NOT a closed subset, since the
$$
\lambda_{\infty} \subset \varepsilon \left( N^4 (2\Gamma (\infty)) - \underset{1}{\overset{\infty}{\sum}} \, h_k \right)
$$
is still lacking to it. We have to add the $\lambda_{\infty}$, like in $(49)$ in order to clinch the closure. We have now a diffeomorphism}
$$
\left[ \left( N^4 (2\Gamma (\infty)) - \underset{1}{\overset{\infty}{\sum}} \, h_k \right) \cup {\rm LAVA} \right]^{\wedge} \underset{\rm DIFF}{=} \left( N^4 (2\Gamma (\infty)) - \underset{1}{\overset{\infty}{\sum}} \, h_k \right)^{\wedge} \underset{\rm DIFF}{=} n \, \# \, (S^1 \times B^3).
\leqno (50)
$$

3) {\it Finally, there is a diffeomorphism, giving our model for $N^4 (\Delta^2)$, from now on, at least until further notice,}
$$
\boxed{
N^4 (2X_0^2)^{\wedge} = \left[ \left( N^4 (2\Gamma (\infty)) - \underset{1}{\overset{\infty}{\sum}} \, h_k \right) \cup {\rm LAVA} \right]^{\wedge} + \sum_1^{\bar n} D^2 (\Gamma_j) \underset{\rm DIFF}{=} N^4 (\Delta^2).
} \leqno (51)
$$

\bigskip

\noindent {\bf Reminder.} In this formula, the 2-handles $D^2 (\Gamma_j)$ are attached directly to $N^4 (2\Gamma (\infty)) - \underset{1}{\overset{\infty}{\sum}} \, h_n$. We have
$$
\left\{ \left[ \left(N^4 (2\Gamma (\infty)) - \underset{1}{\overset{\infty}{\sum}} \, h_n \right) \cup {\rm LAVA} \right]^{\wedge}  + \sum_1^{\bar n} \, D^2 (\Gamma_i) \right\} - \left\{ \left[ \left(N^4 (2\Gamma (\infty)) - \underset{1}{\overset{\infty}{\sum}} \, h_n \right) \cup {\rm LAVA} \right] + \sum_1^{\bar n} \, D^2 (\Gamma_i) \right\} 
$$
$$
= \left\{ \varepsilon  \left(N^4 (2\Gamma (\infty)) - \underset{1}{\overset{\infty}{\sum}} \, h_n \right) + \sum_{\alpha} \varepsilon ({\rm LAVA}_{\alpha}) \ \mbox{(44.1)}  \right\} ;
\leqno (51.1)
$$
here, in the LHS we see twice the same expression, first with a hat, then without one. Clearly also, the $\sum D^2 (\Gamma_i)$ can be erased from that LHS. Finally, the  RHS should be read like in the (51.2) below.

\medskip

Let us make now explicit the topology of the disjoined union occurring in the RHS of (51.1). We start by adding, to each $\bar\varepsilon (X^3 (\alpha))$, the set of endpoints, living at its infinity, i;e. the
$$
\varepsilon (\bar\varepsilon \, X^3 (\alpha)) = \varepsilon \, X^3 (\alpha) \subset \varepsilon ({\rm LAVA}_{\alpha}).
$$
This is glued then to $\varepsilon \left(N^4 (2\Gamma (\infty)) - \underset{1}{\overset{\infty}{\sum}} \, h_n \right)$ according to the surjective map from (45.1)
$$
\varepsilon (X^3 (\alpha)) \twoheadrightarrow \lambda_{\alpha} \subset \varepsilon \left(N^4 (2\Gamma (\infty)) - \sum_1^{\infty} h_n \right) .
$$
This {\ibf is} the topology which (51.1) inherits from the ambient space. And then, in the same spirit as in formulae (48) and (49), we can express now the space in (51.1) as follows, and in the formula below, each $\lambda_{\alpha}$ occurring in $\bar\varepsilon (X^3 (\alpha)) \cup \lambda_{\alpha}$ is to be identified with its counterpart $\lambda_{\alpha} \subset \varepsilon \left(N^4 (2\Gamma (\infty)) - \underset{1}{\overset{\infty}{\sum}} \, h_n \right)$:
$$
\varepsilon \left(N^4 (2\Gamma (\infty)) - \sum_1^{\infty} h_k \right) \cup \sum_{\alpha} (\bar\varepsilon (X^3 (\alpha)) \cup \lambda_{\alpha} ) . \leqno (51.2)
$$

This is exactly what lives at the infinity of the space $\left[ \left(N^4 (2\Gamma (\infty)) - \underset{1}{\overset{\infty}{\sum}} \, h_k \right) \cup {\rm LAVA} \right]^{\wedge}$, i.e. it is exactly

\medskip

$$
\left[ \left(N^4 (2\Gamma (\infty)) - \sum_1^{\infty} h_k \right) \cup {\rm LAVA} \right]^{\wedge} - \left[ \left(N^4 (2\Gamma (\infty)) - \sum_1^{\infty} h_k \right) \cup {\rm LAVA} \right] =
$$

\smallskip

$$
\overline{\left[ \left(N^4 (2\Gamma (\infty)) - \sum_1^{\infty} h_k \right) \cup {\rm LAVA} \right]} - \left[ \left(N^4 (2\Gamma (\infty)) - \sum_1^{\infty} h_k \right) \cup {\rm LAVA} \right] .
$$

\medskip

The various $\lambda_{\alpha} \subset \varepsilon \left(N^4 (2\Gamma (\infty)) - \underset{1}{\overset{\infty}{\sum}} \, h_k \right)$ are, generally speaking, not disjoined from  each other. Moreover, we also find $\lambda_{\infty} \subset \varepsilon \left(N^4 (2\Gamma (\infty)) - \underset{1}{\overset{\infty}{\sum}} \, h_n \right)$.

\bigskip

\noindent {\bf The proofs of Lemmas 7, 8 and 9.} On the boundary of the 4-ball state $(i) = h_i \cup D^2 (C_i)$, we will consider the following two pieces
$$
{\mathcal B}_i^1 \equiv J_i \times D_i^2 \xhookrightarrow{ \qquad } \partial \{{\rm state} \ i \} \xleftrightarrow{ \qquad } {\mathcal B}_i^2 \equiv (\partial J_i \times B_i^*) \cup (I_i \times D_i^*),
\leqno (52)
$$
coming with ${\mathcal B}_i^1 \cap {\mathcal B}_i^2 = \partial J_i \times D_i^2 \subset \partial J_i \times B_i^*$, see here (32). Here, of course $\dim \, ({\rm state} \, (i)) = 4$, $\dim {\mathcal B}_i^1 = \dim \, {\mathcal B}_i^2 = 3$, $\dim \, ({\mathcal B}_i^1 \cap {\mathcal B}_i^2) = 2$. Moreover
$$
\partial ({\rm state} \, (i)) = {\mathcal B}_i^1 \cup {\mathcal B}_i^2 \cup \{ D_i \times \partial D_i^* , \ \mbox{the lateral surface of the $2$-handle}\}.
$$

$$
\includegraphics[width=6cm]{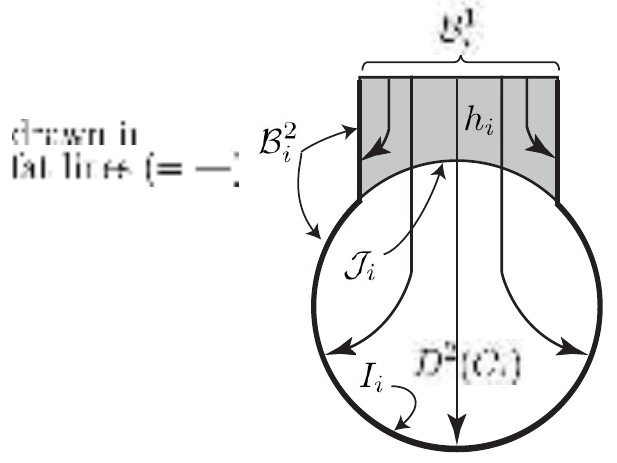}
$$
\label{fig14}
\centerline {\bf Figure 14.}
\begin{quote}
We illustrate here, with two dimension less, the formula (52). The arrows are in red.

With red lines, we have suggested the {\ibf flattenning retraction}
$$
{\mathcal B}_i^1 \subset \{{\rm state} \, i \} \twoheadrightarrow {\mathcal B}_i^2
$$
\vglue-8mm
$$
{\mid\mbox{\hglue -2mm}}_{-\!\!-\!\!-\!\!-\!\!-\!\!-\!\!-\!\!-\!\!-\!\!-\!\!-\!\!-\!\!-\!\!-\!\!-\!\!-\!\!-\!\!-}{\mbox{\hglue -2mm}\uparrow}
$$
\vglue-8mm
$$
\mbox{\footnotesize $\approx$}
$$
This same figure can be read as a profile view of the whole $\{{\rm state} \, i\} = h_i \cup D^2 (C_i)$, with $J_i$ and $\partial D_i = I_i \cup J_i$ displayed, and with profile views of the ${\mathcal B}_i^1 , {\mathcal B}_i^2$ too.
\end{quote}

\bigskip

With one dimension less, formula (52) can be visualized in Figure 14.

\smallskip

Notice that the dumbbell-like connected $3$-manifold from (52) and from the Figure 14, the
$$
{\mathcal B}_i^2 = \{\mbox{attaching zone of $h_i$}\} \cup \{\mbox{(attaching zone of $D^2 (C_i)$)} - \mbox{(lateral surface of $h_i$)}\},
$$
is exactly that part of $\partial ({\rm state} \, i)$ via which, in the absence of incoming arrows $\{ k \to i \} \in {\mathcal M}$, our $4$-cell ${\rm state} \, (i)$ is SPLIT from the rest of the world. And then, those arrows $k \to i$ hit our ${\rm state} \, (i)$ via the piece ${\mathcal B}_i^1$.

\smallskip

There exists a FLATTENING RETRACTION suggested in the Figure 14
$$
{\mathcal B}_i^1 \hookrightarrow \{{\rm state} \, i \} \overset{F=F(i)}{-\!\!\!-\!\!\!-\!\!\!-\!\!\!-\!\!\!-\!\!\!\twoheadrightarrow} {\mathcal B}_i^2 . \leqno (53)
$$
\vglue-7mm
$$
{\mid\mbox{\hglue -2mm}}_{\overset{\approx}{-\!\!-\!\!-\!\!-\!\!-\!\!-\!\!-\!\!-\!\!-\!\!-\!\!-\!\!-\!\!-\!\!-\!\!-\!\!-\!\!-\!\!-\!\!-\!\!-\!\!-\!\!-\!\!-\!\!-\!\!-\!\!-\!\!-\!\!-\!\!-\!\!-}}{\mbox{\hglue -2mm}\uparrow}
$$
\vglue-5mm
$$
\mbox{\footnotesize $F \mid {\mathcal B}_i^1 = {\rm diffeomorphism}$}
$$

Then, there is also a canonical isomorphism (see (32.1))
$$
{\mathcal B}_i^1 = J_i \times D_i^2 \overset{\eta_i}{\underset{\approx}{-\!\!\!-\!\!\!-\!\!\!\longrightarrow}} I_i \times D_i^2 = {\rm Box} \, (i) \quad \mbox{(see (33))} \leqno (54)
$$
and so, via (53), (54), ${\mathcal B}_i^1 , {\mathcal B}_i^2$ and ${\rm Box} \, (i)$ are three distinct diffeomorphic models of the same objet and, in what will follow next we will happily and freely move from one of these models to another, without bothering to change the notation. In this context, in our present easy id $+$ nilpotent context, it can always be assumed that
$$
\lim_{n=\infty} \, {\rm length} \, J_n = \infty \, , \quad \lim_{n = \infty} \, {\rm diam} \, D_n^n = 0, \quad {\rm for} \ n \to \infty \ {\rm in} \ {\mathcal M}. \leqno (55)
$$

It should be stressed that the easy id $+$ nilpotent condition is necessary  at this point. Once that condition is satisfied, we can impose consistently the {\ibf metric} requirement that when in the RED matrix $C \cdot h$ there is an arrow $i \to j$, then we also have that
$$
{\rm length} \ I_i > {\rm length} \ I_j \quad {\rm and} \quad {\rm diameter} \ D_i^2 < {\rm diameter} \ D_j^2  ,
$$
like in the Figure 13. Starting from that, it is not hard to impose the condition (55) too.

\smallskip

Without our easy id $+$ nil, like for instance in the case of the classical Whitehead manifold ${\rm Wh}^3$, these things are not possible. The difficult id $+$ nil does not allow such nice metric conditions. It should be an amusing exercise to figure out the (35) for the Whitehead manifold. With these things, let us go now to Lemma 7, where we fix a Box, the ${\rm Box} \, (i) = I_i \times D_i^2$, and consider all the infinite trajectories of our oriented graph ${\mathcal M}$ (easy id $+$ nil) incoming into $i$. They have the general form below
$$
i \equiv i_0 \leftarrow i_1 \leftarrow i_2 \leftarrow i_3 \leftarrow \ldots . \leqno (56)
$$

With this, we introduce now
$$
\bar\varepsilon \, X^3 (\alpha) \mid {\rm Box} \, (i) \equiv \sum_{\overbrace{\mbox{\footnotesize$i = i_0 \leftarrow i_1 \leftarrow i_2 \leftarrow \ldots$}}} \ \bigcap_{i_n} \ \nu (i_n \to i) (\lambda (i_n \to i) \, I_i \times D_{i_n}^2),
$$
and this formula should make the transversal structure of the lamination ${\mathcal L}_{\alpha}$ transparent. Next, of course, the lamination itself is defined by
$$
\bar\varepsilon \, X^3 (\alpha) = \bigcup_i  \left( \bar\varepsilon \, X^3 (\alpha) \mid {\rm Box} \, (i) \right).
$$
Our lamination is without holonomy, since ${\mathcal M}$ has no closed orbits.

\bigskip

\noindent {\bf Remark.} The easy id $+$ nil implies that our states $i$ may be labelled by positive integers $n \in Z_+$, s.t. if there is an arrow $i_n \to i_m$ in ${\mathcal M}$, then $n > m$. Of course, also, our oriented graph ${\mathcal M}$ is not necessarily a tree, but it certainly has no closed ORIENTED orbits. When the big ${\mathcal M} = \underset{\infty}{\sum} \, {\mathcal M}_{\alpha}$ is being considered, then we may identify
$$
\{\alpha\} \cong \{\mbox{the set of FINAL states, i.e. the states which have no outgoing arrows}\}.
$$

We go back now to a general situation when for our ${\mathcal M}$ we have chosen an order on the states $i_{\alpha}$, s.t. when there is an arrow $i_{\alpha} \to i_{\beta}$ then $i_{\alpha} > i_{\beta}$ (or, schematically, $\alpha > \beta$). We have
$$
\left(N^4 (2\Gamma (\infty)) - \sum_1^{\infty} h_k \right) \cup {\rm LAVA} = \left(N^4 (2\Gamma (\infty)) - \sum_1^{\infty} h_k \right) \cup \ \bigcup_{i_n} \, \{{\rm state} \ i_n\} \leqno (57)
$$
$$
\supset \sum_{i_n} {\mathcal B}_{i_n}^1 , \ \mbox{with this inclusion map being PROPER.}
$$
[The ${\mathcal B}_j^1$'s are well touched by the higher ${\mathcal B}_i^2$'s, but inside our $\left(N^4 (2\Gamma (\infty)) - \underset{1}{\overset{\infty}{\bigcup}} \, h_n \right) \cup {\rm LAVA}$, they are 2-by-2 disjoined. This gives us an injective map, occurring in (57), call it
$$
\sum_{i_n} {\mathcal B}_{i_n}^1 \underset{j}{\longrightarrow} {\rm LAVA}.] \leqno (57.1)
$$

Here we have
$$
\underbrace{\left[ \left(N^4 (2\Gamma (\infty)) - \sum_1^{\infty} h_k \right) \cup \bigcup_{n=1}^{N-1} \{{\rm state} \ i_n \} \right]}_{\mbox{\footnotesize we call this $X^4 (N-1)$}} \cap \ \{{\rm state} \ i_N\} = {\mathcal B}_{i_N}^2 ,
\leqno (58)
$$
and 
$$
\left(N^4 (2\Gamma (\infty)) - \sum_1^{\infty} h_k \right) \cap \ \{{\rm state} \ i_1\} = {\mathcal B}_{i_1}^2.
$$
We will denote
$$
\nu (i_N \to ) \equiv \sum_{{\rm all} \, j < i_N} \{\mbox{the source $\lambda (i_N \to j) \, I_j \times D_{i_N}^2$ (see Figure 13) of the map $\nu (i_N \to j)$ (34)}\} \subset {\mathcal B}_{i_N}^2.
$$
Here, inside our $\underset{j < i_N}{\sum}$ things are not disjoined, but we will not make this explicit.

\smallskip

When we consider the decomposition
$$
{\mathcal B}_{i_N}^2 = ({\mathcal B}_{i_N}^2 - \nu (i_N \to)) \cup \nu (i_N \to)
$$
then, in terms of (58), the ${\mathcal B}_{i_N}^2  - \nu (i_N \to)$ goes to $\partial \left(N^4 (2\Gamma (\infty)) - \underset{1}{\overset{\infty}{\sum}} \, h_k \right)$, while the $(\nu_{i_N} \to)$ goes to the $\left\{\mbox{still free part of} \ \underset{1}{\overset{N-1}{\sum}} \, {\mathcal B}_n^1 \right\} \subset \partial X^4 (N-1)$.

\smallskip

The following map, defined by composing the $F(i)$'s from (53)
$$
\{{\rm state} \ i_N\} \underset{}{\overset{\Phi (i_N) \equiv F(i_1) \circ F(i_2) \circ \ldots \circ F(i_N)}{-\!\!\!-\!\!\!-\!\!\!-\!\!\!-\!\!\!-\!\!\!-\!\!\!-\!\!\!-\!\!\!-\!\!\!-\!\!\!-\!\!\!-\!\!\!-\!\!\!-\!\!\!-\!\!\!-\!\!\!-\!\!\!-\!\!\!-\!\!\!\longrightarrow}} \partial \left(N^4 (2\Gamma (\infty)) - \sum_1^{\infty} h_n \right) , \leqno (59)
$$
is well-defined on ${\mathcal B}_{i_N}^2 \subset \{{\rm state} \ i_N\}$ which it finally flattens on $\partial \left(N^4 (2\Gamma (\infty)) - \underset{1}{\overset{\infty}{\sum}} \, h_n \right)$. Moreover, when we put together the various $\Phi (i_N)$'s then we get the following big, PROPER {\ibf flattening map}
$$
{\rm LAVA} \equiv \bigcup_N \ \{{\rm state} \ i_N\} \underset{\Phi (\infty) \, \equiv \, \underset{N}{\bigcup} \, \Phi (i_N)}{-\!\!\!-\!\!\!-\!\!\!-\!\!\!-\!\!\!-\!\!\!-\!\!\!-\!\!\!-\!\!\!-\!\!\!-\!\!\!-\!\!\!-\!\!\!-\!\!\!-\!\!\!\!\longrightarrow} \partial \left(N^4 (2\Gamma (\infty)) - \sum_1^{\infty} h_n \right) . \leqno (60)
$$
When $j < N$, then $\Phi (i_N)$ factorizes through $\Phi (i_j)$, and hence the two are compatible.

\smallskip

The big $\Phi (\infty)$ is a retraction of LAVA onto the $\delta \, {\rm LAVA} \subset \partial \left(N^4 (2\Gamma (\infty)) - \underset{1}{\overset{\infty}{\sum}} \, h_n \right)$, and we have $\Phi (\infty) \, {\rm LAVA} = \delta \, {\rm LAVA}$.

\newpage
\vglue-10mm
$$
\includegraphics[width=125mm]{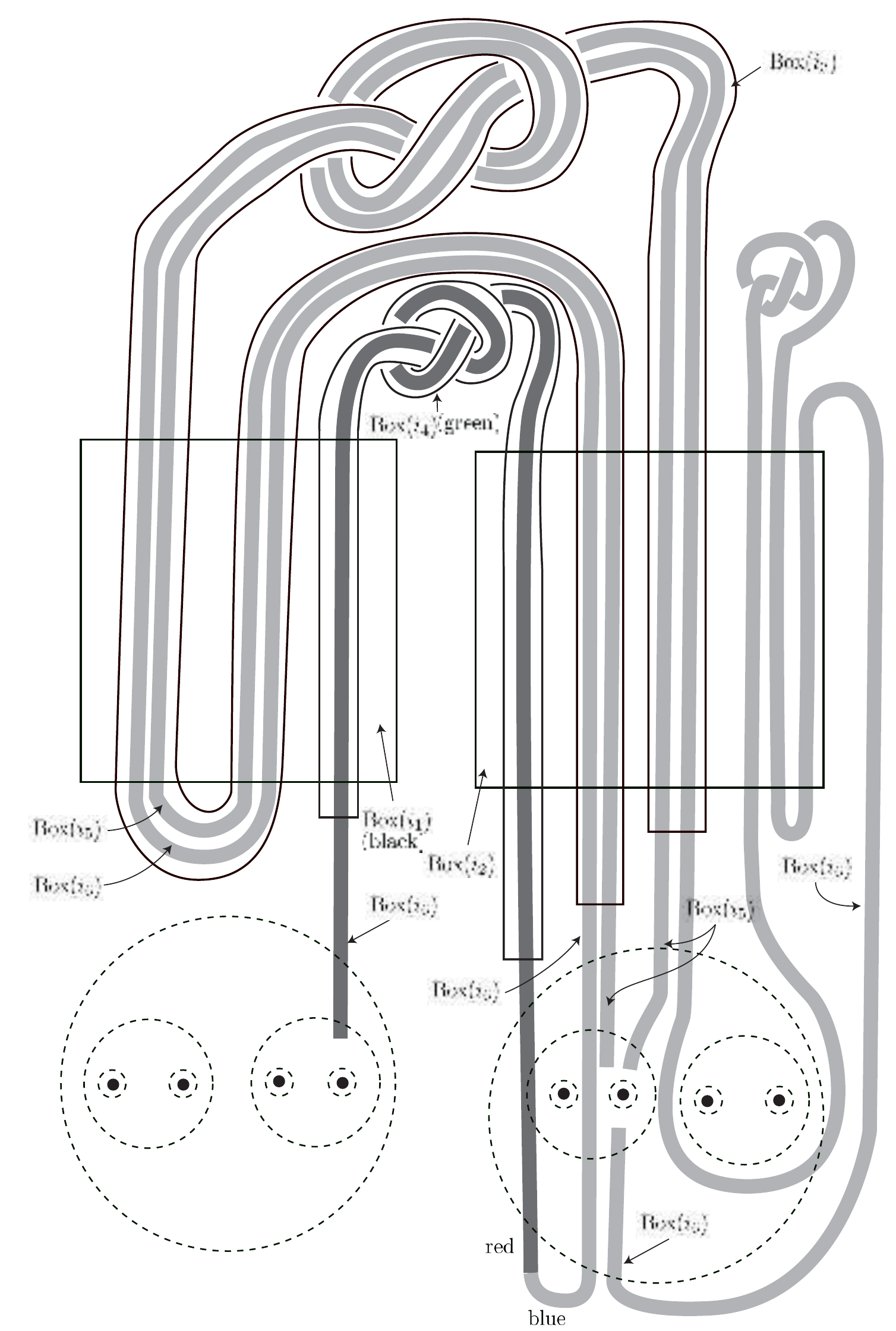}
$$
\label{fig15}
\centerline {\bf Figure 15.} 
\begin{quote}
This figure should illustrate the PROPER embedding $J$ from (41). 

The plane of our figure is supposed to be
$$
\partial \left(N^4 (2\Gamma (\infty)) - \sum_1^{\infty} h_n \right).
$$
The fat points $(= \bullet)$ are in $\varepsilon \left(N^4 (2\Gamma (\infty)) - \underset{1}{\overset{\infty}{\sum}} \, h_n \right)$. One should notice that our ${\mathcal M}$ is generally speaking, NOT a tree, it can have closed non-ordered cycles (not orbits or trajectories). Notice here the occurrence of the following closed curve, contained in ${\mathcal M}$, but which is NOT a cycle
$$
\xymatrix{
i_6 \ar[d]\ar[r] &i_4 \ar[d] \\
i_3 \ar[r] &i_1
}
$$
and the graph becomes more complicated when we introduce the $i_2 , i_5$ too.

LEGEND: ${\rm BLACK} = i_1$ and $i_3$, ${\rm GREEN} = i_4$, ${\rm RED \ and \ BLUE} = i_6$.
\end{quote}

\bigskip

We have a commutative diagram
$$
\xymatrix{
{\rm LAVA} \ar[rdd]^{\Phi (\infty)} &&&\ar@{_{(}->}[lll]^{j \ (57.1)} \ \underset{i_n}{\sum} \, {\mathcal B}_{i_n}^1 \cong  \underset{n}{\sum} \, {\rm Box} \, (i_n) \ar[ld]^{{{\rm \ canonical \, quotient} \atop {\rm map \, projection;}} \atop {\rm and \, here, \, see \, (35).}} \\
&&\underset{\alpha}{\sum} \, X^3 (\alpha) \ar[ld]^{ \ \ J \, ({\rm our \, PROPER} \atop {\rm map \, {\it J} \, (41))}} \\
&\partial  \left(N^4 (2\Gamma (\infty)) - \underset{1}{\overset{\infty}{\sum}} \, h_k \right)
}
\leqno (61)
$$
and here ${\rm Im} \, \Phi (\infty) = {\rm Im} \,  J = \delta \, {\rm LAVA}$.

\smallskip

Figure 15 should give an impressionistic idea of what the map $J$ may look like. Inside each given Box we only see many parallel strands, telescopically embedded inside each other according to the prescription (see the beginning of Lemma 6)
$$
D_{\ell}^2 \xrightarrow{\ \mu (\ell \to m) \ } D_m^2 ,
$$
which is an embedding. Of course, the $J$-image of each individual box is no longer multilinear, like it is in the abstract model.

\smallskip

Very importantly, {\ibf no knotting or linking} ever occurs, inside any given individual box. But, generally speaking, the global $X^3 (\alpha)$ is highly non-simply-connected and knotted too. This $\pi_1 \ne 0$ issue stems from the fact that ${\mathcal M}$ is not necessarily a tree and also from the multiplicities of the given contacts $i \to j$.

\smallskip

Finally, inside the ambient space $\partial \left(N^4 (2\Gamma (\infty)) - \underset{1}{\overset{\infty}{\sum}} \, h_k \right)$ the $J \, \underset{\alpha}{\sum} \, X^3 (\alpha) \approx \underset{\alpha}{\sum} \, X^3 (\alpha)$ is, generally speaking, highly knotted and linked. All these things having been said, one proves (43) by letting each ${\rm LAVA}_{\alpha}$ grow naturally and very slowly, so as to respect our various smoothness conditions, out of its seed $J X^3 (\alpha) \approx X^3 (\alpha)$. The closer it comes to $\varepsilon \left( \partial N^4 (2\Gamma (\infty)) - \underset{1}{\overset{\infty}{\sum}} \, h_n \right)$ the slower LAVA is supposed to grow and, also, the height to which it grows is becoming smaller and smaller. This is suggested in the Figure 16.

\smallskip

The growth process was suggested here as it is happening inside the ambient space $N_1^4 (2X_0^2)^{\wedge} = N^4 (2X_0^2)^{\wedge} \cup (\partial N^4 (2X_0^2)^{\wedge} \times [0,1))$. But then there is also a more intrinsic description. Now, there is no limit for the growth height. But on some pieces the growth process stops at some finite height; the growth there stops at $X^3 (\alpha) \times \{1\} - \bar\varepsilon \, (X^3 (\alpha))$. Then, on the other pieces, the growth continues indefinitely, up to an infinite height, and this defines the
$$
\bar\varepsilon \, (X^3 (\alpha)) \subset X^3 (\alpha) \times \{1\}.
$$
Finally, starting with the diffeomorphism (40.1), which we re-write here:
$$
\left(N^4 (2\Gamma (\infty)) - \sum_1^{\infty} h_k \right)^{\wedge} + \sum_1^{\bar n} D^2 (\Gamma_j) \underset{\rm DIFF}{=} N^4 (\Delta^2) ,
$$
one proves the more serious (51) just like one proves (43), working now with the full $\underset{\alpha}{\sum} \, X^3 (\alpha)$. An impressionistic view of these things is presented in the Figure 16. And remember here that $({\rm LAVA}) \cap \sum \Gamma_j = \emptyset$.

$$
\includegraphics[width=12cm]{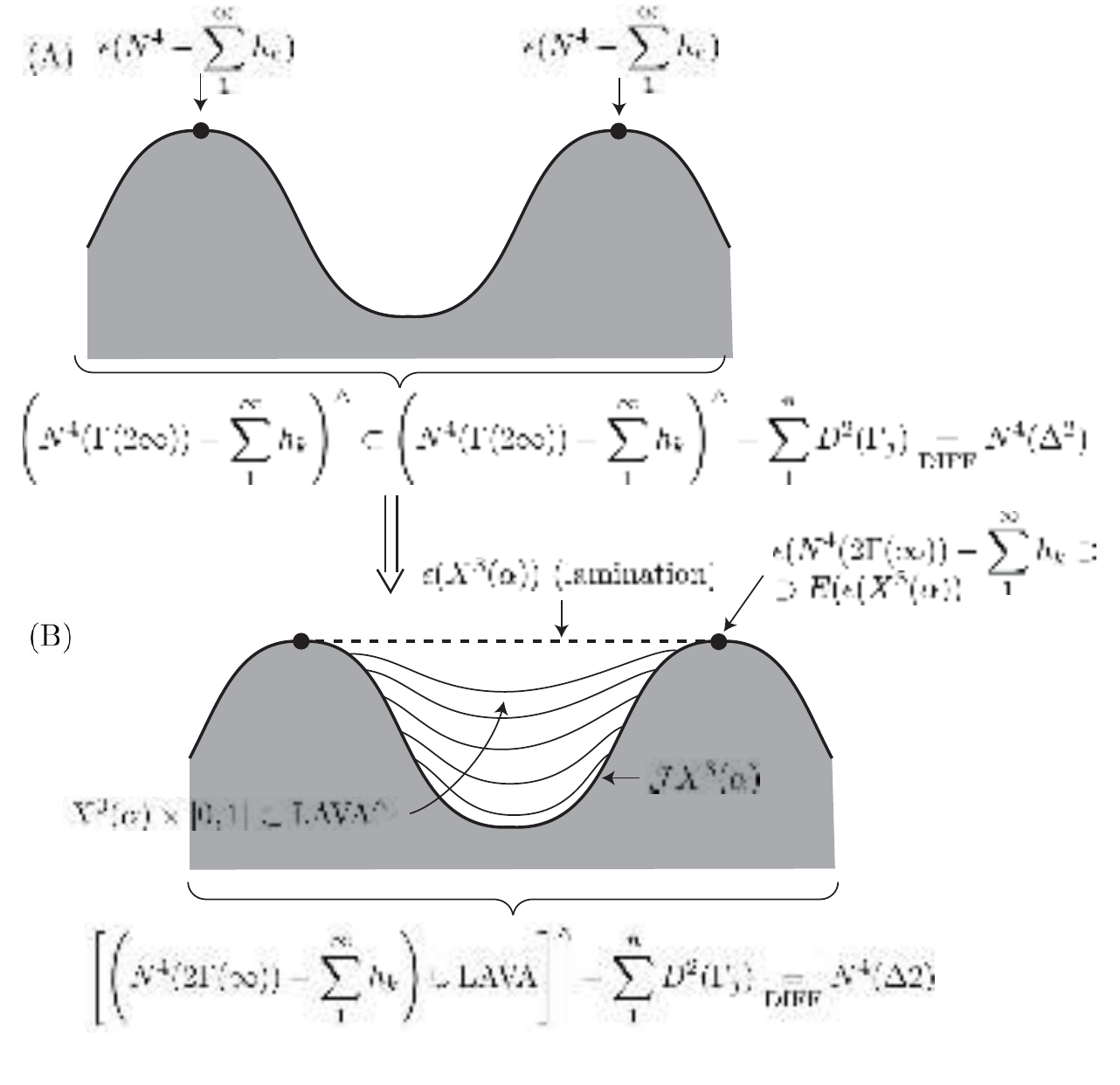}
$$
\label{fig16}
\centerline {\bf Figure 16.} 
\begin{quote}
An impressionistic view of the proof of (58). The step (A) $\Rightarrow$ (B) is a diffeomorphism, where  a whole orchestra of $X^2 (\alpha)$'s grows simultaneously, in a smoothly controlled manner.
\end{quote}

\bigskip

We will move now to the following lighter notation
$$
(L_{\alpha} , \delta L_{\alpha}) = ({\rm LAVA}_{\alpha} , \delta \, {\rm LAVA}_{\alpha})
$$
with which, our previous formulae become
$$
(L_{\alpha} , \delta L_{\alpha}) \underset{\rm DIFF}{=} ( \delta L_{\alpha} \times [0,1] - \bar\varepsilon \,  \delta L_{\alpha} \times \{1\} ,  \delta L_{\alpha} \times \{0\}),
$$
$$
L_{\alpha}^{\wedge} = ( \delta L_{\alpha} \times [0,1]) \underset{\overbrace{\mbox{\footnotesize$\delta L_{\alpha}$}}}{\cup} ( \delta L_{\alpha})^{\wedge} =  (\delta L_{\alpha} \times [0,1]) \cup \lambda_{\alpha} .
$$

\bigskip

\noindent {\bf A complement to the Lemmas 8 and 9; more on the product property of LAVA.}

\medskip

1) {\it There is a PROPER Whitehead collapse
$$
L_{\alpha}^{\wedge} \supset  \delta L_{\alpha} \times [0,1] \xrightarrow{ \ \hat\pi_{\alpha} \ }  \delta L_{\alpha} .
\leqno (62)
$$

\noindent [COMMENT. It is the PROPER collapse, gotten from the $(62)$ above
$$
\left( N^4 (2\Gamma (\alpha)) - \sum h_n \right) \cup \sum_{\alpha} \delta L_{\alpha} \times [0,1] \underset{\underset{\alpha}{\sum} \ \hat\pi_{\alpha}}{-\!\!\!-\!\!\!-\!\!\!-\!\!\!-\!\!\!\!\longrightarrow} \left( N^4 (2\Gamma (\alpha)) - \sum h_n \right),
$$
where $\left( N^4 - \sum h_n) \right) \cup \underset{\alpha}{\sum} \, \delta L_{\alpha} \times [0,1] \subset \left( \left( N^4 - \sum h_n \right) \cup {\rm LAVA} \right)^{\wedge}$, which is behind $(50)$ and $(51)$. This collapse is supposed to proceed {\ibf smoothly}, with smaller and smaller amplitude as we approach the $\varepsilon \Bigl( N^4 (2\Gamma (\infty))$ $- \, \underset{1}{\overset{\infty}{\sum}} \, h_n \Bigl)$, so as to yield a {\ibf diffeomorphism}
$$
\left[ \left(N^4 (2\Gamma (\infty)) - \sum_1^{\infty} h_n \right) \cup {\rm LAVA} \right]^{\wedge} \underset{\rm DIFF}{=} \left(N^4 (2\Gamma (\infty)) - \sum_1^{\infty} h_n \right)^{\wedge} .
$$
The $\underset{\alpha}{\sum} \, \hat\pi_{\alpha}$ above is the inverse of the process via which the $\underset{\alpha}{\sum} \, {\rm LAVA}_{\alpha}$ grows naturally and very slowly out of $\underset{\alpha}{\sum} \, J \times X^3 (\alpha)$. Also, we will denote by $\pi_{\alpha}$ the restriction of $\hat\pi_{\alpha}$ to $L_{\alpha} \subset L_{\alpha}^{\wedge}$.]}

\medskip

2) {\it The $(62)$ will be called the {\ibf strong} product property of LAVA, so as to distinguished from the {\ibf weak} product properties to be developed next, and which is consequence of the strong property. Let us start with a compact bounded, not necessarily connected surface $(S,\partial S)$ coming with an inclusion
$$
(S , \partial S) \xhookrightarrow{ \ \ i \ \ } (\delta L_{\alpha} , \partial \, \delta L_{\alpha}). \leqno (63)
$$
In principle, at least, we will always think of $(S , \partial S)$ as consisting of connected  components, each a copy of $(D^2 , \partial D^2 = S^1)$.

\smallskip

This induces embeddings
$$
\xymatrix{
\left( \hat\pi_{\alpha}^{-1} S , \partial S \times [0,1] \right) \ar[r] &\left( \widehat L_{\alpha} , (\partial \, \delta L_{\alpha}) \times [0,1] \right) \\
\left( \pi_{\alpha}^{-1} S , \partial S \times [0,1] \right) \ar[u] \ar[r] &\left( L_{\alpha} , (\partial \, \delta L_{\alpha}) \times [0,1] \right), \ar[u]
}
$$
with $\hat\pi_{\alpha}^{-1} S = S \times [0,1]$, with $\hat\pi_{\alpha}^{-1} S \cap \bar\varepsilon \, L_{\alpha}$ a {\ibf tame} totally discontinuous (Cantor) subset of $S \times \{1\} \subset \hat\pi_{\alpha}^{-1} S$ and with $\pi_{\alpha}^{-1} S = \hat\pi_{\alpha}^{-1} S - \hat\pi_{\alpha}^{-1} S \cap \bar\varepsilon \, L_{\alpha}$.}

\medskip

3) {\it In this same vein, let 
$$
x \in \Gamma (2\infty) \cup \sum_1^{\infty} D^2 (C_i)
$$
be like in $(28.4)$. Then, in terms of the things just said
$$
N^4 (\mbox{extended cocore of $(x)$}) = \sum_{\alpha} \pi_{\alpha}^{-1} (x) \leqno (64)
$$
and this is a copy of $B^3 = \{\mbox{a {\ibf tame} Cantor set $\subset \partial B^3$}\}$. From now on, we simplify the notation and, unless the opposite is explicitly said, the $(64)$ above will be denoted just by
$$
\{\mbox{Extended cocore $(x)$}\} \subset \left(N^4 (2\Gamma (\infty)) - \sum_1^{\infty} h_n \right) \cup {\rm LAVA} , \leqno (65)
$$
and notice here the capital ``$E$''. The $(65)$ is a PROPER codimension one submanifold, which can be compactified into
$$
\{\mbox{Extended cocore $(x)$}\}^{\wedge} \equiv \sum_{\alpha} \hat\pi_{\alpha}^{-1} (x) \subset \left[ \left(N^4 (2\Gamma (\infty)) - \sum_1^{\infty} h_n \right) \cup {\rm LAVA} \right]^{\wedge} .
$$
There is a diffeomorphism $\{\mbox{Extended cocore $(x)$}\}^{\wedge} \underset{\rm DIFF}{=} B^3$.}

\medskip

4) {\it We have a diffeomorphism of pairs
$$
\left(\left[ \left(N^4 (2\Gamma (\infty)) - \sum_1^{\infty} h_n \right) \cup {\rm LAVA} \right]^{\wedge} , \sum_1^n \{\mbox{Extended cocore $(R_i)$}\}^{\wedge} \right) \leqno (66)
$$ 
$$
\underset{\rm DIFF}{=} n \, \# \, ((S^1 \times B^3) , (*) \times B^3) \ (\mbox{standard}).
$$
}

\bigskip

Here, the $\underset{1}{\overset{n}{\sum}} \, R_i \subset \Gamma (1)$ correspond to the RED $1$-handles of $N^4 (\Delta^2)$.

\smallskip

In what follows next, in the present section we will make more explicit the structure of LAVA. The next items will not be used explicitly, later on, in the proof of our main result, but they may clarify things.

\smallskip

We will call {\ibf elementary} a matrix ${\mathcal M}_0$, \`a la (30), which is both of the easy id $+$ nilpotent type and which has a unique minimal (i.e. final) state $i_1 = h_{i_1} \cup B_{i_1}$. We will assume ${\mathcal M}_0$ to be a {\ibf full submatrix} of our $C \cdot h =$ id $+$ nilpotent, meaning that no state in ${\mathcal M}_0$ receives arrows from outside ${\mathcal M}_0$, when we consider ${\mathcal M}_0$ as a piece of $C \cdot h$. There is a contribution of ${\mathcal M}_0$ to LAVA, a ${\rm LAVA} \, ({\mathcal M}_0)$ or ${\rm LAVA}^{\wedge} \, ({\mathcal M}_0)$ and we will refer to these as being an {\ibf atom}, respectively a {\ibf compact atom}. We will use the following generic notation for the successive states of ${\mathcal M}_0$
$$
i_1 (\mbox{unique minimal state}) \leftarrow i_2 (+ i'_2 + i''_2 + \ldots) \leftarrow i_3 (+ \ldots) \leftarrow \ldots
$$
Let now $S \subset \partial (i_1)$ be a SPOT looking like one of the ${\rm Im} \, (\nu (i \to j_k))$ in Figure 12. Of course our present $S$ cannot be a ${\rm Im} \, (\nu (i_1 \to \ldots))$, since the state $i_1$ is minimal.

\smallskip

At the level of $2X_0^2$ the $\pi ({\mathcal M}_0)^{-1} S$ is a full tree, splitting locally $2X_0^2$, except at its foot $S$. When one goes $4^{\rm d}$, then at the level of $N^4 (2X_0^2)$ our $\pi ({\mathcal M}_0)^{-1} S$ becomes an object which I will denote by $M^3 [S] \underset{\rm DIFF}{=} B^3 - \{$a tame Cantor set $\subset \partial B^3\}$, organized as follows
$$
M^3 [S] = X_1 \, \# \, X_2 \, \# \, X_3 \, \# \, \ldots , \leqno (67)
$$
where $X_1 = \{$the contribution $B^3 (i_1)$ of $i_1\}$, $X_2 = \{$the contribution of $i_2\}$, which is a finite collection of $B^3 (i_2)$'s, with $\#$ a multiple connected sum, a.s.o. Let us say that $X_n = \underset{j=1}{\overset{N(n)}{\sum}} B_j^3 (i_n)$. The Figure 17 suggests the structure of
$$
\left([M^3 (S)]^{\wedge} , [M^3 (S)] \right) = \left( \hat\pi ({\mathcal M}_0)^{-1} (S) , \pi ({\mathcal M}_0)^{-1} (S) \right)
$$
coming with
$$
F_S = [M^3 (S)]^{\wedge} - [M^3 (S)],
$$
which is a Cantor set suggested in Figure 17, by fat points.

\smallskip

Because of typographical reasons, in the Figure 17 we did not represent realistically the organic structure of $[M^3 (S)]^{\wedge}$, which we draw as a square.

\smallskip

For the triplet
$$
\mbox{(compact atom; $\delta$ (atom), atomic lamination)},
$$
where $\delta \, ({\rm atom})$ splits the compact atom from the rest of the universe, we have the explicit description below
$$
\bigl([M^3 (S)]^{\wedge} \times [0,1] ; \left([M^3 (S)] \times \{0\} \right) \cup (S \times [0,1]) \cup \left( [M^3 (S)] \times \{1\} \right), \leqno (68)
$$
$$
F_S \times [0,1]\bigl) = \left({\rm LAVA} ({\mathcal M}_0)^{\wedge} ; \delta \, {\rm LAVA} ({\mathcal M}_0) , \bar\varepsilon ({\mathcal M}_0)\right). 
$$

We move now to a completely general geometric intersection matrix  (30) which is of the easy id $+$ nilpotent type, without being necessarily elementary. We will show now how the ${\rm LAVA}_{\alpha}$ and, more precisely the triple
$$
({\rm LAVA}_{\alpha}^{\wedge} \, ; \ \delta \, {\rm LAVA}_{\alpha} \, , \ \varepsilon (X^3 (\alpha)))
$$
can be reconstructed by glueing together in a locally finite but not necessarily simply-connected manners, compact atoms of type (68).

\smallskip

We start by introducing, on the same lines as $X^3 (\alpha)$ in (35), the following object
$$
X^1 (\alpha) = \left(\sum_{i \, \in \, \alpha} I_i \right) \diagup {\mathcal R}^1 \ (=\{\mbox{the equivalence relation induced by the various $\lambda (i\to j)$}\}).
\leqno (69)
$$
Topologically speaking, this is an infinite, connected, locally finite, highly non-simply-connected graph. But its natural structure is not the standard edge/vertex structure; instead we have a richer kind of structure which we will describe. Roughly speaking $X^1 (\alpha)$ is a {\ibf train-track} with infinite weights, recording the way in which our various $I_i$'s go through each other. Each $I_i$ is here a {\ibf smooth} arc and the maps $\lambda (i\to j)$ are smooth too.
$$
\includegraphics[width=9cm]{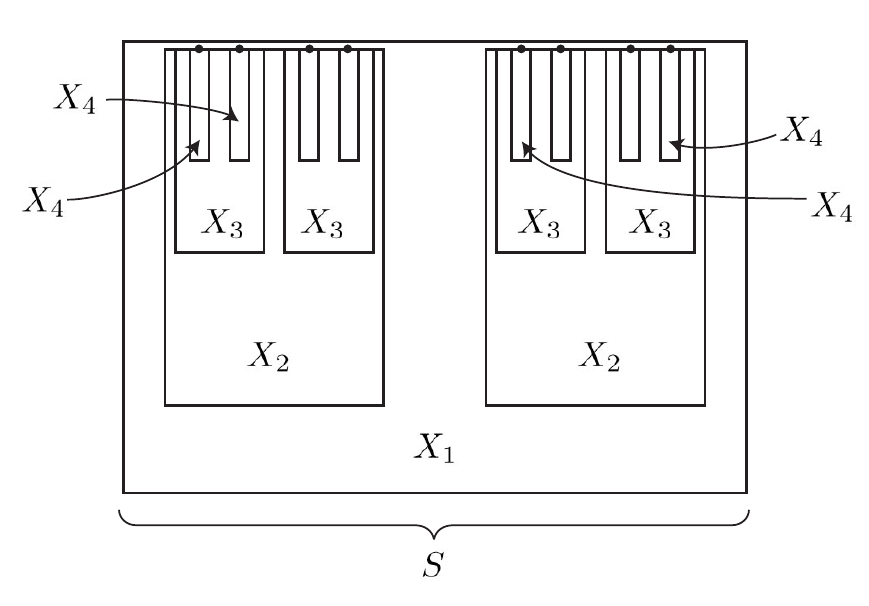}
$$
\label{fig17}
\centerline {\bf Figure 17.} 
\begin{quote}
With one dimension less, this figure should suggest the 3-ball 
$$
[M^3 (S)]^{\wedge} = X_1 \, \# \, X_2 \, \# \, \ldots \, \# \, X_n \, \# \, \ldots \, \cup \{\mbox{Cantor set $F_S$}\},
$$
with $S \subset \partial X_i$. LEGEND: $\bullet =$ points in $F_S$.
\end{quote}

\bigskip

We have obvious maps
$$
X^3 (\alpha) \overset{P}{\underset{{\mbox{\footnotesize ``quotient-space}} \atop \mbox{\footnotesize projection''}}{-\!\!\!-\!\!\!-\!\!\!-\!\!\!-\!\!\!-\!\!\!-\!\!\!-\!\!\!-\!\!\!-\!\!\!-\!\!\!-\!\!\!-\!\!\!-\!\!\!\twoheadrightarrow}} X^1 (\alpha) \underset{{\mbox{\footnotesize inclusion of a}} \atop \mbox{\footnotesize DISCRETE subset}}{\xhookleftarrow{\qquad\qquad\qquad\qquad}} X^0 (\alpha) \underset{\rm def}{=} \sum_i \partial I_i . \leqno (70)
$$

All the {\ibf branching points} of $X^1 (\alpha)$ i.e. the points where, locally speaking $X^1 (\alpha)$ fails to be a smooth line, occur at the $X^0 (\alpha)$ points. For any $p \in X^1 (\alpha) - X^0 (\alpha)$ we can define an elementary ${\mathcal M}_0 (p)$, which should be thought of as the {\ibf weight} of $p$, as follows.

\bigskip

\noindent (71) \quad We will start by introducing $J(p) \equiv \{$the set of $I_i$'s which are such that $p \in I_i$, counted with {\ibf multiplicities}; any given $I_i$ may go finitely many times through the point $p\}$.

\smallskip

By definition, this will be the set of states of ${\mathcal M}_0 (p)$. Then, for $I_i , I_j \in J(p)$, any occurrence $p \in I_j \subsetneqq I_i$, without intermediary $p \in I_j \subsetneqq I_k \subsetneqq I_i$, will be, by definition, an arrow $\{ i \to j \} \in {\mathcal M}_0 (p)$. End of (71).

\bigskip

\noindent {\bf Remarks.} a) Notice that the $I_i$'s, states of $J(p)$, are occurrences $p \in I_i$, not just $I_i$'s. Consider then two distinct occurrences $p \in I_j \subset \overset{\circ}{I}_i$, $p \in I_k \subset \overset{\circ}{I}_i$. We may have here $j=k$ or $j \ne k$. If $j=k$, then these occurrences $p \in I_j$, $p \in I_j$ are distinct states of $J(p)$ (hence so are the $p \in I_i$, $p \in I_i$) and we have two disjoined edges in ${\mathcal M}_0 (p)$. If $j \ne k$, one can show that we have two successive arrows in ${\mathcal M}_0 (p)$. It follows from these things that {\ibf a state in $J(p)$ has at most one outgoing arrow}.

\medskip

b) To each smooth point $p \in X^1 (\alpha)$ as above, corresponds a properly embedded disk $S(p)\subset X^3 (\alpha) = \delta \, {\rm LAVA}_{\alpha}$ and, if we consider
$$
{\rm LAVA}_{\alpha} \xrightarrow{ \ \pi_{\alpha} \ } \delta \, {\rm LAVA}_{\alpha} = X^3 (\alpha) \xrightarrow{ \ P \ } X^1 (\alpha),
$$
then $\pi_{\alpha}^{-1} S(p) = \{$the $\pi ({\mathcal M}_0)^{-1} S(p)$ from (67)$\} \equiv [M^3 (p)]$. This fiber stays constant along connected components of $X^1 (\alpha) - X^0 (\alpha)$ but it {\ibf jumps} at points in $X^0 (\alpha)$. Roughly speaking, {\ibf this is the correct description of} ${\rm LAVA}_{\alpha}$; it will be rendered more precise, below. The typical structure of the train-track $X^1 (\alpha)$ in the neighbourhood of a point $q \in X^0 (\alpha)$ is suggested in the Figure 18.

$$
\includegraphics[width=7cm]{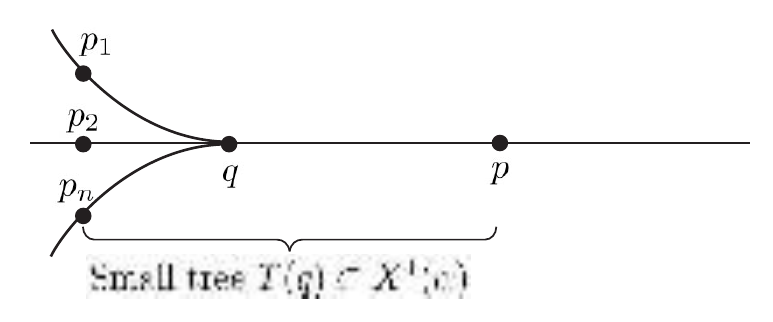}
$$
\label{fig18}
\centerline {\bf Figure 18.} 
\begin{quote}
Here $n=n(q) < \infty$. We see a small piece of $X^1 (\alpha)$, $q \in X^0 (\alpha)$ and $p,p_i$'s are smooth points.
\end{quote}

\bigskip

Together with each $q \in X^0 (\alpha)$ (see Figure 18) comes a smooth PROPER embedding
$$
\sum_1^n [M^3 (p_i)] \xhookrightarrow{ \ \xi (q) \ } [M^3 (p)]
$$
which is such that ${\rm Im} \, \xi (q) = \{\{ X_2 \cup X_3 \cup \ldots \}$ of $p$, (see (67)$\}$, while $\xi (q) (X_j (p_i)) \subset X_{j+1} (p)$, for each $1 \leq i \leq n(q)$.

\smallskip

In terms of the Figure 12, which in real life is of course supposed to be 4-dimensional, the purely $1^{\rm d}$ Figure 18 corresponds to the Figure 19.

$$
\includegraphics[width=110mm]{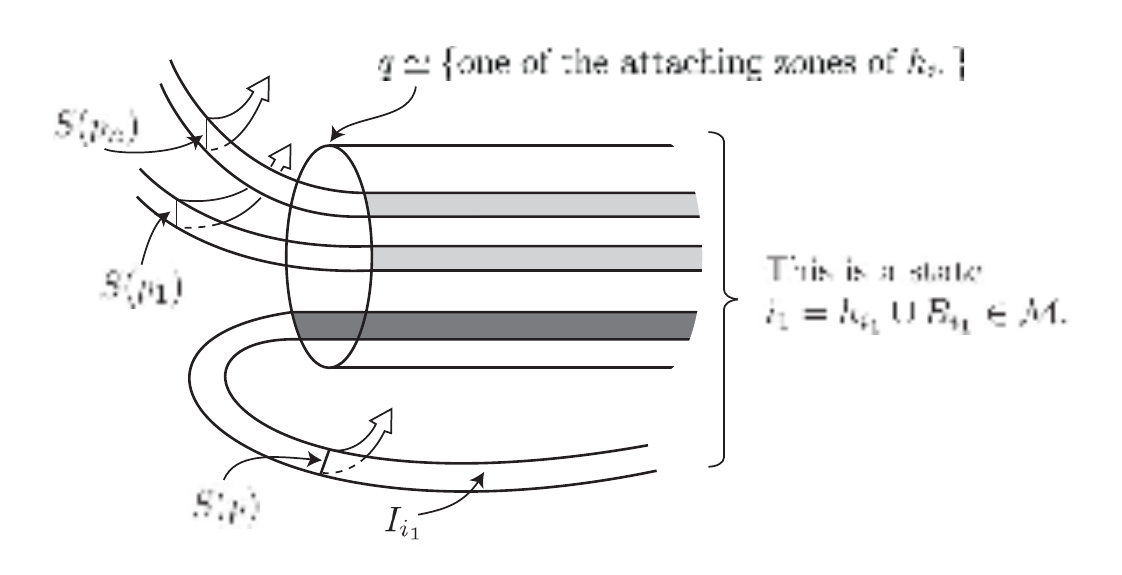}
$$
\label{fig19}
\centerline {\bf Figure 19.} 
\begin{quote}
This is what Figure 18 (which is in $X^1 (\alpha)$) corresponds to, at the level of $X^3 (\alpha)$ and/or ${\rm LAVA}_{\alpha}$. Here $p,p_1,\ldots , p_n$ are smooth points of $X^1 (\alpha)$ (i.e. points of $X^1 (\alpha) - X^0 (\alpha)$). At the level of $X^1 (\alpha)$ there is a unique occurrence $p \in I_{i_1}$ and this is the unique minimal state in $J(p) \subset {\mathcal M}_0 (p)$. We also have
$$
[M^3 (p)] = \{\mbox{extended cocore} \ S(p)\}.
$$
Notice that here $S(p) + S(p_1) + \ldots + S(p_n)$ are properly embedded inside $\delta \, {\rm LAVA}_{\alpha}$. Apart from small scallops coming from the shaded areas, the attaching zone $\partial J_{i_1} \times B^*_{i_1}$ (see (31)) of the 1-handle $h_{i_1}$ is also in $\delta \, {\rm LAVA}_{\alpha}$. So, the discs $S(p) + S(p_0) + \ldots + S(p_n)$ split out of $\delta \, {\rm LAVA}_{\alpha}$ a 3-ball $B^3 (q)$.
\end{quote}

\bigskip

Here is the description of the $(P \circ \pi_{\alpha})^{-1} \, T(q) \subset {\rm LAVA}_{\alpha}$, for a neighbourhood $T(q)$ of $q \in T(q) \subset X^1 (\alpha)$, like in the Figure 18
$$
(P \circ \pi_{\alpha})^{-1} (T(q)) = \sum_{i=1}^{n(q)} [M^3 (p_i)] \times [1,0] \cup [M^3 (p)] \times [0,1], \leqno (72)
$$
where the two pieces are glued together along
$$
\sum_i [M^3 (p_i)] \times \{0\} \approx {\rm Im} \, \zeta (q) \subset [M^3 (p)] \times \{0\}
$$
and where
$$
\sum_i [M^3 (p_i)] \times \{1\} + [M^3 (p)] \times \{1 \}
$$
splits the $(P \circ \pi_{\alpha})^{-1} \, T(q)$ from the rest of the ${\rm LAVA}_{\alpha}$.

\smallskip

Let $m_1 , m_2 , \ldots$ be the minimal states of ${\mathcal M}$. For each $I_{m_i}$ we consider some $I'_{m_i} \subset {\rm int} \, I_{m_i}$. For the generic $p \in {\rm int} \, I_{m_i}$ we have then
$$
(P \circ \pi_{\alpha})^{-1} I'_{m_i} = [M^3 (p)] \times I'_{m_i}, \ \mbox{like in (68)}. \leqno (73)
$$
One can find a family of smooth points ${\mathcal P} \subset X^1 (\alpha)$, s.t.
$$
\{ X^1 (\alpha) \ \mbox{SPLIT along} \ {\mathcal P}\} = \sum_{q \in X^0 (\alpha)} T(q) + \sum_{m_i} I'_{m_i},
$$
yielding

\smallskip

\noindent (74)
\vglue-5mm
\begin{eqnarray}
{\rm LAVA}_{\alpha}^{\vee}& \equiv& \{{\rm LAVA}_{\alpha} \ \mbox{SPLIT along} \ \sum_{p \in {\mathcal P}} \{\mbox{extended cocore} \ (p)\}\} \nonumber \\
&= &\sum_{q \in X^1 (\alpha)} \{ (P \circ \pi_{\alpha})^{-1} (T(q)\} + \sum_{m_i} \{ (P \circ \pi_{\alpha})^{-1} I_{m_i}\}. \nonumber
\end{eqnarray}
Next, glueing back together, in the obvious way, the spare pieces in (74) we finally get our ${\rm LAVA}_{\alpha}$. One can also extract from here a more or less explicit of
$$
{\mathcal L}_{\alpha} \subset {\rm LAVA}_{\alpha}^{\wedge} \supset {\rm LAVA}_{\alpha} \supset \delta \, {\rm LAVA}_{\alpha}.
$$

\section{Constructing exterior discs, which are disjointly embedded}\label{sec4}

We return now to $X^3 ({\rm old}) = \{\mbox{3-skeleton of} \ {\rm int} \, \Delta_1^4\}$ and to its RED $3^{\rm d}$ collapse $X^3 ({\rm old}) \searrow \Delta^3 \equiv \{\mbox{3-spine of}$ $\Delta^4_{\rm Schoenflies}\}$. This collapse kills teh 2-cells $D^2 (\gamma_k^0) \subset X^2 ({\rm old})-\Delta^2 (= 2 \ \mbox{skeleton of} \ \Delta^4)$ (with $\Delta^2 \subset \Delta^3$). More explicitly for every $\gamma_k^0 \subset \{{\rm link}\}$ (9) we have a  simplicial nondegenerate map, with a collapsible source,
$$
B^3 (\gamma_k^0) \xrightarrow{ \ \ F_k \ \ } \overline{X^3 - \Delta^3} , \leqno (75)
$$
with the following features. We have $D^2 (\gamma_k^0) \subset \partial B^3 (\gamma_k^0)$ and we introduce the 2-cell, coming with its non-degenerate map
$$
d_k^2 \equiv \partial B^3 (\gamma_k^0) - {\rm int} \, D^2 (\gamma_k^0) \xrightarrow{ \ f_k \, \equiv \, F_k \mid d_k^2 \ } X^2 ({\rm old}).
$$

With this, our $F_k$ (75) is injective except for possible edge-effects, meaning double points of the $f_k$ above.

\smallskip

Notice that $\gamma_k^0 = \partial d_k^2$. Also, the $F_k (B^3 (\gamma_k^0))$'s are disjoined, except for the following two items:
\begin{enumerate}
\item[$\bullet$)] edge effects,
\item[$\bullet$$\bullet$)] telescopic embeddings $F_p (B^3 (\gamma_p^0)) \subset F_k (B^3 (\gamma_k^0))$ induced by $\gamma_p^0 \subset d_k^2 - \gamma_k^0$.
\end{enumerate}

All this, so far, was in the context of $X^2 ({\rm old})$ and we move now to $X^2 [{\rm new}]$ forgetting about the $3^{\rm d}$ collapse. Every time for our $d_k^2$ above we find $D^2 (\Gamma_j) \times (\xi_0 = 0) \subset f_k \, d^2_k$, we replace it by $\Gamma_j \times [0 \geq \xi_0 \geq -1] \cup (D^2 (\Gamma_j) \times (\xi_0 = -1))$ getting this way a {\ibf new} $(d_k^2 , f_k)$, the only one to be used from now on, and we will not bother to change the notations
$$
\xymatrix{
d_k^2 \ar[rr]^{f_k} \ar[dr] &&2X_0^2 \\
&X_0^2 [{\rm new}] = X_0^2 \times r. \ar[ur]
}
\leqno (76)
$$
One should notice that when we go to the context $X_0^2 [{\rm new}]$, then to the $\gamma_k^0 \subset X^2 ({\rm old})$, we have to add the curves $\Gamma_j \times (\xi_0 = 0)$, $\Gamma_j \subset \Delta^2$.

\smallskip

We will work, in what follows from now in this paper, in the following context which supersedes the (2), for the time being
$$
\xymatrix{
N^4 (\Delta^2) &\!\!\!\!\!\!\!\!\!\!\!\subset &\!\!\!\!\!\!\!\!\!\!\!N^4 (X_0^2 \times r) &\!\!\!\!\!\!\!\!\!\!\!\subset &\!\!\!\!\!\!\!\!\!\!\!N^4 (2X_0^2)^{\wedge} &\!\!\!\!\!\!\!\!\!\!\!\!\!\!\!\!\!\!\!\!\!\!\!\!\subset &\!\!\!\!\!\!\!\!\!\!\!\!\!\!\!\!\!\!\!\!\!\!\!\!N_1^4 (2X^2_0)^{\wedge} \equiv N^4 (2X_0^2)^{\wedge} \underset{\overbrace{\mbox{\footnotesize$\partial N_1^4 (2X_0^2)^{\wedge}$}}}{\cup} [\partial N^4 (2X_0^2)^{\wedge} \times [0,1]]. \\
N^4(\Gamma (1)) \ar[u] &\!\!\!\!\!\!\!\!\!\!\!\subset &\!\!\!\!\!\!\!\!\!\!\!N^4 (\Gamma_1 (\infty)) \ar[u] &\!\!\!\!\!\!\!\!\!\!\!\subset &\!\!\!\!\!\!\!\!\!\!\!N^4 (2\Gamma(\infty)) \ar[u]
}
\leqno (77)
$$
We have
$$
N^4 (X_0^2 \times r) = N^4 (\Gamma_1 (\infty)) + \sum_{C,\Gamma} \{\mbox{the $4^{\rm d}$ 2-handle} \ D^2 (C,\Gamma)\},
$$
where $(C,\Gamma) \in \{$the $\{$link$\}$ [new] from (28.2), with all the $\gamma^0$'s DELETED$\}$. The $N^4 (X_0^2 \times r)$ is the place where a large part of the action of this paper takes place and, when we write $N^4 (\Gamma_1 (\infty))$ rather than just $N^4 (\Gamma (\infty))$, it is because we want the $\Sigma \, b^3 (\beta) \subset \partial N^4 (\Gamma_1 (\infty))$ to be with us.

\smallskip

Let now the set $B[{\rm new}]$ be like in the context of $X^2 [{\rm new}]$, meaning $B[{\rm new}] \equiv B({\rm old}) + \{$For every edge $e_i \subset \Gamma (1) \subset \Delta^2$, there is now a $b_i \in e_i \times (\xi_0 = -1)$, the dual curve of which is the
$$
\eta_i = \partial \, [e_i \times [0 \geq \xi_0 \geq -1]]\}.
$$
At the level of $X^2 [{\rm new}]$, the $B[{\rm new}]$ has an order relation given by the geometric intersection matrix $\eta \cdot B = {\rm id} + {\rm nil}$ of the easy type, of $X^2[{\rm new}]$. This natural BLUE order induced by $\eta \cdot B$ is defined as follows. The off-diagonal terms $\eta_{\ell} \cdot B_k$ ($\ell \ne k$) define an oriented graph of vertices $B[{\rm new}]$ and oriented edges $\ell \to k$. The id $+$ nil condition implies that there are no closed oriented orbits of the graph, and the oriented edge $\ell \to k$ will mean $\ell > q$, generating a partial order relation.

\bigskip

\noindent (77.1) \quad In the ordered set $B[{\rm new}]$, pertaining to the space $X^2[{\rm new}] \supset X_0^2 [{\rm new}] \approx X_0^2 \times r \subset 2X_0^2$, with $X^2[{\rm new}] \subset 2X^2$, the elements
$$
(b_i , \eta_i) = (b_i \times (\xi_0 = -1) , \partial (e(b_i) \times [0 \geq \xi_0 \geq -1])),
$$
are INITIAL elements, without incoming arrows (in the ordered graph) i.e. without higher indices $j > i$ is the ordered set $B[{\rm new}]$. [{\bf Remark}. Contrary to what happened with the $B$ of $X^2$ ($=X^2 ({\rm old})$), in the context of Lemma 4 and the ``Important Digression'' which follows it, $B[{\rm new}]$ is not totally ordered, we only have a partial order relation now. No question now any longer, of indexing the states by $Z_+$ or by $Z$.] End of (77.1).

\bigskip

Very importantly, it is the BLUE $2^{\rm d}$ collapse of $X^2[{\rm new}]$ which concerns us now and that is what (77.1) talks about, and NOT the one of $2X^2$. In the BLUE $2^{\rm d}$ collapse which interests us now, every $b_i \in  \Gamma (1) \times (\xi_0 = -1)$ comes, inside the ordered graph corresponding to $\eta \cdot B = {\rm id} + {\rm nil}$, with a trajectory which, a priori, could have been something like below
$$
\includegraphics[width=6cm]{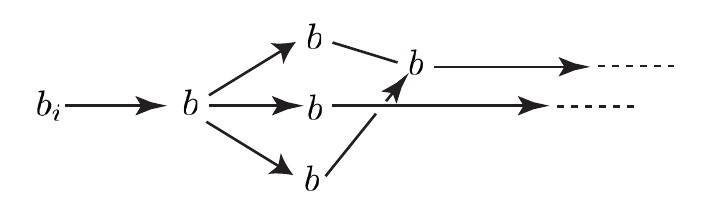}
$$
stopping, in finite time at some dead ENDS.

\smallskip

Actually, the GPS structure which we are using, the condition (5)-(d) and the strategy developed in Figure 35.2 and described afterwards, imposes for $b_i = b_i \times (\xi_0 = -1) \in \Delta^2 \times (\xi_0 = -1)$ a much more restricted form than the one above. There are two cases:

\medskip

I) The trivial case when $b_i \times (\xi_0 = 0) \notin B_0$. Then $b_i = b_i \times (\xi_0 = -1)$ itself is a DEAD END. Theorem~11 for the corresponding $\delta_i^2$ is now trivially verified and we will not dwell longer on this case.

\medskip

II) The $b_i \times (\xi_0 = 0) \subset X^3 \times t_j$ is in $B_0$. Then, to begin with, inside $X^3 \times t_j$, we have a main {\ibf linear} trajectory of $b_i \times (\xi_0 = -1)$, call it $\{$reduced $T(b_i)\}$
$$
b_i = b_i \times (\xi_0 = -1) \to \underbrace{b_i \times (\xi_0 = 0) = b_{i_1} \to b_{i_2} \to b_{i_3} \to \ldots \to b_{i_k}}_{{\rm inside} \, X^3 \times t_j} \ \mbox{(DEAD END)}. \leqno (77.1.{\rm bis})
$$

When $b_{i_k} (77.1.{\rm bis}) \subset e \subset 2X^1 ({\rm COLOUR}) \times t_j$, we might also have additional arrows of type $b_{i_k} \to b(\alpha)$, where $D^2 (\eta (b(\alpha))) \subset [t_j , t_{j \pm 1}]$ and this $b(\alpha)$ is a DEAD END.

\smallskip

So, in its full glory, we have a tree $T(b_i) \supset \{{\rm reduced} \ T(b_i)\}$, representing the complete BLUE trajectory of $b_i \times (\xi_0 = -1)$
$$
\xymatrix{
b_i = b_i \times (\xi_0 = -1) \ar[r] &b_{i_1} \ar[r] \ar[d] &b_{i_2} \ar[r] \ar[d] &b_{i_3} \ar[r] \ar[d] &\ldots \ar[r] &b_{i_k} \ar[d] &\mbox{(DEAD END)} \\
&b(\alpha) &b(\alpha) &b(\alpha) &&b(\alpha)
}
\leqno (77.1.{\rm ter})
$$
$$
\qquad \qquad \qquad \underbrace{\qquad \qquad \qquad \qquad \qquad \qquad \qquad \qquad \qquad \qquad}_{\mbox{\footnotesize some of these may either not exist or be doubles. But anyway, they are all DEAD ENDS.}}
$$

\smallskip

\noindent [The ``double'' here comes from the fact that $D^2 (\eta (b_{i_k}))$ has two vertical edges.]

\bigskip

Let $b_j$ be the generic $b \in \underset{\overbrace{\mbox{\footnotesize$b_i \in (\xi_0 = -1)$}}}{\bigcup} T(b_i)$. To this $b_j$, in $X^2 [{\rm new}]$ corresponds a $D^2 (\eta_j) \subset X^2 [{\rm new}]$. [But then, remember, that at the level of $2X^2 \supset 2X^2_0$ this $D^2 (\eta_j)$ is now a $D^2 (\gamma^1)$ (in the extended sense).]

\smallskip

For our 2-cell $D^2 (\eta_j) \subset X^2 [{\rm new}]$, the BLUE curve $\eta_i$ has a RED counterpart $\Gamma_{\ell}$ or $C_i$ or $\gamma_k^0$ in the link (9) and to this corresponds a $[D^2 (\eta_j)] \subset X^2 [{\rm new}]$ defined as follows: If $\eta_j = \Gamma_{\ell}$ or $C_i$, then $[D^2 (\eta_j)] = \{ D^2 (\Gamma_{\ell})$ or $D^2 (C_i)\}$, i.e. in this case $[D^2 (\eta_j)] = D^2 (\eta_j)$, topologically speaking, as pieces of $X_0^2 [{\rm new}]$. If $\eta_j = \gamma_k^0$, then $[D^2 (\gamma_j)]$ stands just for the collar $\gamma_k^0 \times [0,1]$, where $\gamma_k^0 \times \{1\} \equiv \gamma_k^0 \subset \partial X_0^2 [{\rm new}]$, and, of course, there is no $\{$full $D^2 (\gamma^0)\} \subset X_0^2 [{\rm new}]$. Here $b_j \in e(b_j)$ is like in Figure 7-(II, III), when we are in the generic case, far from the root of the trees $T(b_i)$, respectively like in Figure 7-bis for $b_i \in \Gamma(1) \times (\xi_0 = -1)$, in the non-generic case. To the $b_j$ we attach the following sphere with holes $B_j^2 \subset 2X_0^2$, defined as follows:

\bigskip

\noindent (77.2) \quad When we are in the generic case, then
$$
B_j^2 \equiv \{ ([D^2 (\eta_j)] \times r) \cup (\eta_j \times [r,b]) \cup (D^2 (\eta_j) \times b, \ \mbox{which is always there})\},
$$
and here $e(b_j) \times [r,b]$ is replaced by the surviving collar in the Figure 7-(II, III), and the same for every $e(b_k) \times [r,b]$, where $b_k \subset \eta_j - b_j$, coming clearly with an arrow $b_j \to b_k$ in the oriented graph defined by the $\eta \cdot B$ of $X^2 [{\rm new}]$.

\bigskip

In the case when we consider a $b_i \in \Gamma(1) \times (\xi_0 = -1)$, then $B_i^2$ is defined like above, EXCEPT that now $e(b_i) \times [r,b]$ is like in Figure 7-bis, with all of its interior removed. Of course $e(b_i) \subset \eta_i$.

\smallskip

Since the special treatment of Figure 7-bis is restricted to $e(b_i)$ itself, its effect stays confined at $\xi_0 = - 1$ and, with this, in {\ibf all} cases, the $B_j^2$ is a sphere with holes (see here Figure 19.1 too). This is coming with

\bigskip

\noindent (78) \quad $\partial B_j^2 = c(b_j) + \{$some $c(b_{\ell < j})$ getting arrows $b_j \to b_{\ell}$, see here the order considered in (77.1)$\} +  \{$possibly a $\gamma_k^0$'s, when $\eta_j = \gamma_k^0$ (RED-wise, in $X_0^2 [{\rm new}])\} \subset \partial (2X_0^2)$. [The $\gamma_k^0$ here could be $\Gamma_j \times (\xi_0 = 0)$].

\bigskip

[The BLUE order in $X^2 [{\rm new}]$ is being meant here, organized like in the IMPORTANT DIGRESSION which follows after Lemma 4, adapted afterwards from $X^2$ to $X^2 [{\rm new}]$. Also, in the generic case, the $c(b_j)$ means the smaller of the two boundary curves, of the shaded collar in the case of Figure 7-(II, III), while in the case $b_i \times (\xi_0 = -1) = b_i$, $c(b_i)$ is the boundary of the rectangle in Figure 7-bis. Notice that $X^2 [{\rm new}]$ is used for the BLUE order, while $2X_0^2$ is used for providing the spare parts for the $B_j^2$'s.]

\smallskip

We will denote from now on by $\underset{1}{\overset{M}{\sum}} \, b_i$ the family of the $b_i \in \Delta^2 \times (\xi_0 = -1)$, which is occuring in (77.1) to (77.1-ter).

\bigskip

\noindent {\bf Lemma 10.} 1) {\it For each $b_i \in \underset{1}{\overset{M}{\sum}} \, b_i$ there is a disc $\delta_i^2$, cobounding $c(b_i)$ (i.e. $\partial \delta_i^2 = c(b_i)$), endowed with a non-degenerate map
$$
\delta_i^2 \xrightarrow{ \ \ g_i \ \ } 2X_0^2, \leqno ({\rm 79})
$$
and this $\delta_i^2$ is gotten by the following inductive construction, the induction using here the BLUE {\ibf order} of $X^2 [{\rm new}]$. The construction will only make use of the outgoing trajectories of $\underset{1}{\overset{M}{\sum}} \, b_i$, see here the {\rm (77.1-ter)}.}

\bigskip

\noindent (79-bis) \quad {\it Here is how we CONSTRUCT THE DISC $\delta_i^2$. We start with the disjoined union of all the $B_j^2$'s (see {\rm (77.2)}) for all the $b_j^2 \in \{$reduced $T(b_i)\}$ {\rm (77.1-bis)}. To these, we add now the additional $b(\alpha)$'s from the {\ibf full} $T(b_i)$ {\rm (77.1-ter)}.

\smallskip

These $b(\alpha)$'s, seable in Figure {\rm 35.3-(B)} appear on the boundary of a  vertical $2$-cell $e(b(\alpha)) \times [t_j , t_{j \pm 1}]$ and they also come with an annulus $c(b(\alpha)) \times [0,\varepsilon]$ in the Figure {\rm 7-(II, III)}. Here $c(b(\alpha)) \equiv \{$outer border of the annulus$\} = \{ [e(b(\alpha)) \times r \ (=$ our $e(b(\alpha)))] \cup \partial e(b(\alpha)) \times [r,b] \cup[e(b(\alpha) \times b)]\}$, lives at $t_j$ and bounds a thin $[r,b]$-rectangle, call it $R^2 (b(\alpha))$, normally living at $t=t_j$ too.

\smallskip

In $2X_0^2$, $c(b(\alpha)) = c(b(\alpha)) \times t_j$ bounds the disc
$$
B^2 (b(\alpha)) \equiv (c (b(\alpha)) \times [t_j , t_{j\pm 1}]) \cup R^2 (b(\alpha)) \times t_{j\pm 1}, \leqno ({\rm 79.0})
$$
to which we may or may not add the annulus $c(b(\alpha)) \times [0, \varepsilon]$, shaded in Figure {\rm 7-(II, III)}, according to our convenience, see what is said after the {\rm (79.1)} below.

\smallskip

In order to get the $\delta_i^2$ we start by putting together the $B_j^2$'s from {\rm (77.2)} along the common boundary curves $c(b_j)$. [To be more precise, when we have $b_{j_1} \to b_{j_2}$ in the oriented graph ${\mathcal M} (\eta \cdot B$ of $X^2 [{\rm NEW}])$ then, to begin with, $c(b_{j_2}) \times [0,\varepsilon] = \{$a thin boundary collar, of exterior boundary $c(b_{j_2}) \times \{\varepsilon \}\}$, occurs both in $B_{j_1}$ and in $B_{j_2}$. In $\delta_j^2$ we throw then in $B_{j_1} \cup B_{j_2} - \{$the common $c(b_{j_2}) \times (0,\varepsilon]\}$, rendered now smooth. In this construction, each $c(b_{\ell < j}) \subset \partial B_j^2$ (i.e. with $j \to \ell$), continues inductively with its $B_{\ell}^2$. Moreover, when $\gamma_k^0 \subset \partial B_j^2$ OR $c(b(\alpha)) \subset \partial B_j^2$ occur, then we fill in with $d_k^2$, the one cobounding the $\gamma_k^0$, respectively with $B^2 (b(\alpha))$ (see {\rm (79.0)}). In the context of {\rm (79)}, when we get to the $d_k^2$ we set $g_i \mid d_k^2 = f_k$.} End of (79-bis).

\bigskip

\noindent {\bf Important remark.} In order to construct $\delta_i^2$ we need $2X_0^2$, which provides the necessary spare parts. But the order in which we put things together, is the BLUE order of $X^2[{\rm NEW}]$, {\ibf not} the one of $2X^2$. $\Box$

\bigskip

{\it Given the form of $T(b_i)$ {\rm (77.1-ter)} our $\delta_i^2$ is a disc of boundary $c(b_i)$}

\bigskip

\noindent (79.1) \quad {\it $\delta_i^2 = \{$a main disc with holes ${\mathcal B}_i^2 \equiv \bigcup \ \{$the contribution of the $B_j^2$'s and $B(b(\infty))$'s$\}\} \cup \sum d_k^2$'s, filling in the holes.

\bigskip

[The ${\mathcal B}_i^2$ normally contains the contribution $B^2 (b(\alpha))$ (see {\rm (79.0)}) for each $b(\alpha)$ involved in the full $T(b_i)$. But, for purely expository purposes, we will chose to let only the collar $c(b(\alpha)) \times [0,\varepsilon]$ occur in ${\mathcal B}_i^2$, to begin with, throwing in the new boundary $c(b(\alpha)) \times \varepsilon$ much later than the boundaries $\gamma_k^0 \times \varepsilon$ of ${\mathcal B}_i^2$. This way we avoid for $c(b(\alpha))$ the {\ibf deletions} of the $c(b_{j_2}) \times [0,\varepsilon]$'s mentioned earlier in our construction of $\delta_i^2$.]

\bigskip

We will work from now on with the following big but finite map
$$
\sum_1^M \delta_i^2 \xrightarrow{ \ \ \underset{1}{\overset{M}{\sum}} \, g_i \ \ } 2X_0^2 . \leqno (**)
$$
We denote by $\underset{1}{\overset{N}{\sum}} \, d_k^2$ the $d_k^2$'s occurring in the map above. Here $N=N(M)$.}

\medskip

2) {\it The global map $(**)$ above can be lifted off $2X_0^2$ and changed into a generic immersion (the generic form of maps of $2^{\rm d}$ into $4^{\rm d}$), which is winding tightly around $2X_0^2 \subset N_1^4 (2X_0^2)^{\wedge}$,
$$
\sum_1^M \delta_i^2 \xrightarrow{ \ \  J \ \ } N_1^4 (2X_0^2)^{\wedge} ; \leqno ({\rm 80})
$$
this map $J$ extends the $\underset{1}{\overset{M}{\sum}} \, c(b_i) \subset \partial N^4 (2X_0^2)^{\wedge} \subset N_1^4 (2X_0^2)^{\wedge}$.}

\medskip

3) {\it Generally speaking, the $J$ {\rm (80)} has ACCIDENTS, namely double points $x \in JM^2 (J) \subset N_1^4 (2X_0^2)^{\wedge}$ AND also transversal contacts
$$
z \in J \delta_i^2 \pitchfork 2X_0^2 \subset N_1^4 (2X_0^2)^{\wedge}. \leqno {\rm (81)}
$$

\noindent [{\bf Notations}. For any, generic, map $K \xrightarrow{ \ \psi \ } L$, we use the notations $M_2 (\psi) \equiv \{$set of $x \in K$ s.t. $\# \, \psi^{-1} \psi (r) > 1\} \subset K$. But we can also consider the double points of $\psi$ in the context
$$
K \times K \supset M^2 (\psi) \underset{\mbox{\footnotesize projection on the first factor}}{\xrightarrow{ \qquad \qquad \qquad \qquad \qquad \qquad}} M_2 (\psi) \subset K ,
$$
coming with $\psi M^2 (\psi) = \psi M_2 (\psi) \subset L$.]

\smallskip

The next item is a consequence of {\rm (21.A)}.}

\medskip

4) {\it When it comes to the accidents {\rm (81)} of the form
$$
z \in J \delta_j^2 \pitchfork \Delta^2 \times (\xi_0 = - 1),
$$
then these can only come from the parts $d_k^2 \subset \delta_i^2$ and {\ibf never} from the parts $B_i^2 \subset \delta_j^2$. Here is an immediate consequence of this fact.}

\bigskip

\noindent (81.1) \quad {\it Every transversal contact
$$
z \in J {\mathcal B}_i^2 \cap 2X_0^2 \subset N_1^4 (2X_0^2)^{\wedge}
$$
possesses an $\{$extended cocore $(z)\} \subset 2X_0^2$, since it does not concern $\Delta^2 \times (\xi_0 = -1)$.}

\bigskip

Before going into the proof, notice that in the context of this construction for $c(b_{\ell}) = \partial \delta_{\ell}^2$ we have $c(b_{\ell}) \cdot b_k = \delta_{\ell k}$ so, in some singular sense, the $\underset{1}{\overset{M}{\sum}} \, \delta_{\ell}^2$ are in cancelling position with the 1-handles $\underset{1}{\overset{M}{\sum}} \, b_{\ell}$. We say here ``singular'', because of the fact that there are ACCIDENTS, see above.

\bigskip

\noindent {\bf Proof.} We only need to prove 4). We start by reminding that any edge $e_i \subset \Gamma (1) \times (\xi_0 = -1)$ contains a $b_i \in B [{\rm new}]$. Remember also, that any $\sigma^2 \subset \Delta^2 \times (\xi_0 = -1)$ is, BLUE-wise a $D^2 (\gamma^1)$ of $X^2 [{\rm new}]$; so these $\sigma^2$'s do {\ibf not contribute} to the 2-cells $D^2 (\eta)$ of (77.1). Then, let $b_i \subset e \subset \partial \sigma^2 \subset \sigma^2 \subset \Delta^2 \times (\xi_0 = -1)$. As we know, such a $b_i$ is automatically an INITIAL element in the BLUE order of $X^2 [{\rm NEW}]$, i.e. it has no incoming BLUE arrows. At the level of our $X_r^2 = X_0^2 [{\rm NEW}] \subset X^2 [{\rm NEW}]$, the geometric intersection matrix $\eta \cdot b$ is of the easy id $+$ nil type. The BLUE dual of the $b_i \subset e_i \subset \Gamma (1) \times (\xi_0 = -1)$ is $D^2 (\eta_i) = e_i \times [-1 \leq \xi_0 \leq 0]$, with $c(b_i) = \partial (e_i \times [r,b])$ and, if the edge $e_i \times (\xi_0 = 0)$ contains a $b_j \in B[{\rm NEW}]$, this comes with $j < i$. In this case, the construction of $\delta^2 (b_i)$ does not stop short at $B_i^2$, but it continues further.

\smallskip

All the $B_{\ell}^2$'s present in our construction come from the finite family of those $b \in \underset{1}{\overset{M}{\sum}} \, T(b_i)$ (see (77.1-bis)). If, by any chance, for some $b_i = b_i \times (\xi_0 = -1) \subset \Delta^2 \times (\xi_0 = -1)$ it so happens the $b_i \times (\xi_0 = 0)$, which should be next to our $b_i$ in $T(b_i)$, is NOT in $B[{\rm new}]$, then our $\delta_i^2$ is trivially the $B_i^2 = D^2$ and there are no accidents for it. This really is a very trivial case when everything which we will want to achieve is already automatically there. We will ignore this case, from now on.

\smallskip

The edges $P \times [0 \geq \xi_0 \geq -1]$ never contain $B$'s and so, for our $b_j$ above, occurring as $b^?$ in Figure 19.1, this is the only prospective $b_{j<i} \subset \eta_i$, for the $b_i \in \Gamma (1) \times (\xi_0 = -1)$.

\smallskip

The edge $e = e(b_i) \subset \Gamma (1) \times (\xi_0 = -1)$ gets treated exactly like in 4) from Lemma 5, i.e. like in (21-C) and in Figure 7-bis. The point is that, when we add the $b^3 (\beta)$'s to the Figures 9 which correspond to $\Delta^2 \times (\xi_0 = -1)$ (and each of the two $\beta$'s in the Figure 19.1 corresponds to such a $b^3 (\beta)$), then this does not come (in the figures in question {\ibf with any $2^{\rm d}$ contribution from} $\underset{1}{\overset{M}{\sum}} \, {\mathcal B}_i^2$ (see (79.1)). [Without the strategic decision from 4) Lemma 5, in a figure of type 24 for $\Delta^2 \times (\xi_0 = -1)$ we could find transversal contacts like
$$
B^2 (\beta , z_+ , y_+) \pitchfork \{C(x_- , t_+) \subset \mbox{some curve} \ \Gamma_j \} ,
$$
which would certainly contradict our 4). See here Figure 24, for a concrete illustration. BUT, with our strategic decision 4) in Lemma 5, illustrated in Figure 7-bis, at $\xi_0 = -1$ the $b^3 (\beta)$ only comes with $1^{\rm d}$ attachements, no $2^{\rm d}$ ones, and the bad contact is avoided.]

\smallskip

With these things, clearly when ${\mathcal B}_i^2 \subset \delta^2 (b_i)$ is lifted off $2X_0^2$, then this happens far from $\Delta^2 \times (\xi_0 = -1)$. All the rest of $\delta^2 (b_i)$, source of the map (80) is made out of spare parts $d_k^2$. When it comes to the construction of the ${\mathcal B}_i^2 \subset \delta^2 (b_i)$, this can touch $\Delta^2 \times (\xi_0 = -1)$ only via $B^2 (b_i)$, along the edge $e(b_i)$. The $B^2 (b_i)$ is here like any other $B_j^2$ in the inductive construction from (77.2) and (79-bis). We are here, with $B^2 (b_i)$ in the non-generic case. The edge $e(b_i)$ comes with $\Gamma(1) \times (\xi_0 = -1) \supset e(b_i) \subset c(b_i) \subset \partial B^2 (b_i)$. The $\Delta^2 \times (\xi_0 = -1)$ with its adjacent $B^2 (b_i)$'s, is presented in the Figure 19.1 which should help vizualise our discussion

$$
\includegraphics[width=13cm]{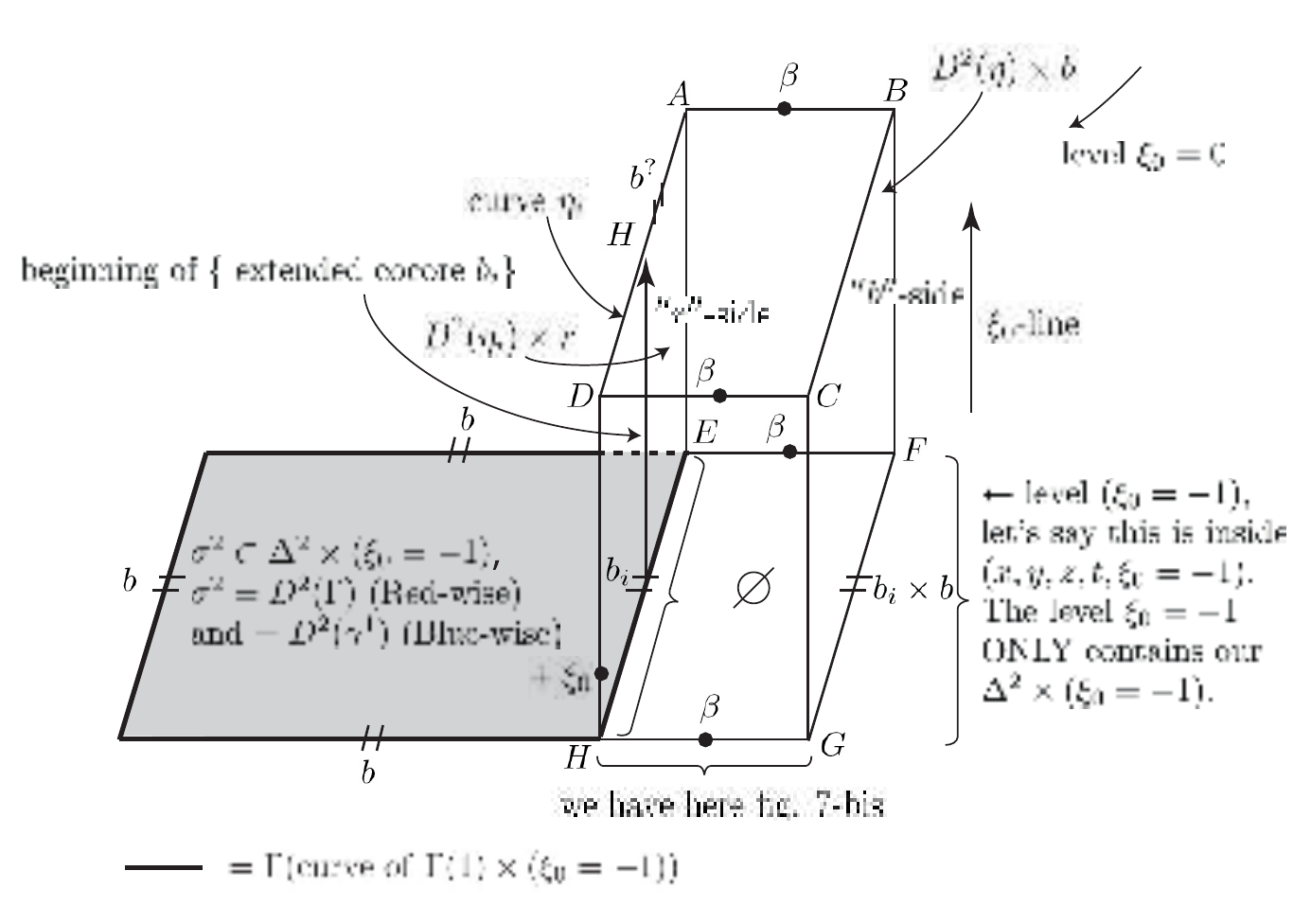}
$$
\label{fig19.1}
\centerline {\bf Figure 19.1.} 
\begin{quote}
The interaction of $\Delta^2 \times (\xi_0 = -1)$ with $B^2 (b_i) = [A,B,C,D,E,F,G,H]$, with the interior of $[E,F,G,H]$ being {\ibf deleted}.
\end{quote}

\bigskip

\noindent {\bf Additional explanations for Figure 19.1.} Without loss of generality, in the context of our figure, we have
$$
JB^2 (b_i) \cap (\Delta^2 \times (\xi_0 = -1)) = e(b_i).
$$

At $\xi_0 = 0$, adjacent to our $b_i$ lives the $b^? = b_{j>i}$ and the question mark is there to suggest that (in the exceptional case), this $b_j$ may be lacking. Anyway, the rectangle $[ABCD]$ adjacent to $b^?$ is like in Figure 7, unlike the lower $[EFGH]$ which is like in Figure 7.bis. Every edge $e \subset \Delta^2 \times (\xi_0 = -1)$ is like our present $e(b_i) = [EH]$, with a $B^2 (b_i)$ like in the Figure 19.1 adjacent to it. End of explanations.

\bigskip

\noindent {\bf Important remark.} All the action in Lemma 10 is concentrated in
$$
N^4 (2X_0^2) \mid ((X_0^2 \times r) \cup (\Gamma (\infty) \times [r,b])). \qquad \Box
$$

\bigskip

The next theorem and its proof will occupy the rest of this section IV.

\bigskip

\noindent {\bf Theorem11. (Existence of external, embedded discs.)}

\smallskip

{\it We can change the map {\rm (80)} into a DIFF embedding, where not only $J$ is changed but the $\underset{1}{\overset{M}{\sum}} \, \partial \delta_i^2$ too
$$
\left( \sum_{i=1}^M \delta_i^2 , \sum_1^M \eta_i ({\rm green}) = \partial \delta_i^2 \right) \underset{\mbox{\footnotesize $J=J$ (embedded)}}{\xrightarrow{\qquad \qquad \qquad \qquad}} (\partial N^4 (2X_0^2)^{\wedge} \times [0,1]  \subset N_1^4 (2X_0^2)^{\wedge} , \partial N^4 (2X_0^2)^{\wedge} \times \{0\}), \leqno {\rm (82)}
$$
for which, eventually, we will have, in terms of the geometric intersection matrix
$$
\eta_j ({\rm green}) \cdot b_i^2 = \delta_{ji}, \ \mbox{for} \ i , j \in \{1,\ldots , M\} , \ \mbox{AND} \ \eta _j ({\rm green}) \cdot \left( B_1 - \sum_1^M b_i \right) = 0. \leqno {\rm (82.1)}
$$

In other terms,  we manage to create a system of external embedded discs (NO MORE ACCIDENTS), in exact cancelling position with $\underset{1}{\overset{M}{\sum}} \, b_i$.}

\bigskip

[Here is what the ``eventually'' in the statement above means. In a first step, realized in this section IV already, we will realize a map
$$
\left( \sum_1^M \delta_i^2 , \sum_1^M \partial \delta_i^2 \equiv C_i (b)(\lambda) \right) \underset{\mbox{\footnotesize $J$}}{\xrightarrow{\qquad \qquad}} (\partial N^4 (2X_0^2)^{\wedge} \times [0,1] , \, \partial N^4 (2X_0^2) \times \{0\}), 
$$
which is completely ACCIDENT free (see (93.1) below), but without having yet the (82.1) satisfied. That will be eventually realized only by Theorem 13 in the next section V.]

\bigskip

\noindent {\bf Remark.} As things stand, right now, there is no handlebody decomposition of $N^4 (\Delta^2)$ having the $\underset{1}{\overset{M}{\sum}} \, b_i$ as its system of 1-handles. Here, notice that $\Gamma(1) - \underset{1}{\overset{M}{\sum}} \, b_i$ is highly disconnected. This means that, Theorem 11, as it stands, cannot be plugged into Lemma 3, our (82.1) not withstanding. Not to speak, also, of the fact that the $\underset{1}{\overset{M}{\sum}} \, c(b_i)$ is not contained in $\partial N^4 (\Delta^2)$ (nor in $N^4 (\Delta^2)$, for that matter).

\bigskip

\noindent {\bf Proof.} Let's stay for the time being at level $X^2({\rm old})$ and consider the purely spatial Figures 4. To them, we will add, in an appropriate way, the $b^3 (t_+)$ OR $b^3 (t_-)$ {\ibf never both simultaneously}. This leads to the Figures 9 where in addition to the already purely spatial smooth sheets, like $(x_+ , y_+ , z_+)$ in Figure 4-B appear now new ones, like $(x_+ , y_+ , t_{\pm})$, in the Figures 9-(B$_{\pm}$). Later $\xi_0$, $\beta$ will occur too.

\smallskip

These figures are considered at the level $X^2 ({\rm old})$, for the time being.

\smallskip

The $f_k \, d^2_k$ (76) consist of unions of local smooth sheets like the $(x_{\pm} , y_{\pm} , z_{\pm})$, $(x_{\pm} , y_{\pm} , t_{\pm})$ above, appearing with possible multiplicities, and not all of them may appear inside $\underset{k}{\bigcup} \, f_k \, d_k^2$. When we move to our enlarged figures from Figures 9, then the simple closed loops in those figures span little 2-cells contained in $S^3 (P) = \partial N^4 (P)$, like $\left(\xy *[o]=<15pt>\hbox{$x_-$}="o"* \frm{o}\endxy , \xy *[o]=<15pt>\hbox{$t_{\pm}$}="o"* \frm{o}\endxy , \xy *[o]=<15pt>\hbox{$x_+$}="o"* \frm{o}\endxy \right)$ or $\left(\xy *[o]=<15pt>\hbox{$x_-$}="o"* \frm{o}\endxy , z_+ , \xy *[o]=<15pt>\hbox{$x_+$}="o"* \frm{o}\endxy \right)$ in Figure 9-A or $\left(\xy *[o]=<15pt>\hbox{$z_+$}="o"* \frm{o}\endxy , \xy *[o]=<15pt>\hbox{$t_+$}="o"* \frm{o}\endxy , \xy *[o]=<15pt>\hbox{$x_-$}="o"* \frm{o}\endxy \right)$ in Figure 9-B$_+$. We will denote there 2-cells, generically by $d^2 (x_- , t_{\pm} , x_+)$ (or, if a letter is without a surrounding circle, we put prentices around (like $d^2 (x_- , (z_+) , x_+)$). Since these $d^2$'s follow the rules of the game established by the paradigmatical Figure 8, we have the following 

\bigskip

\noindent (83) \quad At the level of the Figures 9, inside $\partial N^4 (P)$, the various $d^2$ (space-time, space-time, space-time) are 2-by-2 disjoined except for possible common edges or vertices.

\bigskip

Next we move from $X^2 ({\rm old})$ to the realistic $2X_0^2$, where we will want to take $J\delta^2$ off the $\Delta^2 = \Delta^2 \times (\xi_0 = -1)$. But before we can do that we will describe the changes to be operated at the level of the figures of type 9, when moving from $X^2 ({\rm old})$ to $2X_0^2$.

\bigskip

\noindent (84.1) \quad To our purely space-time figures 9, when we are at level $X_0^2 \times r \subset 2X_0^2$, we add more vertices, namely the following:

\medskip

\noindent $\bullet$) A vertex $\xy *[o]=<15pt>\hbox{$\beta$}="o"* \frm{o}\endxy$ to {\ibf all} Figures 9 $(\subset X_0^2 \times r)$ corresponding to the axis $\zeta$, moving from $r$ towards $b$ (with, remember, $r < \beta < b$ and $\vert \beta - r \vert \ll \vert b-r \vert$. The $[r,b]$ lives somewhere on the $\zeta$-axis.

\medskip

\noindent $\bullet$$\bullet$) A vertex $\xy *[o]=<18pt>\hbox{$+\xi_0$}="o"* \frm{o}\endxy$ (respectively $\xy *[o]=<18pt>\hbox{$-\xi_0$}="o"* \frm{o}\endxy$) for every $P \times (\xi_0 = -1)$, when $P \in \Gamma (1)$, (respectively for every $P \times (\xi_0 = 0)$, $P \in \Gamma (1)$).

\bigskip

\noindent (84.2) \quad This way, new arcs occuring in the REALISTIC FIGURES 9. Every Figure 9 acquires now a $b^3 (\beta) \subset B^3 (+)$. Any Figure 9 acquires a $b^3 (-\xi_0) \subset B^3 (-)$, if it lives at $\xi_0 = 0$, respectively a $b^3 (+ \, \xi_0) \subset B^3 (-)$, if it lives at $\xi_0 = -1$. Here are the {\ibf golden rules} for drawing our realistic Figures 9.

\bigskip

\noindent $\bullet$$\bullet$$\bullet$) For $P \times (\xi_0 = -1) \subset \Delta^2 = \Delta^2 \times (\xi_0 = -1)$, never contradict what is already done in the Figures~9. The ``Figures 9'' are all compatible with the paradigmatical Figure 8, so the present rule makes (22.A) automatically satisfied, provided we proceed like in $X^2 ({\rm old})$, along the edges of $\Delta^2 \times (\xi_0 = -1)$, in particular when the framing for attaching 2-handles are concerned.

\medskip

\noindent $\bullet$$\bullet$$\bullet$$\bullet$) In terms of the basic splitting from (27.6), (28), the Figures 9 live essentially in $B^3 (-)$ (to be more precise in $B^3 (-) \cup B^3 ({\rm in})$), while we have $b^3 (\beta) \subset B^3 (+) (\cup  \, B^3 ({\rm out}))$. Figure 9-bis gives the accurate picture.

\medskip

\noindent $\bullet$$\bullet$$\bullet$$\bullet$$\bullet$) For our $P \times (\xi_0 = 0)$ or $P \times (\xi_0 = -1)$, the $b^3 (- \, \xi_0)$, respectively $b^3 (+ \, \xi_0)$, live inside $B^3 (-)$, like the rest of Figures 9 and, if we think that $t_+$ points towards the observer and $t_-$ away from the observer, then $b^3 (+ \, \xi_0)$ lives higher than $b^3 (t_+)$ and $b^3 (- \, \xi_0)$ lower than $b^3 (t_-)$. The $b^3 (t_+) / b^3 (t_-)$ are mirror images of each other, a condition stemming from Figures 8 and 8-(B).

\smallskip

Now the $\beta$ and the $\pm \, \xi_0$ have nothing to do with $R^4$ (which imposes Figure 8) and no metric restrictions concern them for our $4^{\rm d}$ reconstruction formulae. But then, for reasons of simplicity, we will impose on $b^3 (- \, \xi_0) / b^3 (+ \, \xi_0)$ the same kind of mirror symmetry as for $b^3 (t_-) / b^3 (t_+)$. Other better reasons for this might become clear later on.

\smallskip

[The mirroring is like in the Figure 8-(B) BUT, careful, do not mix up the splittings $B^3 ({\rm out}) \underset{\overbrace{\mbox{\footnotesize $S^2 (P)$}}}{\cup} B^3 ({\rm in})$ and $B^3 (+) \underset{\overbrace{\mbox{\footnotesize $S^2_{\infty}$}}}{\cup} B^3 (-)$. The $S^2 (P)$ is defined by Figure 8-(C) and the relation $S^2 (P) / S_{\infty}^2$ is shown in Figure~9.bis. Approximatively, one may think in terms of $B^3 ({\rm in}) \subset B^3 (-)$.]

\medskip

\noindent $\bullet$$\bullet$$\bullet$$\bullet$$\bullet$$\bullet$) Any crossing of lines involving $\xy *[o]=<15pt>\hbox{$\beta$}="o"* \frm{o}\endxy$, $\xy *[o]=<18pt>\hbox{$\pm\xi_0$}="o"* \frm{o}\endxy$ or space-time lines {\ibf not} at $\Delta^2 \times (\xi_0 = -1)$ are {\ibf acceptable}, as long as (22-B) is also satisfied. We will call this the {\ibf crossing freedom}. This comes with {\ibf local changes in topology}. But, once $C \cdot h = {\rm id} + {\rm nil}$ stays with us, the global topology stays unchanged.

\medskip

\noindent $\bullet$$\bullet$$\bullet$$\bullet$$\bullet$$\bullet$$\bullet$) We do insist that (83) should be satisfied everywhere. [As a consequence of b), e) in the definition (5), in any figure of type 24 there is {\ibf at most one} crossing $\includegraphics[width=8mm]{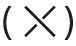}$ of space-time lines. Any change of {\ibf this} one by crossing UP $\Leftrightarrow$ DOWN of the two lines, cannot being about any violation of (83). And then, our (83) does not concern neither $\beta$ nor $\pm \, \xi_0$.] End of (84.2)

\bigskip

\noindent (85) \quad When it comes to the realistic $f_k \, d_k$ (76) it is only the $P \in X_0^2 \times r$ which are involved, and in the corresponding realistic Figures 9, $\xy *[o]=<15pt>\hbox{$\beta$}="o"* \frm{o}\endxy$'s are mute (for $d_k$ but {\ibf not} for $\delta_k^2$).

\bigskip

\noindent {\bf Lemma 11.1.} 1) {\it We can lift the map $\overset{N}{\underset{1}{\sum}} \, d_k^2 \xrightarrow{ \ \Sigma \, f_k \ } X_0^2 \times r \subset 2X_0^2$ to an immersion winding thightly around $X_0^2 \subset N_1^4 (2X_0^2)^{\wedge}$
$$
\sum_1^N d_k^2 \xrightarrow{ \ \ J \ \ } N_1^4 (2X_0^2)^{\wedge} , \leqno (86)
$$
with ACCIDENTS like in {\rm 3)} Lemma $10$. This is part of} (80).

\medskip

2) {\it We can get our LIFT {\rm (86)} by the following two steps process. One starts by lifting $\overset{N}{\underset{1}{\sum}} \, d_k^2$ to a {\ibf not everywhere well-defined} map
$$
\sum_1^N d_k^2 \xrightarrow{ \ \ F \ \ } \partial N^4 (2X_0^2), \leqno (87)
$$
and for describing it one needs a partition, to be made explicit below: $d_k^2 =$ body $d_k^2 + (d_k^2 - \mbox{body} \ d_k^2)$. Here $F \mid (d_k^2 - \mbox{body} \ d_k^2)$ is completely well-defined everywhere and WITHOUT ACCIDENTS. It runs parallel and very close to the $D^2 (\Gamma_j \ \mbox{or} \ C_i) \subset X_0^2 \times r$. More precisely, for each of the $D^2 (\Gamma_j \ \mbox{or} \ C_i)$ which is involved, the corresponding piece of $F (d_k^2 - \mbox{body} \ d_k^2)$ lives on the lateral surface of the corresponding $4^{\rm d}$ handle of index $\lambda = 2$, call it $S^1 \times D^2$, taking the form $(*) \times D^2$.

\smallskip

So, the disc $d_k^2$ is broken into a $\{$disc with holes, called body $d_k^2\}$ of boundary $\{\gamma_k^0 + \mbox{(internal circles)}\}$, while $d_k^2 - \mbox{body} \ d_k^2$ is a disjoined collection of discs $D^2 (C) , D^2 (\Gamma)$, filling the internal circles in question.}

\medskip

3) {\it It is the map
$$
\sum_1^N \mbox{body} \ d_k^2 \xrightarrow{ \ \ F \ \ } \partial N^4 (\Gamma_1 (\infty)) \leqno (88)
$$
which is not everywhere well-defined, since it comes with {\ibf punctures}, i.e. transversal contacts
$$
F(\mbox{body} \ d_k) \pitchfork \{\mbox{attaching zones of the $2$-handles, i.e. the} \leqno (88.1)
$$
$$
\mbox{$\{$curves $\Gamma_j$ or $C_i\}$} \times b^2 (\mbox{framing}) \subset \partial N^4 (2\Gamma (\infty))\}.
$$

These punctures as well as the double points $M^2 (F)$ create the ACCIDENTS of the lifted map {\rm (86)}. Before going on with the {\rm (86)} we will open a}

\bigskip

\noindent {\bf Prentice.} We have a decomposition

\bigskip

\noindent (88.2) \quad $\overset{N}{\underset{1}{\sum}} \, d_k^2 = \biggl\{ \overset{N}{\underset{1}{\sum}} \, d_k^2 \, \Bigl \vert \, X^2 ({\rm old})$ from which we exclude the interiors of the 2-cells in $d_k^2 \cap (\Delta^2 \times (\xi_0 = 0) \biggl\}$ $+$ $\{$a contribution from $(\Gamma (1) \times [0 \geq \xi_0 \geq -1]) \cup (\Delta^2 \times (\xi_0 = -1))$ which requires the deletion just done$\}$.

\bigskip

More explicitly, we change $d_k^2 \mid X^2 ({\rm old})$ by replacing every $\sigma^2 \subset (d_k^2 \mid X^2 ({\rm old})) \cap (\Delta^2 \times (\xi_0 = 0))$ by $(\partial \sigma^2 \times [0 \geq \xi_0 \geq -1]) \cup (\sigma^2 \times (\xi_0 = -1))$.

\smallskip

Once the contribution of $\Delta^2 \times (\xi_0 = 0)$ has been removed, the rest of body $\left(\overset{N}{\underset{1}{\sum}} \, d_k^2  \, \Bigl \vert \,  X^2 ({\rm old})\right)$ stays intact, inside $\partial N^4 (\Gamma_1 (\infty))$ (see here $\bullet$$\bullet$$\bullet$$\bullet$$\bullet$) in (84.2)).

\smallskip

Because of $\bullet$) in (84.2) and because we never mock around with the Figures 9 at $(\xi_0 = -1)$, the body $\left(\overset{N}{\underset{1}{\sum}} \, d_k^2  \, \Bigl \vert \,  (\Delta^2 \times (\xi_0 = -1)) \right)$ comes without any special problems. The body $\left(\overset{N}{\underset{1}{\sum}} \, d_k^2  \, \Bigl \vert \,  (\Gamma (1) \times [0 \geq \xi_0 \geq -1]) \right)$ is a topic to be discussed more in detail, below. End of prentice.

\bigskip

{\it All this ends a first part of our road to $(86)$ and the second part will consist in the $\{$lifting of $(88)$ off from $\partial N^4 (\Gamma_1 (\infty))\} + \{$the much more trivial step of lifting the $F \mid (d_k^2 - \mbox{\rm body} \ d_k^2)\}$. This last piece consists just of parallel copies of the various $D^2 (C_i)$ or $D^2 (\Gamma_j)$. As already said, it is the lift of $(88)$ which creates the ACCIDENTS of $(86)$.

\smallskip

A priori, the punctures in {\rm (88.1)} are accompanied by CLASPS and RIBBONS. Now, at the level of the embellished Figures $9$, if the $b^3 (\pm \, \xi_0)$ are lacking (and see here Figure $24$ for illustration, with all the red part deleted) and if we, moreover, restrict ourselves to the $\overset{N}{\underset{1}{\sum}} \, d_k^2$ part of $\overset{M}{\underset{1}{\sum}} \, \delta_i^2$, and hence forget the $b^3 (\beta)$ (and hence the green part of Figure $24$ too), then it follows from $(83)$ that there are {\ibf no transversal contacts}
$$
d^2 (\mbox{\it space-time, space-time, space-time}) \pitchfork C(\mbox{\it space-time}).
$$

But then, we will have {\ibf transversal contacts}
$$
d^2 (\pm \, \xi_0 , \ \mbox{\it space-time, space-time}) \pitchfork C(\mbox{\it space-time}),
$$
when ``$C$'' may mean actual $C_i$ or $\Gamma_j$. There will be made explicit, completely, and they produce both the {\rm (88.1)} and the
$$
M^2 \left( F \ \biggl\vert \ \sum_{k=1}^N \mbox{\it body} \ d_k^2 \right) .
$$

We will also have
$$
M^3 \left( F \ \biggl\vert \ \sum_{k=1}^N \mbox{\it body} \ d_k^2 \right) = \emptyset = M^3 \left( J \ \biggl\vert \ \sum_1^N d_k^2 \right).
$$
This absence of triple points follows because, at the level of the figures of type $24$, the only conceivable way to produce some $M^3$ would be
$$
d^2 (\pm \, \xi_0 , \ \mbox{\it space-time, space-time}) \pitchfork d^2(\mbox{\it pure space-time}) \pitchfork d^2(\mbox{\it pure space-time}) \leqno (*)
$$
and, once $d^2(\mbox{\it pure space-time}) \pitchfork d^2(\mbox{\it space-time}) = \emptyset$, the $(*)$ is also absent.

\smallskip

When we are {\ibf far from} $\Delta^2 \times [0 \geq \xi_0 \geq -1]$, there are neither CLASPS nor RIBBONS for $\overset{N}{\underset{1}{\sum}} \, d_k^2$. Any CLASP has to involve $d^2 (\pm \, \xi_0 , \ \mbox{space-time})$. Moreover, a direct analysis will show that the CLASPS, occurring all of them in the region $\Delta^2 \times  [0 \geq \xi_0 \geq -1]$, are $2$-by-$2$ disjoined.

\smallskip

Figure $20$ gives a general illustration of the CLASPS of $\underset{k}{\sum} \ d_k^2$ and Figure $21$ of the RIBBONS.

\smallskip

So, our Figure $20$ which lives inside $\partial N^4 (\Gamma_1 (\infty))$ illustrates a clasp of $(88)$ and, a more detailed view is provided by the Figure $22$.}

$$
\includegraphics[width=105mm]{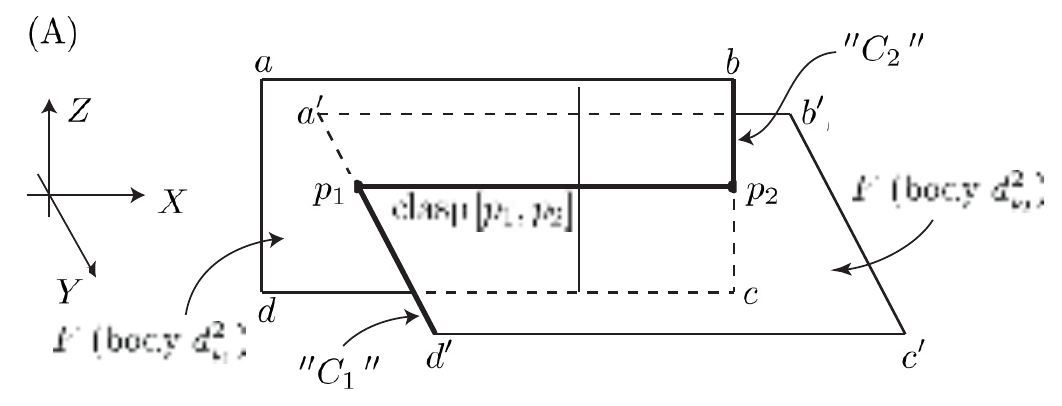}
$$
\begin{quote}
This figure is in $\partial N^4 (\Gamma_1 (\infty))$. The $C_1 = [a' , p_1 , b']$ and $C_2 = [c , p_2 , b]$ are pieces of two curves $\Gamma$ or $C$, attaching zones of 2-handles. The $p_1 , p_2$ are punctures, like in (88.1). Inside the $\partial N^4 (2\Gamma_1 (\infty))$, the pieces of $F ({\rm body} \, d^2)$ seable in this figure, continue beyond the thin contours $[ba] \cup [ad] \cup [dc]$ and $[a'b'] \cup [b'c'] \cup [c'd']$.
\end{quote}

$$
\includegraphics[width=105mm]{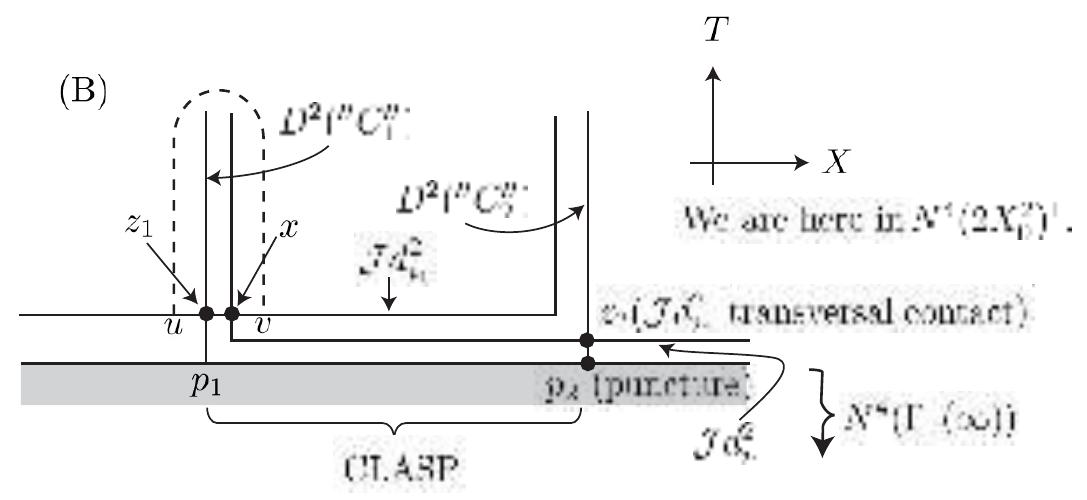}
$$
\label{fig20}
\centerline {\bf Figure 20.} 
\begin{quote}
(A) is generic figure for a CLASP of (88), creating in $4^{\rm d}$ (in (B)) the ACCIDENTS $z_1 , z_2$ and $x$. Figure 23 is a more accurate and detailed version of (B). The process of pushing over the (extended cocore)$^{\wedge}$ of $z_1$, is suggested in dotted lines.
\end{quote}

\bigskip

Here are {\bf some explanations concerning Figure 20.} Beyond the contours $[b,a] \cup [a,d] \cup [d,c]$ and $[a',b'] \cup [b',c'] \cup [c',d']$, body $d_k^2$'s continues, while at the ``$C_1$'' and ``$C_2$'', where the $d_k^2$'s climb on the corresponding $D^2 (\mbox{``$C$''})$'s, the bodies stop. The (B) shows a lift of the $d_{i_1}^2 , d_{i_2}^2$ from (A) into $N_1^4 (2X_0^2)^{\wedge}$, cut by the plane $(X,T)$; this is also transversal to $\partial N^4 (\Gamma_1 (\infty))$, the space of Figure (A). The $[u,v]$ in (B) is the trace of a small disc $D^2 \subset J d_{i_1}^2$, centered at the transversal contact $z_1$ and cutting transversally, both through $D^2(\mbox{``$C_1$''})$ and to $Jd_{i_2}^2$. In the more specific Figure 23, which will supersede the present (B), this disc $D^2$ is put to work. Importantly, the double point $x$ is contained in $D^2$ too. Here $d_{i_1}^2 \subset \{$some $\delta_i^2\}$ and, the $J \, \overset{M}{\underset{1}{\sum}} \, \delta_i^2$ runs close, at some distance $\varepsilon > 0$ from $\overset{M}{\underset{1}{\sum}} \, g_i \, \delta_i^2 \subset 2X_0^2$ (see (79)).

\smallskip

We have here diam $D^2 \gg \varepsilon$, which will make that there is no obstruction for the operation of pushing $J d_{i_1}^2$ OVER the $\{$extended cocore $z_1\}^{\wedge}$ (when it exists), like in the Figure 23.

\bigskip

{\it We continue with}

$$
\includegraphics[width=11cm]{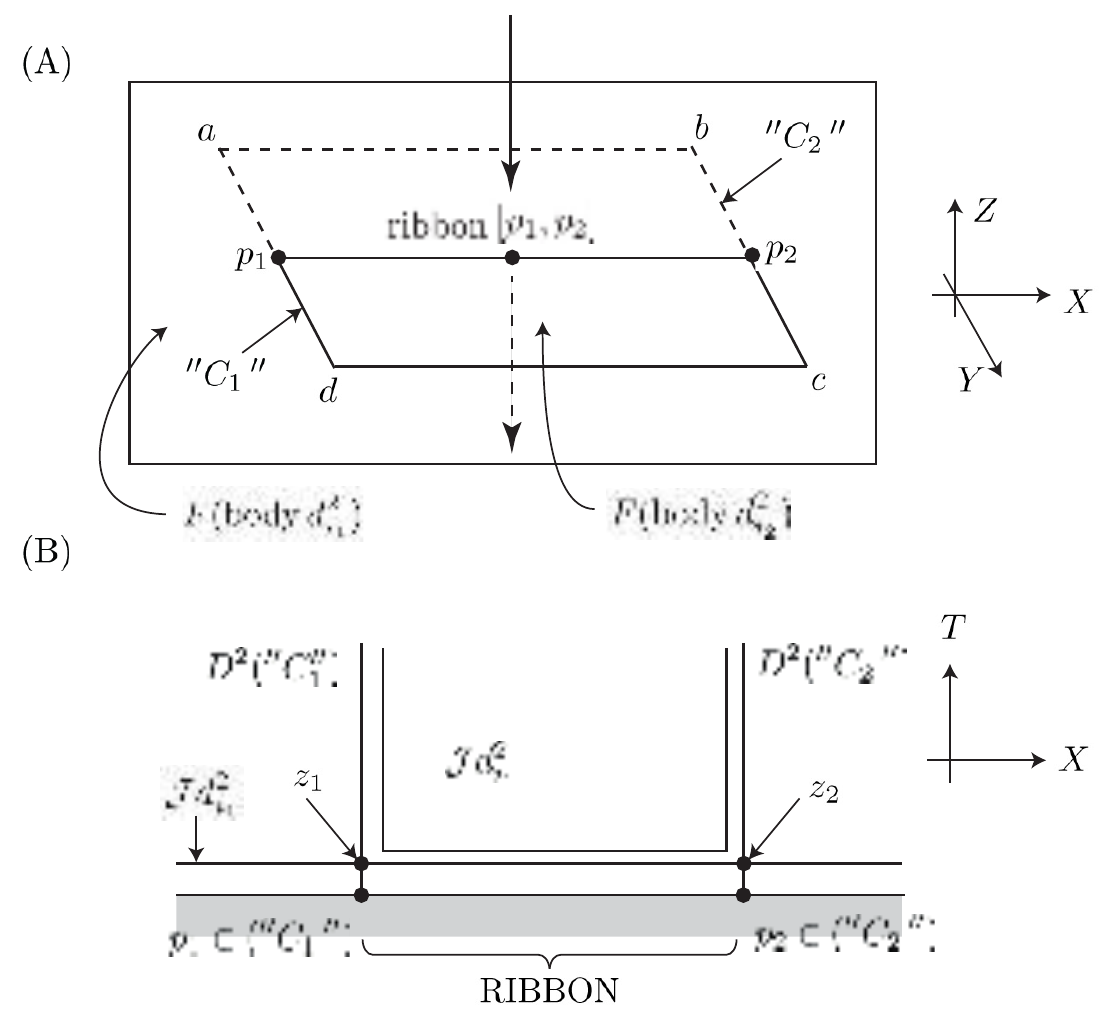}
$$
\label{fig21}
\centerline {\bf Figure 21.} 

\smallskip

\centerline{A RIBBON of the map (88). The vertical arrow, crossing the RIBBON is part of a green arc.}

\bigskip

\noindent {\bf More explanations for the Figures 20 and 21.} One goes from (88) to the corresponding part of (86), by replacing each $F ({\rm body} \ d_k^2)$ by $F ({\rm body} \ d_k^2) \times [0,\varepsilon_k]$ where $F ({\rm body} \ d_k^2) \times \{0\} = F ({\rm body} \ d_k^2)$ and the $[0,\varepsilon_k]$ is outgoing from $\partial N^4 (2X_0^2)$ (direction $[0,1] \subset$ coordinate $T$). The $\varepsilon_k$, depending on $k$ is a {\ibf free choice} which we can make, normally at least, as we please, according to our convenience. Whatever this choice is, at the level of Figure 20-(B) (CLASP), there {\ibf has} to be a double point $x \in JM^2(J)$. But then, in Figure 20-(B) we did make a {\ibf choice} $\varepsilon_{i_2} > \varepsilon_{i_1}$ $(> 0)$ and {\ibf if $z_1$ possesses an extended cocore $(z_1)$} (which is exactly the case when our ``$C_1$'' is a $C_i$ and {\ibf not} a $\Gamma_j$), then as we have suggested in Figure 23 in dotted lines, there is now a procedure for {\ibf killing both $z_1$ and $x$}, by pushing $Jd_{i_1}^2$ over the $\{$extended cocore $(z_1)\}^{\wedge}$. BUT, even if $z_2$ also possesses an extended cocore $(z_2)$, we still cannot apply this procedure both for $z_1$ and $z_2$. The so-called {\ibf ``green arc process''} will be necessary then.

Figure 21-(B) plays the same role with respect to 21-(A) as the 20-(B) with respect to 20-(A). But now, with the choice $\varepsilon_{i_2} > \varepsilon_{i_1}$ which be made (in Figure 21-(B)), there are no double points $x \in JM^2 (J)$ attached to the ribbon. With the choice $\varepsilon_{i_2} < \varepsilon_{i_1}$ there are two such. A priori, both choices are locally possible and it is the global policy of getting rid of accidents, AND of CONFINING the green arcs inside $\partial N_+^4 (\Gamma_1 (\infty))$, which will determine which choices are appropriate and which are not. Notice, also, that while in Figure 20-(B) we have $D^2(\mbox{``$C_2$''}) \subset (Y,T)$, $D^2(\mbox{``$C_1$''}) \subset (Z,T)$, in 21-(B) we find
$$
D^2(\mbox{``$C_1$''}) \subset (Y,T) \supset D^2(\mbox{``$C_2$''}).
$$
End of explanations. $\Box$

\bigskip

$$
\includegraphics[width=135mm]{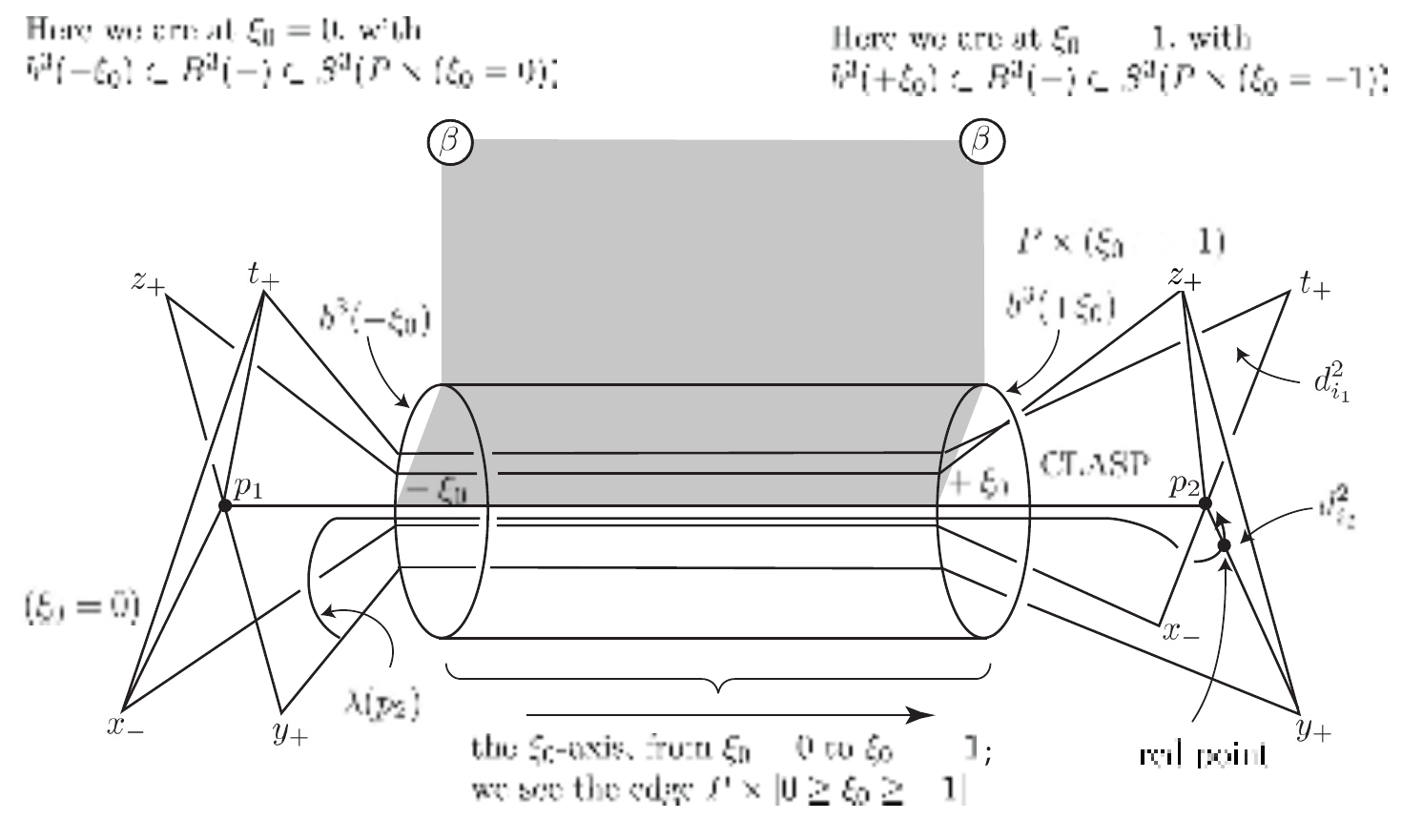}
$$
\label{fig22}
\centerline {\bf Figure 22.} 
\begin{quote}
We see here in the NORMAL situation, the only clasp which (88) may have, occurring between the pair of adjacent vertices $P \times (\xi_0 = 0)$ and $P \times (\xi_0 = -1)$, with the same $P \in \Gamma (1) \subset \Delta^2$. Here to begin with, at $\xi_0 = -1$ for illustration we have Figure 9-(B$_+$) for $P \times (\xi_0 = -1)$, with a $b^3 (+ \,\xi_0)$ (which gets added when we go to $X_0^2 [{\rm new}] \subset 2X_0^2$), living closer, visually, to the observer, i.e. higher, than everything else, seable in the Figures 9. Also, with complete details, the fully embellished Figure 9-(B$_+$) occurs as in Figure 24-(A). This accounts for the situation at $p_2 \in \Delta^2 \times (\xi_0 = -1)$. The situation at $\xi_0 = 0$ is to be discussed in the main text. A priori, making use of our free choices, it would be possible to change this CLASP into a RIBBON. But we could not live with both CLASPS and RIBBONS along the various $P \times [0 \geq \xi_0 \geq -1]$. So, our free choices at $\xi_0 = 0$, will be used to have consistently only CLASPS.

LEGEND: ${\bm -\!\!\!\!\!\bm -\!\!\!\!\!\bm -\!\!\!\!\!\bm -\!\!\!\!\!\bm -\!\!\!\!\!\bm -\!\!\!\!\!\bm -} \, = d^2 (\mbox{space-time, space-time, space-time})$;

$\longrightarrow\!\!\!-\!\!\!-\!\!\!- =$ green arc $\lambda (p_2)$. Here $\lambda (p_2) \subset \{$the $[y_+ , z_+]$ branch which cuts transversally  the $\Gamma_j$ at $p_2\}$ (if we think of it inside body $d^2$). The red intersection point
$$
\lambda (p_2) \pitchfork \{\mbox{RIBBON starting at $p_2$}\}
$$
is harmless, since the ribbon in question will receive the treatment of Figure 28-bis. In particular, at the RHS of our figure, we have $\lambda (p_2) \subset d_{i_2}^2$, coming with $\varepsilon (d_{i_2}^2) < \varepsilon$ $(d_{i_3}^2 \equiv$ the branch of the RIBBON $p_2 \to y_+$). So, the lifted $\lambda (p_2)$ flows under the lifted $d_{i_3}^2$.
\end{quote}

\bigskip

4) {\it We will make use of the following {\ibf crossing-freedom} when we will put up the completely embellished Figures $9$. To begin with, edges starting at $b^3 (\pm \, \xi_0)$ should {\ibf never} be cut by triangles $d^2 (\mbox{\it space-time, space-time,}$ $\mbox{\it space-time})$. Moreover, edges starting at $b^3 (\beta)$ should not be cut by anything at all. Figure {\rm 24-(A)} should illustrate these things.

\smallskip

With this, all the CLASPS and RIBBONS of {\rm (88)}, source of the accidents {\rm (86)}, are of the following types.

\smallskip

[And retain also that, since the $b^3 (\beta)$ is present only for the part ${\mathcal B}_i^2 \subset \delta_i^2$, the {\rm (87), (88)}, which only concern the $d_k^2$-parts, do not contain any $\beta$-contribution.]}

\medskip

4.1) {\it CLASPS internal to $\Delta^2 \times [0 \geq \xi_0 \geq -1]$, occurring exactly along SOME $P \times [0 \geq \xi_0 \geq -1]$'s. Here is the mechanism which creates these CLASPS. To begin with, as a conseqence of our constant use of the GPS structures, in any given Figure $9$ {\ibf  there is at most one crossing of space-time lines}. With this, consider for illustration the crossing $[y_+ , z_+] ({\rm DOWN}) / [x_- , t_+] ({\rm UP})$ in Figure {\rm 9-(B$_+$)} and occurring also in the RHS of Figure $22$. At this point, remember that at $\xi_0 = -1$ we have to add an axis $+ \, \xi_0$ in addition to space-time and at $\xi_0 = 0$ an axis $- \, \xi_0$.

\smallskip

For $+ \, \xi_0$ we have Figure {\rm 24-(A)} as it stands. When it comes to the axis $- \, \xi_0$ at $\xi_0 = 0$, here is what goes on. We certainly have our crossing $[y_+ , z_+] ({\rm DOWN}) / [x_- , t_+] ({\rm UP})$, which {\ibf stays put}, and for $b^3 (- \, \xi_0)$ we chose a position which is the MIRROR IMAGE with respect to $S^2 (P)$ (see Figure {\rm 8-(C)}), of the position of $b^3 (+ \, \xi_0)$ in Figure {\rm 24-(A)}, like for $t_{\pm}$ in Figure {\rm 8-(A)}.

\smallskip

So $b^3 (- \, \xi_0)$ (contrary to $b^3 (+ \, \xi_0)$ in Figure {\rm 24-(A)}) is now lower with respect with anything in space-time. In terms of Figure {\rm 24-(A)} this amounts to the following changes: {\rm i)} a notational change $+ \, \xi_0 \Rightarrow - \, \xi_0$, {\rm ii)} an allowed crossing, performed at each site with a circle drawn as a broken dark line, performed at the level of Figure $24$,
$$
\includegraphics[width=7cm]{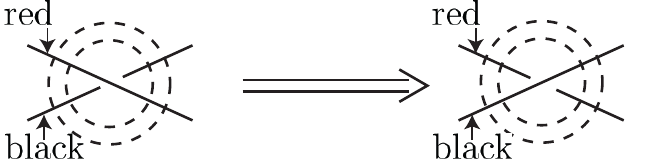}
$$
with RED $= b^3 (- \, \xi_0)$-line, BLACK $=$ space-time line.

\smallskip

It can be checked that this produces the CLASP from Figure $22$.}

\medskip

4.2) {\it All the RIBBONS of $(88)$ occur at $\xi_0 = -1$ AND at $\xi_0 = 0$. In terms of Figure $22$, RIBBONS at $\xi_0 = -1$ go along two neighbouring $p_2$'s. They corresponds to things like $[p_2 , z_+]$, $[p_2 , y_+]$. Dually, RIBBONS at $\xi_0 = 0$ go along two neighbouring $p_1$'s and, again, they correspond to things like $[p_1 , t_+]$, $[p_1 , x_-]$. So, all the accidents of $(88)$ are concentrated at $\Delta^2 \times [0 \geq \xi_0 \geq -1]$.}
End of Lemma 11.1.

$$
\includegraphics[width=125mm]{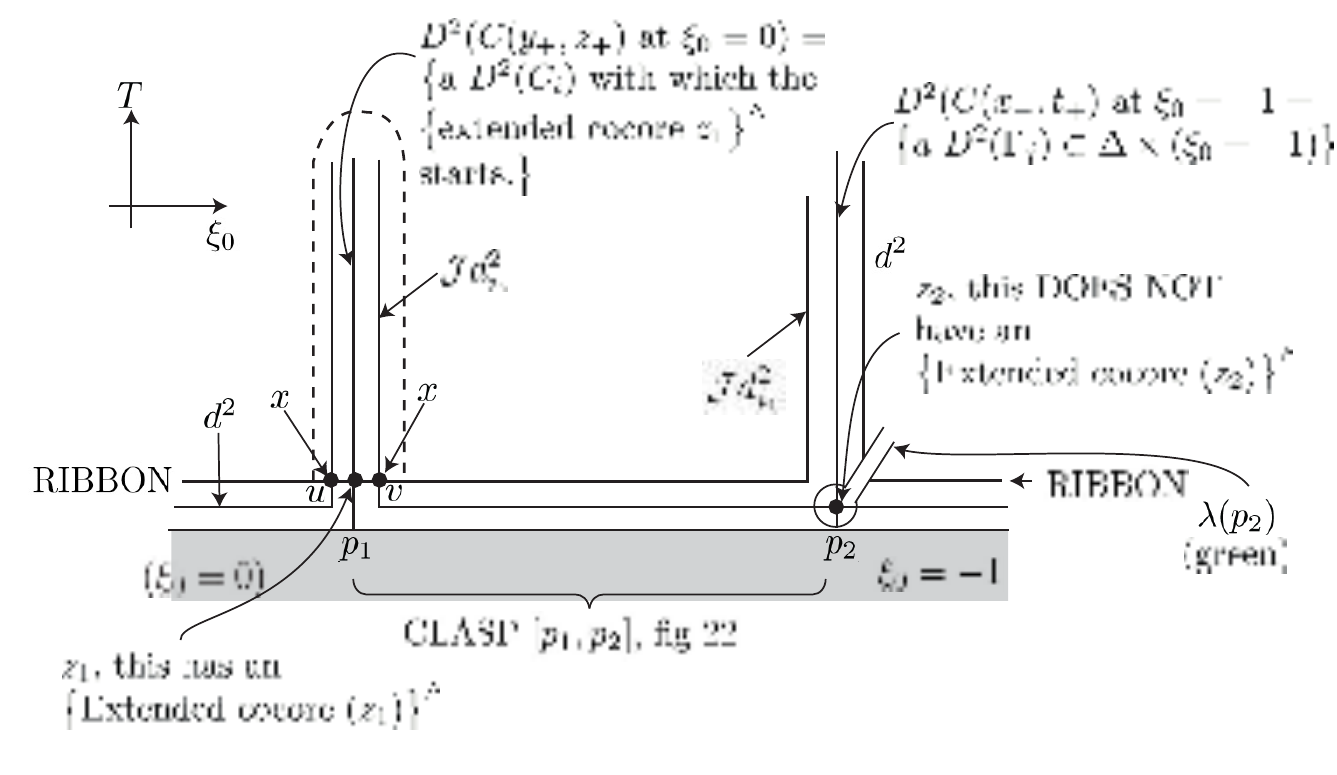}
$$
\label{fig23}
\centerline {\bf Figure 23.} 
\begin{quote}
The real-life effect of Figure 22. All the  ACCIDENTS of $Jd^2$ (86) occur like in this figure, for the various clasps along $[0 \geq \xi_0 \geq -1]$, which are 2-by-2 disjoined, and see here 4.1) in Lemma 11.1, and also 4.2) from the same lemma. The RIBBONS at $\xi_0 = -1$, respectively at $\xi_0 = 0$, form each of them, a system which, IF each $P \times [0 \geq \xi_0 \geq -1]$ would carry a clasp like in Figure 22 would be connected. But, in real life, since there are not necessarily CLASPS along each $P \times [0 \geq \xi_0 \geq -1]$, this breaks into several connected systems, each of them possibly highly non-simply connected. The  clasps along $[p_1 , p_2]$ connect the two multi-systems.

In our present figure, at $p_1$ both $z_1$ and the $x$'s can be demolished by pushing $Jd_{i_1}^2$ over the $\{$extended cocore $(z_1)\}^{\wedge}$. Of course, at $p_2$, such a cocore is not available. But notice that, even if it would be, something else would still be needed. One cannot use a push over the $\{$extended cocores$\}^{\wedge}$, at {\ibf both} ends of a CLASP, even when such extended cocores do exist.
\end{quote}

\bigskip

\noindent {\bf Remarks.} A) In order to define our $\underset{1}{\overset{M}{\sum}} \, \delta_i^2$, we need the BLUE order of $X^2 [{\rm new}]$, and not the BLUE order of $2X^2$. In the context of $2X^2$, on the side $X_0^2 \times r = X_0^2 [{\rm new}]$, all the $2^{\rm d}$ cells are $D^2 (\gamma^1)$'s BLUE-wise, with the real BLUE order transferred on the $b$-side. To define the $\delta_i^2$'s, we need the DOUBLE $2X_0^2$, in particular the (77.2)'s. Also, in this same context, all the $B_i^2$'s occurring inside $\underset{1}{\overset{M}{\sum}} \, \delta_i^2$ are exactly the $B_j^2$'s produced by the BLUE $X^2 [{\rm new}]$-trajectories of the $\underset{1}{\overset{M}{\sum}} \, b_i$; see (77.1).

\medskip

B) So far, we have only looked at the RIBBONS ${\rm body} \, d^2 \pitchfork {\rm body} \, d^2$, but there will also be RIBBONS $\{$body ${\mathcal B}^2$ (long branch) $\pitchfork$ body $d^2$ (short branch)$\}$.

\medskip

C) When it comes to destroying the ACCIDENTS, it will be convenient to have only CLASPS along the $P \times [0 \geq \xi_0 \geq -1]$ and NO RIBBONS going from $P \times (\xi_0 = 0)$ to $P \times (\xi_0 = -1)$. End of Remarks.

\bigskip

Each vertex $P \in 2X_0^2$ comes with an ``embellished Figure 9'', actually Figure 24 and the purely space-time part of these figures (that is what one actually sees in Figures 9), lives inside $B^3 (-)(P)$. But, of course, the complete, real-life, embellished Figures 9 contain the additional items too, namely the following.
\begin{enumerate}
\item[$\bullet$)] A $b^3 (\beta) \subset B^3 (+)$, present for all $P$'s.
\item[$\bullet$$\bullet$)] A $b^3 (+ \, \xi_0$ or $- \, \xi_0) \subset B^3 (-)$ for $P \times (\xi_0 = -1$ or $\xi_0 = 0)$. At $\xi_0 = -1$ this is placed higher than the space-time items, and at $\xi_0 = 0$ lower.
\item[$\bullet$$\bullet$$\bullet$)] In any of the embellished Figures 9 one has to add edges $b^3 (\beta) -\!\!\!-\!\!\!-\!\!\!-\!\!\!- b^3 (\mbox{space-time})$ and at $\{\xi_0 = -1 , \xi_0 = 0\}$ one also adds edges $b^3 (\pm \, \xi_0) -\!\!\!-\!\!\!-\!\!\!-\!\!\!- b^3 (\beta$ or space-time).
\item[$\bullet$$\bullet$$\bullet$$\bullet$)] We have triangles $d^2 (\pm \, \xi_0$, space-time, space-time), $B^2 (\beta, \pm \, \xi_0$, space-time), $B^2 (\beta , \mbox{space-time}$, space-time) in addition to any purely space-time triangle in Figures 9. At $\xi_0 = -1$, the interiors of the triangles $B^2 (\beta , \ldots, \ldots)$ are {\ibf void}, like in the Figures 7-bis and 19.1.
\end{enumerate}

\medskip

Using the {\ibf crossing-freedom}, we may impose the following

\bigskip

\noindent (90) \quad No triangle $d^2 (\ldots , \ldots , \ldots)$ at all cuts through an edge starting at $b^3 (\beta)$. Here $\ldots$ stands for $\{\pm \, \xi_0$ or space-time$\}$. So, there are no CLASPS containing $b^3 (\beta)$ and, when Figure 21 is not connected with $\Delta^2 \times [0 \geq \xi_0 \geq - 1]$, then we have body ``$d_{i_1}^2$'' $\subset {\mathcal B}^2$, and the curve it goes through transversally is $C$ (space-time, space-time) OR $C(\pm \, \xi_0 , \mbox{space-time})$ (see, for this last item the doubly circled crossings in Figure 24).

\bigskip

\noindent (91) \quad There are NO TRIPLE POINTS FOR the map $J$ in (80), i.e.
$$
M^3 \left( J \, \biggl\vert \, \sum_1^M \delta_i^2 \right) = \emptyset .
$$

\noindent [Proof of (91). We know this already for the restriction $J \, \bigl\vert \, \underset{k}{\sum} \, d_k^2$. With everything we have already said, the triple points could only be produced by something of the following type, with possibly the $d^2 (\mbox{space-time})$ replaced by another $B^2 (\beta , \pm \, \xi_0$ and/or space-time)
$$
d^2 (\pm \, \xi_0 , \, \mbox{space-time}) \pitchfork d^2(\mbox{space-time}) \pitchfork B^2 (\beta , \pm \, \xi_0 \ \mbox{and/or space-time}).
$$
In order to conjure these out of existence we will make use, according to the case, either of the STRATEGIC DECISION 4) in the Lemma 5, which comes with the Figures 7-bis and 19.1, OR of the crossing freedom when outside of $\Delta^2 = \Delta^2 \times (\xi_0 = -1)$.

\smallskip

With $b^3 (+ \, \xi_0)$, at $P \times (\xi_0 = -1)$, such triple points do not exist, because the $B^2$-triangles have now no interior; see our Figures 7-bis and 19.1. So we look at $P \times (\xi_0 = 0)$, and now $b^3 (- \, \xi_0)$ is the MIRROR-IMAGE of the $b^3 (+ \, \xi_0)$ from Figure 24-(A). The $b^3 (\beta)$-lines continue to be ``over the table'', just like in Figure 24, while now the $b^3 (- \, \xi_0)$-lines are ``under the table''.

\smallskip

With this, the double-line $d^2 (- \, \xi_0 , x_- , t_+) \cap d^2 (- \, \xi_0 , z_+ , y_+)$ (and there are no other lines $d(-\xi_0 ,\mbox{space-time},$ $\mbox{space-time}) \cap d^2 (- \, \xi_0 , \mbox{space-time}, \mbox{space-time})$) cannot be cut by any $B^2 (\beta , \mbox{space-time}, \mbox{space-time})$. And it can be seen, directly, that it is not cut by $B^2 (\beta , - \, \xi_0 , \mbox{space-time})$ either.

\smallskip

This proves our (91). $\Box$]

\bigskip

\noindent {\bf Complement to Lemma 11.1.} {\it We can extend Lemma {\rm 11.1} from $\sum d_k^2$ to the whole of $\sum \delta_i^2$, where
$$
\delta_i^2 = {\mathcal B}_i^2 \cup \sum d_k^2\mbox{'s, see {\rm (79.1)}.}
$$
To begin with, there is a decomposition
$$
{\mathcal B}_i^2 = \mbox{body} \, {\mathcal B}_i^2 + ({\mathcal B}_i^2 - \mbox{body} \, {\mathcal B}_i^2),
$$
with a not everywhere well-defined map like in {\rm (87)}
$$
\sum_1^n {\mathcal B}_i^2 \xrightarrow{\qquad F \qquad } \partial N^4 (2X_0^2) , \leqno (91)
$$
coming with the analogue of {\rm (88)}, namely with
$$
\sum_1^n \mbox{body} \, {\mathcal B}_i^2 \xrightarrow{\qquad F \qquad } \partial N^4 (\Gamma_1 (\infty)) , \leqno (92)
$$
with $F \mid ({\mathcal B}_i^2 - \mbox{body} \, {\mathcal B}_i^3)$ very much like the $F \mid (d_k^2 - \mbox{body} \, d_k^2)$ at {\rm 2)} in Lemma {\rm 11.1}.}

\medskip 

1) {\it Once nothing cuts through edges of $\partial N^4 (P)$ having $b^3 (\beta)$ at one end (see Figure {\rm 24}), all the punctures of {\rm (92)} concern curves $C (\pm \, \xi_0 , \mbox{space-time})$ or $C(\mbox{space-time, space-time})$. These come with RIBBONS only, more precisely
$$
\mbox{RIBBONS} \, \{\mbox{body} \, {\mathcal B}^2 \} (\mbox{long branch}) \pitchfork \{\mbox{body} \, d^2\} (\mbox{short branch}) \ \mbox{OR RIBBONS} \ \{\mbox{body} \, {\mathcal B}^2\} \pitchfork \{\mbox{body} \, {\mathcal B}^2\}.
$$

All the double points and the transversal contacts when {\rm (91)} is involved, come from the RIBBONS of {\rm (92)} (which has no other double points to which it participates).}

\medskip

2) {\it When we go to $(J(80)) \ \Bigl\vert \ \underset{1}{\overset{M}{\sum}} \, {\mathcal B}_i^2$, this comes with the following ACCIDENTS (part of {\rm 3)} Lemma {\rm 10}): double points $x \in J{\mathcal B}^2 \pitchfork Jd^2$ and transversal contacts $z \in J{\mathcal B}^2 \pitchfork (2X_0^2 - \Delta^2 \times (\xi_0 = -1))$.

\smallskip

To these, a priori, one might have to add double point $x \in J{\mathcal B}^2 \pitchfork J{\mathcal B}^2$ produced by the internal RIBBONS of $\{\mbox{body} \, {\mathcal B}^2\}$. When the detailed treatment of all the RIBBONS will have been reviewed we will see what happens with these potential double points.} End of the complement to Lemma 11.1.

\bigskip

\noindent {\bf Remark.} There is a good a priori reason why $M^2 \left( J \ \Bigl\vert \ \underset{1}{\overset{M}{\sum}} \, {\mathcal B}_i^2 \right) = \emptyset$. Inside each individual Figure 24 (i.e. for each individual $\partial N^4 (P)$), there is a unique $b^3 (\beta)$ spanning a number of distinct $B^2 (\beta , \ldots)$-triangles which, except for common edges (and the common vertex at $b^3 (\beta)$) are disjoined. So, the only potential source for $M_2 \left( J \ \Bigl\vert \ \underset{1}{\overset{M}{\sum}} \, {\mathcal B}_i^2 \right)$ is dry. End of Remark.

\bigskip

What we see now at the level of the Figure 24 below and its likes, happens at the level of (88) and (92). This figure is actually a fully embellished version of Figure 9-$(B_+)$. It exhibits triangles of types $d^2$ and ${\mathcal B}^2$, all being part of the corresponding body ($\ldots$)'s.

$$
\includegraphics[width=105mm]{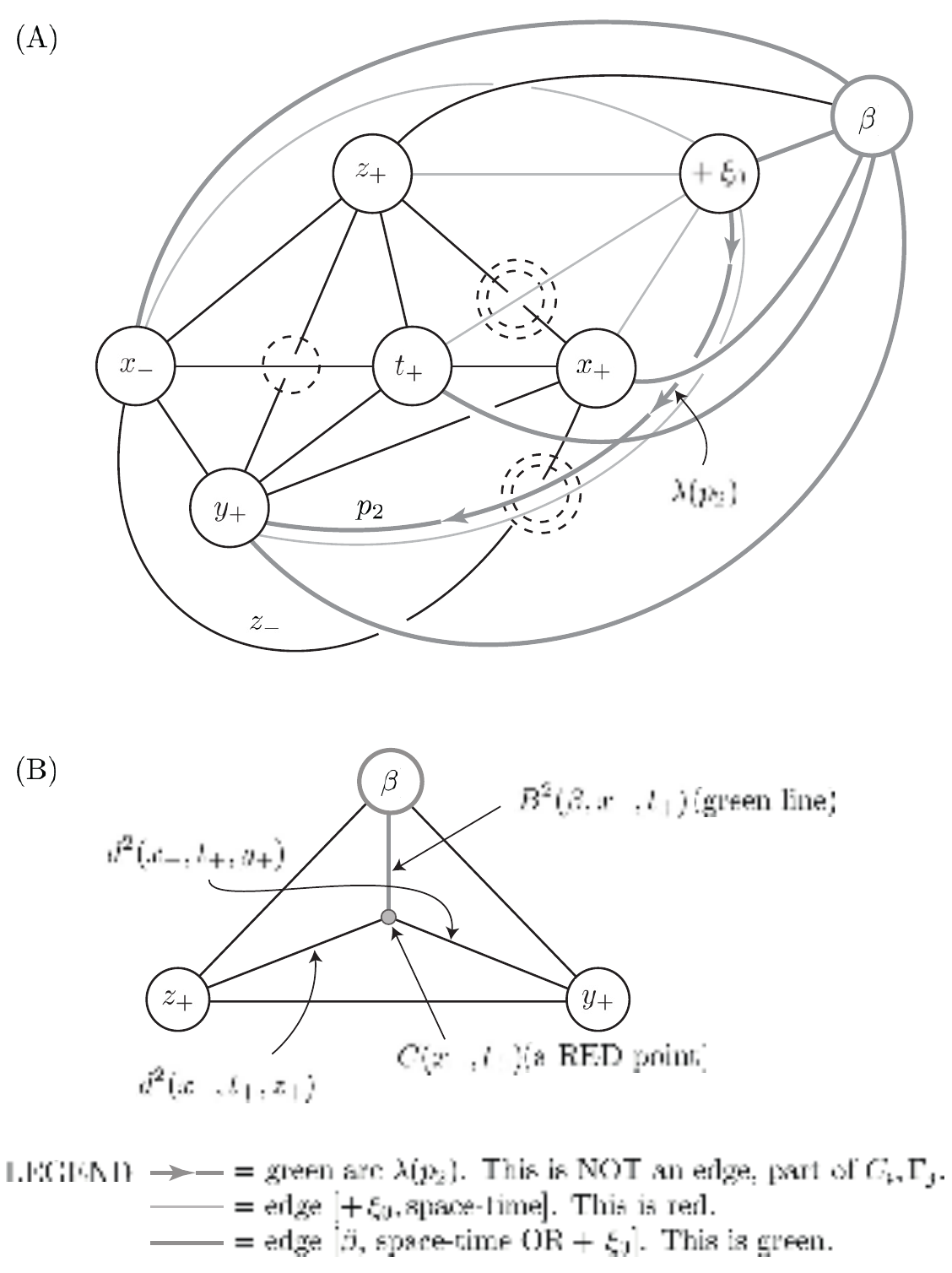}
$$
\label{fig24}
\centerline {\bf Figure 24.} 
\begin{quote}
In (A) we present the completely embellished Figure 9-$(B_+)$, in the context of $P \times (\xi_0 = -1) \in \Delta^2$. Then the $B^2$-triangles have a void interior (see here the Figure 22).

The Figure (B) concerns the context of $P \times (\xi_0 = 0)$ and it presents the $B^2$-triangle $B^2 (\beta , y_+ , z_+)$. [IF this $B^2$ would be with its interior alive, present at $P \times (\xi_0 = -1)$, we would see exactly the same picture; when we go from $\xi_0 = -1$ to $\xi_0 = 0$, it is only $+ \, \xi_0$ which  flips to the mirror-position $- \, \xi_0$.]
\end{quote}

\bigskip

\noindent Here are {\bf some explanations concerning Figure 24.} The {\ibf edges} are occurring, according to the case in black ($=$ pure space-time), red ($=$ one vertex $\pm \, \xi_0$ but no $\beta$), green ($=$ one vertex $\beta$), BUT the green arc marked $\lambda (p_2)$ is {\ibf not} such an edge. The green edges do not cut transversally (like $R \pitchfork R^2 \subset R^3$) through the space-time triangles $d^2 (\mbox{space-time}, \mbox{space-time}, \mbox{space-time})$. The (A) lives in $S^3 (p) = B^3 (-) \underset{S_{\infty}^2}{\cup} B^3(+)$. It completes Figure 9-$(B_+)$, and here is how to read it, starting from the Figure 9-$(B_+)$, in the situation $P \times (\xi_0 = -1)$, or $P \times (\xi_0 = 0)$. Add, inside this same $B^3 (-)$, the $b^3 (+ \, \xi_0)$ higher than anything else (i.e. closer to the observer). This gives an embellished Figure 9-$(B_+)$ which we will flatten (i.e. the space-time, $+ \, \xi_0$ lines) like in a link diagram, with crossings UP/DOWN. Next one adds the $b^3 (\beta) \subset B^3 (+)$ and suspend from it, with green lines, the now conveniently flattened $\{$embellished Figure 9-$(B_+)\}$.

\smallskip

This makes sure that nothing can cut transversally through the $\beta$-lines.

\smallskip

All this should help vizualizing things. Incidentally, I am using here the word ``suspending'' in a rather cavalier manner; I should have said ``coning'', instead. At $P \times (\xi_0 = -1)$ the interiors of the $B^2$-triangles are removed, and we turn to $P \times (\xi_0 = 0)$. Then our $+ \, \xi_0$ becomes $- \, \xi_0$, in a MIRROR-IMAGE POSITION, with respect to what we see in Figure 24-(A).

\smallskip

One may notice that the $B^2$-triangles never cut the lines $C (\pm \, \xi_0 , \mbox{space-time})$. [We only have to consider $P \times (\xi_0 = 0)$ with $b^3 (- \, \xi_0)$. With $b^3 (- \, \xi_0)$ in its MIRROR IMAGE location, ``under the table'', we clearly do not have
$$
B^2 (\beta , \ldots) \pitchfork C(- \, \xi_0 , \mbox{space-time}).]
$$

Figure (B), with a red transversal contact $B^2 (\beta , \ldots ) \pitchfork C(\mbox{space-time}, \mbox{space-time})$, corresponds to the unique crossing of black lines, surrounded by a dotted loop, in (A). And here again, is the move (A) $\Rightarrow$ (B). Imagine the BLACK and RED lines in (A) displayed ``on the table'', like in a link diagram. The green $b^3 (\beta)$ lives at infinity and the $B^2 (\beta , \ldots)$ triangles, like $B^2 (\beta , x_- , t_+)$, $B^2 (\beta , y_+ , z_+)$ are now vertical planks, based on the $[x_- , t_+] , [y_+ , z_+]$ and going to $\infty^{\rm ty}$.

\smallskip

When the ``talk at $\infty^{\rm ty}$'' where $\xy *[o]=<15pt>\hbox{$\beta$}="o"* \frm{o}\endxy$ lives, gets contracted to a point, then we get our Figure (B) and its likes. This is, of course, equivalent to the coning description given earlier.

\bigskip

\noindent {\bf Lemma 11.2.} 1) {\it Using the crossing-freedom we can make that all the transversal contacts $\mbox{body} \, {\mathcal B}_i^2 \cap \mbox{body} \, {\mathcal B}_j^2$ are RIBBONS. Typically, one starts from $B^2 (\beta , z_+ , y_+) \cap B^2 (\beta , x_- , t_+)$ in Figure {\rm 24-(B)}, and then one continues along $P \times [r,b]$, from $r$ to $b$, like in Figure $25$. Here, at the level of the $b$ of $2X_0^2$ we are far from $\Delta^2 \times (\xi_0 = -1)$ and we can make full use of the crossing freedom.}

\medskip

2) {\it From {\rm (90)} it follows that all the contacts $\mbox{body} \, {\mathcal B}_i^2 \cap \mbox{body} \, d_k^2$ are RIBBONS like in Figure {\rm 21-(A)}, coming always with $d_{i_1}^2 (\mbox{long branch}) \Rightarrow {\mathcal B}_i^2$ and $d_{i_2}^2 (\mbox{short branch}) \Rightarrow d_k^2$. RIBBONS like in {\rm 1)} above and like in the present {\rm 2)} are not isolated but they, generally speaking, form highly connected and non-simply connected {\ibf nets}. This means the following FACT:}

\bigskip
 
\noindent (92.1) \quad {\it When, for the map
$$
Z^2 \equiv \sum_1^M \mbox{body} \, {\mathcal B}_i^2 \cup \sum_1^N \mbox{body} \, d_k^2 \xrightarrow { \ \ G \ \ } \partial N^4 (\Gamma_1 (\infty)),
$$
and for some transversal contact ($=$ puncture)
$$
p \in \{ \mbox{Curve $C$ from $\{$link$\}$ {\rm (9)}} \mid X^2 [{\it new}] \} \pitchfork G (\mbox{body} \, {\mathcal B}^2 \cup \mbox{body} \, d^2) ,
$$
we look for an arc $\lambda$, embedded inside the source of $G$, connecting $p$ to the $\partial {\mathcal B}^2$ and free of accidents, then we may be OBSTRUCTED by contacts $\lambda \pitchfork {\it RIBBON}$ (and never by CLASPS, see here Figure {\rm 22}).}

\medskip

3) {\it Let $e = [P_1 , P_2]$ be an edge NOT living at $\xi_0 = -1$. In Figure $7$ we see attached to $e$ a rectangle (possibly with only a boundary collar alive). At the level of $2 X_0^2$, we can always find inside this ``rectangle'' a much thinner one, $e \times [r,\beta]$ (not affected by the deletion in Figures {\rm 7-(II, III)}, BUT non existent (wiped out), with the exception of $(e \times r) \cup \{ P_1 , P_2 \} \times [r,\beta]$ in Figure {\rm 7.bis}, at $\xi_0 = -1$. Changing now the topic, we move next to the source $Z^2$ of $G$, see {\rm (91)}, and we assume that for some $u \in (x,y,z,t,\xi_0)$, in the $B^3 (-) \subset S^3 (P_1)$, in the corresponding Figure $9$, we find a $b^3 (\pm \, u)$. If $P_2$ is a vertex adjacent to $P_1$ in $2X_0^2$, along the $u$-axis, then inside $B^3 (-) \subset S^3 (P_2)$ we have to find $b^3 (\mp \, u)$.}

\bigskip

Moreover, as a reflex of our thin rectangle above, when we are OUTSIDE of $\xi_0 = -1$, then for every such $u$, we find at the level of $Z^2$, source of $G$ (91), living inside $\partial N^4 (\Gamma_1 (\infty))$, a thin rectangle $[\beta (P_1) , \beta (P_2) ; u_{\pm} (P_1) , u_{\mp} (P_2)]$, going along it. (This could be the shaded rectangle in Figure 22, for instance.

$$
\includegraphics[width=145mm]{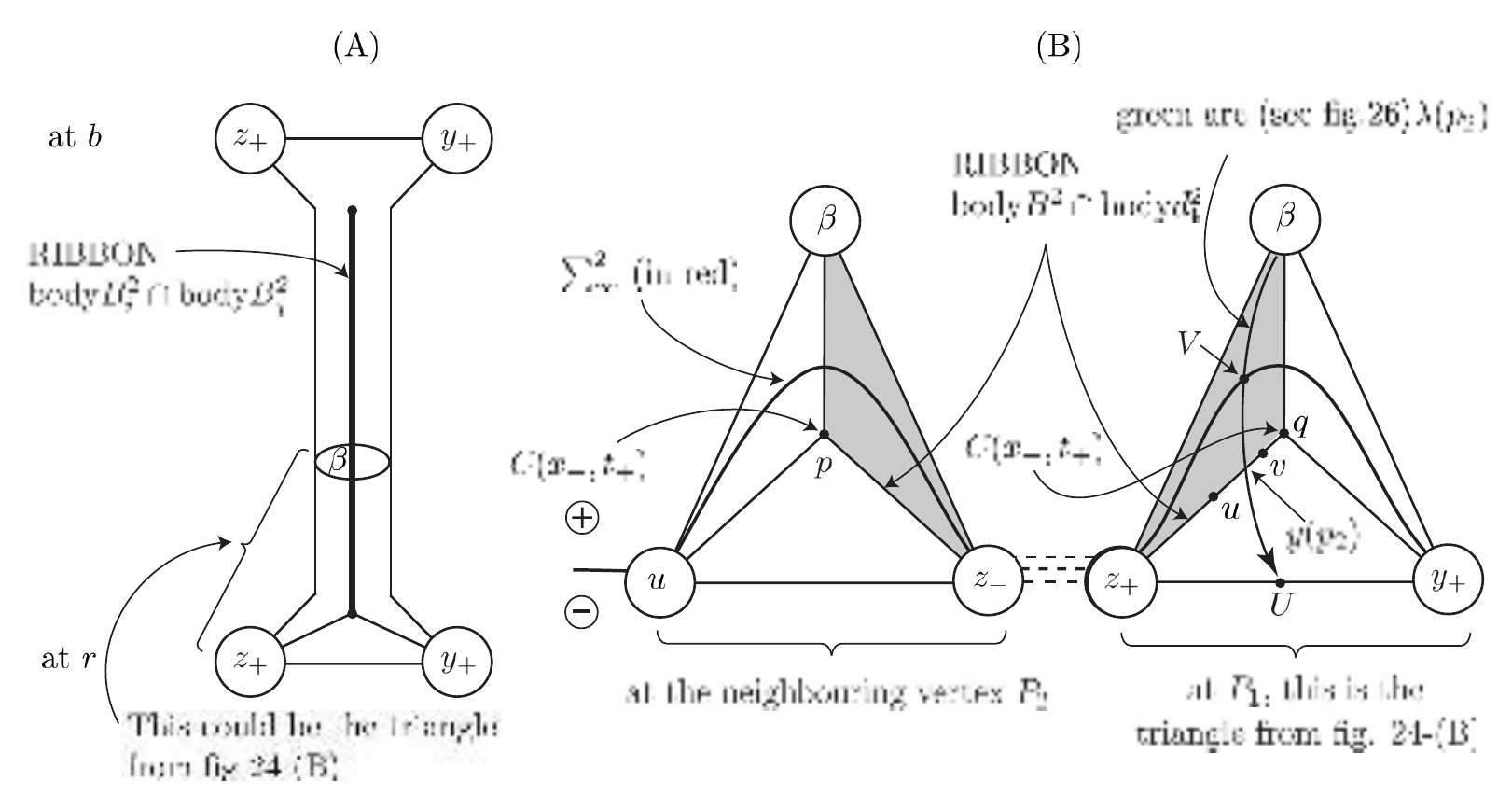}
$$
\label{fig25}
\centerline {\bf Figure 25.} 
\begin{quote}
This figure explains the occurrence of the two types of RIBBONS for $Z^2 \xrightarrow{ \ G \ } \partial N^4 (\Gamma_1 (\infty))$ (91), namely RIBBONS ${\mathcal B}^2 \pitchfork {\mathcal B}^2$ in (A) and RIBBONS ${\mathcal B}^2 \pitchfork d^2$ in (B).

In Figure (B) the two shaded triangles are adjacent to the rectangle $[\beta (P_1) , \beta (P_2) ; z_+ (P_1)$, $z_- (P_2)]$, which communicates with the rest of $\delta^2 \mid X^2 \times r$ via them. This {\ibf rectangle} in question, is like in 3) in the Lemma 11.2. Similarly, in Figure 22 we have shaded the rectangle
$$
[\beta (P \times (\xi_0 = -1)) , \beta (P \times (\xi_0 = 0)) ; + \, \xi_0 \ ({\rm at} \ \xi_0 = -1), - \, \xi_0 \ ({\rm at} \ \xi_0 = 0)].
$$

In (B), the L.H.S. triangle is the natural continuation of the one from the R.H.S., i.e. the one from Figure 24-(B). There is a {\ibf ``principle of continuity''}, meaning that in moving from $P_1$ to $P_2$ the ${\rm body} \, d^2$, ${\rm body} \, {\mathcal B}^2$ has to use the common $2^{\rm d}$ walls in $2X_0^2$, valid both for $P_1$ and for $P_2$. Figure (C) should illustrate this, and we see there the RIBBON $\mbox{body} \, {\mathcal B}^2 \cup \mbox{body} \, d_k^2$, in all its glory.
\end{quote}
$$
\includegraphics[width=65mm]{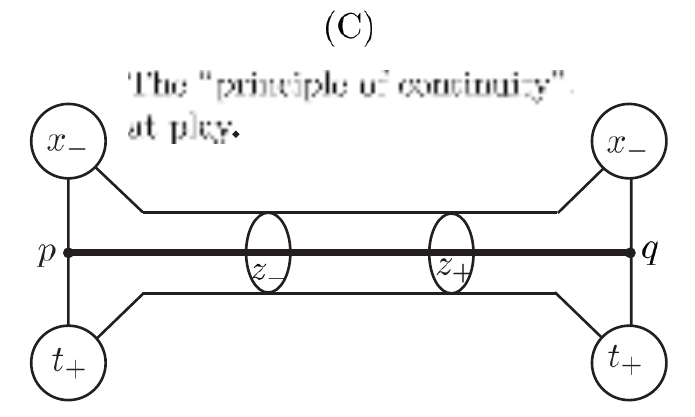}
$$
\label{fig25bis}
\begin{quote}
We see here a piece of $\mbox{body} \, d_k^2$, containing the RIBBON $[p,q] = \mbox{body} \, d_k^2$ $\cap$ ${\rm body} \, {\mathcal B}_i^2$. In $\partial N^4 (\Gamma_1 (\infty))$, where this figure lives, the branch ${\rm body} \, {\mathcal B}_i^2$ cuts transversally the branch $\mbox{body} \, d_k^2$, through $[p,q]$. Finally, Figure (D) should suggest how a triangle like $B^2 (\beta , y_+ , z_+)$ gets generated. In (D) this triangle is represented in green lines.
\end{quote}
$$
\includegraphics[width=9cm]{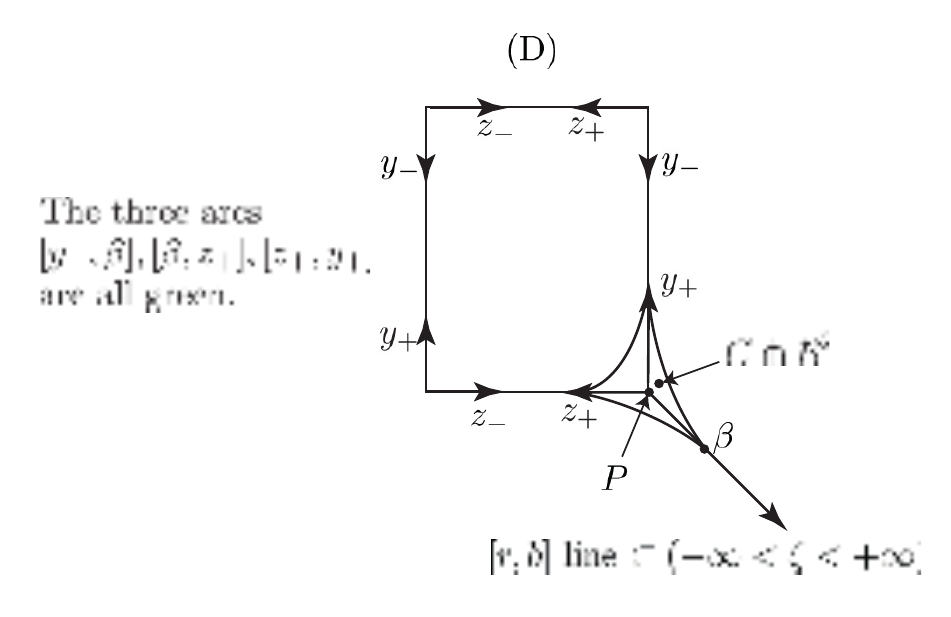}
$$
\vglue -8mm
\label{fig26bis}
\begin{quote}
End of Figure 25.
\end{quote}

\bigskip

The Lemma following next is a ``worst case scenario'', in the sense that we show how to handle all the accident which, with what was done so far in the present paper are a priori possible, not excluded by anything we have said explicitly, yet. It will turn out that the real-life situation will be slightly better than the worst case scenario below.

\bigskip

\noindent {\bf Lemma 11.3. (Destroying the accidents.)} 

\medskip

\noindent 1) {\it There are two procedures for destroying the accidents of $(80)$ which we shall use:}
\begin{enumerate}
\item[I)] {\it The push over the $\{$extended cocore $(z)\}^{\wedge}$ (and when this cocore is not there, our {\rm I)} certainly cannot be used, for instance when $z \in \Delta^2 \times (\xi_0 = - 1)$). Otherwise it is always there. Like in Figure {\rm 20-(B)}, this movement destroys both $z$ and the adjacent double point $x$. But we {\ibf cannot} use {\rm I)} for all the $z$ which come with a cocore.}
\item[II)] {\it The previous description of {\rm I)} was at the level of $N_1^4 (2X_0^2)^{\wedge}$, and we revert now to the $Z^2 \xrightarrow{ \ G \ } \partial N^4 (\Gamma_1 (\infty))$ from {\rm (92.1)} and consider a puncture
$$
p \in Z^2 \pitchfork \{\mbox{a curve} \ C \in \{\mbox{link}\}\} ,
$$
which in $4^{\rm d}$ generates an accident $z \in J\delta^2 \pitchfork 2X_0^2 \subset N_1^4 (2X_0^2)^{\wedge}$. If for $p$ we can find an embedded arc $\lambda \subset Z^2$, joining $\lambda$ to a point $q \in \partial {\mathcal B}^2 \approx \partial \delta^2$, then we can proceed as in Figure $26$, changing $Z^2 \xrightarrow{ \ G \ } \partial N^4 (\Gamma_1 (\infty))$, and accordingly the $(80)$ too, so as to destroy $z$. This modifies, of course $(\delta^2 , \partial \delta^2)$.}
\end{enumerate}
$$
\includegraphics[width=85mm]{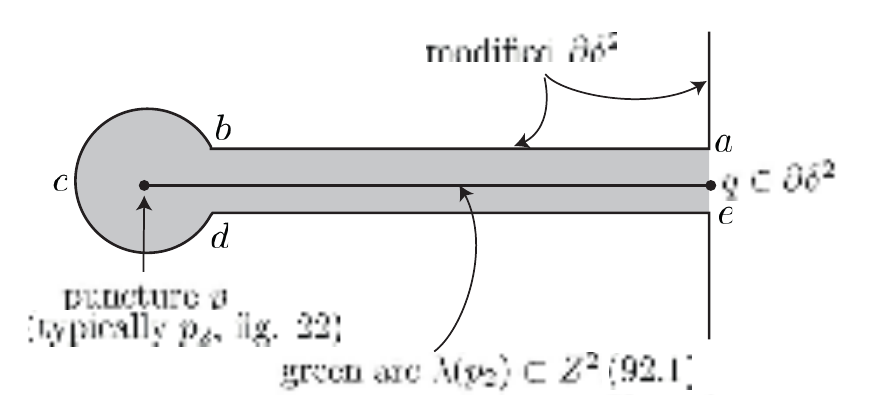}
$$
\label{fig26}
\centerline{\bf Figure 26.}
\begin{quote}
Pushing $\partial \delta^2$ along a green arc $\lambda$. In $4^{\rm d}$, this destroys an ACCIDENT $z$. Then, the closed loop going through $c,b,d$, around $p$, when lifted off $\partial N^4 (2\Gamma(\infty))$, by a multiplication with $+ \, \varepsilon$ (i.e. loop $\Rightarrow$ loop $x (+ \, \varepsilon)$), is linked with $2X_0^2 \subset N_1^4 (2X_0^2)$.
\end{quote}

\bigskip

{\it [With more detail, here is what goes on. Let $A^2 \supset \lambda$, with the $A^2 \subset Z^2$, being a connected piece, {\ibf without} internal accidents, and lifting $Z^2 \to \partial N^4 (\Gamma_1 (\infty))$ into $4^{\rm d}$ means, for $A^2$, an embedded tubular neighbourhood $A^2 \times [0 , \varepsilon = \varepsilon (A^2) > 0] \subset N_1^4 (2X_0^2)^{\wedge}$, with $A^2 \times \{ 0 \} = \{{\rm our} \ A\}$, with the rest $A^2 \times (0,\varepsilon]$ in $N_1^4 (2X_0^2)^{\wedge} - N^4 (2X_0^2)^{\wedge}$, and with $A^2 \{\varepsilon\} \subset JZ^2 \subset J\delta^2$. For our green arc $\lambda$ we consider the shaded neighbourhood $\nu^2 (\lambda) \subset A^2$, from Figure $26$, and its $\nu^2 (\lambda) \times [0,\varepsilon]$. (The $\varepsilon$ may be happily variable here, this does not matter.) The modified $J \delta^2$, with the accident $z$ destroyed, is gotten by deleting $\{$the contribution $\nu^2(\lambda) \times [0,\varepsilon] \subset A^2 \times [0,\varepsilon]\} - \{$its boundary piece, i.e. [the arc $(a,b,c,d,e)$, Figure $26$] $\times [0,\varepsilon]\}$. Of course, for this construction to be possible, it is necessary that on its road from $q \in \partial \delta^2$ to the puncture $p$, the green arc $\lambda (p_2)$ should meet NO OBSTRUCTION; and this issue will be discussed at length. Also we may assume that $\partial Z^2 \supset \partial \delta^2$. With this construction comes a modification of $\partial \delta^2$. The very short arc
$$
q \in [a,e] \subset \partial \delta^2
$$
is to be replaced by the long $[a,b,c,d,e]$. We will denote this change, for $\delta^2 = \delta^2_i$ by
$$
c(b_i) \Longrightarrow c(b_i)(\lambda).]
$$
}

\smallskip

\noindent {\bf Important Remark.} When we consider
$$
\lambda (p_2) \subset GZ^2 \subset \partial N^4 (\Gamma_1 (\infty)),
$$
typically the $N^4 (\Gamma_1 (\infty))$ can see inclusions coming from the $\{h$'s$\} \subset {\rm LAVA}$
$$
(B^3 \times [-\varepsilon , \varepsilon] , S^2 \times [-\varepsilon , \varepsilon]) \subset (N^4 (\Gamma_1 (\infty)) , \partial N^4 (\Gamma_1 (\infty)),
$$
coming with contacts $\lambda(p_2) \cap S^2 \times [-\varepsilon , \varepsilon] \ne \emptyset$ not touching to the $C$'s of the $D^2 (C) \subset {\rm LAVA}$. Let us say that we may a priori see contacts
$$
\lambda (p_2) \cap \{\mbox{extended cocore} \ z\} \ne \emptyset . \leqno (*)
$$
It will turn out that the procedure from Figure 35.4 in the next section V will demolish all the contacts $\{\lambda  (p_2) \cap d_k^2\} \cap \{$cocores $B$ or $R\}$.

\smallskip

On the other hand, the totality of the contacts
$$
\{ \lambda (p_2) \ {\rm or} \ \lambda (q)\} \cap \{ h \in R\}
$$
which will turn out to mean, always $h \in R-B$, will be, all of them, when dangerous for us, part of the parasitical terms in (139), and they will be all taken care of by the CHANGE OF COLOUR, at the end of section VII. So, contacts $\{\lambda (p_2)$ or $\lambda (q)\} \cap R$ should not worry us, but only $\{\lambda (p_2)$ or $\lambda (q)\} \cap B$, and they will be all taken care of too. End of the Important Remark.

\bigskip

2) {\it Because of the lack of extended cocores, the transversal contacts
$$
z_2 \in J d^2 (+ \, \xi_0 , \mbox{space-time, space-time}) \pitchfork (\Delta^2 \times (\xi_0 = -1))
$$
which correspond to the punctures $p_2$ in Figure $22$ (CLASP) obviously need treatment {\rm II}. To the same $p_2$'s are associated RIBBONS
$$
d^2 (+ \, \xi_0 , \mbox{space-time, space-time}) \pitchfork d^2 (+ \, \xi_0 , \mbox{space-time, space-time}),
$$
like the $[p_2 , z_+] , [p_2 , y_+] , \ldots$ Figure $22$. The only punctures at $\xi_0 = -1$ are the $p_2$'s coming with $p_2 \in \mbox{body} \, (d^2 \mid \Delta^2 \times (\xi_0 = -1)) \pitchfork C(\mbox{curve of} \ \Delta^2 \times (\xi_0 = -1)) \subset \partial N^2 (\Gamma_1 (\infty))$.

\smallskip

The Figure $22$ suggests that the clasps $\mbox{body} \, d^2 \pitchfork \mbox{body} \, d^2$ are not an obstruction for $\lambda (p_2)$. For further purpose notice also that
$$
(\lambda (p_2) \mid [P \times [0 \geq \xi_0 \geq -1]]) \cap (B_0 \cup R_0) = \emptyset.
$$
The treatment of the RIBBONS at $\xi_0 = -1$ is suggested in Figure {\rm 28-bis}, while the treatment of the RIBBONS at $\xi_0 = 0$ is suggested in the Figure $29$, WHEN FIGURE $22$ is to be used. Notice how this fits with the LHS of the Figure $23$. It is important here that at $P \times [0 \geq \xi_0 \geq -1]$ we have a CLASP and not a RIBBON.}

\medskip

3) {\it When we join the $p_2$ above with the corresponding $c(b_i) = \partial \delta_i^2$, by an embedded {\ibf long green arc} $\lambda (p_2) \subset Z^2$ $(91)$, and here it should be understood that the various $\lambda (p_2)$'s are two-by-two disjoined, then all the contacts (which are certainly {\ibf in the way} for the process in Figure $26$)
$$
y \in \lambda(p_2) \pitchfork \{\mbox{ACCIDENTS of} \ Z^2 \xrightarrow{ \ G \ } \partial N^4 (\Gamma_1 (\infty))\}
$$
are some intersections with RIBBONS $\mbox{body} \, {\mathcal B}^2 \pitchfork \mbox{body} \, d^2$ (case $\bullet$) below) or RIBBONS $\mbox{body} \, d^2 \pitchfork \mbox{body} \ d^2$ localized at $\xi_0 = 0$ (case $\bullet$$\bullet$$\bullet$$\bullet$) below). The $y = y(p_2)$ occurring here is like in Figure $27$ and RIBBONS $\mbox{body} \, {\mathcal B}^2 \pitchfork \mbox{body} \, d^2$ containing such an obstruction are called {\ibf very special}. Their treatment is explained in Figure $28$. For the obstructing RIBBONS at $\xi_0 = 0$, the treatment is explained in $\{$Figure $29$ with the change ${\mathcal B}^2 \Rightarrow d^2\}$.

\smallskip

Now, in the same connected component of the union of RIBBONS as the very special RIBBONS, there an other ribbons $\mbox{body} \, {\mathcal B}^2 \pitchfork \mbox{body} \, d^2_k$, which we will call {\ibf special} RIBBONS. It may happen that a very special RIBBON contains several $y (p_2)$'s or that very special RIBBONS are adjacent to each other. In order to simplify the exposition, we will assume that none of these things happen, but our procedures extend easily to these more general cases.}

\medskip

4) {\it When a CLASP is lifted to $4^{\rm d}$ then this {\ibf has} to come with one double point $x \in JM^2 (J)$, which can be placed at either end of the clasp. When a RIBBON is lifted to $4^{\rm d}$, then there is a FREE CHOICE which can be made independently for each individual RIBBON: either we create two double points $x_1 , x_2 \in JM^2 (J)$ OR none. }
$$
\includegraphics[width=65mm]{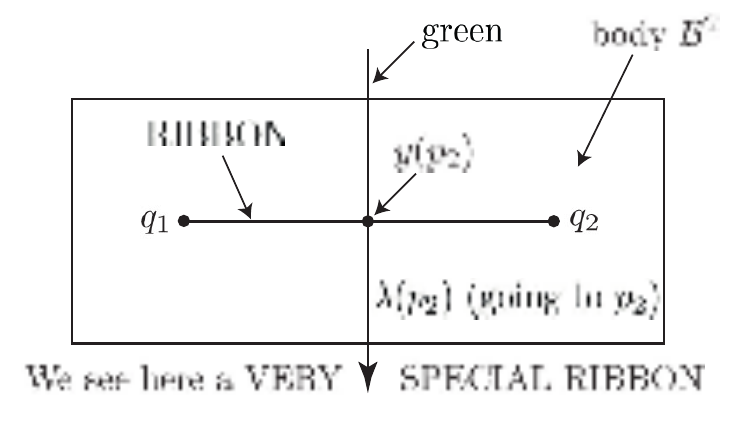}
$$
\label{fig27}
\centerline{\bf Figure 27.}
\begin{quote}
A contact $y(p_2) \in \{$green arc $\lambda (p_2) \cap \{${\ibf very special} RIBBON $[q_1 q_2]\}$, obstructing the $\lambda(p_2)$, like the $\lambda(p_2)$ in Figure 22.
\end{quote}

\bigskip

{\it Here is how this margin of freedom will be made use of. There is, to begin with, the category of RIBBONS called {\ibf generic}, for which the (imposed) choice is {\ibf zero double points}, and these ribbons are the following ones:

\medskip

$\bullet$) The {\ibf very special} RIBBONS and the {\ibf special} RIBBONS coming with them, all of type $\mbox{body} \, {\mathcal B}^2 \pitchfork \mbox{body} \, d^2$. Also, the RIBBONS $\mbox{body} \, {\mathcal B}^2 \pitchfork \mbox{body} \, {\mathcal B}^2$ from Figure {\rm 25-(A)}, when they are in same connected component with them.

\medskip

$\bullet$$\bullet$) The RIBBONS $\mbox{body} \, d^2 \pitchfork \mbox{body} \, d^2$ localized at $\xi_0 = -1$, always having endpoints $p_2$. Since $p_2$ comes with the treatment of the long $\lambda (p_2)({\it green})$, contrary to the case $\bullet$), SHORT GREEN ARCS, like in the Figure $28$, are useless here. Then RIBBONS occur in the RHS of the Figure $23$ and again in Figure {\rm 28-bis}.

\smallskip

[These RIBBONS are like in Figure {\rm  25-(B)}, with the following changes:}
\begin{enumerate}
\item[i)] {\it The $\beta$ becomes $+ \, \xi_0$.}
\item[ii)] {\it The triangles one sees in Figure {\rm 25-(B)} become now the $d^2 (+ \, \xi_0 , \mbox{space-time, space-time})$, $d^2 (+ \, \xi_0 ,$ $\mbox{space-time}$, $\mbox{space-time})$.}
\item[iii)] {\it The $p,q$ become both $p_2$'s.}
\item[iv)] {\it The curves $C$ from Figure {\rm 25-(B)} are now $\Gamma_j$'s, since they live at $\xi_0 = -1$.}
\end{enumerate}

{\it The important fact for the section ${\rm V}$ (CONFINEMENT) and also for the Lemma $12$, is that both at $P_1$ and at $P_2$, to which the two $p_2$'s pertain, the $\Gamma_j$'s which contain them are UP; the Figure {\rm 24-(A)} shows this fact. Notice that, for a puncture $p_2 \in d^2 (+ \, \xi_0 , \mbox{space-time, space-time}) \pitchfork$ ``$C$'' to be there, we have to have a crossing of curves at which our ``$C$'' should be UP.]

\bigskip

We move next to the other category of RIBBONS, called {\ibf exceptional}, for which our choice will be $\# \, JM^2 (J) = 2$. These are:

\medskip

$\bullet$$\bullet$$\bullet$) All the non-special RIBBONS $\mbox{body} \, {\mathcal B}^2 \pitchfork \mbox{body} \, d^2$, as well as the RIBBONS $\mbox{body} \, {\mathcal B}^2 \pitchfork \mbox{body} \, {\mathcal B}^2$ in the same connected component.

\medskip

$\bullet$$\bullet$$\bullet$$\bullet$) The RIBBONS $\mbox{body} \, d^2 \pitchfork \mbox{body} \, d^2$ connected with the puncture $p_1$ at $\xi_0 = 0$, see Figures $22$, $23$. They occur in the LHS of Figure $23$, and the paradigmatic figure for handling them is $29$, for all the exceptional RIBBONS too.

\medskip

In all this story, RIBBONS  $\mbox{body} \, {\mathcal B}^2 \pitchfork \mbox{body} \, d^2$ from Figure {\rm 25-(B)} (or {\rm 25-(C)}) are always treated on the same footing. According to the case, they are treated like in Figure $28$ OR in Figure $29$. Anyway, ${\mathcal B}^2$ is always the long branch and $d^2$ the short one.

\smallskip

More comments will follow now, concerning the $\bullet$) above. So look at the scenario from Figure $28$, explained in $5)$ below.}

\medskip

5) {\it What Figure $28$ presents, is the scenario via which, in the context of our $\bullet$), the long green arc $\lambda (p_2)$ is freed at $y(p_2)$ and hence allowed to continue towards the accident $z_2$ attached to the transversal contact $p_2 \in F \, \mbox{body} \, d_k^2 \pitchfork \Gamma_j \times (\xi_0 = -1)$, in the Figure $23$. The procedure {\rm (II)} (Figure $26$) is being used here at $p_2$. The basic idea is here the following: Once a RIBBON is generic and $({\mathcal B}^2 , d^2)$ are involved, there are {\ibf no} $x_1 , x_2 \in JM^2 (J)$, and we have the inequality
$$
\varepsilon (\mbox{branch} ,\ {\mathcal B}^2) < \varepsilon (\mbox{branch} \, d^2). \leqno (92.2)
$$
This allows us to proceed like in Figure $28$.}

$$
\includegraphics[width=12cm]{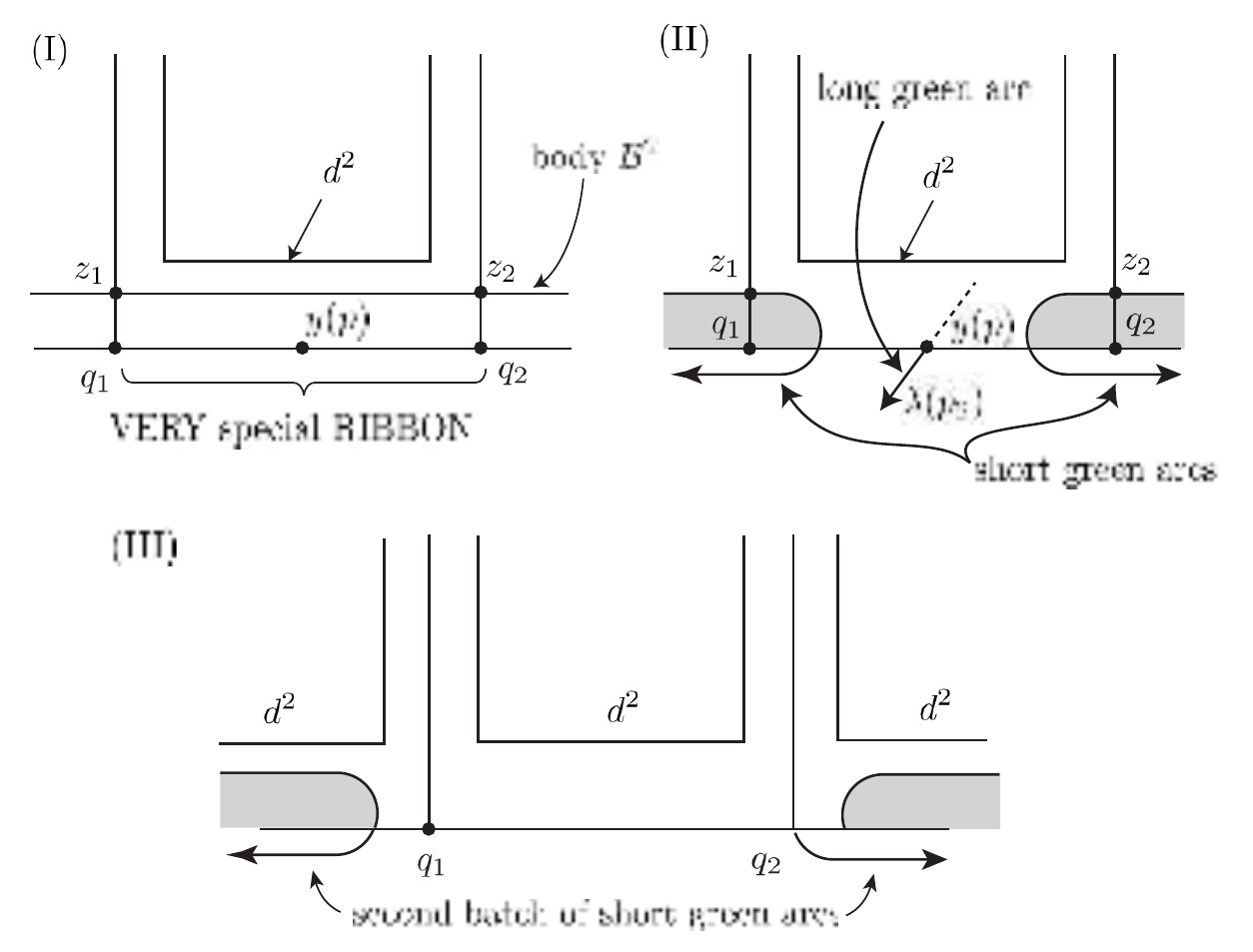}
$$
\label{fig28}
\centerline{\bf Figure 28.}
\begin{quote}
Liberating the main green arcs $\lambda (p_2)$. From $q_1 , q_2$ on, we continue to demolish the accidents in the whole connected component of special RIBBONS, via short green arcs. ONLY procedure II  is present in this figure. The common numerals (I), (II), (III) correspond to the time ordering of the successive steps.
\end{quote}

\bigskip

6) {\it Figure $28$ suggests how the long green arc $\lambda (p_2)$ is able to get to the accident $z_2 \in Jd^2 \pitchfork \Delta^2 \times (\xi_0 = -1)$ occurring at the puncture $p_2$ from the Figure $22$. {\ibf Afterwards}, for the accidents $z_i \in J\delta^2 \pitchfork 2X_0^2$ occurring at the punctures $q_i \in \mbox{body} \, {\mathcal B}^2$ along the whole connected component of very special and special RIBBONS $\mbox{body} \, {\mathcal B}^2 \pitchfork \mbox{body} \, d^2$ from Figure $28$, occurring up-stream from $p_2$, one uses the short green arcs $\lambda(q)$ suggested in Figure $28$. The point here (Figure $28$) is that, with
$$
\varepsilon (\nu^2 (\lambda)) = \varepsilon (\mbox{branch} \, {\mathcal B}) < \varepsilon (\mbox{branch} \, d^2),
$$
the high $Jd^2$ is no longer in the way for $\lambda (p_2)$ which can happily proceed under it, inside the branch ${\mathcal B}^2$. When we get to $\xi_0 = -1$ then we similarly have
$$
\varepsilon (\nu^2 (\lambda)) = \varepsilon (\mbox{branch} \, d^2 (+ \, \xi_0 , \ldots), \mbox{carrying the puncture} \ p_2) < \varepsilon (\mbox{the dual $d^2$-branch}),
$$
and this allows us to proceed like in Figure {\rm 28-bis}, which takes care of $p_2$ (and $z_2$).}

$$
\includegraphics[width=11cm]{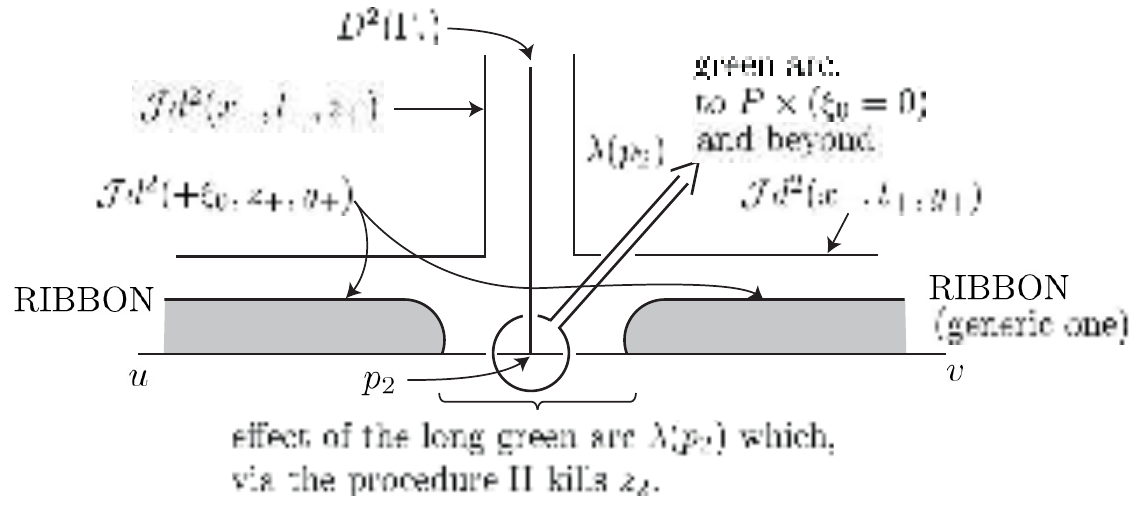}
$$
\label{fig28bis}
\centerline{\bf Figure 28-bis.}
\begin{quote}
We are here at $P \times (\xi_0 = -1)$ and the $\lambda (p_2)$ lives in a plane which is transversal to our figure, cutting through $[u,v]$. See here also the Figures 22 and 23.
\end{quote}

\bigskip

{\it All this takes care of the situation $\bullet$) and, as already said earlier, the $\bullet$$\bullet$) is to be taken care of by Figure~{\rm 28-bis}, which does {\ibf not} need short green arcs.

\smallskip

At this point, the only ACCIDENTS still alive are the transversal contacts $z \in J\delta^2 \pitchfork 2X_0^2$ coupled with double points $x \in JM^2 (J)$ (with $J$ like in $(80)$) coming with the exceptional ribbons from $\bullet$$\bullet$$\bullet)$ $+$ $\bullet$$\bullet$$\bullet$$\bullet)$ at the end of $4)$ above.

\smallskip

[The ribbons $\mbox{body} \, {\mathcal B}^2 \pitchfork \mbox{body} \, {\mathcal B}^2$, like in the Figure {\rm 25-(A)} are either generic, i.e. connected with the special and very special RIBBONS OR exceptional. The generic ones will get the treatment of short given arcs from Figure $28$, while the exceptional ones will be treated in Figure $29$.]

\smallskip

When it comes to the exceptional RIBBONS we make the choice opposite to the one made in connection with Figure $28$, namely now
$$
\varepsilon (\mbox{short branch}) < \varepsilon (\mbox{long branch})
$$
coming with two double points $x \in JM^2 (J)$ for each exceptional RIBBON.

\smallskip

We have now Figure $29$, and we apply the treatment {\rm (I)} which simultaneously kills the $z$ and the $x$'s. Notice that for the global killing of our ACCIDENTS, both procedures {\rm (I)} and {\rm (II)} are needed.}

\smallskip

{\it In the context of Figure $29$ for $\bullet$$\bullet$$\bullet$$\bullet$) we may have possible obstructions
$$
\lambda(p_2) \pitchfork \{\mbox{RIBBONS} \ldots )\}
$$
which are then treated like in Figure {\rm 29-(B)}, without any short arcs.}

\bigskip

With this, the proof of our Lemma 11.3 is finished, but NOT YET the proof of Theorem 11. We have gotten rid of all the ACCIDENTS but the diagonalization condition (82.1), which we will call THE LITTLE BLUE DIAGONALIZATION is not yet realized. That will be done in the next section V.

\bigskip

We end this section with some last COMMENTS. 

$$
\includegraphics[width=115mm]{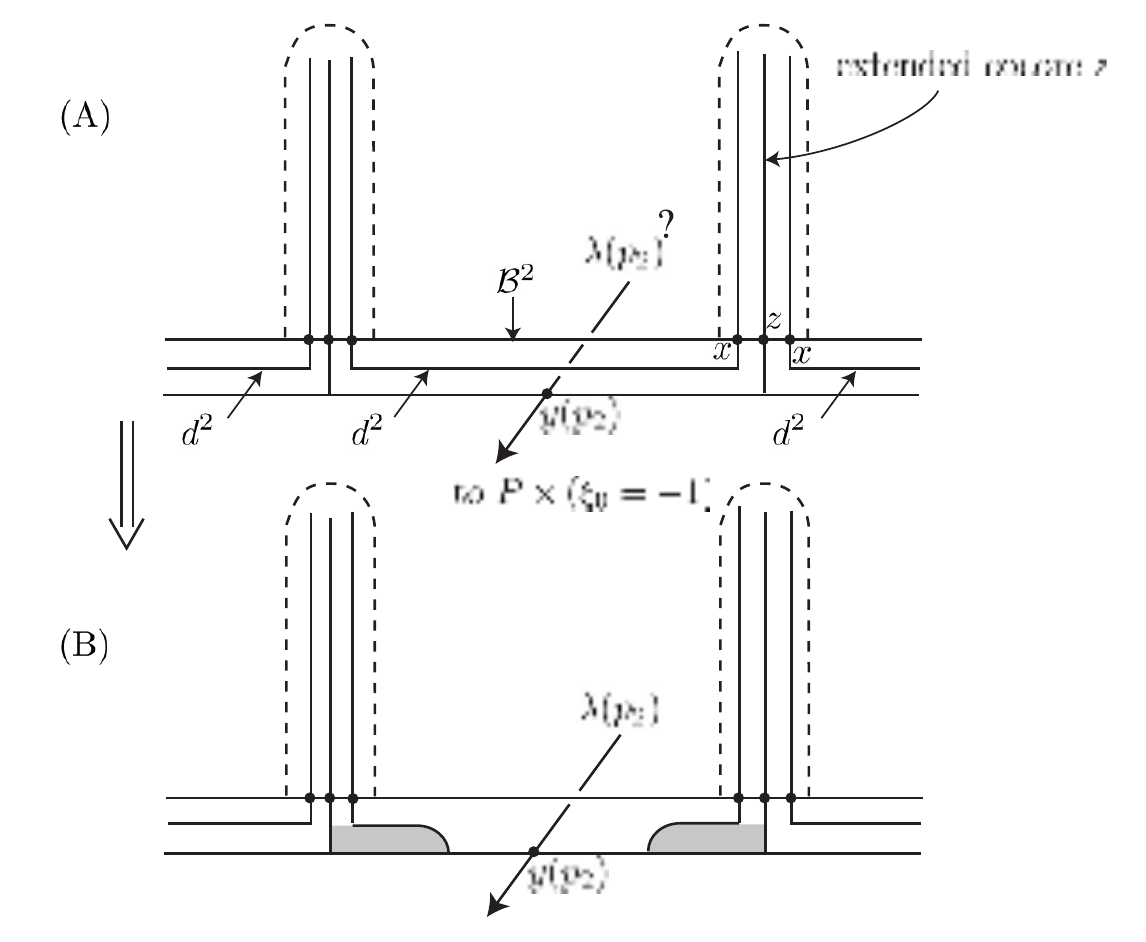}
$$
\label{fig29}
\centerline{\bf Figure 29.}
\begin{quote}
Illustration for the procedure I for killing ACCIDENTS (transversal contacts $z$ and double points $x$, simultaneously). The $\lambda (p_2)$ with the potential obstruction $y(p_2)$ may be present for $\bullet$$\bullet$$\bullet$$\bullet$). [This same figure applies at $\xi_0 = 0$, with ${\mathcal B}^2 \Rightarrow d^2$, IN THE CASE WHEN FIGURE 22 IS USED (CLASP along $[0 \geq \xi_0 \geq -1]$.) It is {\ibf only} then that $y (p_2)$ may be present.]
\end{quote}

\bigskip

We go back now to our family
$$
\sum_1^M b_i = B_1 \cap (\Gamma (1) \times (\xi_0 = -1)) = B_0 \cap (\Gamma (1) \times (\xi_0 = -1)). \leqno (93)
$$
Notice here that each edge $e \subset \Gamma (1) \times (\xi_0 = -1)$ contains its $b_i \in \underset{1}{\overset{M}{\sum}} \, b_i$. On the other hand, since $\Gamma (1) \times (\xi_0 = -1) - \underset{1}{\overset{n}{\sum}} \, R_i \times (\xi_0 = -1)$ is a tree, we have $M \gg n$, i.e. our $\Delta^2 \times (\xi_0 = -1)$ is violently disbalanced, as far as the RED/BLUE balance is concerned. $\Box$

\bigskip

Notice that the destruction of accidents, via I) changes the $J\delta_i^2$ modulo the boundary, which stays unchanged while II) changes the boundary too. So, we will distinguish between the {\ibf original} $C(b_i)$'s and the $C (b_i)(\lambda)$'s $\equiv \{c(b_i)$ MODIFIED by the effect of the green arcs $\lambda (p_2)$ (and the short green arcs)$\}$. It is the
$$
\left( \sum_1^M \delta_i^2 , \sum_1^M \partial \delta_i^2 = C_i (b)(\lambda)\right) \xrightarrow{ \ J \ } (\partial N^4 (2X_0^2)^{\wedge} \times [0,1] , \partial N^4 (2X_0^2)^{\wedge} \times \{0\}) , \leqno (93.1)
$$
which is ACCIDENT-free.

\smallskip

Now, to each $b \in B_0 = \{$the $b \in X^2[{\rm new}] \subset 2X^2$ which live on the $X^2 \times r$ side$\}$, corresponds on the $b$-side, and see here the Figures 7-II, III and 7-bis, a $b \times b \in B_1$, on the $X^2 \times b \approx X_b^2$ side.

\bigskip

\noindent (94) \quad In the geometric intersection matrix, for $1 \leq i \leq M$ we find, BEFORE Lemma 11.3 has been applied, for the original $C(b_i)$,
$$
C(b_i) \cdot b_i = 1 , C(b_i) \cdot (b_i \times b_i) = 1 \ \mbox{and NOTHING else}.
$$

When we move from $C(b_i)$ to $C(b_i)(\lambda)$, then we find

\bigskip

\noindent (94.1) \quad $C(b_i)(\lambda) \cdot b_i = 1$ AND ALSO OFF-DIAGONAL TERMS, both $C(b_i)(\lambda) \cdot \{$the $b_i \times b_i\} = 1$ and $C(b_i)(\lambda) \cdot \{$various $b \in B_0$ touched by the green arcs, long or short$\} \ne 0$.

\bigskip

The off-diagonal terms $C(b_i)(\lambda) \cdot B$ above, are called PARASITICAL. Since $\lambda (p_2)$ never goes along an edge $e \subset \Gamma(1) \times (\xi_0 = -1)$ and since there are no short green arcs at $\xi_0 = -1$, the $\underset{1}{\overset{M}{\sum}} \, b_i$ are never parasitical.

\bigskip

\noindent (94.2) \quad Here is a definition which will be useful. Let $b \in B_0$ and let $b = b_0 \to b_1 \to \ldots \to b_k$ be its normal blue trajectory, defined like in (77.1). This comes with successive $B^2 (b_i)$'s. We will say that $b$ is {\ibf trivial}, if for {\ibf none} of them $B^2 (b_i)$ we have inclusions $D^2 (\gamma^0) \subset B^2 (b_i)$. If, moreover $k=0$, i.e. if $b$ has no outcoming arrows, then we say that it is {\ibf very trivial}.


\section{Confinement and the little blue diagonalization}\label{sec5}

We start by making the following CLAIM, concerning the edges $e(r) \subset \{$curve $c(r)\} \subset \{$link$\}$, displayed in the Figures 7-(I and IV).

\bigskip

\noindent {\bf Claim (95).} \quad Using our margin of freedom for crossings which do not concern $\Delta^2 \times (\xi_0 = -1)$, we can ask that for all the Figures 9 (and/or 24) concerned, and at all its corners, every $e(r) \subset c(r)$ (Figure 7) should be UP.

\bigskip

\noindent {\bf Proof.} We will consider afterwards the case when $e(r) = P \times [0 \geq \xi_0 \geq -1]$. So we look now at $e(r)$ going from $P_1$ to $P_2$ along an axis $u \in (x,y,z,t)$.

\smallskip

Let us say that $e(r)$ goes from $b^3 (u_{\pm})$ at $P_1$ to $b^3 (u_{\mp})$ at $P_2$. The corners of $c(r)$ at $e(r)$ are then $(u_{\pm} , \beta)$ at $P_1$ and $(u_{\mp} , \beta)$ at $P_2$. Figure 24 tells us that all the arcs $[b^3 (u_{\pm}) , b^3 (\beta)]$ are UP at all their crossings, which is what our claim says, concerning them.

\smallskip

Similarly, for the case $e(r) = P \times [0 \geq \xi_0 \geq -1]$ the arcs $[b^3 (\pm \, \xi_0) , b^3 (\beta)]$ are UP too. End of Proof.

\bigskip

Once we have the CLAIM above, without any obstruction the $c(r)$'s can be pushed into $\partial N_+^4 (2\Gamma (\infty))$.

\smallskip

Next, we want to describe a class of {\ibf admissible subdivisions} of $2X_0^2$ which should guarantee, at the BLUE level, the conservation of the following basic feature, part of at the level $X^2 [{\rm new}]$ (and of course at $2X^2$ too)
$$
\boxed{\{\mbox{curves $\Gamma$ at} \ \xi_0 = -1 \} \subset \{\mbox{curves $\gamma^1$}\}} \ .
$$
Remember that, for general subdivisions, the RULES OF THE GAME are
$$
D^2 (\gamma^0) \Longrightarrow \{\mbox{{\ibf one} small $D^2 (\gamma^0)$}\} + \{\mbox{many small $D^2 (C)$'s}\} \ ({\rm RED})
$$
$$
D^2 (\gamma^1) \Longrightarrow \{\mbox{{\ibf one} small $D^2 (\gamma^1)$}\} + \{\mbox{many small $D^2 (\eta)$'s}\} \ ({\rm BLUE}).
$$
That is why we need now special, ``admissible subdivisions''. Our admissible subdivision is presented in Figure 30 and, at $\xi_0 = -1$, we see there a lot of $B({\rm new})$'s only one of which is $B({\rm new}) \cap R$. This way, each edge of $\Gamma (1) \times (\xi_0 = -1)$ continues to have one $B$ and $\Gamma (1) \times (\xi_0 = -1) -R$ continues to be a tree.

$$
\includegraphics[width=14cm]{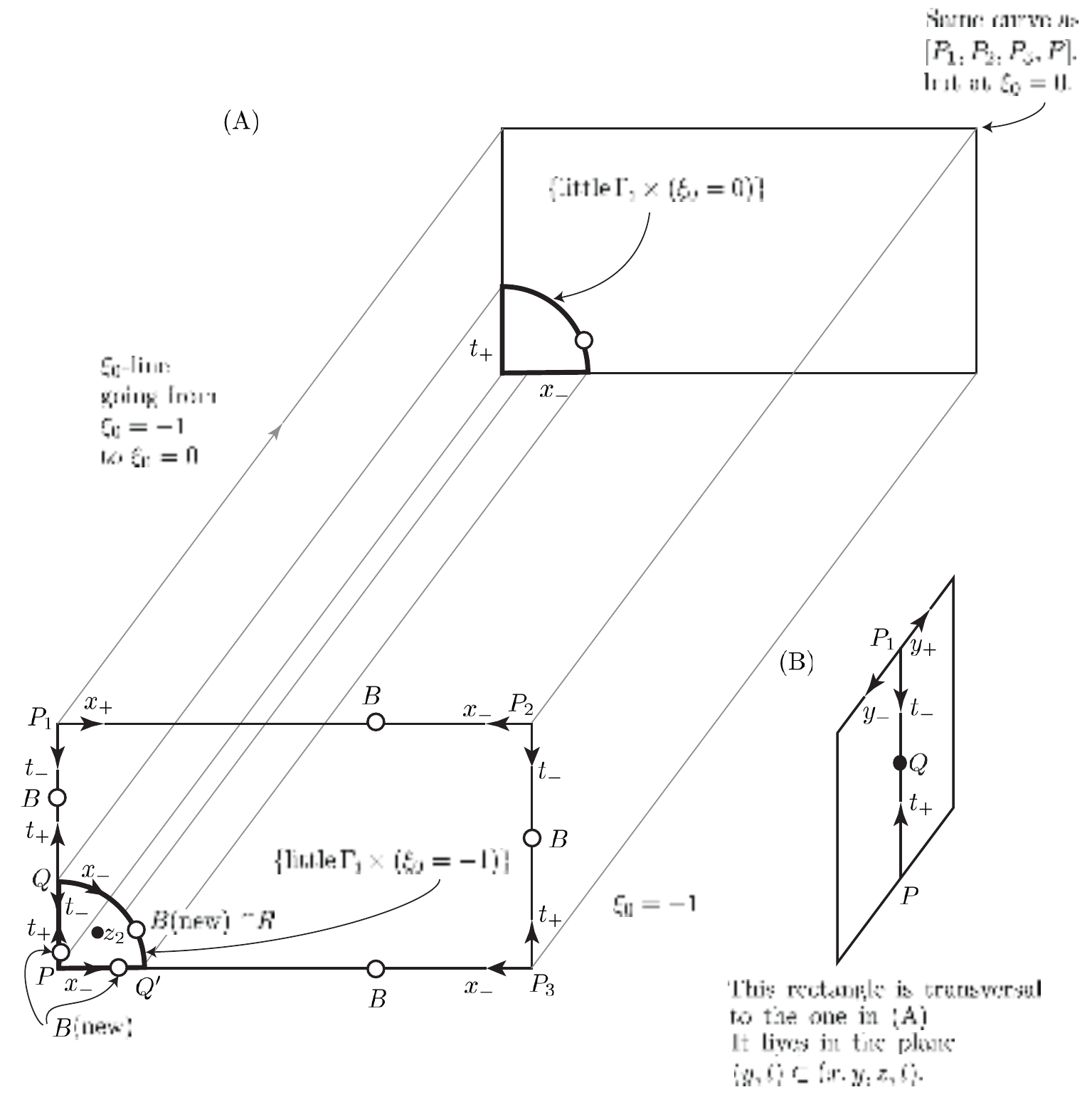}
$$
\label{fig30}
\centerline{\bf Figure 30.}
\begin{quote}
Here $[P,P_1,P_2,P_3]$ is a generic curve $\Gamma_j \subset \Delta^2 \times (\xi_0 = -1)$ with, let us say $P$ the vertex of Figure~22, part of a Figure 9, as displayed in Figure 24. Just like, at $\xi_0 = -1$, the $\Gamma_j \times (\xi_0 = -1)$ breaks (BLUE-wise) into many little $\gamma^1$'s, at $\xi_0 = 0$, the $\Gamma_j \times (\xi_0 = 0)$ also breaks (RED-wise), in many little $\gamma^0$'s (and, of course, BLUE-wise in $\gamma^1$'s too).
\end{quote}

\smallskip

\noindent {\bf Explanations concerning Figure 30.} The $z_2$, placed close to $P$ is an accident $z_2 \in J \delta^2 \pitchfork D^2 (\Gamma_j)$ ($=$ square $[P,P_1,P_2,P_3]$). The $\Gamma_j$ is UP at $P$ (see Figures 22 and 24-A), but certainly not UP everywhere, at all vertices, at least there is no reason for that. We subdivide $[P_1,P_2,P_3,P] \subset (\xi_0 = -1)$ like suggested, BUT then the whole $[P_1,P_2,P_3,P] \times [1 \leq \xi_0 \leq 0]$ too. It is easy to see that this subdivision is now admissible, in the sense defined above. We use this subdivision at every corner $P$ carrying an accident, i.e. a $(p_1 , z_2)$.

\smallskip

This also breaks $\Gamma_j$ into
$$
\Gamma_j = \{\mbox{the {\ibf little} $\Gamma_j$, i.e. the round triangle $[P,Q,Q']$, keeping the accident $z_2$}\} \ +
$$
$$
+ \ \{\mbox{the big $\Gamma_j$ remaining, which is free of accidents}\}.
$$

\bigskip

Figure (B) suggests the environment of $[P_1,P]$ in $\Delta^2 \times (\xi_0 = -1) \subset 2X_0^2$. Also now $z_2 \in J\delta^2 \pitchfork D^2 ({\rm little} \, \Gamma_j)$ and little $\Gamma_j$ is UP at {\ibf all} its corners. This is the aim of the whole thing Figures 24 $+$ 30 show that little $\Gamma_j$ is UP at $P$ and Figure 31 shows that it is UP at $Q$ too. For $Q'$ it is the same.

$$
\includegraphics[width=115mm]{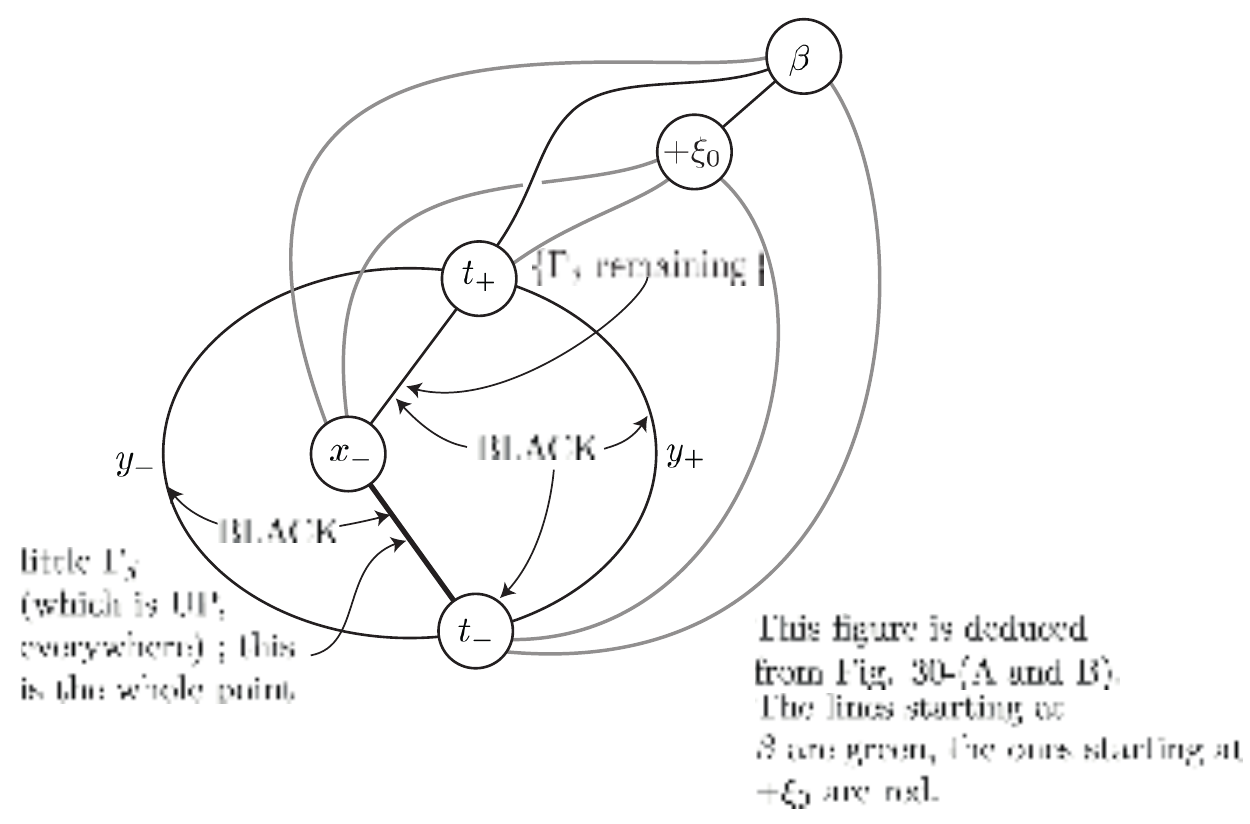}
$$
\label{fig31}
\centerline{\bf Figure 31.}
\begin{quote}
This is the Figure 9 (complete with $+ \, \xi_0$ and $\beta$) for the $Q$, the new vertex introduced in Figure~30. Here only the BLACK arcs are part of curves $\Gamma_j \subset \Delta^2 \times (\xi_0 = -1)$ and at the crossings in which they take part, in this figure, they are always DOWN. The RED and GREEN curves are not part of $\Delta^2 \times (\xi_0 = -1)$ and, with these things, this figure cannot create ACCIDENTS involving $\Delta^2 \times (\xi_0 = -1)$. There is a similar figure at $Q'$ (see Figure 30-(A)). For the colours see the legend of Figure 24.
\end{quote}

\bigskip

Our subdivision introduces an additional $B({\rm new})$'s, one of which is $B \cap R$ which can all be treated like the $b_i$, in the context of Figure 19.1. Since $\Gamma (1)$ increases, a new $R = R_{n+1}$ has to be added to the already existing $\underset{1}{\overset{M}{\sum}} \, R_i$. This is, of course, our $B({\rm new}) \in R \cap B$, Figure 30-(A). Remember at this point that the only function of the set $\Sigma \, R_i$ is to render the $\Gamma (1) - \Sigma \, R_i$ be a tree. It will essentially disappear from our picture.

\smallskip

The construction above extends easily to the accidents carried by $D^2 (C)$'s, outside $\xi_0 = -1$ and which also need green arcs. We need now just a normal subdivision, without worrying about the admissibility condition, which is certainly not required outside of $\Delta^2 \times [0 \geq \xi_0 \geq -1]$. With this we will isolate any transversal contact $z \in J\delta^2 \pitchfork D^2 (C)$, proceeding like in the lower part of Figure 30, and perform, for each $C$, the change (in particular done for $C = \Gamma_j \times (\xi_0 = 0)$)
$$
C \Longrightarrow \{\mbox{the {\ibf little} $C$}\} + \{\mbox{the $C$ remaining}\} .
$$

\bigskip

\noindent {\bf Important Remark.} The accident coming with the $D^2 (\Gamma_j) \times (\xi_0 = 0)$ in Figure 30 certainly does not need a green arc, since we apply for it the push over the $\{$extended cocore$\}^{\wedge}$. Nevertheless, for reasons to become clear later, we still apply to it the procedure from Figure 30. Here $\{{\rm little} \ C\}$ is UP at all its corners and retains the accident $z$.

\smallskip

Concerning now the SPLITTING (27.1), with things as they stood at the end of the preceeding section~IV, we had the normal CONFINEMENT conditions, for the $\{$link$\}$ ($=$ the internal curves, attaching zones of 2-handles if $2X^2$)

\bigskip

\noindent (96) \quad $\underset{1}{\overset{\infty}{\sum}} \, C_i ($of $X_0^2 [{\rm new}]$) $+$ $\underset{1}{\overset{\bar n}{\sum}} \, \Gamma_j (\subset (\xi_0 = -1)) + \{$extended $\gamma^0$'s (18.1)$\} \subset \partial N_-^4 (\Gamma_1 (\infty))$ AND  $\underset{1}{\overset{\infty}{\sum}} \, \eta_i \times b + \sum c(r)$ (see claim (95)) $\subset \partial N_+^4 (2\Gamma (\infty)) \supset \underset{1}{\overset{M}{\sum}} \, c(b_i) (= \partial \delta_i^2)$.

\bigskip

Now, just like in the context of the CLAIM (95), for $b_{i \leq M}$, the $c(b_i)$ is UP at all its corners (since in Figure 24 all the arcs $[\beta , \mbox{space-time}]$ are UP) and so we also have $\underset{1}{\overset{M}{\sum}} \, c(b_i) \subset \partial N_+ (2\Gamma (\infty))$. By contrast, for $\underset{1}{\overset{M}{\sum}} \, c(b_i) (\lambda)$ this is not clear at all, and there are two issues to be settled for the green arcs, in their totality, the long green arc $\lambda (p_2)$ and the short green arcs $\lambda (q)$:

\medskip

{\bf Issue A)} We want the confinement condition
$$
\{\mbox{green arcs $\lambda$}\} \subset \partial N_+^4 (2\Gamma (\infty)) , \ \mbox{to be satisfied}.
$$
This would allow us to add $\underset{1}{\overset{M}{\sum}} \, c(b_i) (\lambda) \subset \partial N_+^4 (2\Gamma (\infty))$ to (96).

\medskip

{\bf Issue B)} We also want to have the following condition satisfied
$$
\{\mbox{green arcs $\lambda$}\} \cdot B \subset \{\mbox{trivial $B$'s (see (94.2))}\}.
$$

\smallskip

These will be essential for getting the (82.1) and hence clinch the proof of the Theorem 11. But notice, that with what we have already done in this section, we can replace (96) with the real life FORCED CONFINEMENT CONDITIONS to be used from now on, until we can get to something even better:

\bigskip

\noindent (97) \quad $\sum \{C_i \ {\rm remaining}\} + \sum \{\Gamma_j \ {\rm remaining}\} + \sum \gamma_k^0 \subset \partial N_-^4 (\Gamma_1 (\infty))$, and then $\sum\{{\rm little} \ C_i\} + \sum\{{\rm little} \ \Gamma_j\} + \sum \eta_i \times b + \sum c(r) \subset \partial N_+^4 (2\Gamma (\infty)) \supset \underset{1}{\overset{M}{\sum}} \, c(b_i)$.

\bigskip

[The $\{{\rm little} \, \ldots\}$ which we want to push into $\partial N_+^4 (2\Gamma (\infty))$, have to be UP at all their corners. And this condition {\ibf is} fulfilled since these $\{{\rm little} \, \ldots\}$ either come from the $D^2 (\ldots)$'s carrying accidents needing green arcs or the $D^2 (\Gamma \times (\xi_0 = 0))$ from Figure 30. For this last one there is no accident but it is UP at the relevant corner, since the $D^2 (\Gamma \times (\xi_0 = -1))$ is, Figure 24-(A).

\smallskip

All the other accidents needing green arcs come from body ${\mathcal B}^2$ and Figure 24-(A) tells us that the $B^2 (\beta , \ldots)$'s are UP at all their corners.]

\bigskip

And once Lemma 12 and also Theorem 13 will be with us, we will be able to improve the RHS of (97), by setting
$$
\partial N_+^4 (2\Gamma (\infty)) \supset \sum_1^M \eta_i ({\rm green}).
$$

\bigskip

[The more precise way in which the confinement condition $\sum c(r) \subset \partial N_+^4 (2\Gamma (\infty))$ should occur, is explained in the Figure 45, which matches well with Figure 32, which displays the confinement
$$
\sum \{ {\rm little} \ C_i \} \subset \partial N_+^4 (2\Gamma (\infty)).]
$$

\bigskip

In a figure  like 45 or 32, we are in $\partial N^4 (2\Gamma (\infty))$, and even if only $2\Gamma (\infty)$ is mentioned explicitly in our drawings (in Figure 45), this is so only for typographical reasons, since the $2\Gamma (\infty)$ is buryed deep inside $N^4 (2\Gamma (\infty))$, far from its boundary. But the curves are there and, in Figure 45-(III), by $c(b)$ we mean $c(b_{i \leq M})$. Forgetting, for the time being, about the issue B), here is how issue A) is taken care of.

\bigskip

\noindent {\bf Lemma 12.} {\it Via a simple isotopic move of the map (see {\rm (92.1)})
$$
Z^2 = \sum_1^M \mbox{body} \ {\mathcal B}_i^2 \cup \sum_1^N \mbox{body} \ d_k \xrightarrow{ \ G \ } \partial N^4 (2\Gamma (\infty)),
$$
we can make that
$$
\sum \{\mbox{green arcs} \ \lambda \} \subset \partial N_+^4 (2\Gamma (\infty)).
\leqno (98)
$$
}

$$
\includegraphics[width=13cm]{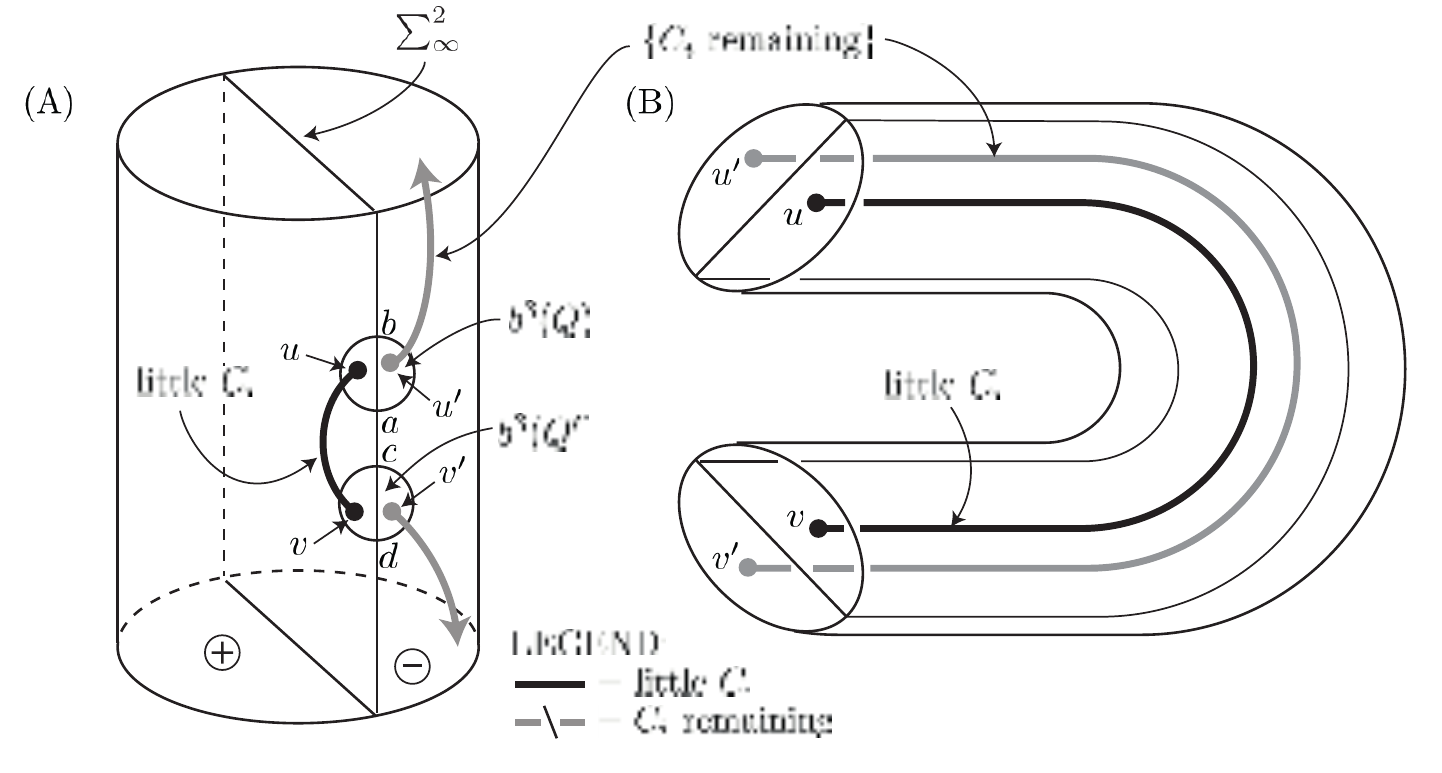}
$$
\label{fig32}
\centerline{\bf Figure 32.}
\begin{quote}
The $4^{\rm d}$ geometry of $\{ {\rm little} \ C_i \} \subset \partial N_+^4 (2\Gamma (\infty))$; (and of course, this applies to $\{ {\rm little} \ \Gamma_i \}$ too). We have a $B^4 (P)$ such that $C_i$ occurs as an arc in $\partial B^4 (P)$. In (A) we have suggested a $B^3 \times [-\varepsilon,\varepsilon]$ cutting through $(B^4 (P),C_i)$ such that $\partial B^3 \times [-\varepsilon,\varepsilon] \subset \partial B^4 (P)$, this is the visible part in Figure 32-(A) and it is far from any other $C$ which may be in the way, between our $C_i$ and $\partial_+ N^4 (2\Gamma (\infty))$. The (B) shows the 1-handle necessary for the transformation
$$
C_i \Longrightarrow \{ {\rm little} \ C_i \} + \{C_i \ {\rm remaining}\}.
$$
In order to see the $\{ {\rm little} \ C_i \}$ close in, we have to put the (A) and (B) together.
\end{quote}

\bigskip

\noindent {\bf Proof of Lemma 12.} We have to look much closer at ${\rm body} \ d^2 \cup {\rm body} \ B^2 \subset \partial N^4 (\Gamma_1 (\infty))$, inside which the green arcs are, anyway, contained. Each $d_k^2 - {\rm body} \ d_k^2$ is a disjoined collection of discs, while $\partial d_k^2 \subset \partial \ {\rm body} \ d_k^2$.

$$
\includegraphics[width=13cm]{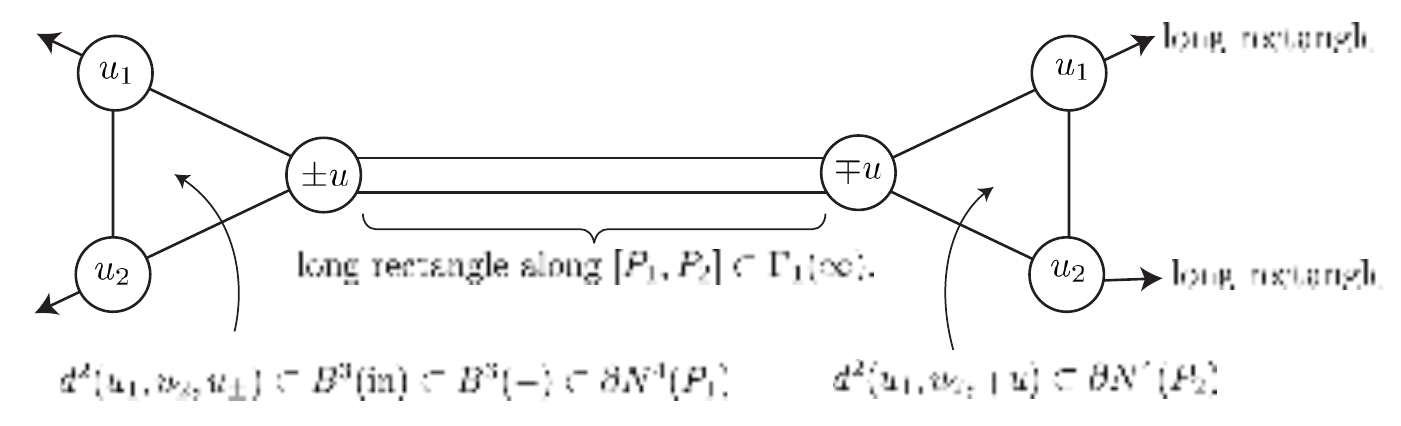}
$$
\label{fig32.1}
\centerline{\bf Figure 32.1.}
\begin{quote}
A piece of ${\rm body} \ d^2_k = U ({\rm triangles}) \cup U ({\rm rectangles})$, continuing along the four arrows. The four sides of the figure are contained in $U C_i \subset \{{\rm link}\}$.
\end{quote}

\bigskip

The disc with many holes, the ${\rm body} \ d^2_k$, is ${\rm body} \ d^2_k = \Bigl(\bigcup \, \{\mbox{individual triangles} \ d^2 (u_1 , u_2 , u_3)$, where $u_i \in (\pm \, x , \pm \, y , \pm z , \pm \, t , \pm \, \xi_0)$, contained inside the $\partial N^4 (P)$'s.They have disjoined interiors and they all occur inside the various embellished Figures 9$\}\Bigl) \cup \sum \{$long rectangles, contained inside the various $\partial N^4 (\Gamma_1 (\infty)) \mid \{$edge $[P_1 , P_2]\}$, and they join two triangles $d^2 (u_1 , u_2 ,u_3)$, like in the Figure 32.1 below$\}$.

\smallskip

In terms of the MODEL FOR $N^4 (\Gamma_1 (\infty))$ given around formulae (27) to (27.5), for $d^2 (u_1 , u_2 ,u_3) \subset B^3 (-)$, there is a vector field along $d^2 (u_1 , u_2 ,u_3)$ normal to it and looking towards $\sum_{\infty}^2$. In terms of the Figure 24, which is the paradigmatical figure of type 9, drawn as a sort of link diagram on its $S^2 (P)$ (and/or $S_{\infty}^2$), this is the direction looking towards the observer, hence, when one is at a crossing, then
$$
\vec v = \{\mbox{direction DOWN $\to$ UP}\}.
$$

When we consider a figure like 32.1, the $\vec v$ extends {\ibf continuously} throughout such a figure, hence continuously along the $[P_1 , P_2]$. Here is how one has to understand what is going on in a complete Figure 32.1. We look again at our paradigmatical Figure 24, and let us say that we want to see what goes on along $P \times [-1 \leq \xi_0 \leq 0]$. Figure 24 displays the $\partial N^4 (P \times (\xi_0 = -1))$ and inside it we see $b^3 (+ \, \xi_0)$. This is met by seven $d^2$-triangles $d^2 (+ \, \xi_0$, space-time, space-time). For expository purposes, we forget here about the fact that at $\xi_0 = -1$ the triangles $B^2 (\beta , \ldots)$ have void interiors. On the 2-sphere $\partial b^3 (+ \, \xi_0)$ our triangles cut a connected graph which is displayed in the Figure 33.1. There, one should be able to read the continuous propagation of $\vec v$ along $P \times [-1 \leq \xi_0 \leq 0]$.

\smallskip

We start by looking now at the main, long green arcs $\lambda (p_2)$ from 3), in Lemma 11.3. These arcs go from some $u \in \{{\rm original} \ c(b_i)\} \subset \partial N_+^4 (2\Gamma (\infty))$ to some $\{{\rm little} \ \Gamma_j \} \subset \partial N_+^4 (2\Gamma (\infty))$, actually to $p_2 \in \{{\rm little} \ \Gamma_j \}$. Figure 22 suggests the intersection $\lambda (p_2) \cap \Delta^2 \cap [0 \geq \xi_0 \geq -1]$. Figure 22 concerns a piece $\lambda (p_2) \cap {\rm body} \ d_k^2$ and, for the time being, we will stay in our proof with these pieces $\lambda (p_2) \cap {\rm body} \ d^2$.

\smallskip

The idea now is to bring $\lambda (p_2) \cap {\rm body} \ d^2$, modulo its end $p_2 \in \partial N_+^4 (\Gamma_1 (\infty))$, fully into $\partial N_+^4 (\Gamma_1 (\infty))$, by pushing it along the vector field $\vec v$. A priori two kinds of obstructions may occur at this point.

\medskip

A) When inside a triangle $d^2 (u_1 , u_2 ,u_3) \subset \partial N^4 (P)$, we may meet a curve $[v_1 , v_2] \subset\{$link$\}$ in the way. This obstruction which is a real one, is illustrated in the Figure 33. Look here at Figure 24-(A) for an explanation.

\smallskip

In Figure 33, the arc $[x_- , t_+]$ (UP) obstructs the local push of $\lambda (p_2)$ into $\partial N_+^4 (2\Gamma (\infty))$. When we are completely outside of $\Delta^2 \times [0 \geq \xi_0 \geq -1]$, which certainly implies far from $\Delta^2 \times (\xi_0 = -1)$, then we can certainly apply the crossing freedom, change UP $\Leftrightarrow$ DOWN around, and in this case, get rid of our obstruction.

\smallskip

We do not push the obstructing curve into $\partial N_+^4 (\Gamma_1 (\infty))$, as a priori we could, since we are not allowed to muck around with the confinement conditions (97). That would perturb the mechanism of section VII below. So we do {\ibf need} here an admissible crossing, coming with its change of local topology.

\medskip

$$
\includegraphics[width=9cm]{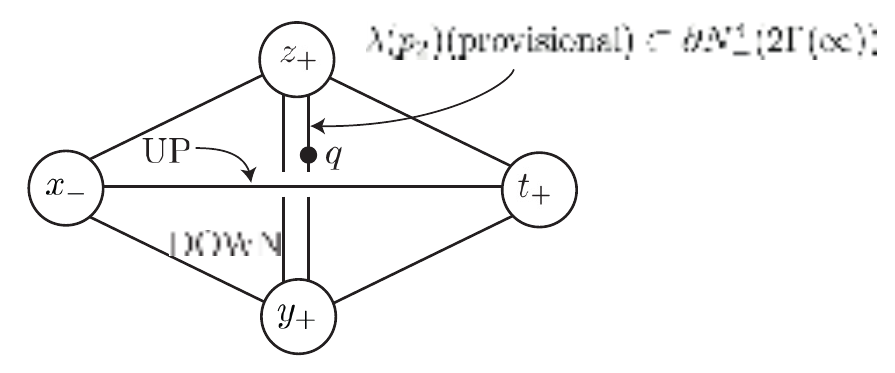}
$$
\label{fig33}
\centerline{\bf Figure 33.}
\begin{quote}
A possible obstruction for pushing the long green arc $\lambda (p_2)$ into $\partial N_+^4 (2\Gamma (\infty))$. The $q$ is here under $p_2$ (Figure 22).
\end{quote}

\bigskip

{\ibf But}, when we are at $\xi_0 = 0$, this admissible crossing would conflict with our policy concerning the LHS of the CLASP from Figure 22 (see here the 4.1) in Lemma 11.1. But then, in terms of Figure 30, when at $P \times (\xi_0 = 0)$, we have $\{ c(x_- , t_+)$ (Figure 33)$\} \subset \{$little $\Gamma_i \times (\xi_0 = 0)\} \subset \partial N_+^4 (2\Gamma (\infty))$, and so there is actually no obstruction to worry about. The {\ibf admissible} subdivision from Figure 30 is essential here as well as the fact that in Figure 30 the $D^2 (\Gamma_j \times (\xi_0 = -1))$ and $D^2 (\Gamma_j \times (\xi_0 = 0))$ get a similar treatment, although the $D^2 (\Gamma_j \times (\xi_0 = 0))$ carries NO ACCIDENT. Next, we look into the possible obstruction when we are at $P \times (\xi_0 = -1)$ and $+ \, \xi_0$ lines could be in the way, in the same manner as the $[x_- , t_+]$ is, in Figure 33. There is a unique CLASP coming with Figure 24; it goes along the $\xi_0$-line, being produced by the crossing (see Figure 22) 

$$
[x_- , t_+] ({\rm UP}) / [z_+ , y_+] ({\rm DOWN}).
$$

\bigskip

With it comes an occurrence of $\lambda (p_2)$, at the level of Figure 24, along $d^2 (+ \, \xi_0 , y_+ , z_+)$. So, let us say that Figure 33 concerns now vertex $P \times (\xi_0 = -1)$, when $p_2$ occurs in Figure 22. We let $\lambda (p_2)$ go from $b^3 (+ \, \xi_0)$ to $b^3 (y_+)$, parallel and close to $[+\, \xi_0 , y_+]$, to begin with, and notice that the line $[\beta , t_+]$ (Figure 24) $\subset \partial N_+^4 (2\Gamma (\infty))$ is NOT an OBSTRUCTION; see here the CLAIM (95).

\smallskip

Then, in the triangle $d^2 (+ \, \xi_0 , y_+ , z_+)$, we let $\lambda (p_2)$ go along the $[y_+ , z_+]$ (Figure 33 at $P \times (\xi_0 = -1))$, now under

$$
[x_- , t_+] ({\rm UP}) \subset \{{\rm little} \ \Gamma_i \times (\xi_0 = -1) \subset \partial N_+^4 (2\Gamma (\infty)),
$$

\medskip

\noindent NO OBSTRUCTION !, all the way up to $p_2$, at the level of the crossing, Figure 33.

\smallskip

All this takes care completely of the potential obstruction A).

\medskip

B) But then, there is also a second OBSTRUCTION (still potential, so far), for pushing $\lambda (p_2) \cap d_k^2$ into $\partial N_+^4 (2\Gamma (\infty))$, namely one may find smooth sheets in body $d^2$ or in body ${\mathcal B}^2$, in the way. But then, these can be pushed into $\partial N_+^4 (2\Gamma (\infty))$, in front of the push of the $\lambda (p_2)$, without any harm. What we see going on in the Figure 33.1 is an instance of the general fact we have just stated.

$$
\includegraphics[width=9cm]{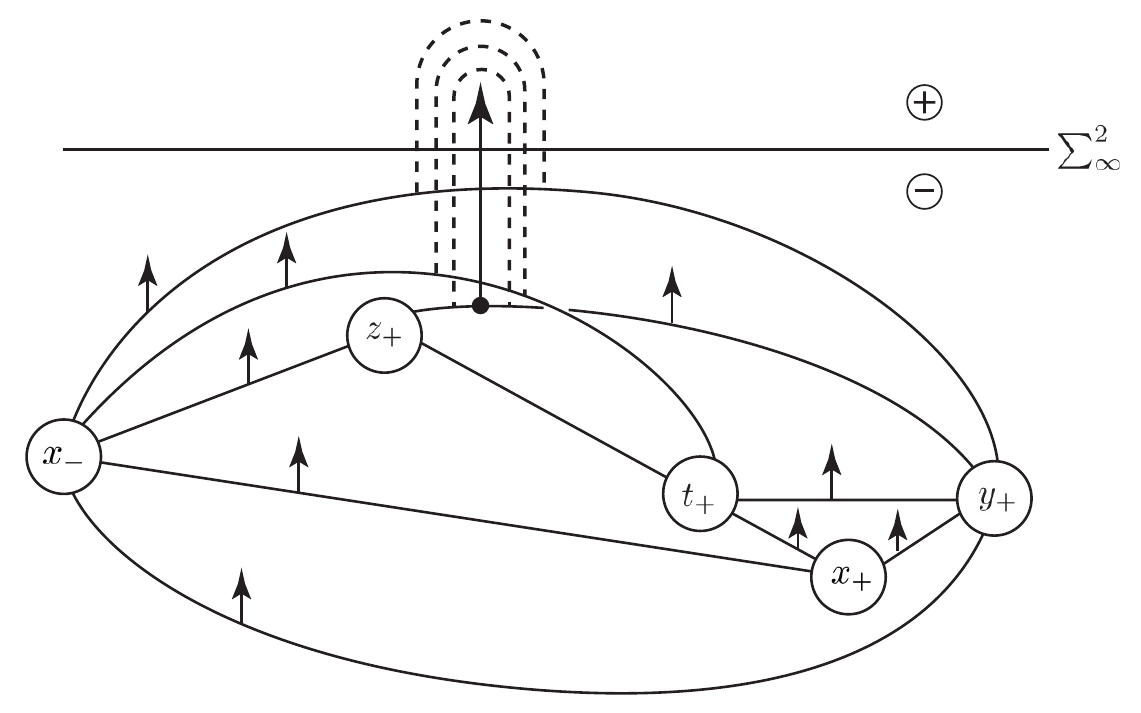}
$$
\label{fig33.1}
\centerline{\bf Figure 33.1.}
\begin{quote}
With $P \times (\xi_0 = -1)$ like in the Figure 24, and with the $b^3 (+ \, \xi_0)$ from that figure, we see here a generic section $S^2 (\xi_0) (0 \geq \xi_0 \geq -1)$ of

$$
\partial N^4 (\Gamma_1 (\infty)) \mid [P \times [0 \geq \xi_0 \geq -1]] = \underbrace{\partial b^3 (\xi_0)}_{\mbox{\footnotesize$S^2 (\xi_0)$}} \times [0 \geq \xi_0 \geq -1].
$$
\end{quote}

\bigskip
\bigskip

\noindent {\bf Some explanations concerning Figure 33.1.}

\medskip

A) (Connection with Figure 24.) The Figure 24 is a projection on the plane of $S_{\infty}^2$ with all the triangles being flat on that plane and the vector field $\vec v$ sticking out of them. What we see here in Figure 33.1, is the way in which the triangles in question cut another 2-sphere, namely the $\partial b^3 (+ \, \xi_0)$. The $S_{\infty}^2$ meets $\partial b^3 (+ \, \xi_0)$ along the equator of $\partial b^3 (+ \, \xi_0)$.

\medskip

B) (How to get this figure from Figure 24.) Imagine $\partial b^3 (+ \, \xi_0)$ as another plane, perpendicular to the $S_{\infty}^2$ from A). This {\ibf is} the plane of the present figure. Each triangle $d^2 (+ \, \xi_0 , u_1 , u_2)$ from Figure 24 corresponds to an arc $[u_1 , u_2]$ here. What Figure 33.1 does, is to give the general idea of how one tackles the second obstruction B.

\smallskip

\noindent End of Explanations.

\bigskip

All this takes completely care of the pushing of $\lambda ({\rm green}) \cap d^2$ into $\partial N^4_+ (\Gamma_1 (\infty))$ and we turn now to body ${\mathcal B}^2$. The body ${\mathcal B}^2$ is made out of disjoined $B^2$-triangles, in the figures of type 9, with
$$
B^2 = B^2 (\beta , \pm \, \xi_0 \ \mbox{or space-time}, \pm \, \xi_0 \ \mbox{or space-time})
$$
and these $B^2$-triangles are joined by rectangles going along the edges. Such a $\beta$-rectangle is displayed in Figure 22 and again in the Figure 34. The $B^2$-triangles can be seen in the Figures 24-(B) and 25-(B). Figure 34 displays a $\beta$-rectangle.

$$
\includegraphics[width=12cm]{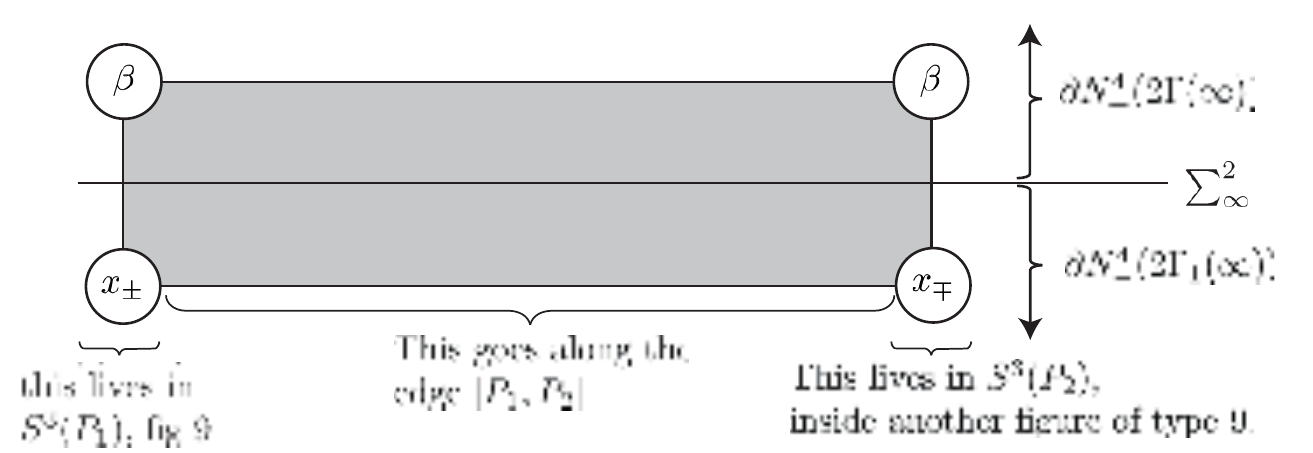}
$$
\label{fig34}
\centerline{\bf Figure 34.}
\begin{quote}
We see here a $\beta$-rectangle going along the edge $[P_1 , P_2 ]$. The body ${\mathcal B}^2$ is made out of $B^2$-triangles and out of rectangles like this one. Each $[\beta , u_{\pm}]$ occurs in some figure of type 9. In Figure 22 we have also displayed a $\beta$-rectangle, shaded, and going from $+ \, \xi_0$ to $- \, \xi_0$. The $[\beta , x_{\pm}]$'s are sides of $B^2$-triangles.
\end{quote}

\bigskip

Of the three sides of a $B^2$-triangle, two are touching to the vertex $\beta$ and, like $\beta$, they are in $\partial N_+^4 (2\Gamma (\infty))$. They are arcs of the curves $C(r \ {\rm or} \ b) \subset \partial N_+^4 (2\Gamma (\infty))$. The third side is contained in an arc [space-time, space-time] (or $[\pm \, \xi_0$, space-time$]) \subset$ curve $C$ (generically $\subset \partial N_-^4 (2\Gamma (\infty))$). This is actually a $\gamma_k^0 = \partial d_k^2$. Figure 25-(B) illustrates the $\Sigma_{\infty}^2 \cap {\mathcal B}^2$ and the transition from ${\mathcal B}^2$ to the $d_k^2$'s, from $(u,z_- , z_+ , y_+)$ down. So, in Figure 25-(B) we have $U \in \partial d_k^2$.

\smallskip

In the figures of type 24, the $B^2$-triangles are higher (i.e. closer to the observer) than the rest. So there is no question of $\lambda ({\rm green}) \cap (B^2$-triangles) to be obstructed like in A).

\smallskip

But there is now another problem, namely we have transversal contacts like in Figure 28
$$
y(p_2) \in \lambda (p_2) \pitchfork \{\mbox{RIBBONS body} \ {\mathcal B}^2 \cap {\rm body} \ d^2 \}.
$$

\bigskip

\noindent {\bf Important Remark.} What we see in the Figure 28 suggests that the transversal contact $y(p_2) \in \lambda (p_2) \pitchfork \{\mbox{very special RIBBON}\}$, is NOT an obstruction for getting to $p_2$. But what we are discussing now is not that, but the potential obstruction for pushing the $\lambda (p_2)$ into $\partial N_+^4 (2\Gamma (\infty))$. End of Remark.

\bigskip

With more details, the contact $\lambda (p_2) ({\rm provisional}) \mid {\mathcal B}^2 \pitchfork \{\mbox{VERY SPECIAL RIBBON}\}$ (Figure 27), is displayed in Figure 25-(B). This concerns now a ($B^2$-triangle) $\subset$ body ${\mathcal B}^2$, where the local piece of $\lambda (p_2)$ which contains $y(p_2)$ is localized, see here the Figure 25-(B).

$$
\includegraphics[width=115mm]{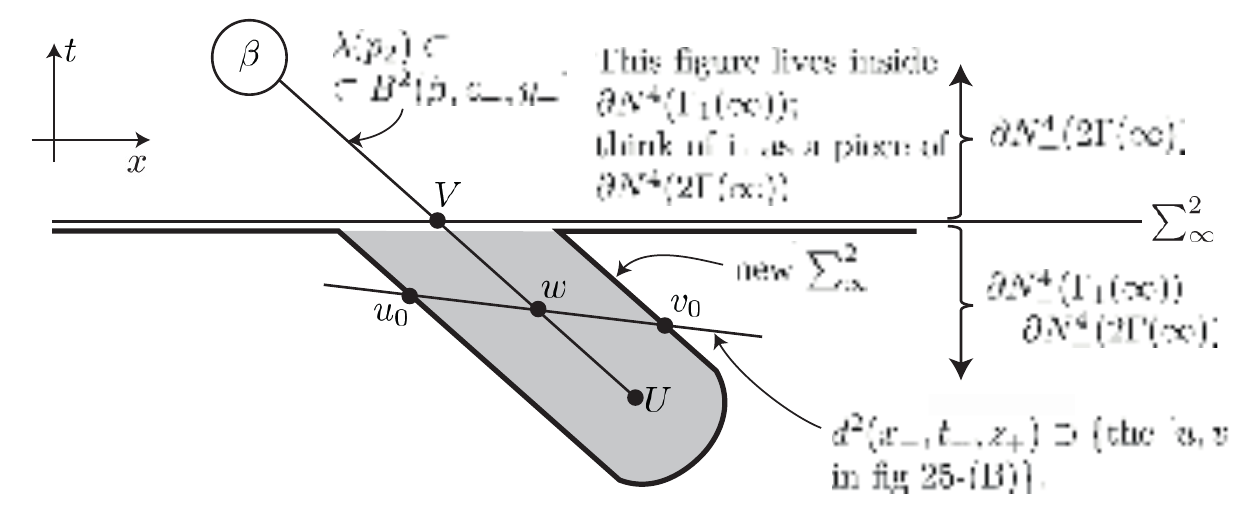}
$$
\label{fig35}
\centerline{\bf Figure 35.}
\begin{quote}
The notations are here like in Figure 25-(B), with the proviso that, the present points $u_0 , v_0$, are the $\{ u,v$ in Figure 25-(B) moved inside $d^2$, until $[u_0 , v_0] \subset {\rm plane} \, (t,x)\}$. The triangle $B^2 (\beta , z_+ , y_+)$ from the 25-(B) lives in a plane which cuts transversally the present plane $(t,x)$, along the green line $[\beta , U]$. Here $w = y(p_2)$.
\end{quote}

Figure 35, where for the time being one should ignore the fat line (which will correspond to a change $\Sigma_{\infty}^2 \Rightarrow$ new $\Sigma_{\infty}^2$), shows a plane $(x,t) \subset \partial N^4 (2\Gamma (\infty))$, which cuts transversally both our triangles $B^2 (\beta , y_+ , z_+) \subset$ body $B^2$ and $d^2 (t_+ , z_+ , x_-) \subset$ body $d^2$, and also the SPLITTING SURFACE $\Sigma_{\infty}^2$. This plane $(t,x)$ contains the line
$$
(\lambda (p_2)({\rm provisional})) \mid [\beta , U] , \mbox{from Figure 25-(B).}
$$
The point here, is that there is no obstruction, meaning NO sacro-sancted confinement conditions, in the way, for performing the following operation: Extend $\partial N_+^4 (2\Gamma (\infty))$ by a $3^{\rm d}$ dilatation {\ibf engulfing} the detail $\left[\xy *[o]=<15pt>\hbox{$\beta$}="o"* \frm{o}\endxy , U\right] \cup [[u_0 , v_0] \subset d^2]$ from Figure 35, with the dilatation in question concentrated around the shaded area in Figure 35.

\smallskip

Notice that we made use here of an engulfing process, rather than of some isotopy of curves. 

\medskip

\noindent [{\bf Remark.} Would we have tried to push isotopically the
$$
[\beta , U] \underset{\overbrace{w}}{\cup} [u_0 , v_0] \ \mbox{from Figure 35}
$$
into $\partial N_+^4 (2\Gamma (\infty))$, then we would have had to drag along appropriate neighbourhoods of it inside $d^2 (\beta , z_+ , y_+)$ $\cup$ $d^2 (x_- , t_+ , z_+)$, keeping at the same time track of boundary conditions. The engulfing is distinctly more economical.]

\smallskip

With this engulfing, which redefines our basic SPLITTING and which hence changes $\Sigma_{\infty}^2$ into $\Sigma_{\infty}^2$ (new), we have now
$$
\left[\xy *[o]=<15pt>\hbox{$\beta$}="o"* \frm{o}\endxy , U\right] \subset \partial N_+^4 (2\Gamma (\infty)).
$$

All this takes completely care of the $\underset{\overbrace{\mbox{\footnotesize $p_2 \in \{$punctures body $\delta^2 \pitchfork \Gamma_i\}$}}}{\sum} \lambda (p_2) \subset \partial N_+^4 (2\Gamma (\infty))$, which has gotten into $\partial N_+^4 (2\Gamma (\infty))$, where we wanted it to be.

\smallskip

So now we have to take care of the $\lambda (q)$'s too.

$$
\includegraphics[width=90mm]{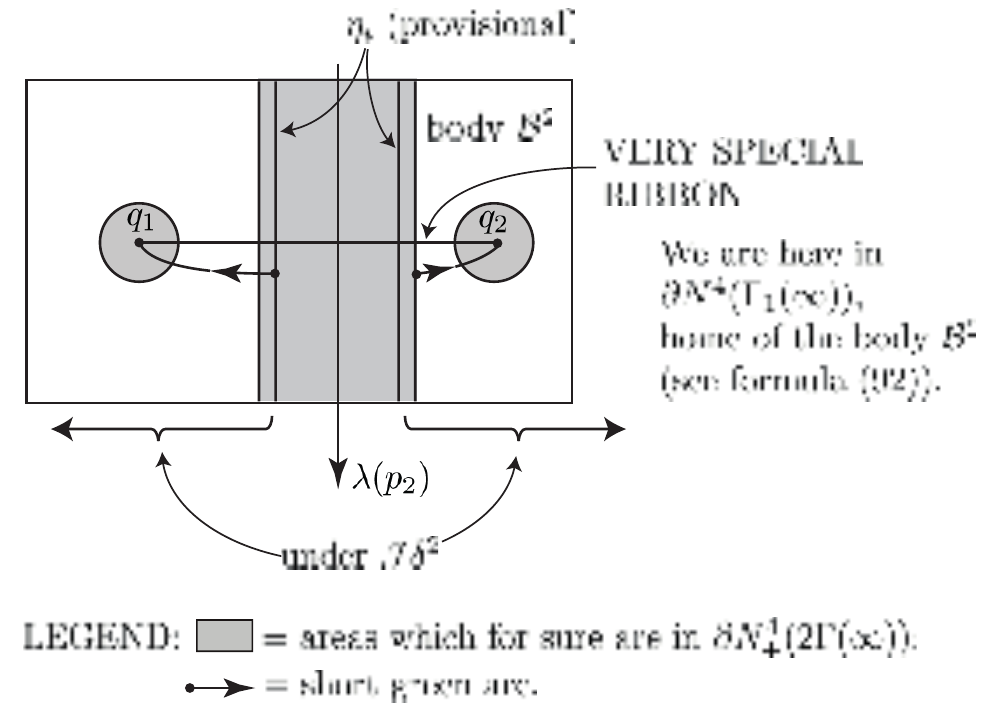}
$$
\label{fig35.1}
\centerline{\bf Figure 35.1.}
\begin{quote}
Compare this with the Figure 28-(II). The $\lambda (p_2)$ has already gone through, by now.
\end{quote}

\bigskip

By combining Figures 28 with 35.1, one sees that there is no obstructions for pushing the short green arcs $\lambda (q_1) , \lambda (q_2)$ into $\partial N_+^4 (2\Gamma (\infty))$. What one sees in Figure 35.1 (and look at 28-(II) too) is what happens along a VERY SPECIAL RIBBON after $\lambda (p_2)$ has gone through $y(p)$. We can go then step by step to the whole component of special RIBBONS (Figure 28) and push all the corresponding short green arcs $\lambda (q) \subset {\rm body} \ {\mathcal B}^2$ into $\partial N_+^4 (2\Gamma (\infty))$; there is no obstruction in the way.

\smallskip

By now we have gotten at level of (93.1), with the $\Sigma \, c(b_i)$ pushed over the $p_2$'s and with the corresponding accidents all destroyed. Let us say that what we have realized is something which goes beyond the (93.1) and where the (98) is now with us too. But we keep the same notations as in (93.1), and our final result is

\bigskip

\noindent (98.1) \quad $\left( \overset{M}{\underset{1}{\sum}} \, \delta_i^2 , \overset{M}{\underset{1}{\sum}} \, \partial \delta_i^2 = c(b_i) ({\rm initial}) \right)$ $\Longrightarrow$ (93.1) $\Longrightarrow$ $\left( \overset{M}{\underset{1}{\sum}} \, \delta_i^2 , \overset{M}{\underset{1}{\sum}} c(b_i)(\lambda) = \partial \delta_i^2 \right)$, when now we also have $\overset{M}{\underset{1}{\sum}} c(b_i)(\lambda) \subset \partial N_+^4 (2\Gamma (\infty))$.

\bigskip

Lemma 12 is now proved.

\bigskip

\noindent {\bf Theorem 13. (The little blue diagonalization)}

\medskip

1) {\it For any long green arc $\lambda(p_2)$ we have a decomposition into two successive pieces $\lambda(p_2) = (\lambda (p_2) \cap {\mathcal B}) \cup (\lambda(p_2) \cap d_k^2)$, with $(\lambda (p_2) \cap {\mathcal B})$ starting at $c(b_i)$ and with $\lambda(p_2) \cap d_k^2$ ending at $p_2$. With this, at the price of complicating the geometric intersections matrices $\Gamma_j \cdot B$ and $\Gamma_j \cdot R$ we can achieve that
$$
(\lambda(p_2) \cap d_k^2) \cdot B = \emptyset . \leqno (99)
$$
This operation does not touch $C \cdot h$ and hence it does not modify the topology of LAVA.}

\medskip

1-bis) {\it The step leading to {\rm (99)} can also be read as an operation which leaves the proper embeddings $B,R \subset N^4 (2\Gamma (\infty))$ unchanged, but  which modifies the $\underset{1}{\overset{\bar n}{\sum}} \, \Gamma_j \subset \partial N^4 (2\Gamma (\infty))$, inside its isotopy class.}

\medskip

2) {\it We also have
$$
(\lambda (p_2) \cap {\mathcal B}) \cdot B \subset \{\mbox{very trivial $B$'s, in the sense of {\rm (94.2)}}\}. \leqno (100)
$$
}

3) {\it When it comes to the short green arcs $\lambda (q)$ we also have, like in $2)$ above,
$$
\lambda (q) \cdot B \subset \{\mbox{very trivial $B$'s, in the sense of {\rm (94.2)}}\}. \leqno (100.1)
$$
}

4) (Reminder) {\it We know already that the system of discs {\rm (98.1)}, which I will re-write here for the convenience of the reader
$$
\left( \sum_1^M \, \delta_i^2 , \sum_1^M \, \partial \delta_i^2 = c(b_i)(\lambda)\right) \xrightarrow{ \ J \ } (\partial N^4 (2X_0^2)^{\wedge} \times [0,1] , \partial N^4 (2X_0^2)^{\wedge} \times \{ 0  \}), \leqno (100.2)
$$
where $\partial N^4 (2X_0^2)^{\wedge} \times [0,1] \subset N_1^4 (2X_0^2)^{\wedge} = N^4 (2X_0^2)^{\wedge} \cup (\partial N^4 (2X_0^2)^{\wedge} \times [0,1])$ has the following features.}
\begin{enumerate}
\item[a)] {\it It is a smooth embedding without any ACCIDENT, i.e. it is a system of discs which is both embedded and external to}
$$
N^4 (2X_0^2)^{\wedge} \subset N^4_1 (2X_0^2).
$$
\item[b)] {\it We have here (and see {\rm (98.1)})}
$$
\sum_1^M c (b_i)(\lambda) \subset (\partial N_+^4 (2\Gamma (\infty))) \cap \partial N^4 (2X_0^2).
$$
\end{enumerate}

5) {\it There is a transformation
$$
\sum_1^M c (b_i)(\lambda) \Longrightarrow \sum_1^M \eta_i (\mbox{green}) \subset \partial N_+^4 (2\Gamma (\infty)) \cap \partial N^4 (2X_0^2) \leqno (101)
$$
which is CONFINED inside $\partial N_+^4 (2\Gamma (\infty))$, and which drags the cobounding $\overset{M}{\underset{1}{\sum}} \, \delta_i^2$ along, so that at the end of $(101)$, with the redefined $\overset{M}{\underset{1}{\sum}} \, \delta_i^2$, we have, like in $(100.2)$
$$
\left( \sum_1^M \, \delta_i^2 , \sum_1^M \, \eta_i (\mbox{green}) = \partial \delta_i^2 \right) \xrightarrow{ \ \mbox{\footnotesize EMBEDDING} \ } (\partial N^4 (2X_0^2)^{\wedge} \times [0,1] , \partial N^4 (2X_0^2) \times \{ 0  \}).
$$
Moreover, now the condition stated in $(82.1)$ is finally satisfied, i.e. we have

\medskip

\noindent $(101.1)$ \quad $\eta_i (\mbox{green}) \cdot b_j = \delta_{ij}$ for $1 \leq i,j \leq M$ AND, at the same time
$$
\eta_j (\mbox{green}) \cdot \left( B_1 - \sum_1^M b_i \right) = 0.
$$
}

6) {\it At this point, we can finally write down the FORCED CONFINEMENT CONDITIONS, superseding from now on the $(97)$, with all the curves, internal and external, presented in their full glory
$$
\sum_i \{ C_i \ \mbox{remaining}\} + \sum_j \{\Gamma_j \ \mbox{remaining}\} + \sum_k \gamma_k^0 \subset \partial N_-^4 (2\Gamma (\infty)), \leqno (101.1\mbox{-{\rm bis}})
$$
AND
$$
\sum_i \{\mbox{little} \ C_i\} + \sum_j \{\mbox{little} \ \Gamma_j\} + \sum \eta_i \times b + \sum c(r) \subset \partial N_+^4 (2\Gamma (\infty)) \supset \sum_1^M \eta_i (\mbox{green}).
$$
}

\bigskip

\noindent {\bf Proof.} We start by reviewing the GPS structure of $X^4$ (see (6)), which extends in an obvious way to $X^2 [{\rm new}]$. And right now we will look into the details of the BLUE collapsing flow.

\smallskip

On the piece $(\Gamma (1) \times [0 \geq \xi_0 \geq -1]) \cup (\Delta^2 \times (\xi_0 = -1))$, things are already settled, since all the $D^2 (\Gamma_i)$'s are $D^2 (\gamma^1)$'s, BLUE-wise, and the BLUE $2^{\rm d}$ collapse proceeds by $\Gamma (1) \times [-1 \leq \xi_0 \leq 0]) \searrow \Gamma (1) \times (\xi_0 = 0)$, a.s.o.

\smallskip

Next, we will take another close look at the BLUE flow on $X^2 ({\rm old}) \subset X^2 [{\rm NEW}]$. The general strategy for this flow is presented in Figure 35.2, and details will follow.

\smallskip

Here is the description of the $2^{\rm d}$ BLUE collapsing flow inside $X^2$ (GPS), as illustrated in Figure 35.2. We insist only on the UPPER-LEFT corner of the figure.

\bigskip

\noindent (101.2) \quad A) On each $X^3 \times t_j$ the trajectories are linear, in the $z$-direction. Any horizontal $2^{\rm d}$ plaquettes ($z =$ const) is a $D^2 (\gamma^1)$ which gets deleted (before any $2^{\rm d}$ collapse can start).

\smallskip

B) Inside $X^3 \times t_i$, the linear $2^{\rm d}$ BLUE collapsing flow, which goes in the $z$-direction, may meet horizontal ($z =$ const) edges $e \subset 2X^1 ({\rm BLACK}) \cup 2X^1 ({\rm ORANGE})$. At this point, the combination
$$
(\mbox{colour of $e$}) \times (\mbox{parity of $i$}), \leqno (*_1)
$$
decides if $e \times [t_i , t_{i \pm 1}] \subset X^2$ (see (6)), and these 2-cells $e \times [t_i , t_{i \pm 1}]$ are declared to be $D^2 (\gamma^1)$'s. Their interiors are to be {\ibf deleted}, like for the $D^2 (\gamma^1)$'s in A) above, before any $2^{\rm d}$ BLUE collapse can start.

$$
\includegraphics[width=11cm]{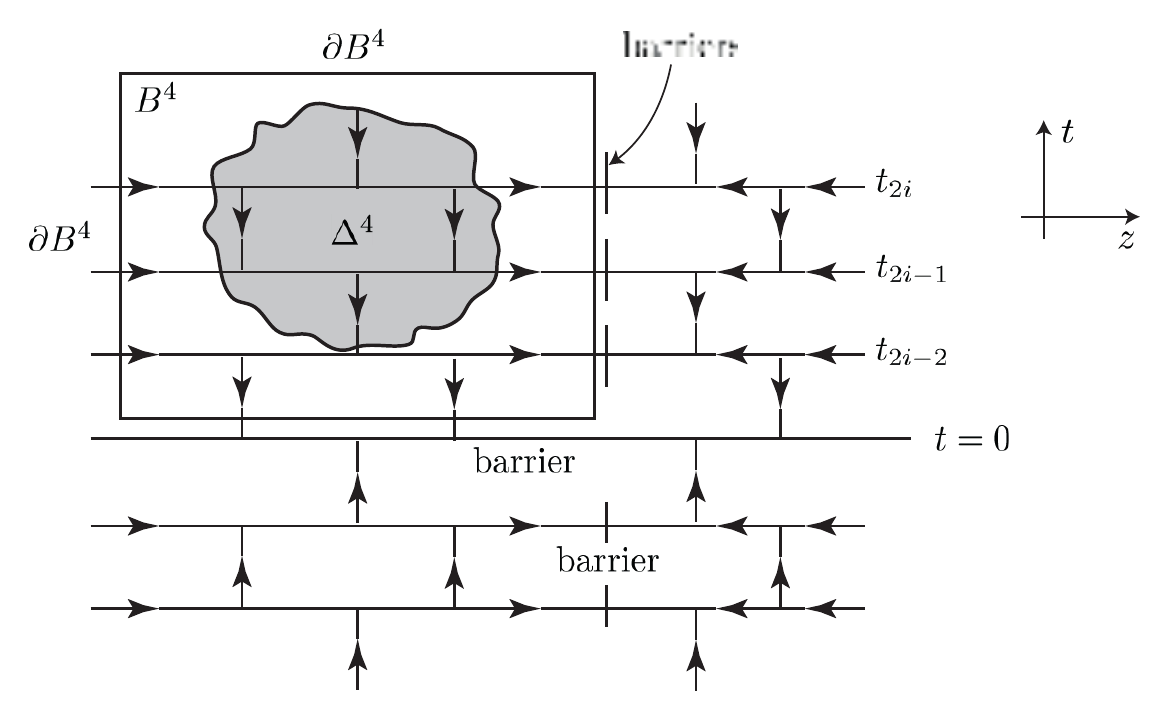}
$$
\label{fig35.2}
\centerline{\bf Figure 35.2.}
\begin{quote}
The general strategy for the $2^{\rm d}$ BLUE collapsing flow on $X^2 ({\rm old})$ (GPS). This is, of course, a very schematical figure.

LEGEND: $\rightarrow \leftarrow = 2^{\rm d}$ flow inside $X^3 \times t_j$; $\uparrow \downarrow \uparrow \downarrow = 2^{\rm d}$ flow inside $[t_j , t_{j+1}]$; $-\!\!-\!\!-\!\!-$, $\mid$ $=$ BARRIERS. The flow can be fixed on the barriers too, but we will do that only if and when needed. Importantly, it is only the upper-left corner, with $\Delta^4_{\rm Schoenflies} \subset B^4$ which really concerns us. The $\partial B^4$ rests on the BARRIERS. Vertical arrows starting at a $t_{\rm odd} > 0$ or $t_{\rm even} < 0$ are ORANGE, all the others are BLACK.
\end{quote}

\bigskip

\noindent [{\bf Remark.} We are a bit cavalier concerning these $D^2 (\gamma^1)$'s. Normally, there should be a $3^{\rm d}$ BLUE collapsing flow, prior to the $2^{\rm d}$ one, which should demolish them. But we are in a purely $2^{\rm d}$ context, so we demolish the $D^2 (\gamma^1)$'s by decree. If it would be there, the $3^{\rm d}$ collapsing flow would also need its $3^{\rm d}$ BARRIERS. We can assume them far from our interesting UPPER-LEFT corner.]

\bigskip

C) After all the ${\rm int} \, D^2 (\gamma^1)$'s are deleted, the BLUE $2^{\rm d}$ collapse takes place, and it leaves us with a residual graph $\Gamma_{\rm residual} \subset (X^1 ({\rm BLACK}) \cup X^1 ({\rm ORANGE})) \mid t_i$. The common part of the two $X^2$'s is generated by things like $a,f,d$, in the Figure 5-(B). With these things, we consider now
$$
\Gamma_{\rm residual} \cap (2X^1 ({\rm BLACK}) \cup 2X^1 ({\rm ORANGE})) \mid t_i . \leqno (*_2)
$$

Any edge $e \in (*_2)$ goes in the $z$-direction and carries a $b_0 \in B_0$; the BLUE $2^{\rm d}$ collapsing flow from $b_0$ goes in the direction $\pm \, t$, along $[t_i , t_{i\pm 1}]$, by the $(*_1)$-dependent rules of Figure 35.2. All these $b_0$'s are very trivial, in the sense of (94.2).

\smallskip

Concerning our same $2^{\rm d}$ BLUE collapsing flow, inside $\Gamma (1) \times [0 \geq \xi_0 \geq -1]$ it follows the direction $+ \, \xi_0$ and all the $D^2 (\Gamma_i) \times (\xi_0 = -1)$'s are $D^2 (\gamma^1)$'s. All this clinches the definition of the $2^{\rm d}$ BLUE flow for $X^2 [{\rm new}]$ (GPS). On the other hand, the RED $2^{\rm d}$ collapsing flow depends on how the $\Delta^4_{\rm Schoenflies}$ is located inside $X^4$, and so it cannot be presented in a similar cavalier and explicit manner. But then, see here the PL-Lemma 3.1 too.

\smallskip

Once the structure of the $2^{\rm d}$ BLUE collapsing flow has been completely unrolled, with its linear {\ibf trajectories}, we will look now at the BLUE $2^{\rm d}$ collapsing trajectory of generic $b_i \in \Gamma (1) \times (\xi_0 = -1)$. This is displayed in Figure 35.3, where one also sees the $\lambda(p_2) \cap {\mathcal B}$.

$$
\includegraphics[width=155mm]{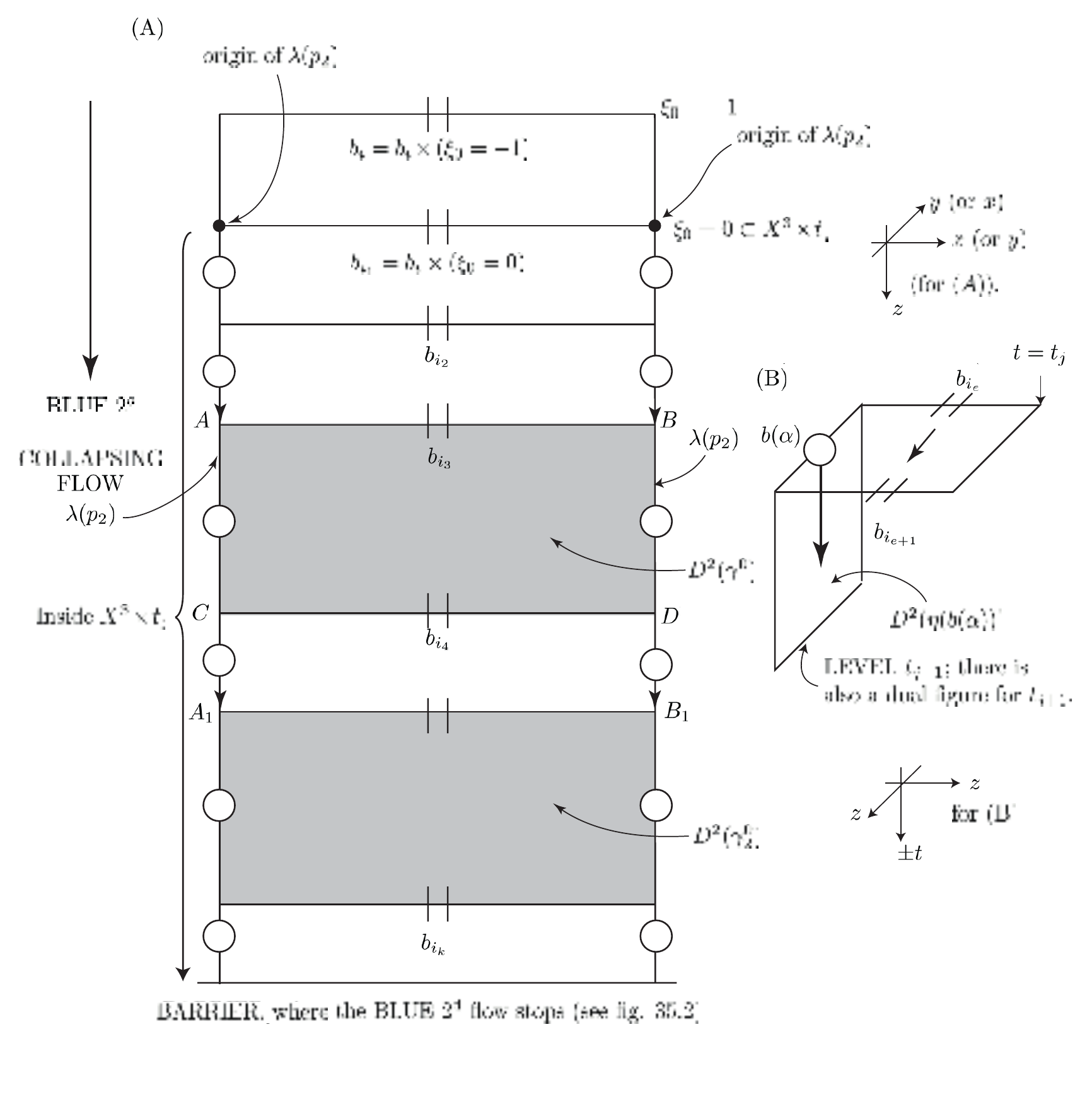}
$$
\label{fig35.3}
\centerline{\bf Figure 35.3.}
\begin{quote}
LEGEND: $\includegraphics[width=8mm]{rectangleclair.pdf} \ = D^2 (\gamma_k^0)$; $\includegraphics[width=3mm]{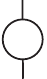} = b(\alpha) \in B_0 (?)$. This $b(\alpha)$ may be or not be in $B_0$, depending on the combination $(*_1)$; $\longrightarrow \ = \lambda (p_2) \cap {\mathcal B}^2$. It has several parallel strands, and at $A,B,A_1 , B_1$, it leaves ${\mathcal B}^2$.
\end{quote}

\bigskip

Forgetting momentarily about the $b(\alpha)$'s, (which are very trivial element of $B$), our trajectory of $b_i = b_i \times (\xi_0 = -1)$ is $b_i \to b_{i_1} \to b_{i_2} \to \ldots$, realized by a vertical, linearly ordered column of $2^{\rm d}$ plaquettes $D^2 (\eta)$. Some of them come with $[D^2 (\eta)] = D^2(\gamma^0)$, shaded in Figure 35.3-(A), and then the piece $\Sigma \, d_k^2 \subset \delta_i^2$ starts.

\smallskip

The $b(\alpha)$'s are the potential sites corresponding to a $b_0 \in B_0$, living on an edge $e \in (*_2)$, produced by the appropriate combination $(*_1)$, see here (101.2). More explicitly, if the corresponding edge $e$ is in $(2X^1 ({\rm BLACK}) \cup 2X^1 ({\rm ORANGE})) \cap \Gamma_{\rm residual}$ then, depending on the combination COLOUR/parity of $j$ (occuring in $t_j$), we have actually a $b \in B_0$ at $b(\alpha)$. I repeat that these $b(\alpha)$'s are very trivial $B$'s (see (94.2)).

\smallskip

On $c(b_i) \cap (\xi_0 = 0)$ we see two green fat points (``origin of $\lambda(p_2)$'') from which green arcs start to each $D^2 (\gamma^0) = [D^2 (\eta)]$. These arcs may have pieces in common, that is OK. These green arcs, seable in Figure 35.3, are $\lambda(p_2) \cap {\mathcal B}$'s.

$$
\includegraphics[width=125mm]{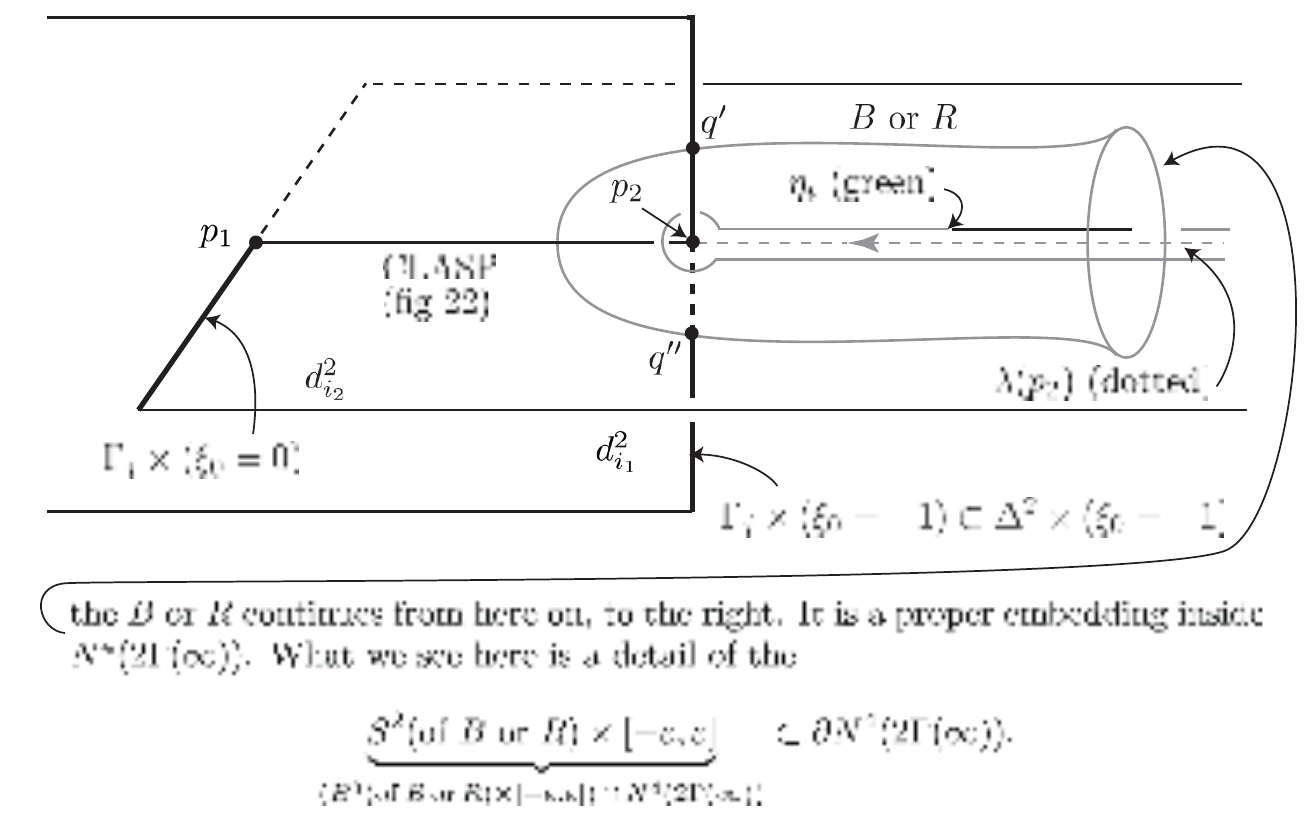}
$$
\label{fig35.4}
\centerline{\bf Figure 35.4.}
\begin{quote}
This figure which lives inside $\partial N^4 (\Gamma_1 (\infty))$ is to be compared to the Figures 20, 22, 23. The contact
$$
(\lambda (p_2) \cap d_k^2) \pitchfork \{{\rm cocore} \ B \ {\rm or} \ R\}
$$
is destroyed by pushing the cocore in question over the puncture $p_2 \in \Gamma_j \times (\xi_0 = -1)$. This will create new intersection points
$$
q' , q'' \in (\Gamma_j \times (\xi_0 = -1)) \cap \{{\rm cocore} \ B \ {\rm or} \ R\} ,
$$
which increase the corresponding geometric intersection matrices.  Here only a piece of the $\{{\rm cocore} \ B \ {\rm or} \ R\}$ is visible. It is a $2$-cell and there are many parallel copies of it, in order to get the full picture around $p_2$. But then each such $2$-cell continues to a whole copy of
$$
S^2 \subset S^2 \times [-\varepsilon , \varepsilon] \subset \partial N^4 (2\Gamma (\infty)),
$$
continuing with a 1-handle $(B^3 (\times [-\varepsilon , \varepsilon]) \subset N^4 (2\Gamma (\infty))$ (BLUE or RED).
\end{quote}

\bigskip

We move now to the point 1) in our Theorem 13. A priori we have transversal contacts
$$
(\lambda (p_2) \cap d_k^2) \pitchfork B \leqno (102)
$$
and these (102) are intermingled along the $\lambda(p_2) \cap d_k^2$ with similar contacts
$$
(\lambda (p_2) \cap d_k^2) \pitchfork R . \leqno (102.1)
$$

The IDEA now is to PUSH THE (102) $+$ (102.1) OVER THE PUNCTURE $p_2 \in \Gamma_j \times (\xi_0 = -1)$. Here are the various noteworthly items concerning this STEP, which is displayed graphically in the Figure 35.4.

\smallskip

A) All the punctures $p_2$ occur on curves $\Gamma_j \times (\xi_0 = -1) \subset \Delta^2 \times (\xi_0 = -1)$, making that our step does not change the matrix $C \cdot h$ and hence it does not touch the LAVA which has its topology stay intact. It does not touch the geometric intersection matrix $\eta \cdot \beta$ either, since both at the levels $X^2 [{\rm NEW}]$ and $2X_0^2$, $D^2 (\Gamma_j \times (\xi_0 = -1))$ is always a $D^2 (\gamma^1)$, never a $D^2 (\eta)$.

\smallskip

B) Like in Figure 35.4, we create contacts $q',q'' \in (\Gamma_j \times (\xi_0 = -1)) \cap \{$cocore of $B$ or $R\}$, actually $q',q'' \in\{$little $\Gamma_j \times (\xi_0 = -1)\} \cap \{$cocore $B$ or $R\} \subset \partial N^4_+ (2\Gamma(\infty))$.

\smallskip

These do NOT OBSTRUCT anything; I mean they leave $\Gamma_j \times (\xi_0 = -1) \subset \partial N^4 (\Gamma_1 (\infty))$ located just as before. The only thing they change are the items $(\Gamma_j \times (\xi_0 = -1)) \cdot B$ and $(\Gamma_j \times (\xi_0 = -1)) \cdot R$ in the geometric intersection matrices, and this is harmless.

\smallskip

C) The way we presented things in Figure 35.4 is to complicate the picture of the $\{$cocore $(B,R)\}$, leaving the curves $\Gamma_j$ in peace. But there is also another way of describing the same step. Leave the triplet
$$
(N^4 (2\Gamma (\infty)); \mbox{cocores of $B$, cocores of $R$)}
$$
completely unchanged; i.e. do not modify Figure 46 which will be crucial for our future COLOUR-CHANGING. BUT then we start by {\ibf shortening} drastically the $\lambda (p_2) \cap d_k^2$, starting at $p_2$, push the $\Gamma_j \times (\xi_0 = -1)$ through our $\{$cocores $R$ and $B\}$, until the $p_2$ curves now very close to the point in $\partial d_k^2 \subset {\rm int} \, \delta_i^2$ where we find $\lambda (p_2) \cap \partial d_k^2$. This clinches the proof of our 1).

\smallskip

The point 2) is readable in the Figure 35.3 when $\lambda (p_2) \cap {\mathcal B}$ is concentrated on the two lateral sides, avoiding the $b_{i_j}$'s and only meeting the very trivial $b(\alpha)$'s. The linearity of the BLUE collapsing trajectories plays here.

\smallskip

We move now to 3) which concerns the short $\lambda(q) \subset {\rm body} \, {\mathcal B}^2 \subset \delta_i^2$. There is here the obvious potential danger, readable in the Figure 25-(B): We may need to propagate $\lambda (q)$ from $p$ to $q$, over the edge $[P_1 , P_2]$ and {\ibf if} this edge is occurring now explicitly in Figure 35.3 and is horizontal, like $[AB]$, the $[P_1 P_2]$ contains one of the non trivial $b_{i_k}$'s from the BLUE collapsing trajectory of $b_i = b_i \times (\xi_0 = -1)$. At this point, we have to look carefully at the precise manner in which the $\delta_i^2$ is being puts up at point 1) in Lemma 10, from the successive pieces $B^2_{i_k}$ (and $B^2 (b(\alpha))$, corresponding to Figure 35.3. Here, when the $B_{i_k}^2$ and $B^2_{i_{k+1}}$ are put together, then the common collar $c(b_{i_{k+1}}) \times [0,\varepsilon]$ gets deleted.

\smallskip

So, let us consider the potentially dangerous situation when $\lambda(q)$ goes along the $[p,q]$ in Figure 25-(B), with (and see here the notations from Figure 25-(B)) $b_{i_{k+1}} \in [P_1 , P_2]$. I CLAIM that, actually, this $\lambda(q)$ is {\ibf not physically present in} $\delta_i^2$, and hence we do not have to worry about it. 

\medskip

[{\bf Proof of the Claim.} We work now with the decomposition $\delta_i^2 = {\mathcal B}_i^2 \cup \Sigma \, d_k^2$'s and the ${\mathcal B}_i^2$ is completely readable from Figure 35.3, which contains all the necessary ingredients. We have

\bigskip

\noindent $\lambda (q) \mid [P_1 , P_2] \subset \{$the two triangles from Figure 25-(B), which live, respectively in $\partial N^4 (P_1), \partial N^4 (P_2)\} \cup \{$The $\beta$-rectangle $[P_1 , P_2]\} \subset\{$the collar zone $c(b_{i_{k+1}}) \times [0,\varepsilon]\}$.

\bigskip

And the collar zone in question gets deleted when we put up $\delta_i^2$. This proves our CLAIM.]

\bigskip

All this shows that the $\lambda (q)$'s finding themselves naturally in the rectangles
$$
\{\mbox{horizontal edge like $[A,B]$, in Figure 35.3}\} \times [r,b]
$$
are not actually in ${\mathcal B}_i \subset \delta_i^2$.

\smallskip

So, the fa\c cade of ${\mathcal B}_i^2$, which is presented explicitly in Figure 35.3, comes with no problem for $\lambda (q) \cdot B$ in 3). We move then to the
$$
\{\mbox{two lateral vertical sides of Figure 35.3}\} \times [r,b], \leqno (**)
$$
which of course is not visible in the figure in question. In that region the only contacts $\lambda(q) \cdot B$ which we may find are with the $b(\alpha)$'s which is OK. [Notice, incidentally, that the procedure from Figure 35.4, via which we have demolished the contacts
$$
(\lambda (p_2) \cap d_k^2) \cap (B \ \mbox{\ibf and} \ R), \mbox{which concerned the} \ \Gamma_j \times (\xi_0 = -1),
$$
could never work for $\lambda (q) \cap R$ for which the COLOUR-CHANGING mechanism is certainly necessary; remember that the $\lambda(q)$'s also contribute to $\eta_i ({\rm green})$.]

\smallskip

Point 3) in our Lemma has by now been proved. Also our $\lambda$'s are all there and point 4) is with us too. So we go to point 5).

\smallskip

Each $c(b(\alpha))$ bounds the disc $B^2 (b(\alpha)) \subset 2X_0^2$, from (79.0). Moreover we have as the only contacts of $c(b(\alpha))$ with $B_1$ the following ones
$$
c(b(\alpha)) \cdot b(\alpha)= 1 = c(b(\alpha)) \cdot (b(\alpha) \times b) . \leqno (102.2)
$$

Along the $b(\alpha)$, the two loops $c(b(\alpha))$ and $c(b_i)(\lambda)$ touch, when $\lambda (p_2) \cap b(\alpha) \ne \emptyset$ and $\lambda (p_2) \subset \delta_i^2$. From here on, we proceed via the following steps.

\medskip

i) Making use of (102.2), the contact $c(b_i)(\lambda) \cdot b(\alpha)$ can be changed into $c(b_i)(\lambda) \cdot (b(\alpha) \times b)$, dragging of course the $\delta_i^2$ along, through an ambient isotopy of $\partial N^4 (2X_0^2)^{\wedge} \times [0,1]$. What we have gained via this step is that now the only off-diagonal contacts $c(b_i) (\lambda) \cdot B_1$ are exactly the following
$$
c(b_i)(\lambda) \cdot (b_i \times b) = 1, \quad c(b_i)(\lambda) \cdot (b(\alpha) \times b) = 1, \leqno (102.3)
$$
both living on the $X_b^2$-side.

\medskip

ii) We use now our BLUE geometric intersection matrix, which has been transported in the $X_b^2$-side,
$$
\eta \cdot B \mid 2X^2 = \mbox{easy id $+$ nilpotent}
$$
and get rid of the off-diagonal terms (102.3). This last step is completely concentrated on the $X_b^2$-side, like i) too. All this also clinches the transformation
$$
c(b_i) \Longrightarrow c(b_i)(\lambda) \Longrightarrow \eta_i ({\rm green})
$$
which is all confined inside $\partial N_+^4 (2\Gamma (\infty))$. Lemma 12 is, of course, used here. End of the Proof of Theorem~13.

\bigskip

\noindent {\bf Remark.} Notice that it is on the $(**)$ above (i.e. the $\{$lateral vertical sides of Figure 35.3$\} \times [r,b]$), that the mechanism from the Figures 28, 29, occurring in the proof of Lemma 11.3 come fully into play.

\section{The balancing of red and blue, and the abstract theory of colour-changing}\label{sec6}

We will consider an infinite connected graph $\Gamma$. In real-life this will be $\Gamma (\infty) = \Gamma (\infty) \times r$ or $\Gamma (2\infty)$, but let us be, for a while, a bit more general. By definition, a discrete set $E \subset \Gamma$ has the $P$-{\ibf property} if $\Gamma - E$ is a tree. An edge $e \subset \Gamma$ contains at most one $x \in E$ and we will often identify $x$ and $e \ni x$, thinking also of $x$ as being a $1$-handle (attached to $\Gamma - E \approx {\rm pt}$).

\smallskip

We are given two discrete subsets $R \subset \Gamma \supset B$ each endowed with the $P$-property. It is also understood that, when $x \in e \ni y$, then $x = y$, leaving us hence with the partitions
$$
R = (R-B) + (R \cap B) , \quad B = (B-R) + (R \cap B).
$$

\bigskip

\noindent {\bf Lemma 14. (The abstract colour-changing process.)} {\it We can get a bijection
$$
R \overset{\Phi}{\underset{\approx}{-\!\!\!-\!\!\!-\!\!\!-\!\!\!\longrightarrow}} B , \ \mbox{which is id on} \ R \cap B , \leqno (103)
$$
via the following INDUCTIVE PROCESS.}

\bigskip

\noindent {\bf Lemma 14. (The abstract colour-changing process.)} {\it There is a bijection
$$
R \overset{\Phi}{\underset{\approx}{-\!\!\!-\!\!\!-\!\!\!-\!\!\!\longrightarrow}} B , \ \mbox{which is id on} \ R \cap B , \leqno (103)
$$
gotten via the following INDUCTIVE PROCESS.}

\medskip

1) {\it Since $R$ has property $P$ for $\Gamma$ and $(\Gamma - R \cap B) - (R-B) = \Gamma - R$, it follows that $R - R \cap B$ has Property $P$ for $\Gamma - R \cap B$ and then, similarly $B-B \cap R$ also has property $P$ for the same $\Gamma - R \cap B$. With this, $B \cap R$ will be just {\ibf mute} in what follows next and we will just work with $\Gamma \equiv \Gamma - R \cap B$, $R-B \cap R$ and $B - B \cap R$. Next, let us pick some arbitrarily chosen $b(1) \in B-R$; since $R \in P$, we have
$$
\Gamma - R = X_1 \cup b(1) \cup Y_1 , \leqno (104)
$$
where $X_1 , Y_1$ are two disjoined trees, joined together by $b(1)$.}

\medskip

2) {\it With all these things, there has to be an $r(1) \in B-R$ which joins $X_1$ and $Y_1$ too. We define then $\Phi (r(1)) = b(1)$ and I claim that the property $P$ is true for $R(1) \equiv R - r(1) + b(1)$. [In this little story, the fact that $\Gamma - R$ is a tree forces the existence of $b(1)$ for $(104)$, while the fact that $\Gamma - B$ is tree, forces the existence of $r(1)$. There is a lot of arbitrariness in the way the map $\Phi$ is constructed, but this will turn out to be OK.]}

\medskip

3) {\it Assume, inductively, that continuing this, for $n = 1,2,\ldots,j$ we have already managed to define a map $r(n) \overset{\Phi}{\longmapsto} b(n)$ (bijection for $n \leq j$), and that the $R(j) \equiv R - \{ r(1) , r(2), \ldots , r(j)\} + \{ b(1) , b(2) , \ldots , b(j)\}$ has Property $P$.

\smallskip

Notice here that $R(j) = R(j-1) - r(j) + b(j)$.

\smallskip

We will pick up some $b(j+1) \in B - R(j)$ and with this, we get a decomposition which is like in $(104)$, but now at a higher level, namely the
$$
\Gamma - R(j) = X_{j+1} \cup b(j+1) \cup Y_{j+1} . \leqno (105)
$$
There has to be some $r(j+1) \in R(j) - B$ connecting $X_{j+1}$ to $Y_{j+1}$. Then we set, by definition $\Phi (r(j+1)) = b(j+1)$ and this kind of process continues now indefinitely. This {\ibf is} our INDUCTIVE PROCESS.}

\medskip

4) {\it And this process gives us a map
$$
R \supset B \cap R + \sum_1^{\infty} r(n) \overset{\Phi}{-\!\!\!-\!\!\!-\!\!\!-\!\!\!\twoheadrightarrow} B
$$
which is surjective on $B$ and also injective on its domain of definition.}

\medskip

5) {\it But I {\ibf claim} that we also have $R \cap B + \underset{1}{\overset{\infty}{\sum}} \, r(n) = R$, and this clinches then our $(103)$. The $\Phi$ above is a bijection.}

\bigskip

With this our statement of Lemma 14 is ended, but let us notice, right away, that there is a slightly different way of formulating 4) $+$ 5). So, let us restart from the end of 3) when the INDUCTIVE PROCESS has already been unrolled.

\smallskip

From this on, our lemma says the following things (and whenever appropriate, the $R,B$'s may mean here $R-B , B-R$).

\medskip

6) This process can be continued until all of $B$ is exhausted. When that point has been reached we have gotten, for $j = \infty$, a set
$$
R(\infty) = R - \{\mbox{a certain subset} \ R_0 \subset R \} + B .
$$

7) We have $R(\infty) = B$. Hence, also, $R_0 = R$ and at the grand end of the INDUCTIVE PROCESS, all the red elements of $R$ {\ibf have turned} BLUE. 

\smallskip

This ENDS the alternative way of stating Lemma 14. $\Box$

\bigskip

\noindent {\bf Remarks.} A) Of course, so far all this is a purely ABSTRACT story, without any geometrical counterpart. But that geometric counterpart will come soon. And, in four dimensions, via handle-sliding and LAVA MOVES, we will be able to play a sufficiently high truncation of the $R(\infty) = B$ from 7) above, so that we should get the BIG BLUE DIAGONALIZATION
$$
\eta_j ({\rm green}) \cdot \{\mbox{extended cocore} \ (b_i)\}^{\wedge} = \delta_{ij} , \quad {\rm for} \ 1 \leq i,j \leq M.
$$

B) In the process, the $C \cdot h = {\rm id} + {\rm nilpotent}$ will get VIOLATED, but the LAVA MOVES will be such that the PRODUCT PROPERTY OF LAVA will always stay with us, and it is via that property that the extended cocores will be defined. End of Remarks.

\bigskip

\noindent {\bf Proof of Lemma 14.} In the Figure 36 we have represented, schematically, the decomposition (104). In the same figure, we introduce the single closed loop $\partial \, C_1^2 \to \Gamma$, boundary of the purely abstract 2-cell $C_1^2$, which is such that
$$
\partial \, C_1^2 \cdot r(1) = 1 = \partial \, C_1^2 \cdot b(1) , \quad {\rm and} \quad \partial \, C_1^2 \cdot (R-r(1)) = 0 = C_1^2 \cdot (R(1) - b(1)).
$$

$$
\includegraphics[width=8cm]{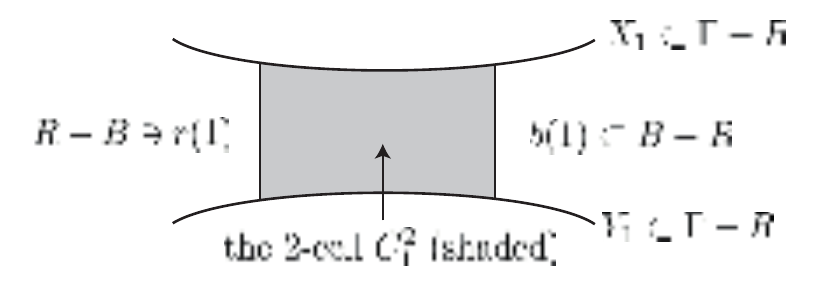}
$$
\label{fig36}
\centerline{\bf Figure 36.}
\begin{quote}
Illustration for (104). $\Gamma$ really means here $\Gamma - R \cap B$. The $\partial \, C_1^2$ is a simple loop and $C_1^2$ (shaded) is abstract. The position of the subsets $\{ r(1) + b(1)\} \cap X_1$ and $\{ r(1) + b(1)\} \cap Y_1$ determines how, in our drawing, we orient $X_1 , Y_1$ with respect to each other.
\end{quote}

\bigskip

We will prove now that $R(1) \in P$ and the same kind of arguments will work for $R(j) \in P$ too. The set $R$ is a free basis for $H_1 (\Gamma , {\rm rel} (\Gamma - R)) = H_1 \, \Gamma$. We introduce the abstract 2-cell $C_1^2$ like in Figure 36. Then, in
$$
\Gamma -b(1) \subset \Gamma \xhookrightarrow{ \ p_1 \ } \Gamma \cup C_1^2 \equiv \Gamma \underset{\partial \, C_1^2}{\cup} C_1^2
$$
\vglue-5mm
\hglue 62mm$
{\mbox{\hglue -5mm}{\mid\mbox{\hglue -2mm}}}_{\overset{f_1}{-\!\!-\!\!-\!\!-\!\!-\!\!-\!\!-\!\!-\!\!-\!\!-\!\!-\!\!-\!\!-\!\!-\!\!-\!\!-\!\!-\!\!-}}{\mbox{\hglue -2mm}\uparrow}
$

\medskip

\noindent the homology class $[r(1)] \in H_1 \, \Gamma$ gets killed (exactly) by $p_1$, while $f_1$ is a homology equivalence. (It is, of course, a homotopy equivalence already, but that is immaterial for our arguments.) Since $C_1^2$ is far from $R-r(1) = R(1) - b(1)$, the following restriction of $f_1$ is also a homology equivalence
$$
\Gamma - R(1) = (\Gamma - b(1)) - (R(1) - b(1)) \underset{f_1}{\xrightarrow{ \qquad }} \Gamma \cup C_1^2 - \underbrace{(R(1) - b(1))}_{= \, R - r(1)} \, ;
$$
this fact is graphically displayed in Figure 36, where $R-r(1)$ is not at all represented; let us say that it is not physically present in the figure.

\smallskip

Since $\Gamma - R$ is connected, so is also the $\Gamma - (R-r(1))$, and this implies that the following is connected too
$$
\Gamma \cup C_1^2 - (R -r(1)) = \Gamma \cup C_1^2 - (R(1) - b(1)). \leqno (106)
$$

It follows that $\Gamma - R(1)$ is connected. Next I CLAIM that $H_1 (\Gamma - R(1)) = 0$.

\medskip

\noindent [{\bf Proof of the Claim.} By going from $\Gamma$ to $\Gamma \cup C_1^2$, we kill $[r(1)] \in H_1 (\Gamma) = H_1 (\Gamma \, {\rm rel} (\Gamma - R))$. When we further delete $R(1) - b(1)$ in the RHS of the formula (106) and hence get the $\Gamma - R(1)$, we also kill all the other $[x] \in H_1 \, \Gamma$; one uses here $R(1) - b(1) = R-r(1)$.] Hence $R(1) \in P$.

\bigskip

By now, the first step in our induction is completely taken care of.

\smallskip

For the step $j-1 \Rightarrow j$, we replace Figure 36 by its higher analogue, Figure 37. In the context of this figure, we find that
$$
\partial \, C_j^2 \cdot b(j) = 1 = \partial \, C_j^2 \cdot r(j) \quad {\rm and} \quad C_j^2 \cdot (R(j-1) - r(j)) = 0 = \partial \, C_j^2 \cdot (R(j) - b(j) ). \leqno (107)
$$

The argument used above can be extended to show that $R(j) \in P$. This clinches the inductive process, and we are left with proving point 5); the 4) should be clear already, by now.

$$
\includegraphics[width=8cm]{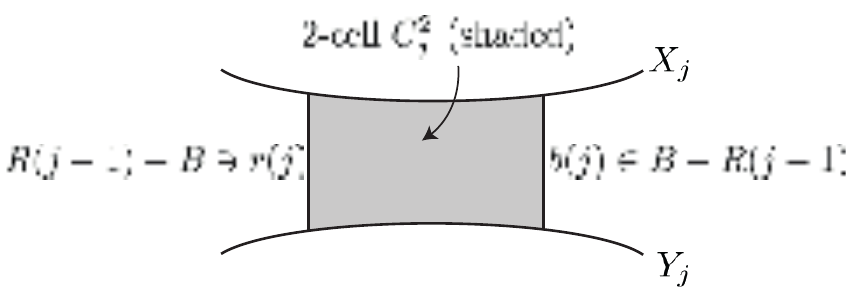}
$$
\label{fig37}
\centerline{\bf Figure 37.}
\begin{quote}
Illustration for (105.j), i.e. for $\Gamma - R(j-1) = X_j \cup b(j) \cup Y_j$, the higher stage of our intersection. Here
$$
R(j-1) = R - \sum_{\ell = 1}^{j-1} r(\ell) + \sum_1^{j-1} b(\ell) \supset \sum_{j+1}^{\infty} r(\ell).
$$
This means that we have
$$
\partial \, C_j^2 \cdot \sum_{j+1}^{\infty} r(\ell) = 0 ;
$$
see the main text.
\end{quote}

\bigskip

We have
$$
R \underset{\mbox{\footnotesize inclusion}}{\xhookleftarrow{ \qquad \Psi \qquad }} B \cap R + \sum_1^{\infty} r(n) \underset{\mbox{\footnotesize bijection}}{\xrightarrow{ \qquad \Phi \qquad }} B ,
$$
and we need to prove that
$$
T \equiv \Gamma - \left[ (B \cap R) + \sum_1^{\infty} r(n) \right]
$$ 
is a tree (i.e. that $H_0 \, T = Z$, $H_1 \, T = 0$); in that case $\Psi$ has to be bijective. [Of course, we could also have worked here with $R \cap B$ deleted everywhere, i.e. with $R - B \cap R \hookleftarrow \underset{1}{\overset{\infty}{\sum}} \, r(n) \to B-B\cap R$, a.s.o.] And since, anyway all our construction has proceeded with the $B \cap R$ deleted, we certainly have $\partial \, C_j^2 \cdot (R \cap B) = \emptyset$ and this allows us to introduce the space $(\Gamma - B \cap R) \cup \underset{1}{\overset{\infty}{\sum}} \, C_j^2$. Here, there are two maps, and in the formula below, the $\ell$ means some arbitrary generic level:
$$
T \xhookrightarrow{ \ \ \zeta_1 \ \ } (\Gamma - R \cap B) \cup \sum_1^{\infty} C_j^2 \xleftarrow{ \ \ \zeta_2 \ \ }(\Gamma - R \cap B) - \sum_{n=1}^{\ell} \ell (n) \cup \sum_{\ell + 1}^{\infty} C_j^2 .
$$
For every $j$, we have $\partial \, C_j^2 \cdot r(j) = 1$ and $\partial \, C_j^2 \cdot \underset{j+1}{\overset{\infty}{\sum}} \, r(n) = 0$, i.e. the matrix $\partial \, C_j^2 \cdot r(i)$ is id $+$ nil (of the easy type). This implies that, at the level of $\Gamma$, the following two operations: deleting $\underset{1}{\overset{\infty}{\sum}} \, r(n)$ OR adding $\underset{1}{\overset{\infty}{\sum}} \, C_j^2$ are, homologically speaking equivalent. This shows that $\zeta_1$ is a homology equivalence.

\smallskip

Next, start with $\underset{n=1}{\overset{j-1}{\sum}} \, b(n) \subset R(j-1) - r(j)$, and then notice that, for all $j$'s, $\partial \, C_j^2 \cdot b(j) = 1$ and $\partial \, C_j^2 \cdot \underset{1}{\overset{j-1}{\sum}} \, b(n) =0$. [This last fact follows from
$$
\partial \, C_j^2 \subset \underbrace{\Gamma - R(j-1)}_{\subset \, \Gamma - \underset{1}{\overset{j-1}{\sum}} \, b(n)} \ + \ r(j).]
$$
It follows that all the inclusion maps below are homology equivalences too
$$
\left( \Gamma - R \cap B - \sum_1^{\ell} b(n) \right) \cup \sum_{\ell + 1}^{\infty} C_j^2 \subset \left( \Gamma - B \cap R - \sum_1^{\ell-1} b(n) \right) \cup \sum_{\ell}^{\infty} C_j^2 \subset \ldots \subset (\Gamma - B \cap R) \cup \sum_1^{\infty} C_j^2.
$$
[{\bf Explanation.} The inclusion
$$
\left( \Gamma - R \cap B - \sum_{n=1}^{\ell} b(n) \right) \cup \sum_{\ell + 1}^{\infty} C_j^2 \subset \left( \Gamma - R \cap B - \sum_{n=1}^{\ell-1} b(n) \right) \cup \sum_{\ell}^{\infty} C_j^2 
$$
consists in: i) leaving $b(\ell)$ in place, undeleted, and ii) adding $C_{\ell}^2$. In the context of the map above, $C_{\ell}^2$ goes once through $b(\ell)$ and it also meets various $b \, (L > \ell)$. So it exactly cancells the increase of $H_1 (\ldots)$ represented by leaving $b(\ell)$ in place, undeleted.]

\medskip

This implies that $\zeta_2$ is a homology equivalence, just like $\zeta_1$.

\bigskip

Since $\zeta_1$ is a homology equivalence and since we also have that $(\Gamma - R \cap B) \cup \underset{1}{\overset{\infty}{\sum}} \, C_j^2$ is connected, the $T$ is connected too. We have to show now that $H_1 \, T = 0$. So, we start with the inclusion $T \overset{i}{\subset} \Gamma$. For any, arbitrary $\ell$, we have a commutative diagram where for the first equality $(H_1 \, T = \ldots)$ we use the fact that both $\zeta_1$ and $\zeta_2$ are homology equivalences.

$$
\xymatrix{
\mbox{\small$H_1 \, T = H_1 \left(\left( \Gamma - R \cap B - \underset{1}{\overset{\ell}{\sum}} \, b(n)\right) \cup \underset{\ell + 1}{\overset{\infty}{\sum}} \, C_j^2 \right)$} \ar[rr]^-{i_*} \ar[dr]_{I_{\ell}}^{\approx} &&\mbox{\small$H_1 \, \Gamma = H_1 (\Gamma \, {\rm rel} \, (\Gamma -B))$} \ar@{->>}[dl]^{\eta_{\ell}} \\
&\mbox{\large${H_1 (\Gamma \, {\rm rel} (\Gamma -B)) \diagup \atop \left\{ \underset{1}{\overset{\ell}\sum} \, b(n) , \, R \, \cap \, B , \, \underset{\ell + 1}{\overset{\infty}{\sum}} C_j^2 \right\}.}$}
} \leqno (108)
$$

In the diagram above, the map $i$ is an inclusion of graphs $T \xhookrightarrow{ \ i \ } \Gamma$, making that, automatically, $i_*$ is injective. Also, quite trivially $\eta_{\ell}$ surjects and $I_{\ell}$ is bijective.

\medskip

\noindent [{\bf Explanations.} Start with the inclusion $\Gamma = R \cap B - \underset{1}{\overset{\ell}{\sum}} \, b(n) \subset \Gamma$ and with the induced homology map $H_1 \left(\Gamma - R \cap B - \underset{1}{\overset{\ell}{\sum}} \, (b(n)) \right) \underset{\alpha_{\ell}}{\xrightarrow{ \quad } } H_1 (\Gamma) = H_1 (\Gamma \, {\rm rel} \, \Gamma - B)$. Here ${\rm Im} \, \alpha_{\ell} = H_1 (\Gamma \, {\rm rel} \, \Gamma - B) / \left\{ R \cap B ,\underset{1}{\overset{\ell}{\sum}} \, b(n)\right\}$. When to $\Gamma - R \cap B - \underset{1}{\overset{\ell}{\sum}} \, b(n)$ one adds $\underset{\ell + 1}{\overset{\infty}{\sum}} \, C_j^2$, then at the level of $H_1$ one gets the further quotient.
$$
{\rm Im} \, \alpha_{\ell} \, \diagup \sum_{\ell + 1}^{\infty} C_j^2 , \ \mbox{a.s.o.}]
$$

Assume now that there exists some $0 \ne x \in H_1 \, T$. From (108) it follows that $I_{\ell} (x) = \eta_{\ell} \, i_* (x) \ne 0$. But then, when we kill
$$
[R \cap B] + \sum_{j=1}^{\ell} \, [b(j)] + \sum_{\ell +1}^{\infty} C_j^2 \subset H_1 (\Gamma \, {\rm rel} \, (\Gamma - B)),
$$
for $\ell \to \infty$, all of $H_1 (\Gamma \, {\rm rel} \, (\Gamma - B))$ gets eventually killed. This means a contradiction. Hence $H_1 \, T = 0$ and our lemma is proved. $\Box$

\bigskip

In all the little theory above, $\Gamma$ was some generic, abstract graph. We move back now to our $\Gamma (1) \subset \Gamma (\infty) = \Gamma (\infty) \times r \subset 2\Gamma (\infty)$ and to the two families $R_0 , B_0 \subset \Gamma(\infty) \times r = \Gamma (X^2 [{\rm new}])$ and the larger $R_1 , B_1 \subset 2\Gamma (\infty)$, each of them four having the property $P$ in their appropriate context. It is understood that the {\ibf subdivisions (like in Figure {\bf 30})}, are incorporated in these definitions. The $\Gamma (1)$ is supposed to incorporate them too. The (93) will be understood in this extended context, with a new meaning for $n,M$. The $\Gamma (1)$ is $\Gamma (1) \times (\xi_0 = -1)$ and it has an excess of BLUE, in the sense that $\Gamma (1) - B_0 \subset \Gamma (\infty) - B_0$ is not connected, in the tree $\Gamma (\infty) - B_0$; it is actually violently disconnected, each vertex corresponds to a  connected component. So, $B_0 \cap \Gamma (1) \subset \Gamma (1)$ does not have property $P$.

\smallskip

I will show now a {\ibf preliminary construction} proceeding for the time being inside the $\Gamma$ (actually our $\Gamma (\infty)$) which is now fixed once and for all.

\smallskip

In our $\Gamma (\infty) - B_0$ we can certainly find a geodetic arc $g_1 \subset \Gamma (\infty) - B_0$ which joins two distinct components of $\Gamma (1) - B_0$ and which only touches $\Gamma (1) - B_0$ via its ends.

\smallskip

We can iterate this process inside $\Gamma (\infty) - B_0$, until we get
$$
\Gamma (2) = \Gamma (1) \cup g_1 \cup g_2 \cup \ldots \cup g_{\chi} \subset \Gamma , \ \sum_1^{\chi} g_i \subset \Gamma - B_0 ,
\leqno (109)
$$
which is now such that $\Gamma (2) - B$ is a {\ibf tree}. In this situation, we will say that $\Gamma (2)$ is {\ibf well-balanced} for BLUE. Also, the formula (109) {\ibf defines} for us the quantity $\chi$.

\smallskip

Now, the $\Gamma (1)$ was well-balanced for RED, but the $\Gamma (2)$ is NOT. We will express the BALANCING features just mentioned, by $\vert \Gamma (1) \cap R_0 \vert = b_1 (\Gamma (1))$ (good balance for RED, the notation here being $\vert \Gamma (1) \cap R_0 \vert = \# \, \{{\rm of} \ R_0 \subset \Gamma (1)\}$). Next, $\vert \Gamma (2) \cap B_0 \vert = b_1 (\Gamma (2))$, but for $\Gamma (2) - R_0 \subset \Gamma (\infty) - R_0$ we get a disconnected $\Gamma (2) - R_0$, coming with $\vert \Gamma (2) \cap R_0 \vert \geq b_1 \, \Gamma (2) = \vert \Gamma (2) \cap B_0 \vert$, hence (potentially at least), a disbalanced for RED.

\medskip

In terms of the various Betti numbers, and with the quantity $\chi$ {\ibf defined like in {\bf (109)}}, (do not mix it up with the Euler characteristic), we have $\chi = b_0 (\Gamma (1) - B_0) - 1$, by construction.

\smallskip

Then, as a direct consequence of the structure of (109)
$$
\chi = b_1 (\Gamma (2)) - b_1 (\Gamma (1)) = \vert \underbrace{\Gamma (2) \cap B_0}_{= \, \Gamma (1) \, \cap \, B_0} \vert - \vert \underbrace{\Gamma (1) \cap R_0}_{\varsubsetneqq \, \Gamma (2) \, \cap \, R_0} \vert .
$$

We will use the notations
$$
R_0 \cap \Gamma (2) \supset R_0 \cap \Gamma (1) = \sum_{i=1}^n H_i^r \, , \quad R_0 - \sum_1^n H_i^r = \{ h_1 , h_2 , \ldots \},
$$
where the order is chosen such that $C \cdot h = {\rm id} + {\rm nilpotent}.$ Here $C \equiv \Sigma \, \{C_i \ {\rm remaining}\} + \Sigma \, \{{\rm little} \ C_i\}$, see (101).

\bigskip

We will present now the real CONSTRUCTION which will realize the BALANCING of RED and BLUE, but this time 1-handles will slide and the ambient infinite graph $\Gamma (\infty)$ (or $\Gamma (2\infty)$) will be changed. To give an idea of what is going on, we start for illustration with the case when in (109) we have $\chi = 1$ and also $\vert \Gamma (2) \cap R_0 \vert > b_1 \, \Gamma (2)$. Our $(\Gamma (1) - B_0) \cup g_1$ is now a tree and, since $b_1 \left( \underbrace{\Gamma (1) \cup g_1}_{\mbox{\footnotesize now, this is $\Gamma (2)$}} - \underset{i=1}{\overset{n}{\sum}} \, H_i^r \right) > 0$, the $\Gamma (1) \cup g_1 - \underset{i=1}{\overset{n}{\sum}} \, H_i^r$ cannot be housed inside $\Gamma (\infty) - R_0$, and we find that $g_1 \cap R_0 = \{ x_1 , x_2 , \ldots , x_p ; y_1 \} \ne \emptyset$. The notation is chosen here such that, in the RED order $y_1 > \{$the $x_i$'s$\}$.

\smallskip

We slide now $y_1$ over the $x$'s, as it is suggested to do in Figure 38. This defines a transformation $\Gamma (\infty) \Rightarrow \Gamma (\infty)$ [balanced], where $\chi = 1$.

\smallskip

We have assumed here that $\vert \Gamma (2) \cap R_0 \vert > b_1 \, \Gamma (2)$ and, when it so happens that $\vert \Gamma (2) \cap R_0 \vert = b_1 \, \Gamma (2)$, then $\{ y;x \} = \emptyset$ and our transformation is mute. Notice the following facts.

\bigskip

\noindent (110.1) \quad Both $\Gamma (\infty)$ [balanced] $- \, R_0$ and $\Gamma (\infty)$ [balanced] $- \, B_0$ are now trees.

\bigskip

\noindent (110.2) \quad In our special case $(\chi = 1)$, the $\Gamma (1) - B_0$ has consisted of two trees; joined via $g_1$. It follows, that when we move from $\Gamma (1) \subset \Gamma (\infty)$ to $\Gamma (1) \cup g_1 ({\rm new}) \subset \Gamma (\infty)$ [balanced], then $\Gamma (1) \cup g_1 ({\rm new})$ is now {\ibf well-balanced}, both for RED and for BLUE.

\bigskip

\noindent (110.3) \quad Our transformation brings the following new trajectory for the RED oriented graph $C \cdot h$:
$$
\xymatrix{
x' \ar@{-}[d] \ar[rr]_{{C \cdot h \atop \mbox{\footnotesize old trajectory}}} &&y_1 \ar[rr]_-{\rm SLIDE} &&x_i \, \mbox{(Fig. 38)} \ar[rr]_-{{C \cdot h \atop \mbox{\footnotesize old trajectory}}} &&x'' \\
{ \ } \ar@{-}[rrrrrr]_{\mbox{\footnotesize new trajectory}} &&&&&&{ \ } \ar[u]
}
$$

But this does not change the RED order of the $h$'s.

$$
\includegraphics[width=12cm]{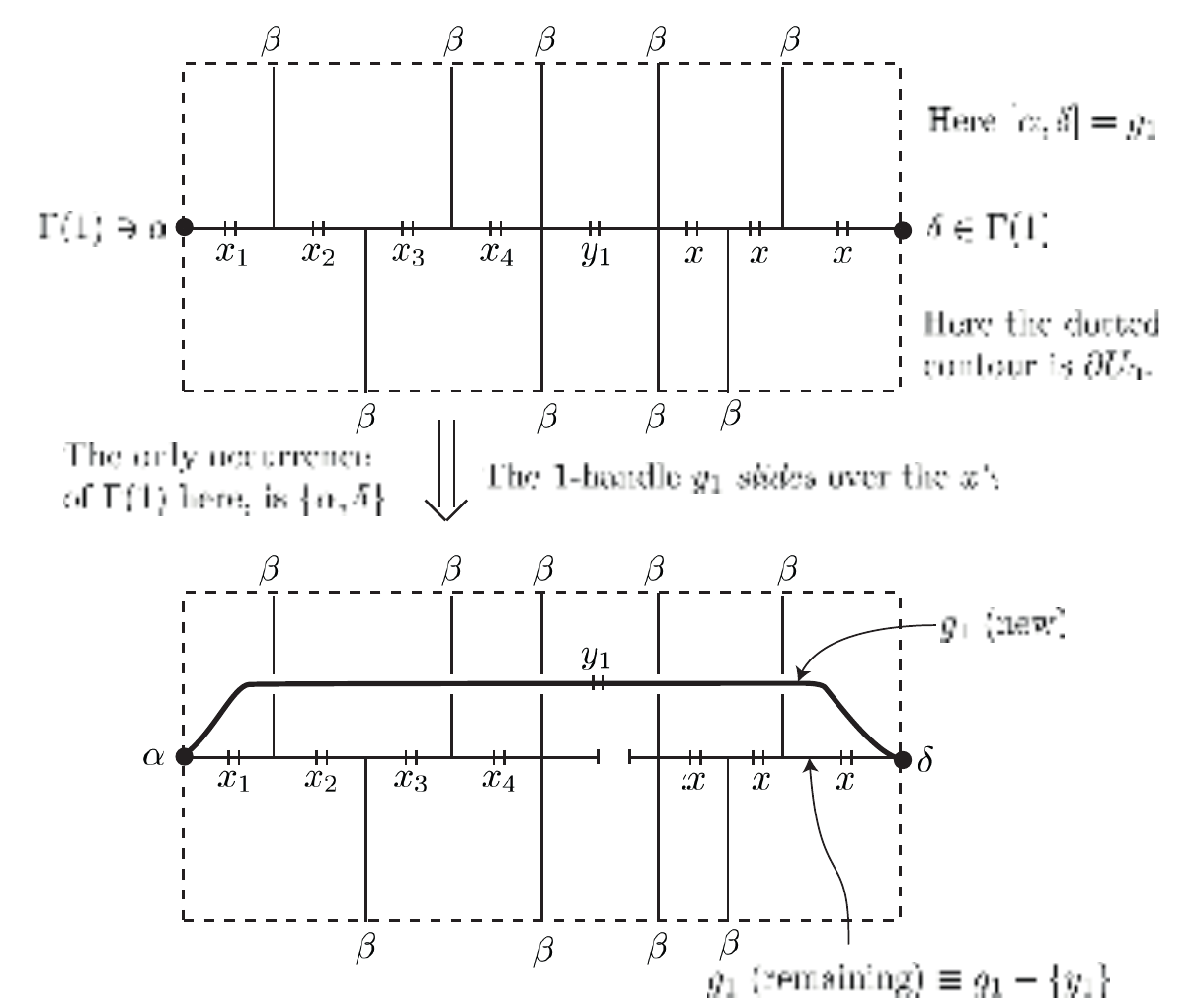}
$$
\label{fig38}
\centerline{\bf Figure 38.}
\begin{quote}
Illustration for the transformation $\Gamma (\infty) \Rightarrow \Gamma (\infty)$ [balanced]. In the upper figure, inside the dotted contour, we can see a neighbourhood $g_1 = [\alpha , \delta] \subset U_0 \subset \Gamma (1)$, which connects with the outer world through the sites $\{ \alpha , \delta , \beta \}$. It is understood here that $\{ x , y_1 \} \subset R_0 - B_0$.
\end{quote}

\bigskip

\noindent {\bf The abstract balancing Lemma 15.}

\smallskip

1) {\it We consider now the general case and introduce the quantity
$$
\chi \equiv \vert \Gamma (2) \cap B_0 \vert - \vert \Gamma (1) \cap R_0 \vert = M-n, \leqno (111)
$$
already mentioned above.

\smallskip

The case $\chi = 1$ has already been dealt with when we discussed Figure $38$. In general, we have a disjoined partition into connected components, with $\chi + 1 \geq 1$
$$
\Gamma (1) - B_0 = \sum_{\alpha = 1}^{\chi + 1} X_{\alpha} \subset \Gamma (\infty) - B_0 .
$$
In $\Gamma (\infty) - B_0$ there is a geodetic arc $g_1$ connecting two connected components of $\Gamma (1) - B_0$ and meeting $\Gamma (1) - B_0$ only through its end-points. I CLAIM, at this point, that we have $g_1 \cap R_0 \ne \emptyset$.

\medskip

Here is the {\ibf proof of this claim}. Assume that $g_1 \cap R_0 = \emptyset$. Then with $(\Gamma (1) - R_0) \cup g_1$ hooked now between two trees
$$
\Gamma (1) - R_0 \subset (\Gamma (1) - R_0) \cup g_1 \subset (\Gamma (\infty) - R_0),
$$
and with the $g_1$ which rests with both its ends on $(\Gamma (1) - R_0)$, we have $b_1 ((\Gamma (1) - R_0) \cup g_1) > 0$, which is a contradiction. The CLAIM is proved.

\medskip

We have $g_1 \cap R_0 = \{ y_1 ; x_{11} , x_{12} , \ldots , x_{1n}\}$, with $y_1 > \{ x_{11} , \ldots , x_{1n}\}$, (the RED order is being meant here), and we can treat now the $g_1$ like in the Figure $38$. This leads to a FIRST TRANSFORMATION, where $y_1$ slides over the $\underset{i}{\sum} \, x_{1i}$, call it
$$
(\Gamma (\infty) , \Gamma (1)) \Longrightarrow (\Gamma (\infty)_1 , \Gamma (1)_1 \equiv \Gamma (1) \cup g_1 (new)). \leqno (112)
$$

Here $\Gamma (1)_1 - B_0$ has one connected component less than $\Gamma (1) -B_0$. So we have 
$$
\Gamma (1)_1 - B_0 = \sum_{\alpha=1}^x \, X_{\alpha} (1) \subset \Gamma (\infty)_1 - B_0 .
$$

There is now a next geodetic arc $g_2 \subset \Gamma (\infty)_1 - B_0$ connecting, just like before two connected components of $\Gamma (1)_1 - B_0$, in a clean manner. Just like for $(112)$ we have now
$$
g_2 \cap R_0 = \{ y_2 ; x_{21} , \ldots , x_{2m} \}.
$$
This means a transformation, analogous to $(112)$
$$
(\Gamma (\infty)_1 , \Gamma (1)_1) \Longrightarrow (\Gamma (\infty)_2 , \Gamma (1)_2 \equiv \Gamma (1)_1 \cup g_2 (\mbox{new})).  \leqno (113)
$$
After $\chi$ steps we get the final composite transformation
$$
(\Gamma (\infty) , \Gamma (1)) \Longrightarrow (\Gamma (\infty)[\mbox{balanced} \,] \equiv \Gamma (\infty)_{\chi} \, , \ \Gamma (3) \equiv \Gamma (1)_{\chi}),  \leqno (114)
$$
which is our BALANCING RED/BLUE CONSTRUCTION.

\smallskip

Our transformation
$$
(\Gamma (\infty) , \Gamma (1)) \Longrightarrow (\Gamma (\infty)[\mbox{balanced} \,] , \ \Gamma (3)),
$$
which changes $\Gamma (1)$ into $\Gamma (3)$, does not change the connections which $\Gamma (1)$ had, already, with the outside world.}

\medskip

2) {\it At the end of the construction above, both $\Gamma (\infty) [\mbox{balanced} \,]$ and $\Gamma (3)$ are well-balanced, both for BLUE and for RED. More explicitly, we find the following items

\bigskip

\noindent $(115)$\vglue-10mm\begin{eqnarray}
\Gamma (3) \cap B_0 &= &\Gamma (1) \cap B_0 = \sum_1^M b(i) \ \mbox{(see $(93)$) and} \nonumber \\
\Gamma (3) \cap R_0 &\equiv &\sum_1^M H_i^r \equiv \Gamma (1) \cap R_0 + \{\mbox{the $y$'s}\}, \ \mbox{so that now} \nonumber \\
b_1 (\Gamma (3)) &= &\vert \Gamma (3) \cap R_0 \vert = \vert \Gamma (3) \cap B_0 \vert = M . \nonumber
\end{eqnarray}
}

3) {\it We consider $B_0 \equiv \{ b_1 , b_2 , \ldots , \ldots , b_M , b_{M+1} , \ldots \}$ with $\{ b_1 , b_2 , \ldots , b_M \}$ like above, written in the increasing BLUE order. I claim it is possible to order $R_0$ starting with $\Gamma (3) \cap R_0$,
$$
R_0 = \{ H_1^r , H_2^r , \ldots , H_M^r ; H_{M+1}^r , \ldots \} ,
$$
so that the following things should happen: one can apply the infinitistic Lemma $14$ to $R_0 \subset \Gamma (\infty)[\mbox{balanced} \,] \supset B_0$ in such a way that $\Gamma (3)$ is invariant for the ABSTRACT COLOUR-CHANGING PROCESS, and that moreover
$$
\Phi (H_i^r) = b_i \ \mbox{for all} \ 1 \leq i \leq M. \leqno (115.1)
$$
}

\smallskip

\noindent {\bf Proof.} The only item not already proved in the statement itself is the 3). So, since $\Gamma (3)$ is already well-balanced for like RED and BLUE, we can start with a first application of Lemma 14 to the finite graph $\Gamma (3)$, so as to achieve just (115.1) for $i \leq M$. This fixes then an order on the $\underset{1}{\overset{M}{\sum}} \, H_i^r$. Next, we turn to our Lemma 14, this time in earnest, using its full infinite glory, and this time for the $\Gamma (\infty)$[balanced], but without touching to $\Gamma (3)$ any longer. Our ABSTRACT BALANCING LEMMA 15 is now proved. $\Box$

\bigskip

\noindent {\bf Important Remark.} In the action of $\Gamma (\infty) \Longrightarrow \Gamma (\infty)$[balanced], Figure 38, the $R_0 \cap B_0 \subset B_0$ gets deleted, from the very beginning. Hence one has not have to worry about $R \cap B$, neither in the context of the abstract Lemma 15 nor when it comes to its geometric implementation $(t=0) \Longrightarrow \left( t = \frac12 \right)$ explained below, leading to lemma 17. $\Box$

\bigskip

Retain, at this point, that by now all the four graphs are {\ibf trees}:
$$
\Gamma (\infty) [{\rm balanced}] - R , \ \Gamma (\infty) [{\rm balanced}] - B , \ \Gamma (3) - R , \ \Gamma (3) - B.
$$

Also our little abstract theory is now finished and the next section will show how to implement it, geometrically, in $4^{\rm d}$.

\smallskip

Let us go back now to the pre-Lemma 15 family
$$
R = \{R_1 , R_2 , \ldots R_n ; h_1 , h_2 ,\ldots , h_p , \ldots \} = \mbox{(with a different notation)} = \{ \underbrace{H_1^r , H_2^r , \ldots , H_n^r}_{{\rm this \ is } \ R_1 , R_2 , \ldots , R_n} ; H_{n+1}^r , H_{n+2}^r , \ldots \} .
$$
In terms of notations from (8) and of the link (99.1), we have $C \cdot h = {\rm id} + {\rm nil}$, of easy type, establishing the duality $\{ C \} \approx \{ h \}$.

\bigskip

\noindent {\bf Fact} (116) \quad In terms of this last duality, these $h$'s which are dual to the $\{$little $C\}$'s are all in $R \cap B$ and hence they can {\ibf never} occur among the $\underset{1}{\overset{M-n}{\sum}} \, y_1^2$ from Lemma 15. The $\{$little $C\}$ is here like in (97).

\smallskip

At this point, I will make a, for the time being {\ibf purely notational PROMOTION} for the 1-handles and the curves of our $\{$link$\}$, at a time before we start implementing geometrically the Lemmas 14, 15 and which I will call the time $t=0$. So here is the

\bigskip

\noindent {\bf Promotion Table:}
$$
\sum_1^n H_i^r \Longrightarrow \sum_1^M H_i^r \equiv \sum_1^n H_i^r + \sum_1^{M-n} y_j ,
$$
with the $y_j$'s promoted, honorifically at least, as 1-handles of $\Delta^4_{\rm Schoenflies}$; next
$$
\sum_1^{\bar n} \Gamma_j \Longrightarrow \sum_{1}^{\overset{=}{n}} \Gamma_j \equiv \sum_1^{\bar n} \Gamma_j + \left\{\mbox{the $C$'s dual to the} \ \sum y_j = \sum_{n+1}^M H_i^r \right\} .
$$
Here the promotion of the $C$'s, as attaching zones of $2$-handles of the $\Delta^4_{\rm Schoenflies}$ is really only {\ibf honorific}. The corresponding 2-handles are no longer directly to $\Gamma (1)$ (with $\Gamma (1) = \Gamma (1) \times (\xi_0 = -1)$), nor even to $\Gamma (3)$. The promotion will stay purely honorific, even when Lemma 15 will have been implemented geometrically. The next promotion are:

\smallskip

The mute promotion $\underset{1}{\overset{M}{\sum}} \, b_i \underset{\approx}{\Longrightarrow} \underset{1}{\overset{M}{\sum}} \, b_i$, where remember that our $b_i$'s are the $(\Gamma (1) \times (\xi_0 = -1)) \cap B$. And, after this trivially mute promotion, we have
\begin{eqnarray}
\sum_1^{\infty} h_i &\Longrightarrow &\sum_1^{\infty} h_i - \sum_{n+1}^{M} H_j^r , \nonumber \\
\sum_1^{\infty} C_i &\Longrightarrow &\sum_1^{\infty} C_i - \sum_{n+1}^{M} \Gamma_j . \nonumber
\end{eqnarray}

\noindent (Now $\bar n \geq n , \overset{=}{n} = \bar n + (M-n)$.) $\Box$

\bigskip

By definition, after this promotion we are at time $t=0$, at the level of $2X_0^2$ and we continue to have (at $t=0$)
$$
C \cdot h = \mbox{easy id $+$ nilpotent.}
$$

In the process leading to our final desired resul in the next section VII, there will be other successive times after $(t=0)$, let us say they are
$$
t=0 , \ t = \varepsilon , \ t=2 \, \varepsilon , \ldots , t = \frac12 , t + \frac12 + \varepsilon , \ldots , t=1 .
$$

With the new smaller families $\underset{i}{\sum} \, C_i$, $\underset{k}{\sum} \, h_k$, the new smaller
$$
{\rm LAVA} \, (t=0) \equiv \sum_1^{\infty} h_j \cup D^2 (C_j) \quad \mbox{(AFTER PROMOTION),} \leqno (117)
$$
continues to have the STRONG PRODUCT PROPERTY and, outside the sites living in the $\underset{1}{\overset{\overset{=}{n}}{\sum}} \, {\rm int} \, D^2 (\Gamma_j)$, the extended cocores continue to make sense.

\smallskip

We have $\underset{1}{\overset{\overset{=}{n}}{\sum}} \, \Gamma_j \cap {\rm LAVA} (t=0) \ne \emptyset$, since the   promoted $\Gamma_j$'s can happily touch to LAVA ($t=0$). They are not directly attached to $\Gamma (1)$ not to $\Gamma (3)$ for that matter. But $\left( \underset{1}{\overset{\bar n}{\sum}} \, \Gamma_j \right) \cap {\rm LAVA} (t=0) = \emptyset$,  while generally speaking, we also have $\underset{1}{\overset{M}{\sum}} \, H_i^r \cap {\rm LAVA} (t=0) \ne \emptyset$ and, using the LAVA ($t=0$) one can build the extended cocore of the $H_i^r$'s ($i \leq M$). When it so happens that $H_i^r \cap {\rm LAVA} (t=0) = \emptyset$, then we have extended cocores of $H_i^r = H_i^r$ itself, or if one wants, the cocore $H_i^r$ itself. This is a trivial case, on which we will not dwell.

\smallskip

Coming back to our successive times $t$ from above, they parametrize the times of successive geometric constructions happening in $4^{\rm d}$, inside $N_1^4 (2X_0^2)^{\wedge}$, and they will be carefully described in the next section.

\smallskip

There will be a time $t=\frac12$, when the $R/B$-BALANCING will have been geometrically realized and then a final time $t=1$ when the COLOURS WILL HAVE CHANGED sufficiently so as to make the following GRAND BLUE DIAGONALIZATION possible, with $1 \leq i,j \leq M$
$$
\boxed{\eta_j ({\rm green}) \cdot \{\mbox{extended cocore} \ H_i^b\}^{\wedge} = \delta_{ij}} \ .
$$

The point is that we will change the diffeomorphic model of $\Delta^4_{\rm Schoenflies}$, more precisely of $N^4 (\Delta^2)$ so that the $\{\mbox{extended cocore} \ H_{i \leq M}^b\}^{\wedge}$ will be its 1-handles, with the $H_j^r$'s completely out of the picture.

\smallskip

But then, besides $t = \frac12$ and $t=1$ there will be other intermediary times too. At any time, things like $\underset{1}{\overset{M}{\sum}} \, H_i^r \cap {\rm LAVA} (t\ne 0) \ne \emptyset$ will always be with us; they are unavoidable. With them, a RED diagonalization would conflict with the sacro-sancted CONFINEMENT. So, we will {\ibf never} try to diagonalize things like
$$
\sum_1^M \eta_i ({\rm green}) \cdot \sum_1^M H_j^r \ \mbox{nor, of course} \ \sum_1^M \Gamma_i \cdot \sum_1^M H_j^r .
$$

Finally, both the BLUE and the RED flow have to be used for our constructions. Both have to be controlled. Also, even after we will make sense geometrically of $N^4 (\Gamma (3))$ the $\underset{1}{\overset{M}{\sum}} \, D^2 (\Gamma_i)$ are still NOT directly attached to it and to make sense of the promoted things as handles of index one and two of $\Delta^4$ we need to move to {\ibf an infinite set-up}. [Possibly a finite high truncation would  do, but the infinite set up is, finally less messy and more manageable.]

\smallskip

Let us define, as a first $4^{\rm d}$ step
$$
{\mathcal Z}^4 (t=0) \equiv N^4 (2\Gamma_1 (\infty)) - \sum_1^{\infty} h_i \supset \sum_1^M H_i^r . \leqno (118)
$$
This can be easily compactified into
$$
\hat{\mathcal Z}^4 (t=0) \equiv {\mathcal Z}^4 (t=0) \cup \varepsilon ({\mathcal Z}^4 (t=0)) = \overline{{\mathcal Z}_0^4 (t=0)} \subset N^4_1 (2X_0^2)^{\wedge}. \leqno (119)
$$

\medskip

\noindent {\bf Easy fact} (119.1) \quad The pair $\left(\hat{Z}^4 (t=0) , \underset{1}{\overset{M}{\sum}} \, H_i^r \right)$ is {\ibf standard}, i.e.
$$
\left(\hat{\mathcal Z}^4 (t=0) , \sum_1^M H_i^r \right) \underset{\rm DIFF}{=} \left( M \, \# \, (S^1 \times B^3) , \sum_1^M (*) \times B^3 \right) .
$$

\bigskip

\noindent {\bf Lemma 16. (The time $t=0$ compactification)}

\medskip

1) {\it (Reminder) Modulo an appropriate, not necessarily ambient isotopy, inside our ambient space $N_1^4 (2X_0^2)^{\wedge}$, for each $\alpha$ (see $(38)$), we have
$$
\mbox{LAVA}_{\alpha} (t=0)^{\wedge} \equiv \{ \delta \, \mbox{LAVA}_{\alpha} (t=0) \times [0,1]\} \underset{\overbrace{\mbox{\footnotesize $\delta \, \mbox{LAVA}_{\alpha} (t=0) \times \{0\}$}}}{\cup} \{ (\delta \, \mbox{LAVA}_{\alpha} (t=0) \times \{0\}) \cup \lambda_{\alpha}\} = \leqno (120)
$$
$$
= \, \overline{\mbox{LAVA}_{\alpha} (t=0)} \subset N_1^4 (2X_0^2),
$$
where this last ``$ \, =$'' is an equality of subsets inside $N_1^4 (2X_0^2)^{\wedge}$.

\smallskip

More globally, we have
$$
\overline{\mbox{LAVA} (t=0)} = \overline{\sum_{\alpha} \mbox{LAVA}_{\alpha} (t=0)} = \left( \bigcup_{\alpha} \mbox{LAVA}_{\alpha} (t=0)^{\wedge} \right) \cup \lambda_{\infty} \quad \mbox{(see $(46)$)}.
$$
}

\medskip

2) {\it The following object is a smooth compact $4^{\rm d}$ handlebody of genus $M$
$$
\hat Z^4 (t=0) \equiv {\mathcal Z}^4 (t=0) \cup \bigcup_{\alpha} \mbox{LAVA}_{\alpha} (t=0) \cup \lambda_{\infty} ,
$$
and we have
$$
\left(\hat Z^4 (t=0)) , \sum_1^M \{\mbox{extended cocore} \ H_i^r \}^{\wedge} \right) \underset{\rm DIFF}{=} \left( M \, \# \, (S^1 \times B^3) , \sum_1^M (*) \times B^3 \right). \leqno (121)
$$
}

\medskip

3) {\it There is a diffeomorphism for $\Delta^4 = \Delta^4_{\rm Schoenflies}$, and of course, like usual, by $\Delta^4$ we actually mean $N^4 (\Delta^2)$,
$$
\Delta^4 \underset{\rm DIFF}{=}  \hat Z^4 (t=0) + \sum_1^{\overset{=}{n}} D^2 (\Gamma_i) \subset N_1^4 (2X_0^2)^{\wedge}  \underset{\rm DIFF}{=}  \Delta^4 \cup \partial \Delta^4 \times [0,1] , \leqno (122)
$$
and this is how $\Delta^4$ and $(1)$ will be conceived from now on, at least as long as we stay at time $t=0$. Via a not necessarily ambient isotopy in the ambient space $N_1^4 (2X_0^2)^{\wedge}$, can be assumed that}
$$
\hat Z^4 (t=0) + \sum_1^{\overset{=}{n}} D^2 (\Gamma_i) = \overline{N^4 (2X_0^2)} = N^4 (2X_0^2)^{\wedge}. 
$$

\bigskip

\noindent {\bf Proof.} Using the product property of LAVA one gets a collapse
$$
\hat Z^4 (t=0) \xrightarrow{ \ \ \pi \ \ } \hat{\mathcal Z}^4 (t=0), \leqno (122.1)
$$
as follows. Start by collapsing away the
$$
\sum_1^M \{\mbox{extended cocore} \ H_i^r \}^{\wedge} \cap \left( \sum_{\alpha} \delta \, {\rm LAVA}_{\alpha} \times [0,1] \right) ,
$$
and next collapse the $\underset{\alpha}{\sum} \, \delta \, {\rm LAVA}_{\alpha} \times [0,1]$ too. This leads to a first diffeomorphism
$$
\left( \hat Z^4 (t=0) , \sum_1^M \{\mbox{extended cocore} \ H_i^r \}^{\wedge} \right) = \left( \hat{\mathcal Z}^4 (t=0) , \sum_1^M H_i^r \right) ,
$$
which allows us to deduce our (121) from (119.1).

\smallskip

Next, we exhibit a big collapse
$$
\hat Z^4 (t=0) + \sum_1^{\overset{=}{n}} D^2 (\Gamma_i) \xrightarrow{ \ \ \{{\rm BIG} \, \pi\} \ \ } N^4 (\Gamma (1)) + \sum_1^{\bar n} D^2 (\Gamma_i),
$$
which is all we need for (122). Here is now how to produce the Big $\pi$. 

\smallskip

After our PROMOTION, we have for our $\Delta^4_{\rm Schoenflies}$ the 1-handles $H_i^r$, with $i=1,2,\ldots ,n,n+1, \ldots ,M$ and attaching curves $\Gamma_j$, with $j=1,2,\ldots , \bar n , \bar n + 1 , \ldots , \overset{=}{n}$, when $\overset{=}{n} - \bar n = M-n$. All this is purely honorific, abstract story, of course. We may assume here, concerning the ordering of the $H_i^r$'s, that, in the RED order, if $n < j < i \leq M$, then $H_i^r > H_j^r$, so that
$$
\Gamma_{\bar n + i} \cdot \{\mbox{extended cocore} \ H_j^r\} = 0 \quad {\rm and} \ \Gamma_{\bar n + i} \cdot \{\mbox{extended cocore} \ H_i^r\} = 1. \leqno (123)
$$
Remember, also, that $M-n = \overset{=}{n} - \bar n$.

\smallskip

The increasing RED order means here that in the oriented graph gotten from the off-diagonal part of the RED geometric intersection matrix, the arrows go like this: $H_M^r \longrightarrow H_{M-1}^r \longrightarrow \ldots \longrightarrow H_{n+1}^r$.

\smallskip

To get our $\{{\rm BIG} \, \pi\}$ we proceed via the following steps.

\medskip

a) The (123) allows us to collapse away, first LAVA $(t=0) \cap \{$extended cocore $H_M^r\}$, next $H_M^r \cup D^2 (\Gamma_M)$. 

\medskip

b) Next, using again (123) we can proceed similarly, in order, for $H_{M-1}^r , H_{M-2}^r , \ldots , H_{n+1}^r$. All these collapses are compatible with the $\pi$ (122.1) and, since we know already that $\underset{1}{\overset{\bar n}{\sum}} \, \Gamma_i \cap {\rm LAVA} = \emptyset$, all the rest of LAVA $(t=0)$ can be collapsed away too, leaving us with a naked, i.e. lava-free, last collapse
$$
\left( {\mathcal Z}^4 (t=0) - \sum_{n+1}^M H_i^r \right) \cup \sum_1^{\bar n} D^2 (\Gamma_i) \searrow N^4 (\Gamma (1)) \cup \sum_1^{\bar n} D^2 (\Gamma_i).
$$
This ends the Proof of Lemma 16. $\Box$

\bigskip

For the $\underset{1}{\overset{M}{\sum}} \, b_i = B \cap \Gamma (1)$, with the indexing from (115.1), we will also use from now on the notation $\underset{1}{\overset{M}{\sum}} \, H_i^b$. And here, just like for $\underset{1}{\overset{M}{\sum}} \, H_i^r$, we have the
$$
\sum_1^M \{{\rm extended} \ H_i^b \}^{\wedge} \subset \hat Z^4 (t=0), \leqno (*)
$$
defined via the RED flow $C \cdot h$.

\smallskip

Eventually it will be $(*)$ which will be the system of 1-handles of our $\Delta^4$ (122).

\smallskip

We also have now our system of embedded exterior discs from Theorem 13, i.e. the
$$
\sum_1^M (\delta_i^2 , \eta_i ({\rm green}) \equiv \partial \delta_i^2) \xrightarrow{ \ {\rm embedding} \, J \ } \bigl(\partial N^4 (2X_0^2)^{\wedge} \times [0,1] \subset \leqno (124)
$$
$$
N_1^4 (2X_0^2)^{\wedge} \ (\mbox{our ambient space}) , \partial N^4 (2X_0^2)^{\wedge}\bigl). 
$$

\smallskip

In lieu of the simple-minded $\eta_i ({\rm green}) \cdot b_j$, which was diagonal, we consider now the more complex geometric intersection matrix

\bigskip

\noindent (125) \quad $\eta_i ({\rm green}) \cdot \{$extended cocore $H_j^b \}^{\wedge} = \delta_{ij} + \{$a {\ibf parasitical} off-diagonal term $\eta_i ({\rm green}) \cdot h_k$, where $h_k \in R_0 - B_0$ and, also $h_k \subset \{$extended cocore $H_j^b \}\}$. A finite set of $k$'s is involved here.

\bigskip

Here $1 \leq i,j \leq M$. Moreover, it is the little blue diagonalization from (101.1), in its full glory, which makes that, in the context of (125) we have that the parasitical $h_k$'s are $h_k \in R_0 - B_0 = R_1 - B_1$.

\smallskip

The presence of these parasitical $h_k \in \eta_i ({\rm green})$, $i \leq M$, mean that, at this stage in the game, the exterior discs
$$
\sum_1^M (\delta_i^2 , \partial \delta_i^2 = \eta_i ({\rm green}))
$$
are certainly NOT in cancelling position with the 1-handle $\underset{1}{\overset{M}{\sum}} \, \{$extended cocore $H_i^b\}^{\wedge}$.

\smallskip

The aim of the next section is to exhibit a handlebody decomposition for our $N^4 (\Delta^2)$ where the  $\underset{1}{\overset{M}{\sum}} \, \{$extended cocore $H_i^b\}^{\wedge}$ {\ibf are} the 1-handles.

\smallskip

And then, via the COLOUR-CHANGE, we will put the $\underset{1}{\overset{M}{\sum}} \, \delta_i^2$ in cancelling position with them. This will allow us, eventually, to appeal to Lemma 3 and show that $\Delta^4_{\rm Schoenflies} \in $ GSC.

\smallskip

The aim stated above will be achieved by realizing geometrically in $4^{\rm d}$ the ABSTRACT LEMMAS 14 and 15.



\section{Balancing and change of colour, done now geometrically {\small(The $4^{\rm d}$ geometry of the passage $(t=0) \Rightarrow (t=1/2) \Rightarrow (t=1)$ and the compactifications at times $t=\frac12$ and $t=1$.)}}\label{sec}

We consider now Lemma 15 and its proof and we will concentrate on the typical step (112) which, keeping in mind the notations from the Figure 38, we write now as
$$
(\Gamma (2\infty) , \Gamma (1)) \Longrightarrow (\Gamma (2\infty)_1 [\mbox{balanced}] , \Gamma (1) \cup g_1 ({\rm new})). \leqno (126)
$$
This is the first of the $\chi$ steps to the $\Gamma (2\infty) [\mbox{balanced}]$ in (114) and here $g_1 \cap R_0 = \{y ; x_1 , x_2 , \ldots , x_n \}$, with

\bigskip

\noindent (126-bis) \quad $y_1 = H_{n+1}^r > \{x_{11} , \ldots , x_{1m} \}$ (in the RED order of $2X_0^2$); the $y_1$ is PROMOTED and, after this promotion $\{x_{11} , \ldots , x_{1m} \} \subset \overset{\infty}{\underset{1}{\sum}} \, h_i$.

\bigskip

We want to realize now this step (126), geometrically inside $N_1^4 (2X_0^4)^{\wedge}$, as
$$
(N^4 (2\Gamma (\infty)) , N^4 (\Gamma (1)) \Longrightarrow (N^4 (2\Gamma (\infty))_1 [\mbox{balanced}], N^4 (\Gamma (1) \cup g_1 ({\rm new})). \leqno (126\mbox{-ter})
$$

So now, in real life, and not just abstractly as in the last section, an {\ibf active} 1-handle $y (\mbox{non-LAVA}) = (B_a^3 \times [-\varepsilon , \varepsilon] , \partial B_a^3 \times [-\varepsilon , \varepsilon]) \subset (N^4 (2\Gamma (\infty)) , \partial N^4 (2\Gamma (\infty)))$ is sliding over a passive 1-handle $x ({\rm LAVA}) = (B_p^3 \times [-\varepsilon , \varepsilon] , \partial B_p^3 \times [-\varepsilon , \varepsilon]) \subset (N^4 , \partial N^4)$. Really very schematically, what we talk about here is suggested in Figure~39.

$$
\includegraphics[width=16cm]{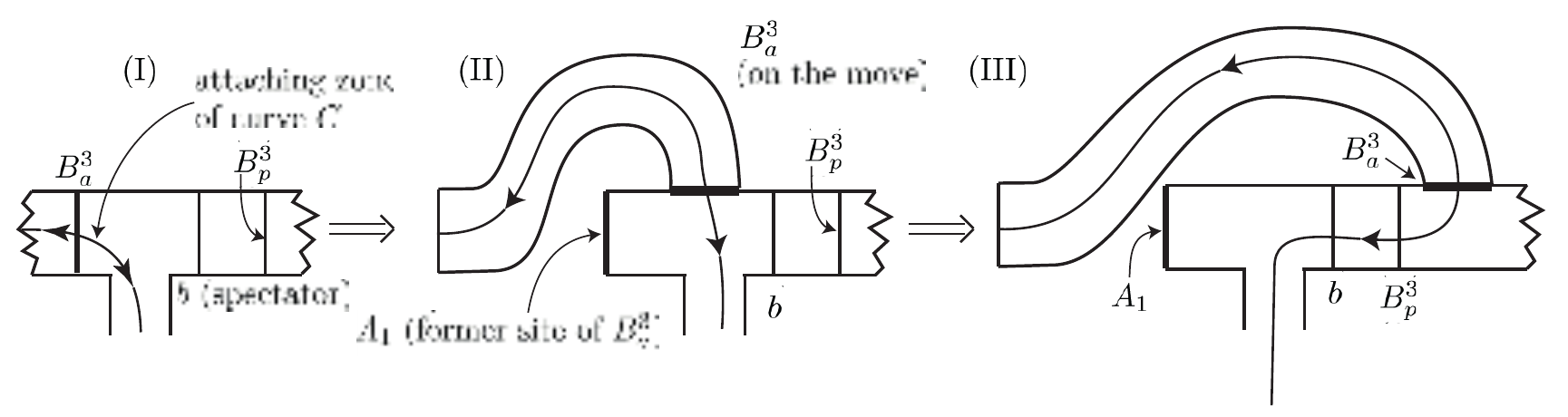}
$$
\label{fig39}
\centerline{\bf Figure 39.}
\begin{quote}
A highly schematical $2^{\rm d}$ representation of the sliding of an {\ibf active} handle, with $B_a^3$ as cocore, over another RED handle, the $B_p^3$ ({\ibf passive}). Besides the RED 1-handles $B_a^3 , B_p^3$, there is here also a {\ibf spectator} BLUE 1-handle $b$. The scenario for our step is I $\Rightarrow$ II $\Rightarrow$ III. 

LEGEND: $\leftrightarrow \, =$ this suggests the attaching zone of a 2-handle, resting on $B_a^3$, which gets dragged via covering isotopy. It is painted green.
\end{quote}

\bigskip

But a much more detailed description than the one provided by Figure 39 is actually necessary.

\bigskip

\noindent (126.1) \quad We want to preserve the STRONG PRODUCT property of LAVA, and, during our step (125) the $C \cdot h =$ (easy) id $+$ nil might get violated. In order to take care of this we will need to add to our LAVA both {\ibf LAVA bridges} and {\ibf LAVA dilatations}. 

\bigskip
\bigskip

\noindent [{\bf A philosophical comment.} In this paper, there are some sacro-sancted things never to be violated, while others may be.

\smallskip

Here are our two sacro-santed principles which can never be violated nor trespassed:
\begin{enumerate}
\item[I)] The CONFINEMENT condition, inside $\partial N_+^4 (2\Gamma (\infty))$.
\item[II)] The STRONG PRODUCT PROPERTY of LAVA.
\end{enumerate}

Now when LAVA was first introduced, II) was assured by the condition $C \cdot h = $ id $+$ nilpotent. But then, when $1$-handles will start sliding over each other, dragging the 2-handles along, then

\medskip

III) The property $C \cdot h =$ id $+$ nilpotent {\ibf might get violated}.

\bigskip

The II) will be maintained by internal LAVA operations. Also, the II) is enough for defining the $\{$extended cocore $x\}^{\wedge}$, which will be with us, even with the violation from III.]

\bigskip

The violation of $C \cdot h =$ id $+$ nil will actually come later. Since right now $y > x$ (RED order) there is yet not danger.

\bigskip

Remembering that LAVA is growing out of
$$
\delta \, {\rm LAVA} \subset \partial \left(N^4 (2\Gamma (\infty)) - \sum_1^{\infty} h_i\right) \equiv \partial M^4 ,
$$
and for building up a LAVA bridge, we start by adding, inside $\partial M^4$ a $3^{\rm d}$ 1-handle $\delta {\mathcal H} \subset \partial M^4$ to $\delta \, {\rm LAVA}$, and then we let it grow into a $\delta {\mathcal H} \times [0,1]$, glued in an obvious way to the rest of LAVA. The LAVA dilatations are defined similarly.

\bigskip

\noindent (126.2) \quad  We also want to respect the basic SPLITTING $\partial N^4 (2\Gamma (\infty)) = \partial N_-^4 (\Gamma_1 (\infty) ) \underset{\overbrace{\mbox{\footnotesize $\Sigma_{\infty}^2$}}}{\cup} \partial N_+^4 (2\Gamma (\infty))$ (even with the change of $2\Gamma (\infty)$ in (126)). And, with this, the forced CONFINEMENT conditions (101.1-bis) should be respected too.

\bigskip

At all times in the process (I) $\Rightarrow$ (II) $\Rightarrow$ (III), Figure 39, $B^3_{(a \, {\rm or} \, p)}$ should be SPLIT as
$$
B^3 = \frac12 B^3 (+) \underset{\overbrace{\mbox{\footnotesize $\sigma^2 = \Sigma_{\infty}^2 \cap B^3 \, (= 2$-cell)}}}{\cup} \frac12 B^3 (-).
$$

$$
\includegraphics[width=180mm]{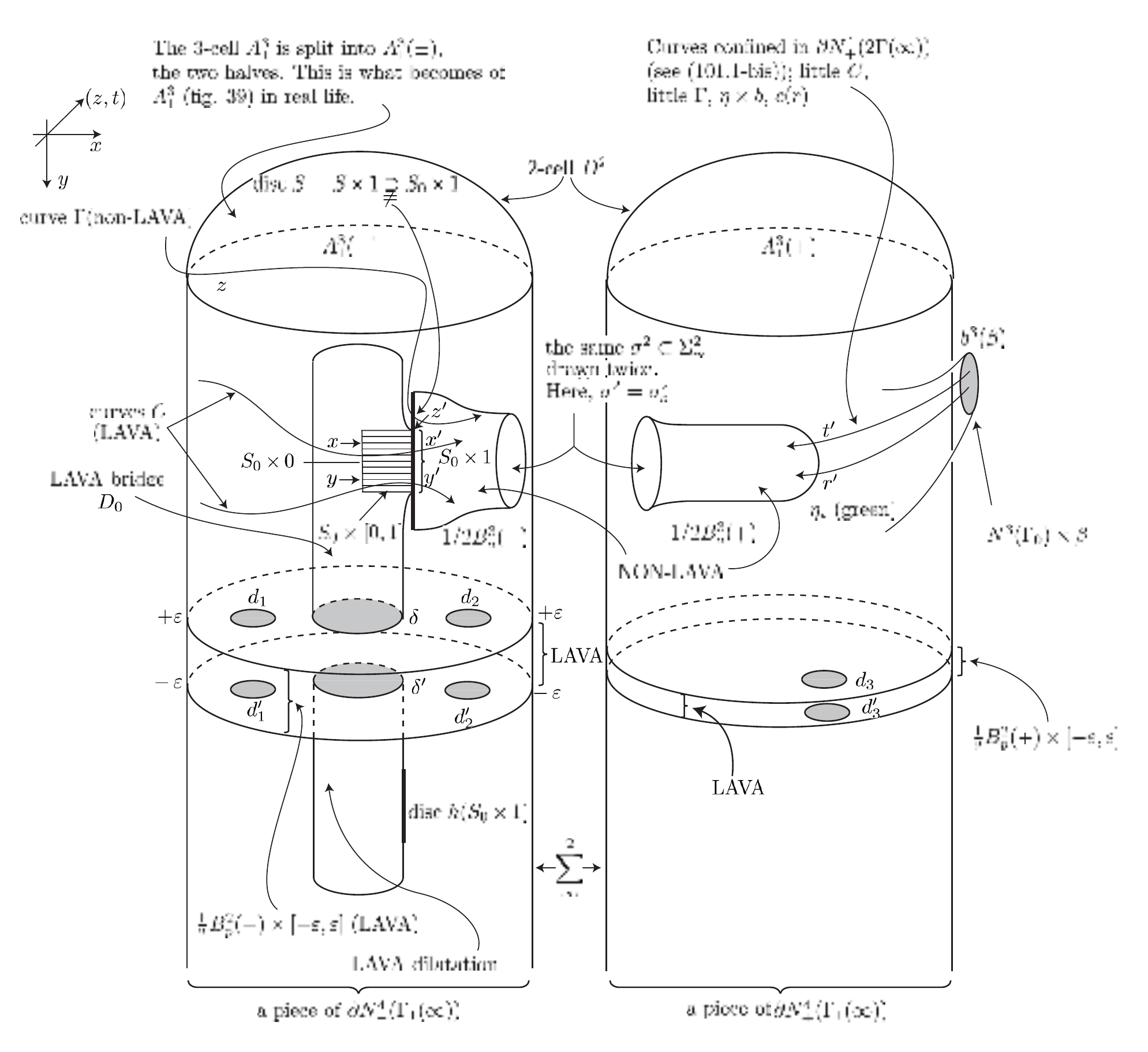}
$$
\label{fig40}
\centerline{\bf Figure 40.}
\begin{quote}
A realistic view of Figure 39-(II). Separated into two distinct $\pm$ halves, this figure suggests two $3^{\rm d}$ pieces out of which a small part of the bigger $\partial N^4 (2\Gamma (\infty))$ can be reconstructed. For typographical reasons, several geometric ingredients occur twice in this figure. This is, for instance, the case for $D^2 , \sigma^2 , \Sigma_{\infty}^2$.

We also have here $S_0 \times [0,1] \subset \delta \, {\rm LAVA}$, with $S_0$ a disc. The LAVA just mentioned is actually LAVA {\ibf dilatation}.

We can reconstruct here also $\partial B_p^3 = \frac12 \, B_p^2 (-) \cup \frac12 \, B_p^2 (+)$ and the whole of $B_p^3 \times [-\varepsilon , \varepsilon]$ is LAVA. Of course, ${\rm int} \, B_p^3$, living inside $N^4 (2\Gamma (\infty))$ is not visible here. But more LAVA is present here than what is drawn: through the pairs of shaded discs $d_1 + d'_1$, $d_2 + d'_2$, $d_3 + d'_3$ go attaching zones $C_i \times D_i^*$ of 2-handles belonging to LAVA. These things are formally 1-dimensional, HENCE THEY ARE NOT IN THE WAY for the $B_a^3$ moving from its initial position where it is glued to LAVA (with bridges and dilatation included, along ${\rm LAVA} \supset S_0 \times 1 \subsetneqq S \times 1 \subset \partial B_p^3$.) The isotopy
$$
\xymatrix{
S_0 \times 1 \ar@{^{(}->}[d] \ar@{=>}[r] &h (S_0 \times 1) \ar@{^{(}->}[d] \\
S \times 1 \ar@{=>}[r] &h(S \times 1)
}
$$
happens along an arc we call $\lambda$, which is essentially our {\ibf long cigar}
$$
(\mbox{LAVA bridge}) + \mbox{lava dilatation} \subset \partial N_-^4 (2\Gamma (\infty)).
$$
We have $\lambda \cap (C_i \times D_i^*) = \emptyset$.
\end{quote}

\bigskip

\noindent {\bf Further comments and explanations concerning the Figure 40.}

\smallskip

At the points $x' , y'$, and $t' , v'$, the various curves (attaching zones of 2-handles) climb on the moving 1-handle, as the arrows indicate.

\smallskip

As an immediate consequence of the inequality in (126-bis) we have that
$$
\{\mbox{extended cocore} \ B_a^3\} \cap B_p^3 = \emptyset . \leqno (126.3)
$$
Here is the $\delta \, {\rm LAVA}$ which is explicitly present in Figure 40:
$$
\mbox{The pieces} \ \{D_0 = \delta \, \mbox{LAVA bridge}\} + \{\delta \, \mbox{LAVA dilatation}\} + \frac12 B_p^3 (\pm) \times \{ -\varepsilon , +\varepsilon\}.
$$
These three pieces are glued together via $\delta , \delta'$. Moreover, we have
$$
S_p^2 \times (-\varepsilon , +\varepsilon) \subset \partial \, \mbox{LAVA } - \delta \, \mbox{LAVA} .
$$

We have, also
$$
S_p^2 = \partial B_p^3 = \partial \left( \frac12 B_p^3 (+) \underset{\sigma^2}{\cup} \frac12 B_p^3 (-)\right) = \frac12 B_p^2(+) \cup \frac12 B_p^2 (-).
$$
The $B_p^3 \times [-\varepsilon , \varepsilon]$ is hidden from our sight, in the fourth dimension, and $B_p^3 \times \{-\varepsilon , +\varepsilon\} \subset \delta \, {\rm LAVA}$. BUT, things like $d_i ({\rm shaded}) \times [-\varepsilon , +\varepsilon]$, groing from $d_i$ to $d'_i$ are in $\partial \, {\rm LAVA} - \delta \, {\rm LAVA}$. Here, some $D^2 (C)({\rm LAVA})$ gets attached to $N^4 (2\Gamma (\infty))$. Same for $\delta \times [-\varepsilon , +\varepsilon] \subset \partial \, {\rm LAVA} - \delta \, {\rm LAVA}$.

\smallskip

Inside the $\frac12 B_p^2 (-)\times [-\varepsilon , \varepsilon]$ our cigar melts into LAVA and disappears from the visual field.

\smallskip

The $A_1^3 = A_1^3 (-) \underset{\overbrace{\mbox{\footnotesize $D^2$}}}{\cup} A_1^3 (+)$ from Figure 40 corresponds to the site $A_1$ visible in the Figures 39-(II, III). It is the former position of $B_a^3$, before the sliding starts (see Figure 39-(I) too).

\bigskip

\noindent (126.4) \quad Very importantly, there is no $B$, certainly no $B \cap R$, present in the context of the step from Figure~40. This fact stems from the following item: in the context of the Abstract Lemma 15, when it comes to the arcs $\underset{1}{\overset{\chi}{\sum}} \, g_i$ which houses the $\{ y ; x_1 , x_2 , \ldots , x_n \}$'s, we automatically have
$$
\left( \sum_1^{\chi} g_i \right) \cap B = \emptyset.
$$
End of (126.4).

$$
\includegraphics[width=155mm]{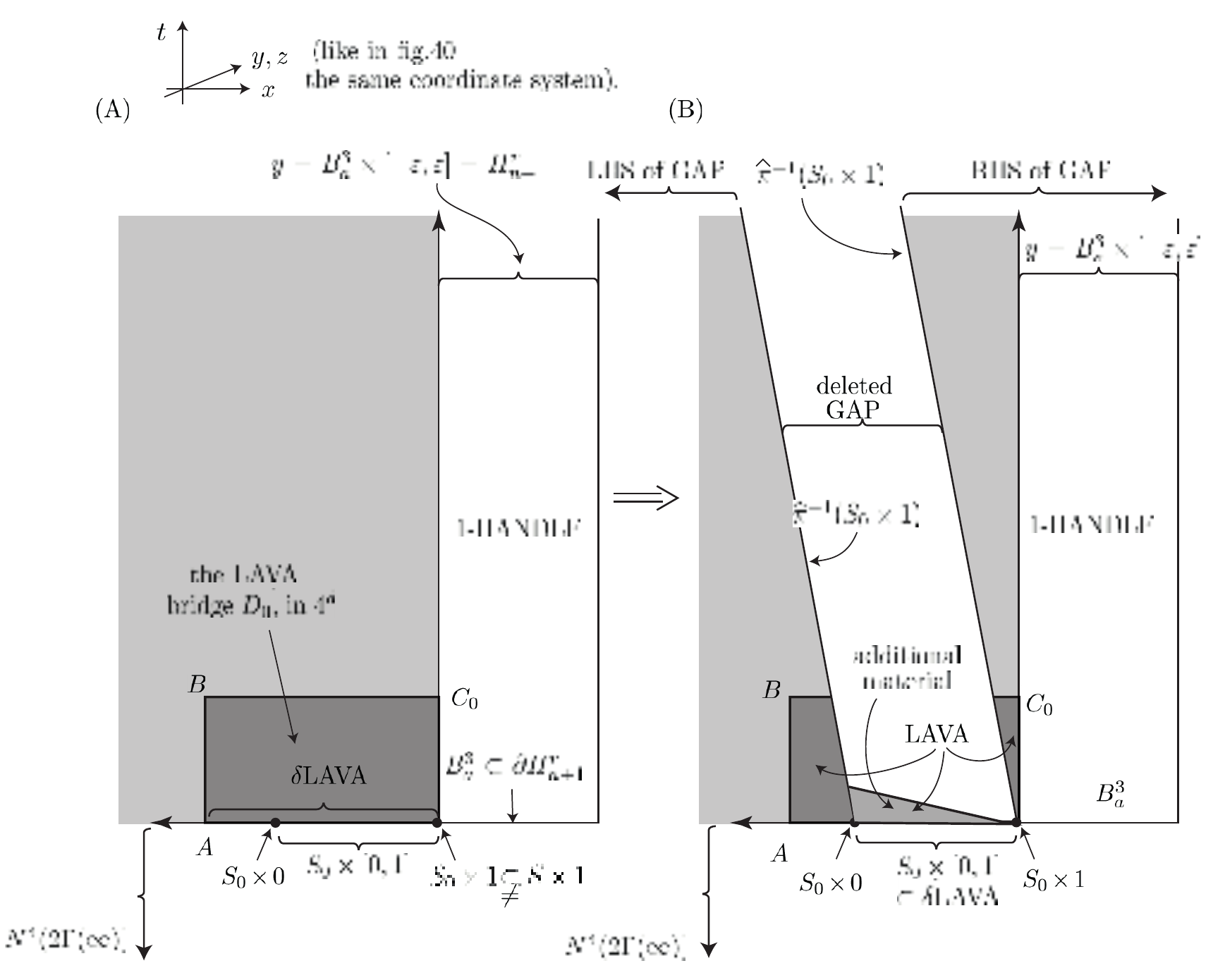}
$$
\label{fig41}
\centerline{\bf Figure 41.}
\begin{quote}
Figure (A) is drawn at the level of Figure 40 which here corresponds just to the line on the bottom ($=$ the bottom floor). The $S_0 \times [0,1]$ here is the one in Figure 40 and the doubly shaded area corresponds to the LAVA bridge $D_0$. We are supposed to be here in the plane $(x,t)$ of one of the $D^2 (C)$, like the $C$ going through $[x,x']$ in Figure 40. So, everything shaded (simply or doubly, is here LAVA. The $[A , S_0 \times 0 , S_0 \times 1]$ corresponds to the $\delta$ LAVA(bridge), and it is located at the bottom floor. The (B) will be explained later, in the main text. And the $S \times 1$ spreads in the $(y,z)$ direction, surrounding $S_0 \times 1$.

LEGEND: $\includegraphics[width=1cm]{rectanglefonce.pdf} = $ LAVA dilatation, in $4^{\rm d}$, $\includegraphics[width=1cm]{rectangleclair.pdf} = D^2(C)$, $\leftrightarrow \, = C$ (attaching curve of the 2-handle $D^2(C)$. The $[ABC_0]$-line is part of it. It CLIMBS OVER $D_0$; in Figure 40 we see only the projection of this on the {\ibf floor} $\partial N^4 (2\Gamma (\infty))$ of the present figure.
\end{quote}

\bigskip

Finally, in our Figure 40, enough things have been added to LAVA so that, in its move $S \times 1 \subset \partial B_p^3$ should stay glued to it, via the $S_0 \times 1 \subset S \times 1$. The reason why $\Gamma_j$ is treated differently than the $C_i$'s, in Figure 40 is that, for it, the kind of extended cocore producing the GAP in Figure 41-(B) is not available to us.

\smallskip

\noindent End of the comments and explanations concerning Figure 40.

\bigskip

BUT there is also the following big problem with our elementary step, as in Figure 40 may suggest it, namely the following item:

\bigskip

\noindent (127) \quad In the move along the isotopy $S \times 1 \Rightarrow h (S \times 1)$, the action handle $B_a^3$ certainly stays glued to LAVA (this is the purpose of the LAVA bodies). BUT, at some point in this move, the $S \times 1 \subset \partial B_a^3$ has to leave the $\partial (\delta \, {\rm LAVA})$, and cover a piece of the naked lateral surface of the passive 1-handle $B_p^3 \times [-\varepsilon , \varepsilon]$, i.e. a piece of
$$
\partial B_p^3 \times [-\varepsilon , \varepsilon] \subset \partial \, {\rm LAVA} - \delta \, {\rm LAVA}.
$$
End of (127).

\bigskip

This is a potential danger for our sacro-sancted STRONG PRODUCT PROPERTY of LAVA. In what follows next we will show how to circumvent this potential difficulty.
$$
\includegraphics[width=75mm]{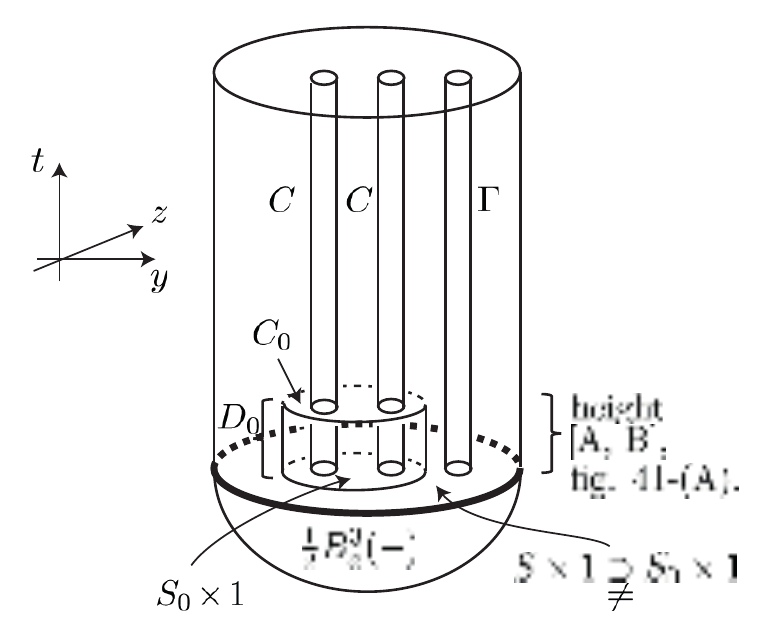}
$$
\label{fig42}
\centerline{\bf Figure 42.}
\begin{quote}
A section $x = {\rm const}$; here $x=x(S_0 \times 1)$, in terms of Figure 41-(A). We see the lateral surface of the 1-handles $H_{n+1}^r = B_a^3 \times [-\varepsilon , \varepsilon]$, Figure 41-(A) and the lateral surface of the cigar $D_0$, glued to each other. The $(({\rm curve} \ C) \cap D_0)$ is here an $x$-projection. The $\Gamma$, unlike the $C$, does not ride on the LAVA bridge $D_0$.
\end{quote}
$$
\includegraphics[width=115mm]{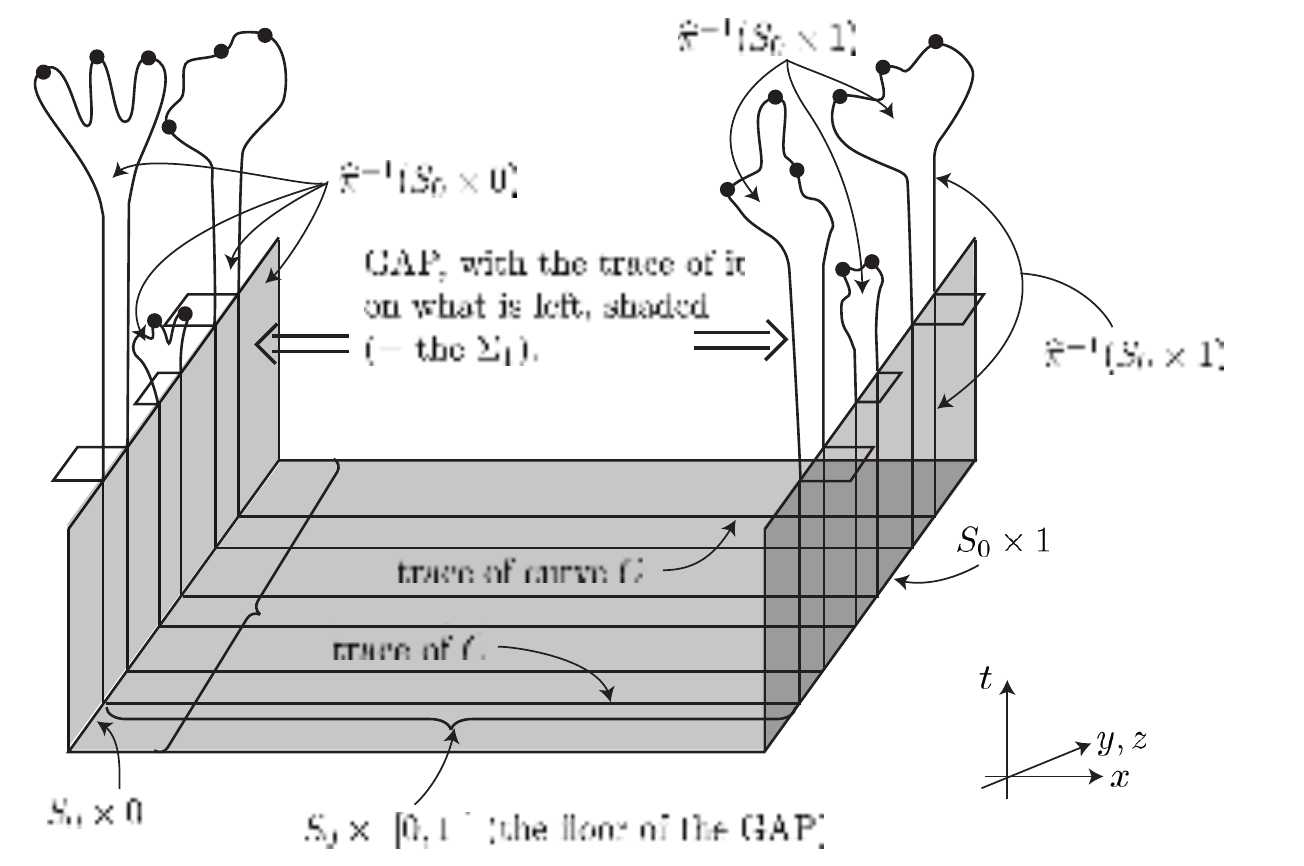}
$$
\label{fig43}
\centerline{\bf Figure 43.}
\begin{quote}
The structure of the GAP, from the Figure 41-(B). The fat points stand for the lamination, living at infinity; see the subsection COMPACTIFICATION, in section III. The $\Gamma$ will runs parallel to the ghostly $C$'s, but staying always {\ibf outside de GAP}. The splitting hypersurface $\Sigma_1$ from (131) is also visible here. Do not mix it up with the $\Sigma_{\infty}^2$. The GAP spans between two arrows signs.
\end{quote}

\bigskip

As said above, the slide of $y = B_a^3$ over the $x = B_p^3$ (corresponding to $S \times 1 \Rightarrow h(S \times 1))$, if conducted in the obvious na{\"\i}ve manner (which Figure 40 may suggest), is not clearly preserving the product property of LAVA. So we will present another transformation, going from the same initial $S \times 1$ to the same final $h(S \times 1)$, but in a manner where, in a manifest way, the PRODUCT PROPERTY is preserved at {\ibf every} intermediary moment. And, as far as we are concerned, it is only the final stage, at $t=1/2$, which interests us, and not the details of $(t=0) \Rightarrow \left( t = \frac12 \right)$, when the PRODUCT PROPERTY OF LAVA $\left( t = \frac12 \right)$ is concerned.

\smallskip

Our product property is encapsulated in the nice retractions $\hat\pi, \pi$ below
$$
\xymatrix{
{\rm LAVA}^{\wedge} \ar[r]_{\hat\pi} &\delta \, {\rm LAVA} \ar@{<->}[d]^{\rm id} \\
{\rm LAVA} \ar@{^{(}->}[u] \ar[r]_{\pi} &\delta \, {\rm LAVA},
} \leqno (128)
$$
which come with collapsible fibers. We will make big use now of (128), as well as of its extension to the $L_0^4 \supset {\rm LAVA} (t=0)$, to be introduced next. We define
$$
\delta L_0^4 \equiv \delta \{{\rm LAVA} \, (t=0) + \mbox{bridges and dilatations}\} \supset S_0 \times [0,1],
$$
which is visible in Figures 40 and 41-(A). We tilt next the $\hat\pi^{-1} (S_0 \times [0,1])$, with $S_0 \times [0,1] \subset (t=0)$ (and this means the value zero of our coordinate axis $t$, in the coordinate system from Figures 40 to 42, and not the time $t=0$), like the ``GAP'' from Figure 41-(B), which displays the tilting.

\smallskip

What follows next, starting with Figure 41-(B) will be a purely internal LAVA affair, concerning $L_0^4$, not reflected outside of LAVA, in particular not reflected for the time being on the $D^2 (\Gamma)$ (see Figure 40, where both $C$'s and $\Gamma$'s climb on $B_a^3 \times [-\varepsilon , \varepsilon] = y$) and $D^2 (\Gamma)$ will be relocated only in the end. With this, we introduce the following object:
$$
{\rm GAP} \equiv \hat\pi^{-1} (S_0 \times [0,1]) \subset \hat L_0^4 , \quad \mbox{see Figure 41-(B),} \leqno (129)
$$
which is split from the rest of $\hat L_0^4$ by
$$
\Sigma \equiv \hat\pi^{-1} (S_0 \times 0) + \hat\pi^{-1} (S_1 \times 1), \leqno (130)
$$
and we will also consider
$$
\Sigma_1 = \Sigma \cup (S_0 \times [0,1]) \subset \partial \, {\rm GAP}, \leqno (131)
$$
separating the GAP from the rest of the world. All this little story above can be vizualized in the Figures 41-(B) and 43. [But then, at the level of Figure 40, where it was introduced, our $S_0 \times [0,1]$ was a part of the $3^{\rm d}$ $\delta$ (LAVA bridge) and it looks that, when we delete the GAP, as we will do, it will totally disappear. We do not want this to happen, and so in Figure 41-(B) we have left a bit of $4^{\rm d}$ $\{${\ibf additional material}$\} \subset$ LAVA bridge (shaded and shielding $S_0 \times [0,1]$), which should stay put (i.e. undeleted) when the GAP is taken away. But we will not stress this in our notations.] As said, these things are vizualizable in Figures 41-(B), 43. The point is here the following

\bigskip

\noindent (132) \quad If we consider the following {\ibf closed} subsets of ${\rm LAVA} \, (t=0)^{\wedge} \equiv \hat L_0^4$, namely $\hat\Lambda^4 \equiv \{\hat L_0^4$ with that GAP deleted, modulo the splitting hypersurface $\Sigma\} \subset \hat\Lambda_1^4 \equiv \hat\Lambda^4 \cup (S_0 \times [0,1]) \cup \{$the additional material, shaded in Figure 41-(B)$\}$ and resting on $S_0 \times [0,1]$, where the $\hat\Lambda^4$ is a manifold and $\hat\Lambda^4_1$ not (even if we leave in place the $4^{\rm d}$ additional material just mentioned) then, at the {\ibf local} level of $S_0 \times 1$, in Figure 41-(B), the two pieces of $\hat\Lambda_1^4$, namely the LHS of the GAP $\supset S_0 \times 0$ and the RHS of the GAP $\supset$ the 1-handle $y$, {\ibf only} hang together, at level $\hat\Lambda_1^4$, via $S_0 \times 1$, which is codimension two. End of (132).

\bigskip

Of course $\hat\Lambda_1^4$ is connected, globally, and dim $(S_0 \times 1) = 2$, as said. As a reflex of (128), we have two PROPER Whitehead collapses, both infinite of course, $\hat\pi (\hat\Lambda^4)$ and $\hat\pi ({\rm GAP})$, occurring in the formulae below
$$
\hat L_0^4 \underset{\mbox{\footnotesize DELETION}}{-\!\!\!-\!\!\!-\!\!\!-\!\!\!-\!\!\!-\!\!\!-\!\!\!-\!\!\!-\!\!\!-\!\!\!\longrightarrow} \Lambda^4 \xrightarrow{ \ \ \hat\pi (\hat\Lambda^4) \ \ } (\delta \, {\rm LAVA}) \cap \hat\Lambda^4 , \leqno (132.1)
$$
and
$$
{\rm GAP} \xrightarrow{ \ \ \hat\pi ({\rm GAP}) \ \ } \Sigma_1 \leqno (132.2)
$$
out of which we can reconstruct $\hat\pi$ (of course $\hat L_0^4 = \Lambda^4 \cup {\rm GAP})$, here just like ${\rm LAVA}^{\wedge}$ in (128), the $\hat L_0^4$ also has the PRODUCT PROPERTY, expressed via the $\hat\pi (\hat L_0^4)$.

\smallskip

We can perform now in order, the following steps

\bigskip

\noindent (132.3) \quad a) Extend the first line in (132.1) to the infinite collapse
$$
\{ \hat\Lambda_1^4 \ \mbox{without the ($\{$additional material$\}$} - S_0 \times [0,1])\} \underset{\hat\pi (\hat\Lambda_1^4)}{-\!\!\!-\!\!\!-\!\!\!-\!\!\!-\!\!\!\longrightarrow} \delta \, {\rm LAVA}.
$$
b) Next, perform the compact dilatation
$$
\{ \hat\Lambda_1^4 \ \mbox{without the ($\{$additional material$\}$} - S_0 \times [0,1])\} \longrightarrow \Lambda_1^4 \ \mbox{(132) (which includes the additional material)}.
$$
The two items a), b) together, express the PRODUCT PROPERTY (albeit a singular one) for $\hat\Lambda_1^4$. This is, schematically,
$$
\hat\Lambda_1^4 \underset{\hat\pi (\hat\Lambda_1^4)}{-\!\!\!-\!\!\!-\!\!\!-\!\!\!-\!\!\!\longrightarrow} \delta \, {\rm LAVA}.
$$
c) Finally, refine the (132.2) to

\noindent GAP $\underset{\hat\pi ({\rm GAP})}{-\!\!\!-\!\!\!-\!\!\!-\!\!\!-\!\!\!\longrightarrow} \{\Sigma_1$ REDEFINED, by changing its bottom $S_0 \times [0,1]$ into the roof of the additional material, resting on $S_0 \times [0,1]$, which is shaded in the Figure 41-(B)$\}$.

\smallskip

With these things, I claim now that the pair
$$
\left([\hat\Lambda_1^4  \cup \{\mbox{additional material (Fig. 41-(B))}\}] \underset{\overbrace{\mbox{\footnotesize $\Sigma_1$ (REDEFINED)}}}{\cup} {\rm GAP} , \delta L_0^4 \right),
$$
has the PRODUCT PROPERTY.

\smallskip

All these various last items are the ingredients to be used in the next moves (133.I), (133.II). End of (132.3)

\bigskip

For the little story above, one should keep in mind that
$$
\delta \, {\rm LAVA} = [(\delta \, {\rm LAVA}) \cap \Lambda^4] \cup [S_0 \times [0,1]].
$$

$$
\includegraphics[width=145mm]{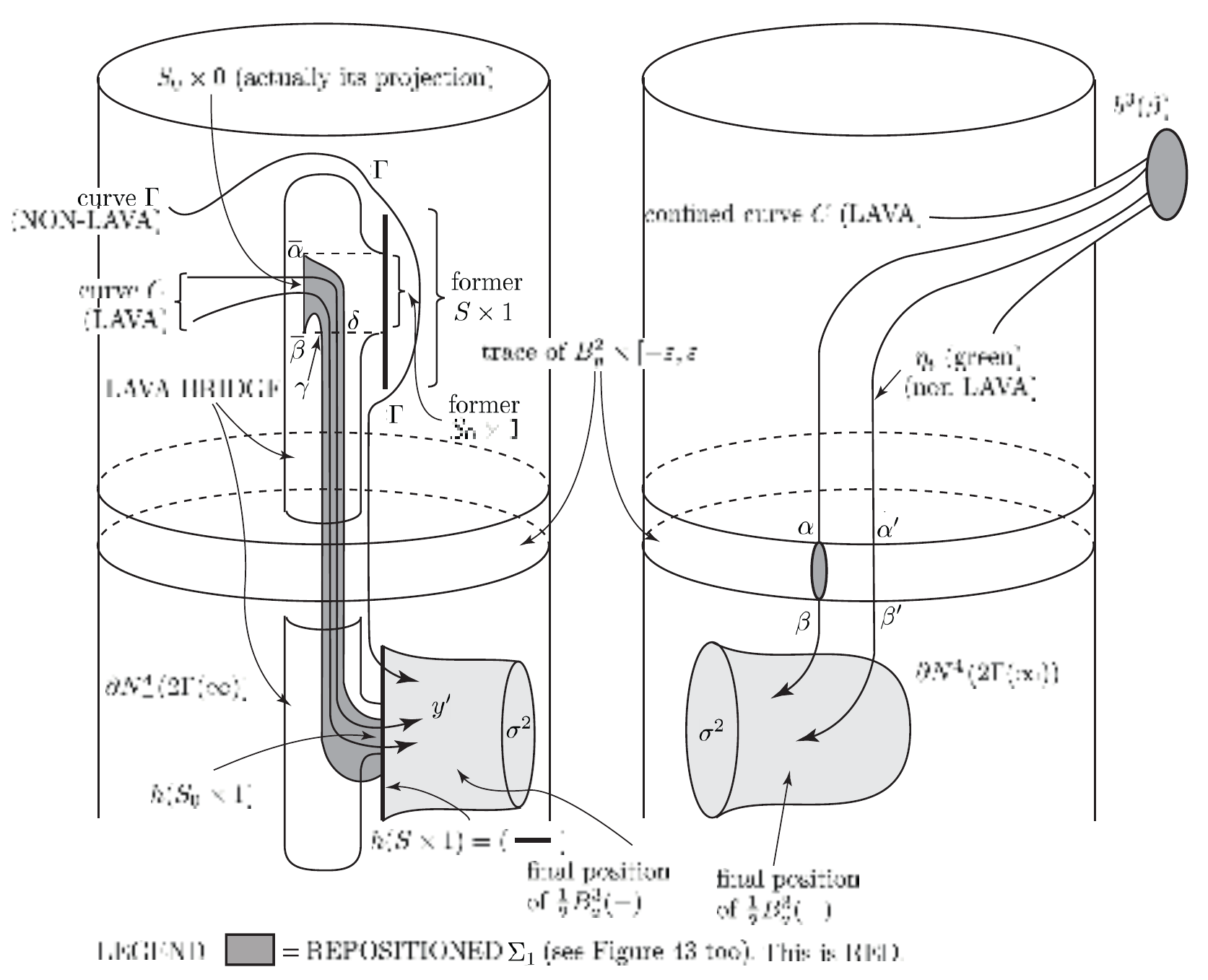}
$$
\label{fig44}
\centerline{\bf Figure 44.}

\bigskip

\noindent {\bf Explanations for the Figure 44.}

\smallskip

This figure, which should be compared to 40, suggests the effect of the step (133.I), on the site of the Figure 40. Many obvious notations, which should be just like in the Figure 40 have not been reproduced again here. The strongly shaded area is part of the REPOSITIONED $\Sigma_1$. We have
$$
\{{\rm final} \ B_a^3 \} \supset h(S \times 1) \supsetneqq h (S_0 \times 1) \subset {\rm LAVA}.
$$
The REPOSITIONING brings the
$$
\{\mbox{RHS of the GAP, Fig. 41-(B)}\} \supset \{\mbox{the active 1-handle} \ B_a^3 \times [-\varepsilon , + \varepsilon]\},
$$
into a correct position, so that it is glued now to $h(\delta \times 1)$.

\smallskip

Here, the $\Sigma_1$ is first REDEFINED like in (132.3) and next REPOSITIONED, like above. 

\smallskip

End of the explanations.

\bigskip

With these things, there is now an alternative procedure for performing our sliding of $y = B_a^3$ over $x = B_p^3$ in a way in which the initial and final stages are exactly as before, but where the PRODUCT property is manifest at all intermediary stages.

\bigskip

\noindent (133.I) \quad Taking advantage of (132) we relocate on top of $N^4 (2\Gamma (\infty)) - B_a^3 (y)$, in the presence of $D_0$, and CORRECTLY, the
$$
\{\mbox{GAP, white in Fig. 41-(B), the additional material stays in place}\} \cup B_a^3 \times 
[-\varepsilon , + \varepsilon]\},
$$
and correctly means here glued to $h (S \times 1)$ from the Figure 40. At the same times, the $\{$LHS of the GAP, same Figure 41-(B)$\}$ is left in place. In view of the PRODUCT PROPERTY of $\hat\Lambda_1^4$, as expressed by (132.3), this step, as performed, preserves the global PRODUCT PROPERTY of LAVA.

\smallskip

The Figure 44 should suggest how we relocate CORRECTLY the bottom $S_0 \times [0,1]$ of $\Sigma_1$. In terms of the pairs described in the formula (132.3), we have a
$$
\Sigma_1 ({\rm REPOSITIONED}) \equiv \Sigma_1 (\mbox{first REDEFINED and then RELOCATED}).
$$

\bigskip

\noindent (133.II) \quad Next, we {\ibf fill the GAP}. The RED band in the Figure 44 corresponds to the $\Sigma_1 ({\rm REPOSITIONED})$ above (which might be suggested by the Figure 43 too), which at the time of the Figure 44 in question is just a DRAWING ON TOP OF LAVA $\hat\Lambda_1^4$. This drawing is, partially, on top of the $\{$additional material, which is neither deleted nor moved$\}$. In the Figure 44, this part of the ride is the area $[\bar\alpha , \bar\beta , \gamma , \delta]$.

$$
\includegraphics[width=11cm]{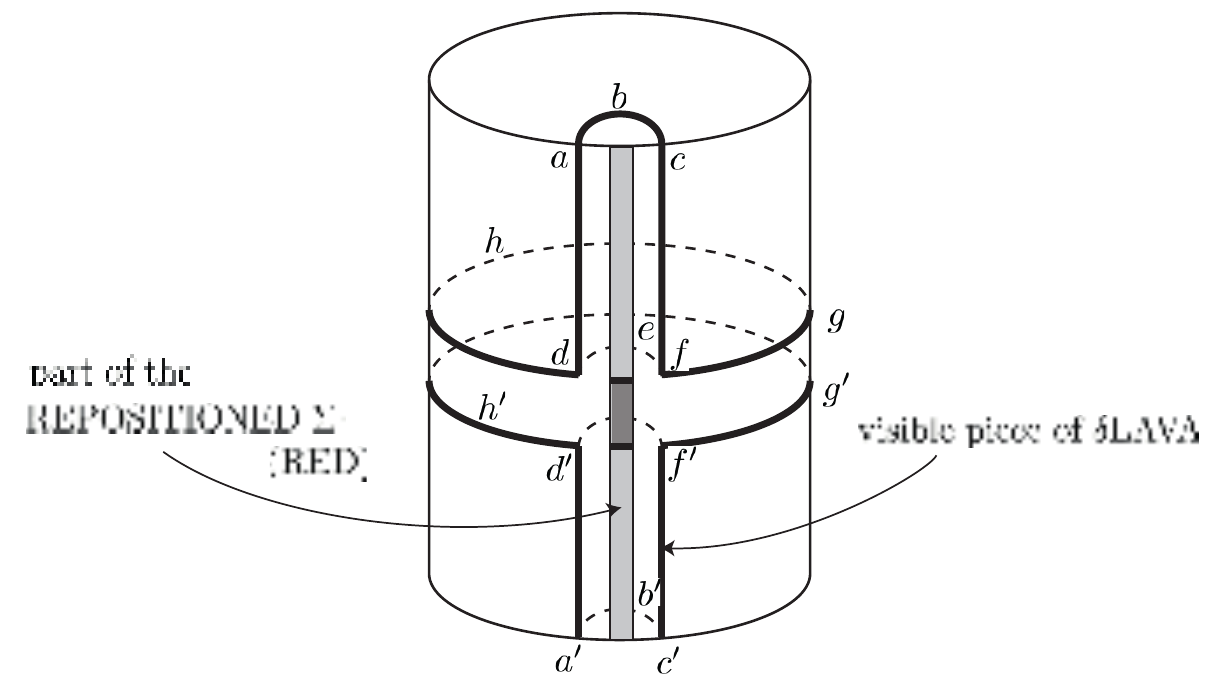}
$$
\label{fig441}
\centerline{\bf Figure 44.1.}
\begin{quote}
A detail of Figure 44 is being presented here, with one dimension less. We see the red band riding on top of LAVA $\hat\Lambda_1^4$, precisely over $\partial \hat\Lambda_1^4 - \delta$ LAVA (the {\ibf same} $\delta$ LAVA as in Figure 40, not affected by our DELETIONS and REPOSITIONINGS). Here $\delta \, {\rm LAVA} \equiv \delta L_0^4$.
\end{quote}

\bigskip

The point here is also that our red drawing (the $\Sigma_1$ (REPOSITIONED)) does not touch the
$$
\delta L_0^4 = \delta ({\rm LAVA} (t=0) , \ \mbox{with bridges and dilatations included}),
$$
which keeps its integrity intact during the DELETIONS and RELOCATIONS, {\ibf except} at the singular spot $S_0 \times 1$, where it glues to the active 1-handle.

\smallskip

So, with our filling of the GAP, in the present step (133.II) we have a pair
$$
([\Lambda_1^4 \cup \{\mbox{additional material}\}] \underset{\overbrace{\mbox{\footnotesize$\Sigma_1$ ({\rm REPOSITIONED})}}}{\cup} {\rm GAP} , \delta L_0^4)
$$
which, in the one hand, manifestly has the PRODUCT property and on the other hand is also isomorphic to the pair which actually interests us $\{({\rm LAVA} , \delta \, {\rm LAVA})$ at the end of the slide of the $y = B_a^3 \times [-\varepsilon , \varepsilon]$ which is NON-LAVA over $x = B_p^3 \times [-\varepsilon , \varepsilon]$, which is LAVA$\}$, i.e. LAVA $\left( t = \frac 12 \right)$.

\smallskip

[The formula above is to be compared with a very similar decomposition occurring inside (132.3), but the two $\Sigma_1$'s are not quite the same.]

\smallskip

Once all this has been achieved, we can also bring $\Gamma_j$ to its correct final position, the one suggested in Figure 44.

\bigskip

\noindent {\bf Additional explanations.} A) In the Figure 44.1, the fat contours $[a,b,c,f,g,h,d,a]$ and $[a',b',c',f',g',h',d',$ $a']$ bound $\delta L_0^4$, which, in the context of our figure is 2-dimensional. It splits the presently 3-dimensional LAVA $\hat\Lambda_1^4$, on top of which the red repositioned $\Sigma_1$ rides, from the rest. The doubly shaded red zone in Figure 44.1 corresponds to the place where the repositioned $\Sigma_1$ in Figure 44 goes through $B_a^3 \times [-\varepsilon , \varepsilon]$, i.e. through $\partial \, {\rm LAVA} - \delta \, {\rm LAVA}$.

\smallskip

B) At the fat site $[\alpha \beta]$ on the $\oplus$ side, in Figure 44, we have a contact
$$
(\delta \, {\rm LAVA} \mid X_b^2) \cap {\rm LAVA} \mid X_r^2
$$
where the two things glue together. This is OK (for the product property), since the big collapse ${\rm LAVA}^{\wedge} \to \delta \, {\rm LAVA}$ goes (in terms of $2X_0^2$) from the BLUE side to the RED side. (We will come back to this later.) The $\eta_i$ (green) is the boundary of an exterior disc, so it does not ``stick'' at $[\alpha' \beta']$, where there is no problem. End of explanations.

\bigskip

\noindent {\bf Lemma 17. (The time $t = \frac12$ compactification.)} 1) {\it Proceeding with a finite number of steps like the one described above, starting from the schematical Figure $39$, we can realize geometrically the abstract balancing $(2\Gamma (\infty) , \Gamma (1)) \Longrightarrow (2\Gamma (\infty) [{balanced}], \Gamma (3))$ and we will use now the notation $2\Gamma (\infty) [balanced] \equiv 2\Gamma (\infty)$ (time $t = \frac12$). So inside $N_1^4 (2X_0^2)^{\wedge}$ we have a transformation proceeding via 1-handle slidings respecting {\rm (126.1), (126.2)}
$$
(N^4 (2\Gamma (\infty)) , N^4 (\Gamma (1)) \Longrightarrow (N^4 (2\Gamma (\infty)) \left( t = \frac12 \right), N^4 (\Gamma(3)). \leqno (134)
$$
This transformation brings us from time $(t=0)$ to time $(t=\frac12)$, as far as the $1$-skeletons are concerned. It also drags along, via covering isotopy, the link $(97)$, hence the transformed family of curves $\underset{1}{\overset{\infty}{\sum}} \, C_i$ as well as the $\underset{1}{\overset{\infty}{\sum}} \, D^2(C_i)$ which cobounds it. This brings about the LAVA $(t=\frac12)$, which has the STRONG PRODUCT PROPERTY just like LAVA $(t=0)$.}

\smallskip

2)  {\it The transformation from $(134)$ leaves the $\underset{i=1}{\overset{M}{\sum}} \, H_i^b = \underset{1}{\overset{M}{\sum}} \, b_i$ intact. From now on, this is {\ibf the} family of $1$-handles of $N^4 (\Gamma (3))$ and we will forget about the $\underset{1}{\overset{M}{\sum}} \, H_i^r$, completely, in the rest of this paper.

\smallskip

We have a diffeomorphism
$$
\left( N^4 (\Gamma (3)) , \sum_1^M H_i^b \right) \underset{\rm DIFF}{=} \underset{i=1}{\overset{M}{\#}} (S_i^1 \times B_i^3 , (*) \times B_i^3), \leqno (135)
$$
i.e. the pair in the LHS of $(135)$ is standard.}

\smallskip

3) {\it As a consequence of $(122)$ there is a next diffeomorphism which, until further notice, supersedes the $(122)$
$$
\Delta^4 = N^4 (\Delta^2)_{\rm Schoenflies} \underset{\rm DIFF}{=} N^4 (2X_0^2)^{\wedge} \left( t = \frac12 \right) = \Biggl[ \left((N^4 (2\Gamma(\infty)) \left( t=\frac12 \right) = \sum_1^{\infty} h_i \left( t=\frac12 \right) \right) \cup \leqno (136)
$$
$$
\cup \, \mbox{LAVA} \left( t=\frac12 \right)\Biggl]^{\wedge} + \sum_1^{\overset{=}{n}} D^2 (\Gamma_j) \subset N_1^4 (2X_0^2)^{\wedge} = \Delta^4 \cup \partial \Delta^4 \times [0,1].
$$
[In order to deduce the $(136)$ from $(122)$, the passage $(t=0) \Rightarrow (t = \frac12)$ has to be conceived in the obvious na{\"\i}ve manner, and NOT like in the context of the Figures $41$ to $44.1$ (that was necessary just for making sure of the PRODUCT PROPERTY of LAVA $(t=\frac12)$), and the covering isotopy theorem can then be applied.]}

\smallskip

4) {\it We also have a diffeomorphism
$$
\left\{ \left[
\underbrace{\left(N^4 (2\Gamma (\infty))\left( t = \frac12 \right) - \sum_1^{\infty} h_i \left( t = \frac12 \right) \right) \cup \mbox{LAVA} \left( t = \frac12 \right)^{\wedge}}_{call \, this \ \hat Z^4 \left( t=\frac12 \right)}
\right] , \sum_1^M \{\mbox{extended cocore} \ H_i^b\}^{\wedge} \right\} =
$$
$$
\underset{\rm DIFF}{=} \underset{i=1}{\overset{M}{\#}} (S^1_i \times B_i^3 , (*) \times B_i^3), \mbox{i.e. the pair in the LHS is standard.}
$$

So, the $\underset{1}{\overset{M}{\sum}} \, \{\mbox{extended cocore} \ H_i^b\}^{\wedge}$ {\ibf is} the system of $1$-handles of $N^4 (2X_0^2)^{\wedge}$ from} (136).

\smallskip

5) {\it In the same context of $(136)$ we have the system of exterior discs, which is smoothly embedded
$$
\left( \sum_1^M (\delta_i^2 , \partial \delta_i^2 = \eta_i (\mbox{green})) \right) \xrightarrow{ \ J \ } (\partial \Delta^4 , \partial \Delta^4 \times [0,1]). \leqno (138)
$$
}

6) (Reminder) {\it We have $\eta_i (\mbox{green}) \cdot \left( B - \underset{1}{\overset{M}{\sum}} \, H_i^b \right) = 0$ but the geometric intersection matrix of interest for us is now the

\bigskip

\noindent $(139)$ \quad $\eta_i (\mbox{green}) \cdot \{\mbox{extended cocore} \ H_j^b \}^{\wedge} = \delta_{ij} + \{${\ibf parasitical} RED terms $\eta_i (\mbox{green}) \cdot h_k (t=\frac12)$, where $h_k (t=\frac12) \subset\{$exterior cocore $H_j^b \}\}$, here, of course, $1 \leq i,j \leq M$, and $h_k \in R-B$.}

\bigskip

\noindent {\bf Sketch of proof.} Enough has been said, since the beginning of the present section VII, so as to make point 1) completely clear. Concerning (135), notice the following, at the purely abstract level of the BALANCING LEMMA 15. Every edge $e \subset \Gamma (1) = \Gamma(1) \times (\xi_0 = -1)$ carries a $b_i \in B_0$, and this is the family $\underset{1}{\overset{M}{\sum}} \, b_i$.

\smallskip

So $\Gamma (1) - \underset{1}{\overset{M}{\sum}} \, b_i = \Gamma (1) - B_0 = \underset{\alpha = 1}{\overset{\chi + 1}{\sum}} \, \chi_{\alpha}$ is a collection of very small graphs $\chi_{\alpha}$, each with a single vertex and each with two ends, corresponding to two distinct $b_i ,b_j$, with $1 \leq i,j \leq M$. When the mechanism which proves Lemma 15 is applied to this situation, then we get a tree $\Gamma (1) \cup \underset{1}{\overset{\chi}{\sum}} \, g_i ({\rm new}) - B_0$, the vertices of which come organized in pairs, each pair corresponding to the endpoints of one of the BLUE handles in $\underset{1}{\overset{M}{\sum}} \, H_i^b$, like in the drawing below:
$$
\includegraphics[width=11cm]{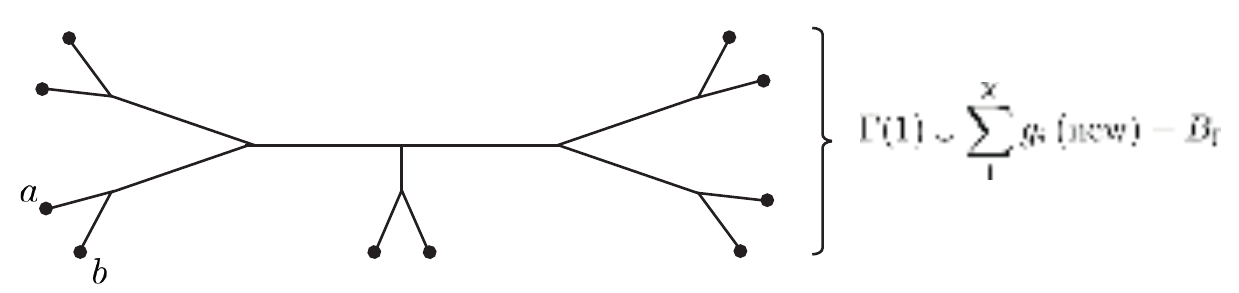}
$$
Here, $a,b$ are endpoints of some given $b_i$ and in real life they do not have to be as close to each other as in the drawing above, but rather like in the next drawing:
$$
\includegraphics[width=7cm]{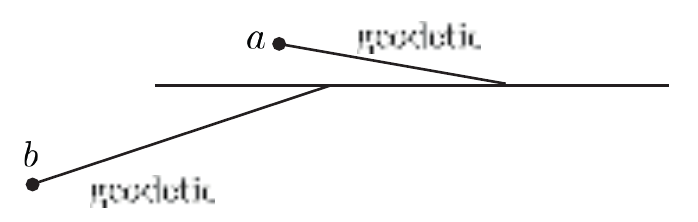}
$$

The (135) immediately follows from this. Then, as already said, the (136) in 3) follows from $\{$(122) from $(t=0)\}$, by appealing to the covering isotopy theorem, and making use of the PRODUCT property of LAVA $(t=\frac12)$.

\smallskip

The basic fact
$$
\left( \hat Z^4 \left( t = \frac12 \right) \cup {\rm LAVA} \, \left( t = \frac12 \right)^{\wedge} , \sum_1^M \{\mbox{extended cocore} \ H_i^b \}^{\wedge} \right) \underset{\rm DIFF}{=} \underset{i=1}{\overset{M}{\#}} (S_i^1 \times B_i^3 , (*) \times B_i^3)
$$
follows by combining (135) with the strong product property of LAVA $(t=\frac12)$. The rest of the Lemma 17 is essentially a reminder, and our proof ends here. \hfill $\Box$

\bigskip

The rest of this paper will be devoted to the elimination of the parasitical red terms in (139). We will call this last step
$$
\left( t=\frac12 \right) \xLongrightarrow{ \ \ \mbox{\footnotesize THE COLOUR-CHANGING} \ \ } (t=1). \leqno (140)
$$

Here are, to begin with, two remarks concerning the parasitical terms $\{ k_k (t=\frac12)$ from (139)$\}$.

\medskip

$\bullet$) Since the total length of $\underset{1}{\overset{M}{\sum}} \, \eta_i ({\rm green})$ is finite, i.e. length $\left(\underset{1}{\overset{M}{\sum}} \,  \eta_i ({\rm green})\right) < \infty$, there are, to begin with, only {\ibf finitely} many parasitical RED terms in (139).

\medskip

$\bullet\bullet$) Since $\eta_i ({\rm green}) \cdot \left( B_1 - \underset{1}{\overset{M}{\sum}} \, H_i^b \right) = \emptyset$, these terms are all in $R_0 - B_0 = R_1 - B_1$.

\bigskip

We go back for a while now to $t=0$ and we will denote by $h(r)$ the $h_i \subset X_0^2 \times r (\approx X_0^2 [{\rm new}]) \subset 2X_0^2$ and by $h(b)$ the others, i.e. the $h_i \subset X_b^2 \subset 2X_0^2$ (see here (16)). This gives a partition
$$
\sum_1^{\infty} h_i = \sum_1^{\infty} h_{\ell} (r) + \sum_1^{\infty} h_k (b) \leqno (141)
$$
inducing, by the $C \cdot h = {\rm id} + {\rm nil}$ duality, a second partition
$$
\sum_1^{\infty} C_i = \sum_1^{\infty} C_{\ell} (r) + \sum_1^{\infty} C_k (b).
$$
In terms of the Figures 7, 7.bis, all the $R_1 \cap B_1$'s located in $\Gamma (\infty) \times b$, are $h_k (b)$'s. Not withstanding what the notation may suggest, because of the duality 1-handles $\underset{\approx}{\longleftrightarrow}$ curves suggested in Figure 7, the curves $c(r)$ from the Figures 7-(I and IV) are among $c_k (b)$'s, while the curves $c(b)$ from Figures 7-(II, III), and 7.bis, are NOT part of $\underset{1}{\overset{\infty}{\sum}} \, C_i$ at all, since the corresponding $D^2 (c(b)) \subset 2X^2 - 2X_0^2$ (and the RED flow is confined to $2X_0^2$).
$$
\includegraphics[width=160mm]{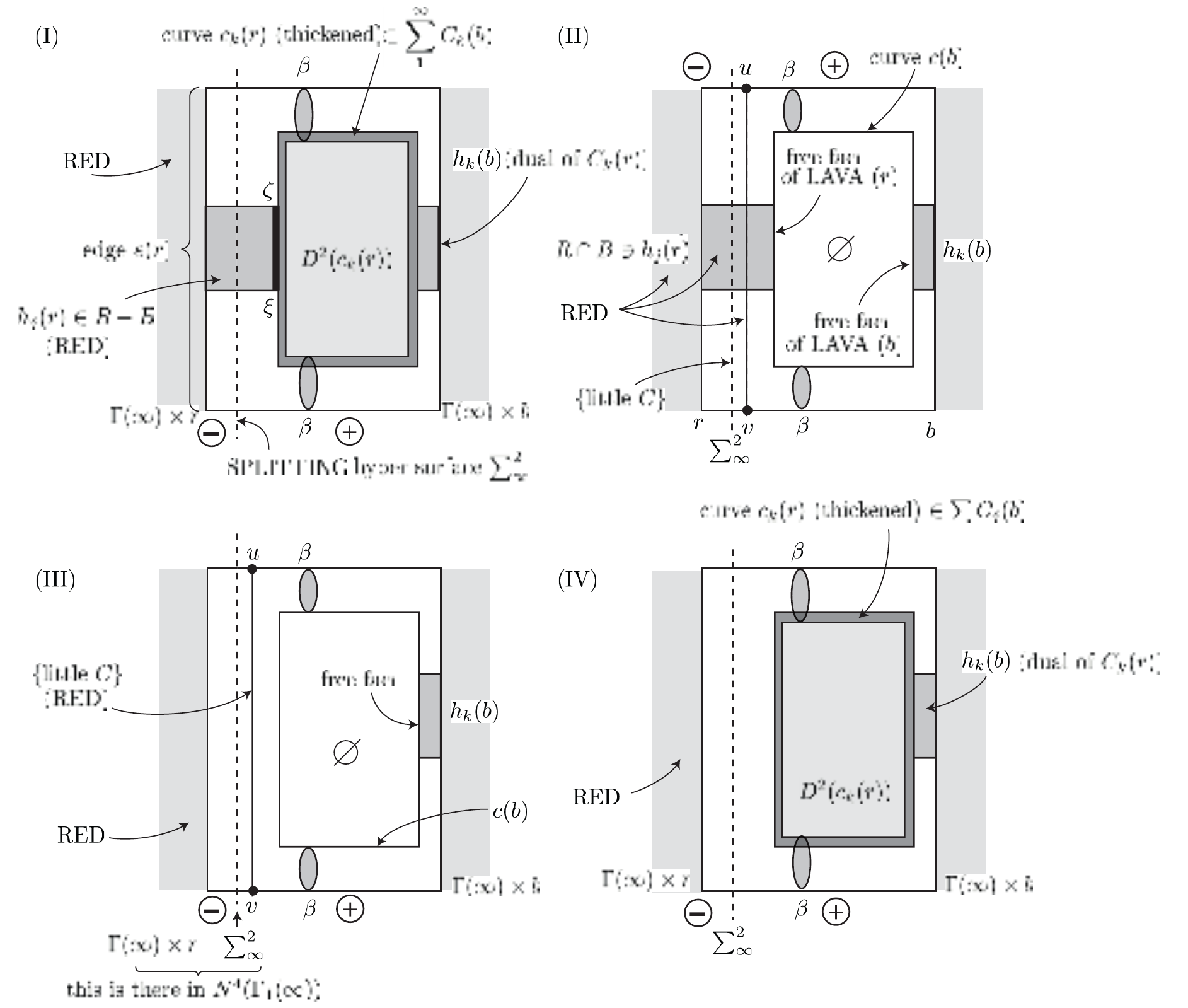}
$$
\label{fig45}
\centerline{\bf Figure 45.}
\begin{quote}
This is Figure 7, embellished. In particular, we see here LAVA $(r)$ and LAVA $(b)$. Do not mix up notations like $c(r)$ (going back to the rectangles of Figure 7) and $C(r)$ occurring in (141). For instance, the $c(r)$ in (I) is a $C_k (b)$ because its dual $h$ is $h_k (b)$. The curves dual to the $h_k (b)$'s in Figures (II), (III) are $C_k (b)$'s living in $X_b^2$. The $h_k (b)$ in (IV) is the dual of $c_k (r) \in \underset{1}{\overset{\infty}{\sum}} \, C_j (b)$.

The Figure 7-bis becomes here like (II) or (III), EXCEPT that: (i) There is no longer any LAVA $(r)$ contribution and, moreover, (ii) There is no longer any added $2^{\rm d}$ meat alive in the purely 1-dimensional $2\Gamma (\infty) \mid$ (Figure 7-bis). Very importantly too, at the free LAVA $(b)$ site ($=$ face in (II), (III)), the internal collapse ${\rm LAVA} \, (b) \searrow \delta \, {\rm LAVA} \, (b)$ can start. The sign ``$\emptyset$'' in (II), (III) means, of course, that the corresponding disc is not there in $2X_0^2$ (but only in $2X^2$, where it is present).

LEGEND: $\includegraphics[width=1cm]{rectangleclair.pdf}$ (RED) $= {\rm LAVA}_r$, $\includegraphics[width=1cm]{rectangleclair.pdf}$ $= {\rm LAVA}_b$, $\zeta {-\!\!\!-\!\!\!-\!\!\!-\!\!\!-} \, \xi$ (occurring {\ibf only} in (I)) $=$ the GLUEING TERM, explained in the main text. 

\smallskip

In the drawing (I') below, a detail of (I) is a bit more accurately rendered.
$$
\includegraphics[width=155mm]{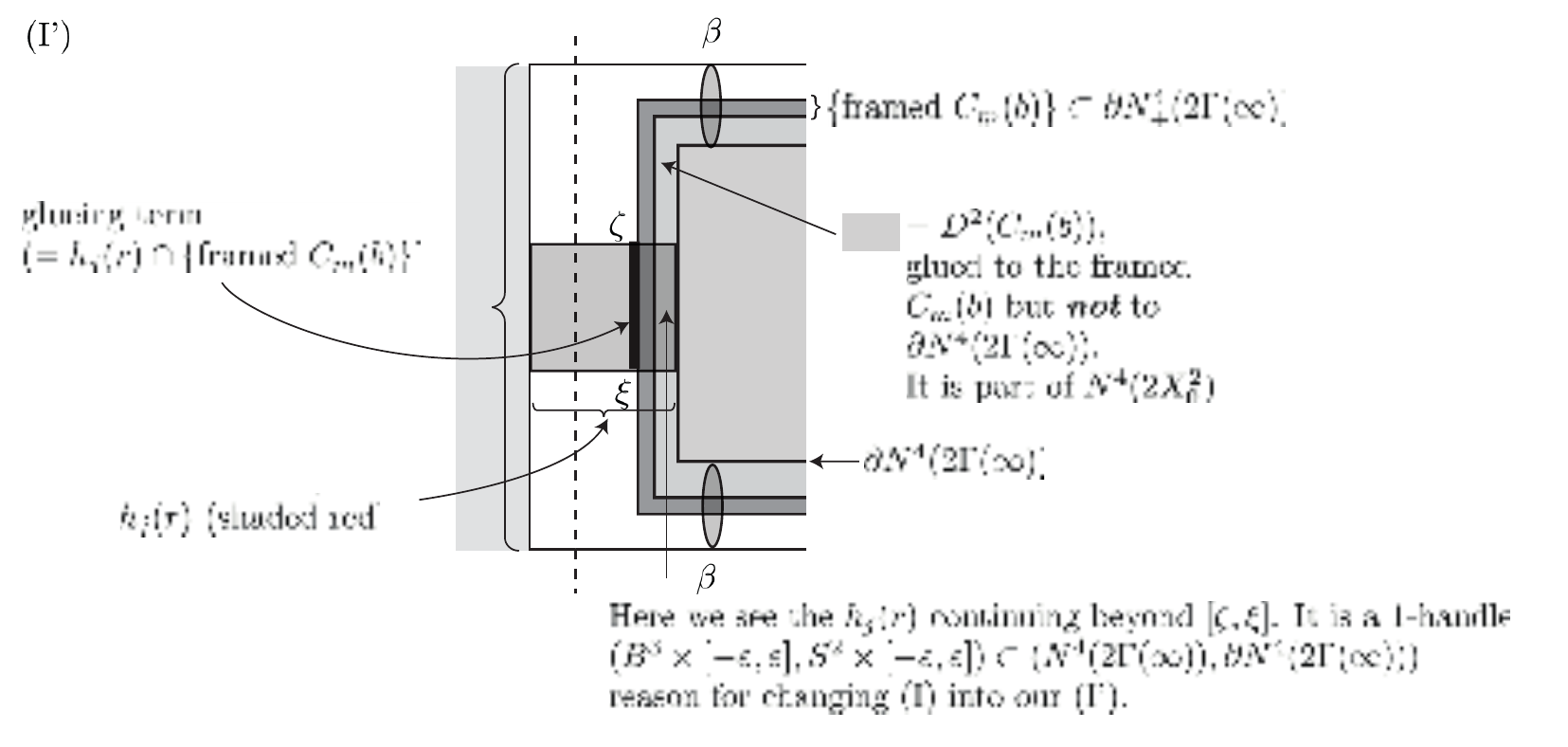}
$$
\label{fig45bis}
\end{quote}

We will enlarge now the purely 1-dimensional $2\Gamma (\infty)$ with some $2^{\rm d}$ meat, as it is suggested in Figure 45, and all this is eventually part of $N^4 (2\Gamma(\infty))$. With the forced confinement conditions (101.1-bis), we have set that $\Sigma c(r) \subset \partial N_+^4 (2\Gamma (\infty))$ and Figure 45 should clarify this. We actually have
$$
\sum_1^{\infty} \{{\rm little} \, C_i ({\rm or} \ \Gamma_j) \} + \sum_1^{\infty} \eta_i \times b + \sum_1^{\infty} c_k (r) \subset \partial N_+^4 (2\Gamma (\infty)) ,
\leqno (142)
$$
but then, $\Sigma \{$little $C_i\} \subset \underset{1}{\overset{\infty}{\sum}} \, C_{\ell} (r)$, while
$$
\sum_1^{\infty} \eta_i \times b + \sum_1^{\infty} c_k (r) \subset \sum_1^{\infty} C_n (b) .
$$
On the other hand, since the duals of the $\{$little $C\}$'s are in $R \cap B$, they will not contribute much to the present story. Finally, we have
$$
\sum_1^{\infty} C_n (r) = \underbrace{\sum_1^{\infty} \{C_i ({\rm remaining})\}}_{{\rm in} \ \partial N_-^4 (2\Gamma (\infty))} + \underbrace{\sum_1^{\infty} \{{\rm little} \ C_i \}}_{{{\rm in} \ \partial N_+^4 (2\Gamma (\infty)) \atop {\rm (by \, our \, construction)}}} .
$$
We will denote by $\{$little $h_i\} \subset \underset{1}{\overset{\infty}{\sum}} \, h_{\ell} (r)$ the dual of $\{$little $C_i\}$. Since both $\Gamma (\infty) - R$ and $\Gamma (\infty) - B$ have to be trees, this forces the $\{$little $h_i\}$ to be in the family $R_0 \cap B_0$. So, in Figure 45-(II) it could (sometimes) happen that the RED $h_j (r)$ we see is a $\{$little $h_i\}$. The dual curve $\{$little $C_i\}$ occurs then as two lines of red colour both inside the $\oplus$ part, in the two contiguous Figures (II) and (III), so that the two lines should melt in to a single closed curve. See here the (II), (III) in our Figure 45, with the points $u$ and $v$ common to both. In terms of Figure 32, the (II) corresponds to the little 1-handle and (III) to the main body. The obsesrvant reader will have already figured out that we need two adjacent Figures (III) and three red lines in all, but that does not change our story. But before saying more concerning the pair $(\{{\rm little} \, C_i \} , \{{\rm little} \, h_i\})$, I will introduce the pairs
$$
({\rm LAVA}_r \, (t=0) , \delta \, {\rm LAVA}_r  \, (t=0)) = \Biggl(\sum_{\ell = 1}^{\infty} (h_{\ell} (r) \cup C_{\ell} (r)) , \leqno (143)
$$
$$
\delta \, {\rm LAVA}_r \, (t=0) \equiv (\delta \, {\rm LAVA}  \, (t=0)) \cap \partial \, {\rm LAVA}_r  \, (t=0)\Biggl) = \partial \, {\rm LAVA}_r  \, (t=0) \cap \partial \left[ N^4 (2\Gamma (\infty)) - \sum_1^{\infty} h_{\ell} (r) \right] \Biggl),
$$
and
$$
({\rm LAVA}_b  \, (t=0) , \delta \, {\rm LAVA}_b  \, (t=0)) = \Biggl(\sum_{k=1}^{\infty} (h_k (b) \cup C_k (b)) , \delta \, {\rm LAVA}_b  \, (t=0) \equiv
$$
$$
\equiv (\partial \, {\rm LAVA}_b  \, (t=0)) \cap \partial \left[ N^4 (2\Gamma (\infty)) - \sum_1^{\infty} h_{\ell} (b) \right] \Biggl).
$$
Very importantly, there is also a glueing term between the two lavas, namely
$$
\mbox{GLUEING TERM} \equiv \delta \, {\rm LAVA}_b  \, (t=0) \cap \partial \, {\rm LAVA}_r  \, (t=0), \leqno (144)
$$
illustrated in the Figure 45-(I), when it occurs as $[\zeta , \xi]$.

\bigskip

\noindent {\bf Lemma 18.} 1) {\it Individually, each of the two lavas, $( {\rm LAVA}_r  \, (t=0) , \delta \,  {\rm LAVA}_r  \, (t=0))$ and $( {\rm LAVA}_b $ $(t=0) , \delta \,  {\rm LAVA}_b  \, (t=0))$ have the STRONG PRODUCT PROPERTY. That property, for the full $( {\rm LAVA}$  $(t=0) , \delta \,  {\rm LAVA}  \, (t=0))$ is a succession of two infinite collapses ${\rm LAVA}_b  \, (t=0) \longrightarrow \delta \, {\rm LAVA}_b  \, (t=0)$ and ${\rm LAVA}_r  \, (t=0) \longrightarrow \delta \, {\rm LAVA}_r  \, (t=0)$ starting from the free faces of $\partial \, {\rm LAVA}_r  \, (t=0)$, which include now the GLUEING TERM.}

\medskip

2) {\it We have
$$
{\rm LAVA}_b \, (t=0) \subset \{\mbox{the $\oplus$ side of the SPLITTING of $N^4 (2X_0^2)$ by $\Sigma_{\infty}^2$}\}.
$$
}

\bigskip

The bulk of ${\rm LAVA}_r \, (t=0)$ lives on the $\ominus$ side of the SPLITTING by $\Sigma_{\infty}^2$. But then, for obvious geometrical reasons, each $h_k (r) = B_k \times [-\varepsilon , \varepsilon]$ which is definitely ${\rm LAVA}_r$ {\ibf has} to have its $\ominus$ half and its $\oplus$ half. For instance, in Figure 45-(I) it is through its $\oplus$ half that the corresponding $h_j(r)$ gets on top of its GLUEING TERM. See also Figure 49 when each $B^3 (r(j))$, $B^3 (b(j))$ has a $\ominus$ half and a $\oplus$ half. This has nothing to do with the $1$-handle in question being ${\rm LAVA}_r$ or ${\rm LAVA}_b$.

\smallskip

Concerning the $(\{{\rm little} \, C_i\} , \{{\rm little} \, h_i\})$ we have the following

\bigskip

\noindent {\bf Lemma 19.} 1) {\it Since $\{{\rm little} \, h_i\} \in R_0 \cap B_0$, it is not concerned by the active part of the COLOUR CHANGE (Lemma $14$ and its geometric realization, in the rest of this section).}

\medskip

2) {\it We have $(\{{\rm little} \, h_i\} \cup D^2 (\{{\rm little} \, C_i\})) \subset {\rm LAVA}_r \, (t=0)$, and also $\underset{i}{\sum} \, (\{{\rm little} \, h_i\}) \cup D^2 (\{{\rm little} \, C_i\}) \cap {\rm LAVA}_b \, (t=0) = \emptyset$ (no corresponding glueing terms exist).}

\medskip

3) {\it The $(\{{\rm little} \, h_i\} \cup \{{\rm little} \, C_i\})$ is an endpoint for the RED collapsing flow of ${\rm LAVA}_r \, (t=0)$, in the sense that there are no outgoing flow-lines from $i$.}

\medskip

All these things which we have just said have concerned the ${\rm LAVA} \, (t=0)$ and, of course, together with the process $(t=0) \Longrightarrow \left( t = \frac12 \right)$ which we have described, comes also a transformation ${\rm LAVA} \, (t=0) \Longrightarrow {\rm LAVA} \left( t = \frac12 \right)$ which involves items $(141)$ to $(144)$, in a manner which
should be more or less automatic. For instance, in the context of the Figure 40, in the sliding of $y$ over $x$, all the curves which we see on the side of $\partial N_+^4 (2\Gamma (\infty))$ are $C(b)$'s, meaning $\eta \times b$'s, or $c(r)$'s, or of course $\eta_{\ell} ({\rm green})$'s OR pieces of $\{{\rm little} \, C\}$, in the worst case. Then $h_i$ is actually RED (precisely an $R \cap B$). So these $\{{\rm little} \, C\}$'s are certainly  mute in the context of Figures 40 to 49, when they have not been explicitly drawn. They may be assumed far from the moving $B_a^3 = y$.

\smallskip

Also, the $[\alpha , \beta]$'s in Figure 44 are all new GLUEING TERMS, and there is a large margin of freedom concerning their positioning. We move now to the next geometric step
$$
\left( t = \frac12 \right) \underset{\mbox{\footnotesize THE BIG COLOUR-CHANGE}}{\xLongrightarrow{ \qquad \qquad\qquad \qquad \qquad \qquad \qquad }} (t=1),
$$
the aim of which is to realize the BIG BLUE DIAGONALITY condition
$$
\eta_i ({\rm green}) \cdot \{\mbox{extended cocore} \ H_j^b \}^{\wedge} = \delta_{ij} \ \mbox{AND} \ \eta_i ({\rm green}) \cdot \left( B - \sum_1^M H_j^b \right) = 0.
$$
What this means is that we both want to get rid of those finitely many RED parasitical terms $h_k = h_k \left( t=\frac12 \right) \in R-B$ which occur in (139) and, at the same time, {\ibf not} create any additional $\eta ({\rm green}) \cdot B$. For this purpose, we start by constructing a finite, sufficiently high truncation $\Phi_1$ of the isomorphism $\Phi$ (103), a truncation which is such that $\Phi_1$ is gotten by stopping the inductive process in Lemma 14 at some finite level which we denote by $j_1$. This $j_1$ is high enough so that we should have

\bigskip

\noindent (144.1) \quad The $j_1$ is sufficiently close to $\infty^{\rm ty}$ so that $R(j_1)$ should contain all the $\Phi$ ($h_k$, parasitical term in (139)). Put differently, we want to have $\{$All the $h_k$ which are $R-B$'s touched by $\underset{1}{\overset{M}{\sum}} \, \eta_{\ell} ({\rm green})$, off-diagonally$\} \subset \{$The family $r(1) , r(2) , \ldots , r(j_1)\} \subset R-B-R(j_1)$. These $h_k$'s are all, automatically LAVA $(t)$, for all $t$'s until they change their colour, and then we have
$$
\Phi (r(i)) \in B-R , \ {\rm for} \ i \leq j_1 .
$$
So, {\ibf all} the parasitical $h_k$'s will have been already turned BLUE, at time $j_1 = (t=1)$. Their place inside the LAVA will have been taken by the BLUE $\Phi (h_k ({\rm parasitical}))$, which are {\ibf untouched} by $\underset{\ell}{\sum} \, \eta_{\ell} ({\rm green})$. Here are somme additional explanations concerning the family
$$
\mbox{``LHS of (144.1)''} \equiv \{\mbox{the parasitical $h_k \in R-B$, touched by} \ \sum_{\ell} \eta_{\ell} ({\rm green})\}.
$$

$\bullet$) Since length $\underset{\ell}{\sum} \, \eta_{\ell} ({\rm green}) < \infty$, $\#$ (``LHS of (144.1)'') $< \, \infty$ too.

$\bullet$$\bullet$) Our LHS of (144.1) corresponds to the (139).

$\bullet$$\bullet$$\bullet$) In what follows next, in this paper, the family $\{$LHS of (144.1)$\}$ will only {\ibf decrease} when $t$ increases.

\bigskip

\noindent (144.2) \quad At any time $t$, here is the list of constituent pieces of LAVA $(t)$: $1$-handles $h_i$ or $\Phi (h_i)$, RED or BLUE, which are LAVA and the corresponding $D^2 (C_i)$ (dragged over $\Phi (h_i)$ when $h_i \Rightarrow \Phi (h_i)$, and this is part of LAVA$_r$), then also LAVA BRIDGES and LAVA DILATATIONS.

\smallskip

The total contribution of these additional items will be compact. The additional pieces insert themselves, when that is needed, in order to preserve the PRODUCT PROPERTY OF LAVA, when $C \cdot h = {\rm id} + {\rm nil}$ is contradicted, between a $D^2 (C_i)$ and a $h_j$ where we have $D^2 (C_i) \cdot h_j \ne 0$ (i.e. $i \to j$) in the geometric intersection matrix.

\smallskip

Any $h \in R-B$ might be replaced at some time $t_0 (h)$, and then for all times $t \geq t_0 (h)$, by the corresponding $\Phi_1 (h) \in B-R$, after which the $\Phi_1 (h)$ replaced $h$ inside LAVA, while $h$ leaves the scene for ever.

\smallskip

This will be the GEOMETRIC $4^{\rm d}$ REALIZATION OF THE COLOUR-CHANGE. End of (144.2)

\bigskip

In the meanwhile, during the step $\left( t = \frac12 \right) \Longrightarrow (t=1)$, the condition $C \cdot h = {\rm id} + {\rm nilpotent}$ will get violated and when we will have to get the $\{$extended cocores$\}$'s, in particular the all-important $\{$extended cocore $b_{i \leq M}\}$'s, we will have to appeal directly to the PRODUCT PROPERTY of LAVA.

\smallskip

With this, both the contribution of the population of $h$'s and of the LAVA bridges inside the extended cocores, may increase in time; but see here also what is said concerning the parasitical $h_{\ell}$'s (see the (139) for these parasites).

\smallskip

The LAVA bridges will be living inside $\partial N_-^4 (2\Gamma (\infty))$ (more precisely their piece of $\delta (\mbox{LAVA bridges}) \subset \partial N_-^4 (2\Gamma (\infty))$, far from $\overset{M}{\underset{1}{\sum}} \, \eta_{\ell} ({\rm green}) \subset \partial N_+^4 (2\Gamma (\infty))$.

\smallskip

On the other hand, the total population of these $h_{\ell}$ {\ibf touched} by $\overset{M}{\underset{1}{\sum}} \, \eta_{\ell} ({\rm green})$ (and they will all be certainly LAVA), will constantly {\ibf decrease} in time, during the ministeps of the $(t=j-1) \Rightarrow (t=j)$ occurring in the (144.3) below. They will eventually all {\ibf turn BLUE}, via the COLOUR-CHANGING process.  So, at time $t=1$, as a consequence of the little blue diagonalization already achieved at (95) above, for any $h_i$ still alive inside LAVA, we will find
$$
\sum_{\ell = 1}^M \eta_{\ell} ({\rm green}) \cdot h_i = \emptyset.
$$

Now, in order to achieve the (144.1), we will realize, for all the $j \leq j_1$, geometrically, in $4^{\rm d}$ inside $N_1^4 (2X_0^2)^{\wedge}$, the abstract steps $R(j-1) \Rightarrow R(j)$ from Lemma 14. This means that we will have a whole collection of intermediary times, between $t = \frac12$ and the final $t=1$
$$
\left( t = \frac12 \right) (\mbox{which we also call} \ (t=(j=0)), t=(j=1) , t=(j=2) , \ldots , (t=j-1) , \leqno (144.3)
$$
$$
(t=j) , \ldots , (t=j_1) \equiv (t=1). 
$$
At all the times $t$ (144.3) the $N^4 (2\Gamma (\infty))$ will be unchanged, with us, but then there will be fluid times $j-1 < t < j$, when $N^4 (2\Gamma (\infty))$, at the intermediary moments between the ones in (144.3), will change. It will become then a fluid, time-dependent object denoted $Z(t)$. But, I insist, at all the integral moments $t$ in (144.3), we always have $Z(t) = N^4 (2\Gamma (\infty))$.

\smallskip

We will describe now the main inductive step, namely

\bigskip

\noindent {\bf The geometry of the step $(t=j-1) \Longrightarrow (t=j)$, $j \leq j_1$ in (144.3)} .

\smallskip

We insist, like for $t=0 \Rightarrow t=\frac12$, that (126.1) and (126.2) should be satisfied. The $\Gamma_1 (\infty) \subset \Gamma (2\infty)$ (to be changed into $N^4 (\Gamma_1 (\infty)) \subset N^4 (2\Gamma (\infty))$), at the initial stage $t=j-1$, is shown in the Figure 46, which is a vastly more elaborate form of the Figure 37. But, before we start looking into the figure, here comes the

\bigskip
\bigskip

\noindent {\bf Important comment (144.4) (concerning (126.1)).} We have the following general formula, valid at all times
$$
{\rm LAVA} = {\rm LAVA}_r \underset{\overbrace{\mbox{\footnotesize GLUEING TERMS}}}{\cup} {\rm LAVA}_b .
$$

\medskip

Independently of each other, both ${\rm LAVA}_r$ and ${\rm LAVA}_b$ share the strong product property. With this, however we place our GLUEING TERMS, the {\ibf global} lava has automatically the STRONG PRODUCT PROPERTY, with a collapsing flow starting in ${\rm LAVA}_b$, reaching ${\rm LAVA}_r$ at the GLUEING TERMS, now free faces, and continuing through the ${\rm LAVA}_r$. When we want to check (126.1) at various stages of our construction, we should keep this in mind.

$$
\includegraphics[width=14cm]{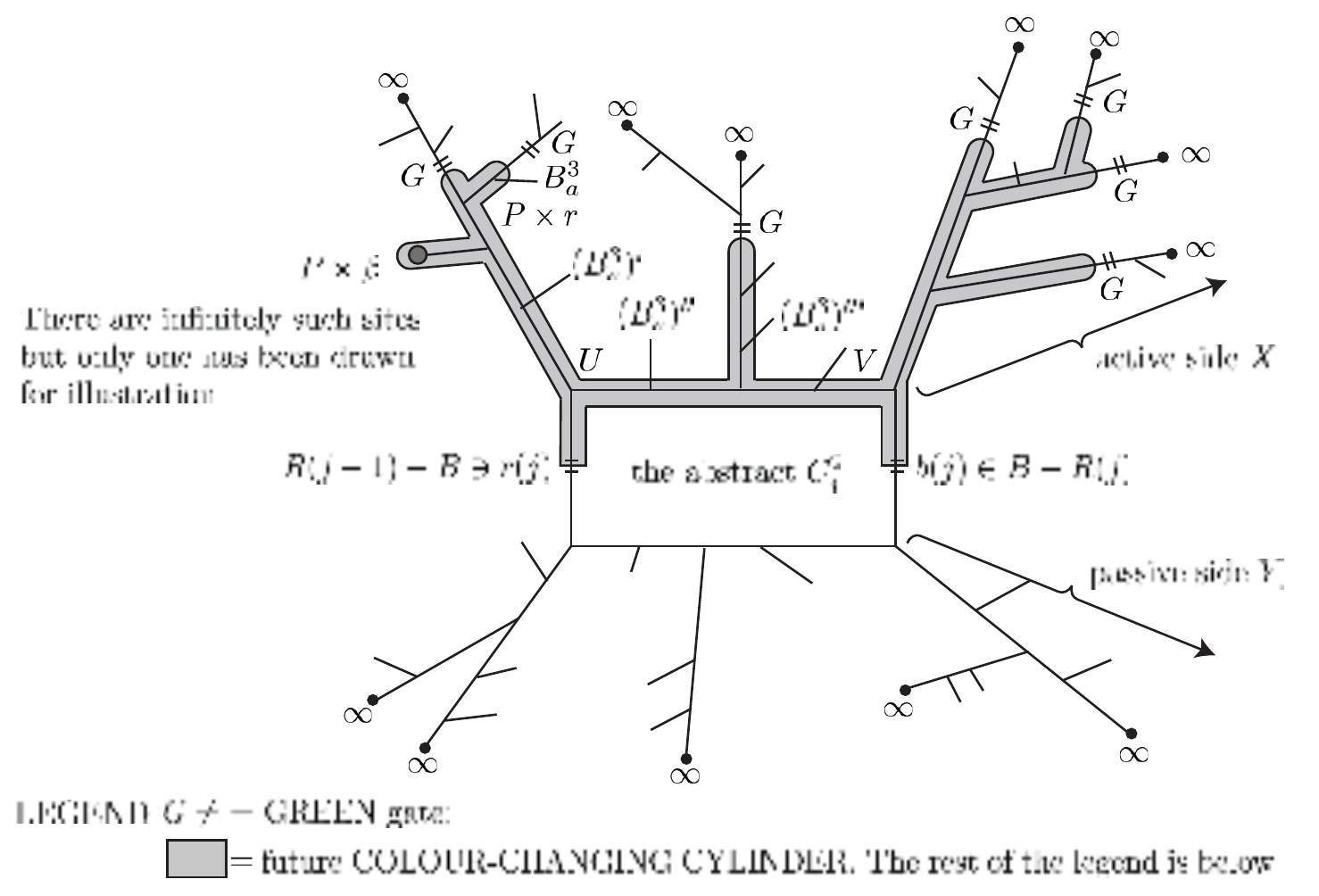}
$$
\label{fig46}
\centerline{\bf Figure 46.}
\begin{quote}
We see here the $\Gamma_1 (\infty) \left( t = \frac12 \right) = \Gamma_1 (\infty) [{\rm balanced}] \subset 2\Gamma (\infty) \left( t = \frac12 \right)$, which we will just denote by $\Gamma$ for typographical simplicity. At stage $j-1$ of the inductive process in the ABSTRACT COLOUR-CHANGING Lemma 14, realized now geometrically in $4^{\rm d}$, this is what we will find, just before $j-1 \Rightarrow j$ starts. So, what we have displayed, in the drawing above, is
$$
\{\mbox{the tree} \ \Gamma -R(j-1) = X_j \cup b(j) \cup Y_j \} \cup \{\mbox{the edge} \ r(j)\}.
$$
Remember the $R \cap B$ is mute, in particular so are the $\{{\rm little} \ h_i\}$'s.

So we need not worry about $R \cap B$ in the process
$$
\left( t = \frac12 \right) \Longrightarrow (t=1) \leqno (*)
$$
which implements geometrically the COLOUR CHANGE, leading to Lemma 20.

The $\Gamma (3)$, result of the BALANCING, is by now already with us and the $\Gamma (3) - R(j-1)$ lives deep inside the passive part $Y_j$, which is untouched by the action in $(*)$. It is exactly this presence of $\Gamma (3)$ inside one (and only one) of the $X$ or $Y$, which decides which is the active $X_j$ and which is the passive $Y_j$. The LEGEND continues now: ${-\!\!\!-\!\!\!-}$ ($=$ short RED arc) $=$ loose end of red arcs in $R(j-1) - \{r(j)\}$. Each such arc comes with two loose ends and by joining all these pairs together, we can reconstruct our $\Gamma (\equiv \Gamma_1 (\infty)[{\rm balanced}]$ at stage $j-1$); $\parallel\!\!\!\!\relbar \, =$ green GATE, to be explained in the main text. One has also drawn as a sample a $P \times [r,\beta]$, but there are many more such; they open the road going from our $\Gamma_1 (\infty)$ to the whole of $\Gamma (2\infty)$.
\end{quote}

\bigskip

Here are some additional EXPLANATIONS for Figure 46. We have finitely many sites denoted $G =$ ``GREEN GATE''; these are not in $R \cup B$, reason for a different COLOUR, but they should be compared with the RED and/or BLUE $1$-handles.

\smallskip

Here is what the finite system of GATES $\{G\} \subset X_j$ does for us.

\bigskip

\noindent (145.1) \quad It separates a compact part of $X_j$ from $\infty^{\rm ty}$. This compact part is the one touching to $b(j) , r(j)$ and the intermediary connecting part between them. Inside it, we will create the COLOUR-CHANGING CYLINDER.

\bigskip

\noindent (145.2) \quad The system $\{G\}$ is sufficiently close to $\infty^{\rm ty}$ so that we should have $\{$extended cocore $G\} \cap r(j) = \emptyset$.

\bigskip

One should use here the following general property which, for illustration, I state in terms of $N^4 (2\Gamma (\infty)) \subset N_1^4 (2X_0^2)^{\wedge}$. Let $S_n \subset 2\Gamma (\infty)$ be a site which is such that $S_n$ possesses an $\{$extended cocore $S_n\} \subset N^4(2X_0^2)$ and such that $\underset{n = \infty}{\lim} S_n = \infty$, in $\Gamma (2\infty)$. Then, we also find that
$$
\lim_{n = \infty} \{\mbox{extended cocore} \ S_n \} = \infty , \ {\rm in} \ N^4 (2X_0^2).
$$

Moreover, if one adds here all the appropriate quantifyers, there is then also a UNIFORMITY property too, and this yields what we want above. 

\smallskip

Let us be even a bit more precise. In the same vein as above, we can ask that, whenever there is a contact $\{G\} \cap C_{\ell} \ne \emptyset$ with, let us say, $C_{\ell} = C_{\ell} (r)$, then in the RED order one has $h_{\ell} > r(j)$. Next, we have

\bigskip

\noindent (145.3) \quad With $b(j)$ and $r(j)$ being given, the distinction between $Y_j$ being passive and $X_j$ being active depends on the $\Gamma (3) - R(j-1)$ being inside $Y_j$, which fixes the distinction. Then $\{G\}$ is located in $X_j$. With $r(j) + b(j)$ deleted, $X_j$ and $Y_j$ are no longer connected with each other.

\bigskip

\noindent (145.4) \quad We can choose our family $\{G\}$ sufficiently close to $\infty^{\rm ty}$ so that, we not only should have the (145.2) above, but also the following feature
$$
\{G\} \cap \Sigma \, \eta_{\ell} ({\rm green}) = \emptyset .
$$

We use here the fact that $\Sigma \, \eta_{\ell} ({\rm green})$ is finite. Also, each $G$ is a 1-handle not in $R \cup B$, of its new colour ``GREEN'', so that, like for the $B_a^3$ in the Figure 40 we have
$$
B^3 (G) = \frac12 \, B^3 (G) (+) \underset{\overbrace{\mbox{\footnotesize $\sigma^2$}}}{\cup} \frac12 \, B^3 (G)(-)
$$
and, a priori, we could have $\partial N_+^4 (2\Gamma (\infty)) \supset \frac12 \, N^3 (G)(+) \cap \eta_{\ell} ({\rm green}) \ne \emptyset$, something which we want to avoid, and {\ibf will avoid}, by chosing $\{G\}$ very close to $\infty^{\rm ty}$

\bigskip

\noindent (145.4-bis) \quad The $\{G\}$ induces a finite partition
$$
X_j = X_j^0 \underset{G}{\cup} X_j^1 \underset{G}{\cup} \ldots \underset{G}{\cup} X_j^{\rho}
$$
where each $X_j^{\varepsilon}$ is a tree, $X_j^0$ being finite and the others infinite.

\bigskip

This family $\{G\}$ serves to allow for a {\ibf finitistic} step $j-1 \Rightarrow j$; wihout it, the step in question would require an infinite process, all by itself. Now, in order to perform our inductive step $j-1 \Rightarrow j$, on the road from $t=\frac12$ to $t=1$, we will have to break it into three successive parts, namely the STEPS I, II and III, to be described below.

\bigskip

\noindent {\bf Important comments (145.4-ter).} At the level of Figure 46 all the 1-handles $R \cap B$ have been cut, like in the 1) of the ABSTRACT COLOUR-CHANGING LEMMA 14. This means that there is no contribution of $\{$Figure 45-(II) with $\{$little $C\}$ present$\}$, to our figure. The $\{$Figure 45-(III) with $\{$little $C\}$ (but without the $h_j (r) \in \{$little $h\} \in B \cap R\}$ is present, and we see its reflex in the upper side of the Figure 49. But, in Figure 46, the $R \cap B$'s, including our $\{$little $h\}$'s are mute.

\bigskip

\noindent {\bf The Step I.} In terms of the Figure 46, but now geometrically in $4^{\rm d}$, this consists in the following slidings of 1-handles over other $1$-handles, in succession:

\smallskip

a) We start with a first succession of elementary steps, for which we consider the various
$$
B_a^3 ({\rm LAVA}) \in \{\mbox{the short red arcs sticking out of $X_j^0$ (145.4-bis)}\}.
$$
These $B_a^3$'s are a finite family of 1-handles presenting themselves like th $B_a^3$ in the Figure 39-(II), except that contrary to that figure, these $B_a^3$'s are now LAVA.

\smallskip

The $B_a^3$'s have to slide over $B_p^3 = r(j)$ which is LAVA too.

\smallskip

b) Next, comes a second succession of elementary steps, where each individual $G = B_a^3$ (non LAVA) slides over the same $B_p^3 = r(j)$ (LAVA), as above.

\smallskip

Since the case b) is somewhat simpler than the case a), we will describe it first.

\bigskip

\noindent {\bf The slide of $G$ (non-LAVA) over $r(j-1)$ (LAVA).} This is very much like in the Figures 40 to 44, which can be used again, with the obvious change $y \Rightarrow G =$ our new $B_a^3$, $x \Rightarrow r(j) = $ our new $B_p^3$. Condition (124) which implied that $\{$extended cocore $y\} \cap x = \emptyset$ is replaced now by the (145.2). Otherwise, with our obvious changes of notation, this step is to be described by the same Figures 40 to 44 as before, and we leave it at that, as far as the detailed description of the step is concerned, except for the following COMMENT (Caveat!). The (126.4) is now no longer with us, so that in the Figure 40 for $G$, which we have not drawn, there could happily be BLUE spectators, like in the schematical Figure 39. But, because of (145.4) this {\ibf cannot} come with unwanted contacts $\eta_{\ell} ({\rm green}) \cdot \left( B - \underset{1}{\overset{P}{\sum}} \, H_i^b \right)$ which would contradict our little BLUE DIAGONALIZATION. [The $G$'s can be chosen sufficiently close to $\infty^{\rm ty}$ so that $(G \cup (X_j^1 \cup \ldots \cup X_j^{\rho})) \cap \Sigma \, \eta_{\ell} ({\rm green}) = \emptyset$.]

\bigskip

We will describe now the sliding of the finitely many $B_a^3 ({\rm LAVA}) \equiv \{$any of the small RED arcs shooting out of $X_j^0 \supset [U,V]$, in the Figure 46$\}$. Keep in mind that this slide of $B_a^3$ (LAVA) over $B_p^3 ({\rm LAVA}) = r(j)$ is to take place before the sliding of the Gates $G$ (non LAVA), over the same $B_p^3$ (LAVA). It is only for expository purposes that the slide of the $G$'s was discussed first.

\bigskip

\noindent {\bf The slide of $B_a^3$ (LAVA) over $B_p^3$ (LAVA) $= r(j)$.} At this point there is a big complication occurring in our story. Both $B_a^3$ and $B_p^3$ are  members of the family $\{ h \} \cap (R-B)$ and in the RED order, both cases
$$
B_p^3 > B_a^3 \qquad {\rm OR} \qquad B_a^3 > B_p^3
$$
are possible. Now, when $B_p^3 > B_a^3$, then the sliding of $B_a^3$ over $B_p^3$ brings about a VIOLATION of our so far sacro-sancted feature $C \cdot  h = {\rm id} + {\rm nilpotent}$, which we will have to drop now. Actually, from now on, we will have to live without it.

\smallskip

The condition $C \cdot  h = {\rm easy \ id} + {\rm nilpotent}$, is actually less sacro-sancted than confinement, a condition which will never get violated. We still want to perform the present step so that the conditions (126.1) and (126.2) should be respected.

\smallskip

The point here is that, as long as lava continues to have its STRONG PRODUCT COMPANY, we can define extended cocores using that, without invoking $C \cdot  h = {\rm id} + {\rm nilpotent}$.

$$
\includegraphics[width=135mm]{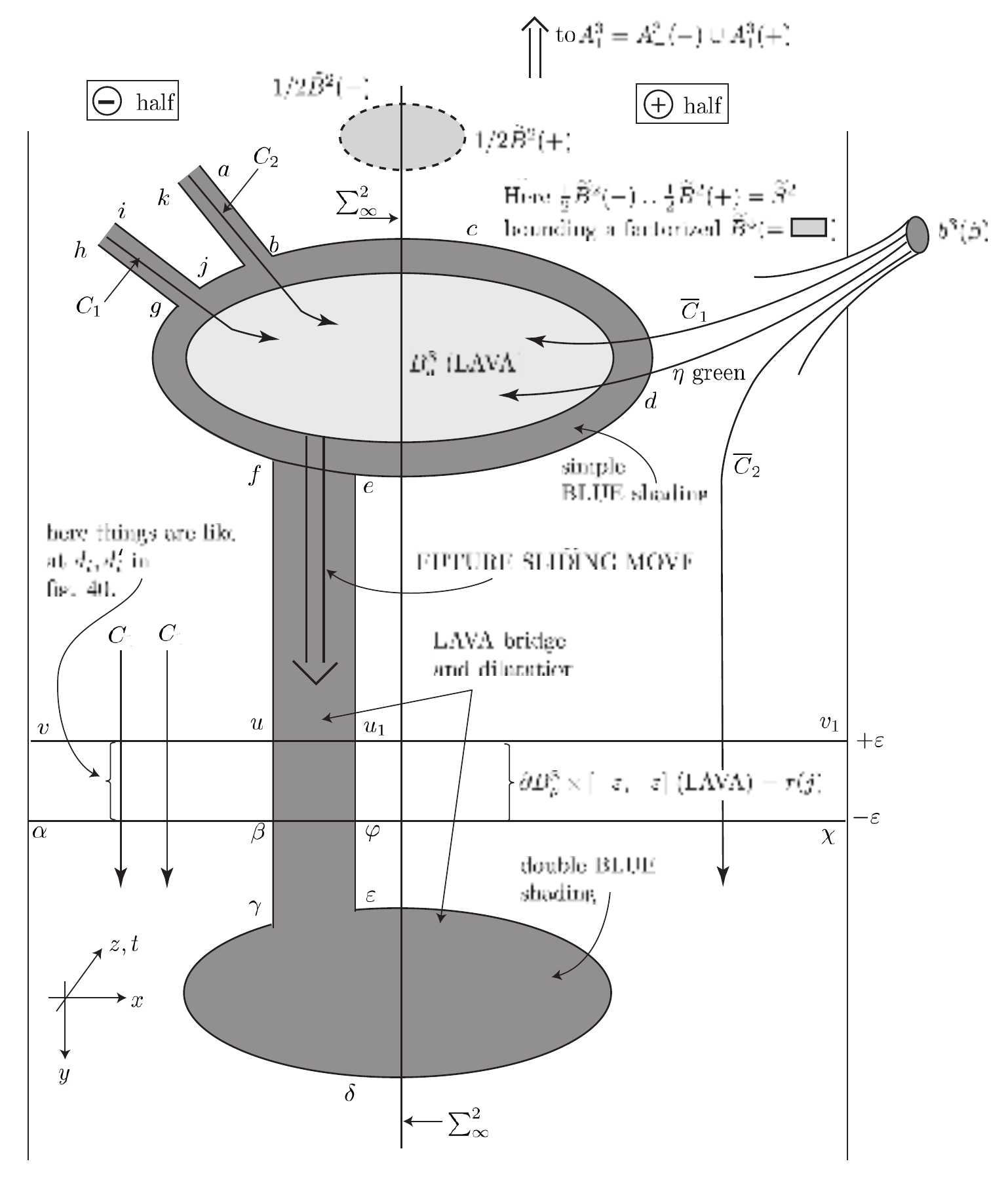}
$$
\label{fig47}
\centerline{\bf Figure 47.}
\begin{quote}
We see here the beginning of the slide of $B_a^3$ (LAVA) over the $r(j) = B_p^3$ (LAVA). For the sake of typographical simplicity we do not use here the realistic dimensions, like in Figure 40, with which this figure should be compared, but we proceed more schematically, with one dimension less. This figure is supposed to stand for a relevant part of $\partial N^4 (2\Gamma (\infty))$. More will be said concerning it in the main text. The $B_a^3 ({\rm LAVA}) \cup C_1 \cup C_2$ are drawn in RED.
\end{quote}

\bigskip

For our present sliding of $B_a^3$ (LAVA) we will proceed like in the Figure 40, but paying special attention to LAVA, which includes now $B_a^3$. Figure 40 is  replaced by the Figure 47. And the only things in this figure which are supposed to be alive, before our action starts, are the RED $B_a^3$ (LAVA), the passive $B_p^3 ({\rm LAVA}) \times [-\varepsilon , \varepsilon]$ and the various curves riding over them (only the ones riding over $B_a^3$ are actually drawn). Except for the $\eta ({\rm green})$'s all the other curves living on the $\ominus$ side (respectively on the $\oplus$ side) are part of ${\rm LAVA} (r) \left( t= \frac12 \right)$, respectively of ${\rm LAVA} (b) \left( t= \frac12 \right)$.

\medskip

Then a preliminary action starts (and we are here still at the level of Figure 47):

\medskip

a) The original red LAVA extends to the shaded blue area $[a,b,c,d,e,f,g,h,i,j,k]$ with the RED $B_a^3$ sitting now on top of it, dragging along the various curves like $C_1 , C_2 , \bar C_1 , \eta ({\rm green})$ which climb on $B_a^3 \times [-\varepsilon , +\varepsilon]$. At this point our contour $[a,b,\ldots ,k]$ above is in $\partial \, \delta \, {\rm LAVA}$, with $\delta \, {\rm LAVA}$ under the whole shaded area (with the shaded red ($\includegraphics[width=1cm]{rectangleclair.pdf}$) {\ibf sitting on top of} the shaded blue ($\includegraphics[width=1cm]{rectanglefonce.pdf}$)).

\medskip

[Here is the scenario for the upper part of Figure 47. In the beginning there is only the RED part. Then the blue $[j,b,a,d,e,g,j]$ is glued to it; and the red blob climbs on it, sliding. Then, finally along the curve $C_1 , C_2$ and UNDER THEM, we add the blue thin shaded zones all along their free parts. Importantly, all that red part is now disconnected from $\delta \, {\rm LAVA}$.]

\medskip

b) Next, like in Figure 40 a LAVA BRIDGE and LAVA dilatation grows out of the a) above. This is doubly shaded ($\includegraphics[width=1cm]{rectanglefonce.pdf}$) and when this doubly shaded area meets $[\partial B_p^3 \times [-\varepsilon , \varepsilon] ({\rm LAVA})$], it melts into it, like in Figure 40. To be very explicit, by now $\partial (\delta \, {\rm LAVA})$ has become the disjoined union
$$
\{ [a,b,c,d,e, u_1 , v_1] \cup [k,j,i] \cup [h,g,f,u,v]\} + \{[\alpha , \beta , \gamma , \delta , \varepsilon, \varphi , \chi]\}.
$$

Notice that there are no curves $\Gamma_j$ in the Figure 47. Normally they are in $Y_j$, out of the universe of our figure. But one way or another, this is not something to bother us. A piece of curve $\Gamma_j$ surviving inside the $X_j$ of Figure 46 and climbing on the $B_a^3$ (Figure 47), can be treated now on par with the $C_i$'s; there is now no GAP like in Figure 41-(B), which had forced a different treatment of the $\Gamma_j$'s with respect to the $C_i$'s, before.

\medskip

Now the real action of sliding $B_a^3 ({\rm LAVA})$ over $B_p^3 ({\rm LAVA}) = r(j)$, can take place, and Figure 48 refers to this slide.

\medskip

The important point is that, from the viewpoint of $({\rm LAVA} , \delta \, {\rm LAVA})$, the sliding of the $B_a^3 ({\rm LAVA})$ over $B_p^3 ({\rm LAVA})$ is now an {\ibf isotopic internal} to LAVA. So, although we have lost the feature $C \cdot h = {\rm id} + {\rm nilpotent}$, the LAVA retains its STRONG PRODUCT PROPERTY. Moreover, the (144.2) stays with us too.

\medskip

Very importantly, when moving from Figure 47 to Figure 48, the $\delta \, {\rm LAVA}$ stays unchanged. Very importantly too, the $\{$extended cocores $S\}$ (for ``site'' $S$), in particular the all-important $\{$extended cocore $H_{j\leq M}^b\}$, are {\ibf defined} from now on, not by appealing to $C \cdot h = {\rm id} + {\rm nil}$ (which is, in principle, no longer with us), but {\ibf directly to the} PRODUCT PROPERTY OF LAVA.

\medskip

Here is an additional explanation, concerning our Figures 47, 48.

\medskip

In the context of Figure 47 we may happily have BLUE spectators $B$, and contrary to that had happened for $G$, this may come (at the level of Figure 48) with {\ibf unwanted contacts}
$$
(\mbox{BLUE spectator $B$}) \cap \eta_{\ell} ({\rm green}) \ne \emptyset , \leqno (145.5)
$$
which are a temporary violation of the little BLUE DIAGONALIZATION, and to which we will have to come back later on.

\medskip

It is the STEP III which will take care of this potentially dangerous complication. In the same vein, our move may also introduce new contacts, very much unwanted too (see Figure 48-(A))
$$
(B_p^3 = r(j)) \cap \eta_{\ell} ({\rm green}) \ne \emptyset , \leqno (145.6)
$$
exactly as consequence of the $B_a^3 \cap \eta_{\ell} ({\rm green})$ from Figure 47.

$$
\includegraphics[width=125mm]{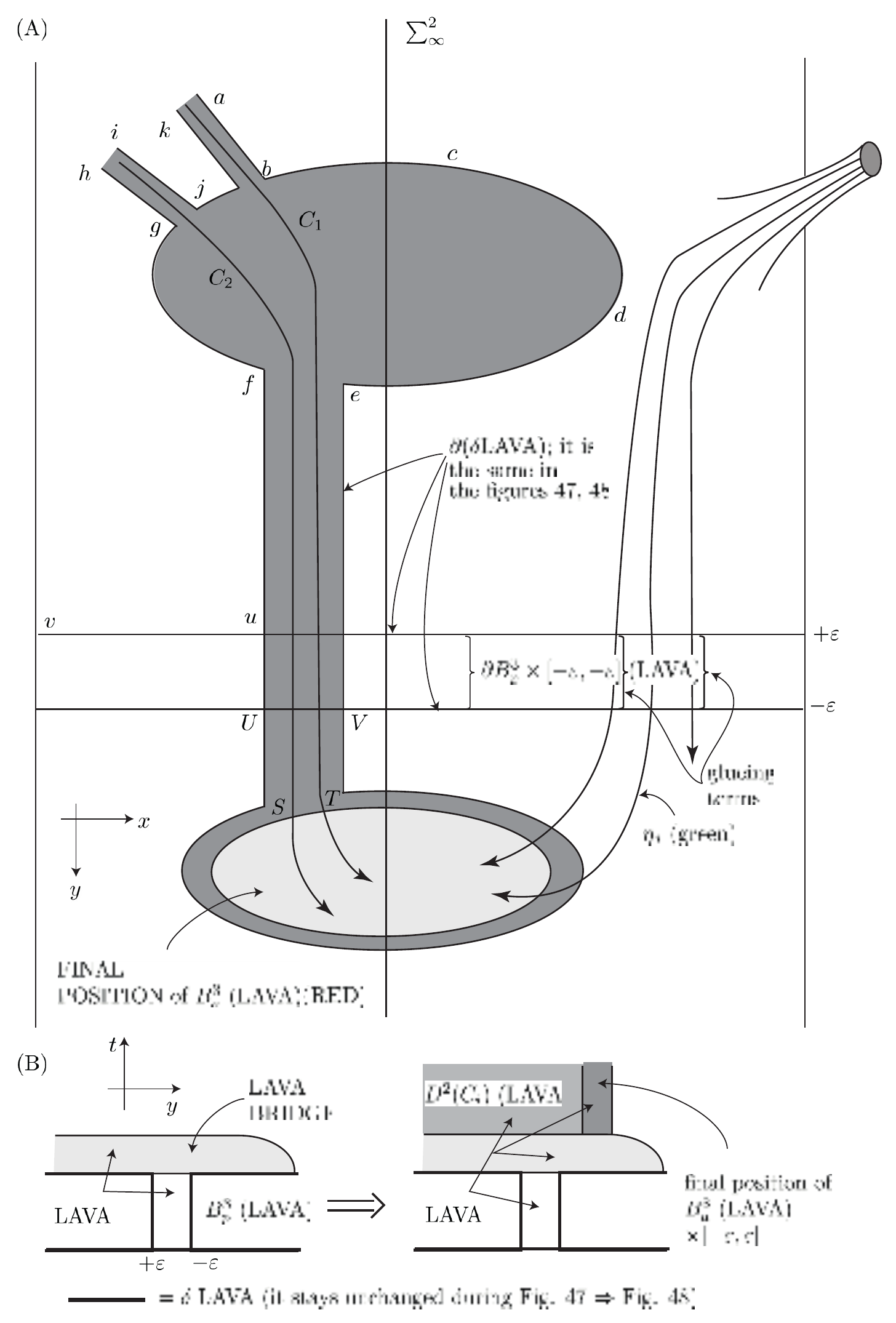}
$$
\label{fig48}
\centerline{\bf Figure 48.}
\begin{quote}
The end of the slide of $B_a^3 ({\rm LAVA})$ over $r(j) = B_p^3 ({\rm LAVA})$. This is an internal LAVA isotopic-move, making that the PRODUCT PROPERTY is automatically preserved. The $B_a^3 ({\rm LAVA})$, $D^2 (C_i)$ are all RED here.
\end{quote}

\bigskip

All this will be taken care of below. [The BLUE spectators mentioned above, and which are like the $b$'s in the Figure 39, come exactly from $B \cap X_j^0 - \{ b_j\}$, see here Figure 46 and formula from (145.4-bis).]

\bigskip

\noindent {\bf Step II.} As a result of STEP I, our $N^4 (2\Gamma (\infty))$ has been changed into a fluid object, the time-depending $Z^4 (t)$, where $j-1 \leq t \leq j$. The net result of our STEP I (first part of $j-1 \Rightarrow j$), is that at the time
$$
\{ t = j-2\varepsilon \} \equiv \{\mbox{the time $t$ when the STEP I above has been completed}\},
$$
the situation is the following. Once the $B_a^3$'s, i.e. the small RED arcs which stick out of $X_j^0$ in the Figure 46, and the green gates $G$ too, have slided over $r(j)$, one can put together a COLOUR-CHANGING CYCLINDER
$$
B^3 \times [r,b] \subset Z^4 (t=j-2\varepsilon), \leqno (145.7)
$$
connecting $B^3 (r(j))$ to $B^3 (b(j))$ and which is factorized too, i.e. divided naturally into $\xy *[o]=<12pt>\hbox{$\pm$}="o"* \frm{o}\endxy$-halves. Figure~49 presents the lateral surface of this cyclinder
$$
S^2 \times [r,b] = \left( \frac12 \, B^2(+) \times [r,b] \right) \underset{\overbrace{\mbox{\footnotesize$\Sigma_{\infty}^2$}}}{\cup} \left( \frac12 \, B^2(-) \times [r,b] \right). \leqno (146)
$$
In the Figure 46, at its schematical $1^{\rm d}$ level, our (145.7) corresponds to that piece of the graph, resting on $[r(j) , U , V , b(j)]$ and stretching all the way to the gates $G$. To make things completely clear, this piece, producing the COULOUR-CHANGING cylinder is surrounded by the shaded area with its boundary $\includegraphics[width=3cm]{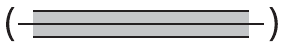}$. After the sliding of the $B_a^3$'s and the $G$'s, this is separated now from the rest of what has become out of the graph $\Gamma$ in Figure 46 at time $t=j-2\varepsilon$, by exactly $r(j)$ and $b(j)$. The red arcs which in Figure 46 stick out of the shaded area, are the $B_a^3 ({\rm LAVA})$'s which have already slided over $r(j)$.

\smallskip

Our COLOUR-CHANGING CYLINDER is made exactly out of the following spare parts, which in a $1^{\rm d}$ version are vizualizable in Figure 46.

\bigskip

\noindent {\bf List of spare parts of the colour-changing cylinder:}

\smallskip

i)  The edge $[r(j) , U]$; ii) The edge $[U,V] \cup \{$all the part of $X_j$ between the $[U,V]$ and the GATES $G$ (see Figure 46)$\}$, iii) The edge $[V , b(j)]$. What we have just unrolled corresponds to the $X_j^0$ in (145.4-bis) and here both the $B_a^3 ({\rm LAVA})$'s and the GATES $G$ have already slided away, leaving behind them ghostly terms $\# \, B^3 (G)$, $\# \, B_a^3$, suggested in the Figure 49. There are very much like the $A_1^3 (-) \underset{D^2}{\cup} A_1^3 (+)$ in Figure 40.

\smallskip

With this our list of spare parts in now closed.

\bigskip

The $\# \, B^3 (G)$, $\# \, B_a^3$ mentioned above are factorized as follows
$$
\# \, B^3 ({\rm ghost}) = \frac12 \, B^3 ({\rm ghost})(-) \cup \frac12 \, B^3 ({\rm ghost}) (+) ,
$$
and see here the Figure 49. Here are some facts to be kept in mind.

\bigskip

\noindent (146.1) \quad By construction, the COLOUR-CHANGING cylinder $B^3 \times [r,b]$ does NOT communicate with the outer world of $Z^4 (t=j-2\varepsilon)$ through the ghostly $\# \, B^3$'s. And then, in the same vein,

\bigskip

\noindent (146.2) \quad NO curves in the cylinder stay hooked at the $\# \, B^3 ({\rm ghost})$ (something which would be quite an inconvenience for performing our next STEP II (the actual COULOUR-CHANGE)).

\bigskip

\noindent {\bf Explanation:} Look at Figure 39. At the level of Figure 46 we are like in the Figure 39-(I), with a green $C$ going through $B_a^3$, which is like our present moving $1$-handles $G$ (NON-lava) or $B_a^3$ (LAVA), moving over $r(j)$. Once the 1-handle slides, and we move from Figure 39-(I) to the Figures 39-(II and III), the $A_1$ is now like our $\# \, (G \ {\rm or} \ B_a^3)$ (ghost), and one can see in these last two Figures that there is now no curve going through the ghostly sites $A_1$. Anyway, Figure 39 can serve here as an illustration. The green curve is not hooked at $A_1$ in (II), (III). A more complex choregraphy of sliding handles and curves dragged along could be put up, corresponding to the real-life situation of all the $\# \, B^3$ (ghost). End of (146.2).

\medskip

\noindent (147) \quad But the colour changing cylinder $B^3 \times [r,b]$ does communicate with the outside world, not only through its ends, but also through the following sites: $b^3 (\beta)$ which is {\ibf not} factorized, and which lives completely on the $\oplus$ side, the $b^3 (Q,Q')$ which are factorized, but only the $\oplus$ is explicitly drawn in the Figure 49. Then $b^3 (Q,Q')$'s are attaching zones of short $1$-handles like in Figure 32, cut along their $h \in R \cap B$ which are mute in our very present story. The $\{{\rm little} \ C_i\}$ which is seable in the Figure 49 should be compared with the red line in Figure 45-(III). On the $\oplus$ side we find the $\underset{j}{\sum} \, c(r_j) + \underset{i}{\sum} \, \{{\rm little} \ C_i\} \subset \frac12 \, B^3 (+) \times [r,b]$, which can get out of the colour-changing cylinder only through the sites of type $b^3 (\beta)$, $b^3 (Q \ {\rm or} \ Q')$. End of (147).
$$
\includegraphics[width=140mm]{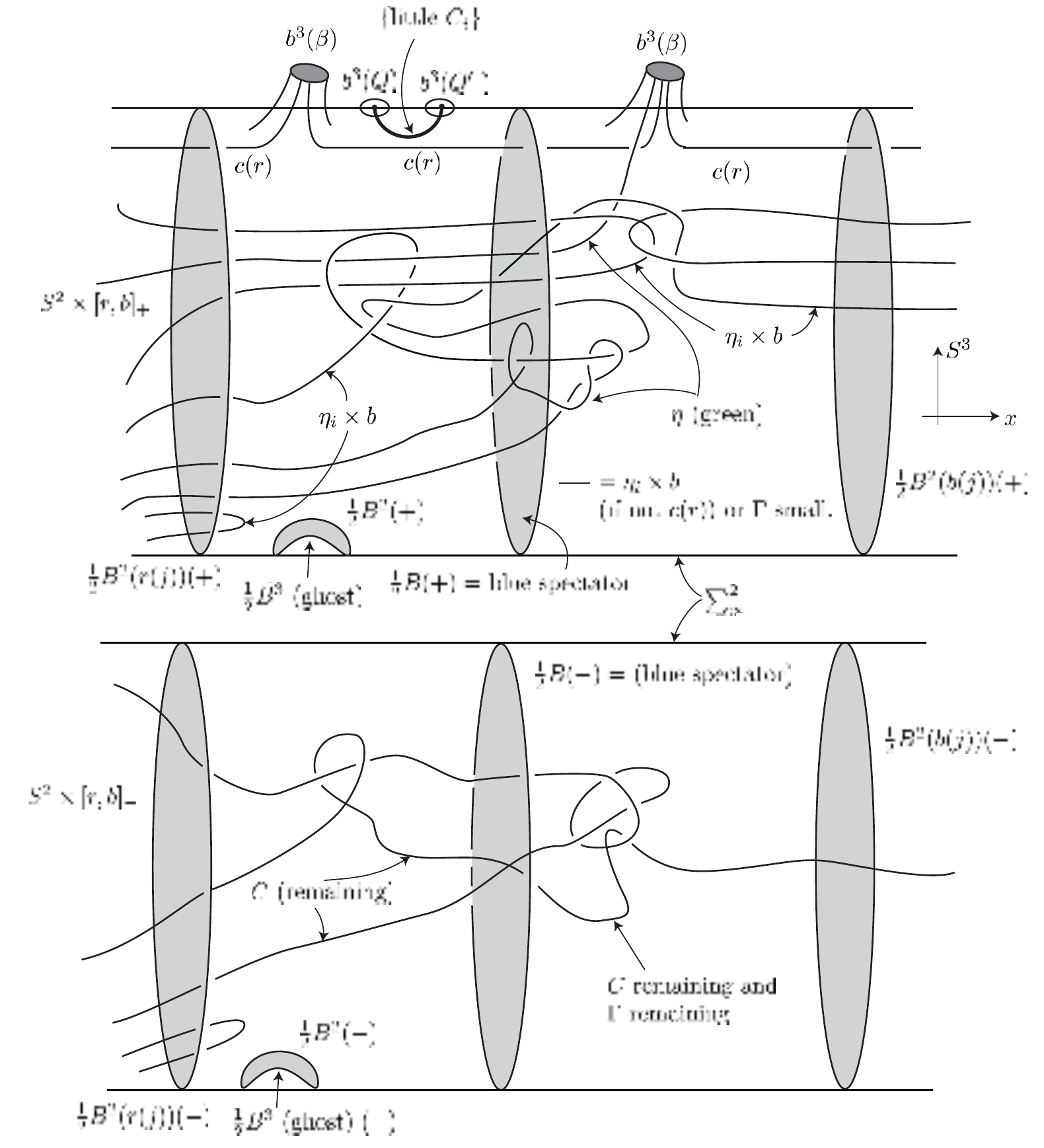}
$$
\label{fig49}
\vglue-7mm
\centerline{\bf Figure 49.}

\begin{quote}
The COLOUR-CHANGING CYLINDER. In this cylinder the curves in the $(\pm)$-halves $S^2 \times [r,b]_{\pm}$ are {\ibf not entangled} with each other, and therefore they can be treated independently. The contacts $(\Sigma \, c(r) + \Sigma \, \eta_i \times b)$ $\cap \, \frac12 \, B^2 (r(j)(+))$ are GLUEING TERMS, like in (144.3). [The $h_j (r) \approx B^3 (r(j))$ occurs in Figure 45-(I), where we can also see its GLUEING TERM with $c(r)$. The glueing term with $\eta_i \times b$ has been drawn for pedagogical purposes.] The GLUEING TERMS may be severed and reglued somewhere else on $\partial \, {\rm LAVA}_r$, without any harm.

Also, inside $S^3 \times [r,b] (-)$ we have LAVA bridges, which we have not drawn here. They will follow the general fate of the curves. Here $1/2 \, B^2 (b(j))(-) \cup 1/2 \, B^2 (b(j))(+) = \partial B^3 (b(j))$.
\end{quote}

\bigskip

In terms of the coordinate system from Figure 49, our STEP II, which is part of the bigger step $j-1 \Rightarrow j$, is a translation move to the right, expressible as $x \longmapsto x + [r-b]$, which, in a factorized manner, meaning operating in such a way that the $\pm$-SPLITTING should be respected, {\ibf superposes} in Figure 49 the $\frac12 \, B^3 (r(j))_{\pm}$ on top of $\frac12 \, B^3 (b(j))_{\pm}$, i.e. $S^3 (r(j))$ on top of $S^3 (b(j))$. At the level of $Z^4 (j-2\varepsilon) \Rightarrow Z^4 (j-\varepsilon)$ (with $t=j-\varepsilon$ the time when the COLOUR-CHANGING step is completed), the equality of sets $r(b) = b(j) \subset Z^4 (j-\varepsilon)$ has been established. The ``site'' $r(j)$ will be forgotten from now on. In more detail, we proceed as follows as far as this ``isotopic move of 1-handle cocores'' is concerned.

\smallskip

If the word ``isotopic'' occurs here between quotation marks, it is that during the STEP II, the $B^3 (r(j))$ which moves, inside the COLOUR-CHANGING CYLINDER, into the position $B^3 (b(j))$, sweeps through $b^3 (\beta)$, $b^3 (Q,Q')$ and the $B^3 ({\rm ghost})$. So, when all the rest of the world of $\Gamma (2\infty)$ is taken into account, our STEP II is certainly NOT  {\ibf  an isotopy}, of anything. It is only a non-isotopic geometric transformation corresponding to the abstract equality $\Phi (r(j)) =  b(j)$ in the Lemma 14. What our STEP II achieves is to force the equality of sets $B^3 (r(j)) = B^3 (b(j))$, keeping the $B^3 (b(j))$ in place as it is. Next, we also have the following items.

\medskip

\noindent $\bullet$) We do not budge anything else on the $S^3 \times [r,b] (+)$ side, where all the curves stay put, as they are, and the ${\rm LAVA}_b$ too. The only move on the $\oplus$ side is the $r(j) \Rightarrow b(j)$ itself.

\medskip

\noindent $\bullet$$\bullet$) [{\bf Important Remark:} Of course, our transformation $r(j) \Rightarrow b(j)$ is, from the viewpoint of the mere COLOUR-CHANGING CYLINDER, Figure 49 is simple-minded isotopy of 1-handle cocores, moving $B^3 (r(j))$ to the right until it is superposed on $B^3 (b(j))$. But since this {\ibf sweeps} through the sites $b^3 (\beta)$, $\# \, B^3 ({\rm GHOST})$ it cannot be in any way an isotopy. It is a way to put flush and bones on Lemma~14, actually on the $\Phi (r(j)) = b(j)$, and we will {\ibf decree} that $B^3 (r(j))$ leaves LAVA, while $B^3 (b(j))$ enters it. But the metaphor of translation move to the right, is very convenient for describing the fate of the objects in the $\pm$ sides of the colour-changing cylinder.] With this, on the $S^3 \times [r,b] (-)$ side, all the $\{$curves$\} + (\mbox{LAVA BRIDGES} \ \subset {\rm LAVA}_r)$ {\ibf move solidarily} with $\frac 12 \, B^3 (r(j))(-)$ and hence get glued to $\frac12 \, B^3 (b(j))(-)$ at the end of the process, and they stick there from there on. This movement is made possible, without any obstructions, because the ghostly $B^3 (-)$'s are untouched by the curves; this has been explained in (146.2). Notice, also, that the conjunction of $\bullet$) and $\bullet$$\bullet$) is made possible by the combination
$$
\mbox{SPLITTING $+$ CONFINEMENT}.
$$

The topologies of both the individual LAVA$_r$ and LAVA$_b$ stay intact, but we may sever GLUEING TERMS like $((\eta_i \times b) + c(r)) \cap \frac12 \, B^3 (r(j))(+)$, where $\frac12 \, B^3 (r(j))(+) \subset {\rm LAVA}_r$, {\ibf decreeing} that $((\eta_i \times b) + c(r)) \cap \frac12 \, B^3 (b(j))(+)$ are now NEW GLUEING TERMS. Also, now the $B^3 (b(j)) = \frac12 \, B^3 (b(j)(-)) \cup \frac12 \, B^3 (b(j)(+))$ becomes LAVA$_r (t=j-\varepsilon)$, by decree too. Here, of course, as already said, $\{ t = j- \varepsilon \} \equiv \{$The time when STEP II inside $j-1 \Rightarrow j$ has been completed$\}$.

\smallskip

In our move on the $S^3 \times [r,b](+)$ side, we do not worry about severing connections ${\rm LAVA}_b \cap B^3 (r(j))$, these are just freely movable GLUEING TERMS ${\rm LAVA}_b \cap {\rm LAVA}_r$.

\smallskip

It should be stressed that the transformation $r(j) \Longrightarrow b(j)$ operates at the level of the $1$-handle cocores themselves, and not at the level of $\{$extended 1-handle cocores$\}$.

\smallskip

Here the $r(j)$ may be the one in Figure 48 and one should also keep in mind that in the context of the Figure in question, the $C \cdot h = {\rm id} + {\rm nilpotent}$ might be bumped. But then, also, the $r(j)$ drags along all its $\{{\rm LAVA}_r\} + \{$the corresponding LAVA bridges$\}$, and puts them on top of $b(j)$ which, by degree, becomes itself LAVA.

\smallskip

The ${\rm LAVA}_r$ stays unchanged, up to homeomorphism and hence it keeps its PRODUCT structure intact. The $\frac12 \, B^3 (r(j))(+)$ has never contributed to the ${\rm LAVA}_b$ and hence the $\frac12 \, B^3 (b(j))(+)$ will not contribute to it either. But these $B^3$'s contribute to ${\rm LAVA}_r$ of course, at the appropriate times. [And, at the level of Figure 49, inside the $S^2 \times [r,b]_+$, the only contributions to ${\rm LAVA}_r$ can only come from $\frac12 \, B^3 (r(j))(+)$ OR from $\frac12 \, B^3 (b(j))(+)$; they do not occur simultaneously, of course.] 

\smallskip

The ${\rm LAVA}_b$ stays intact and then, because of Lemma 18, the global LAVA continues to have its product structure too. Here are other features of our step, as described above.

\bigskip

\noindent (147.1) \quad  No new contacts $b(j) \cdot \underset{1}{\overset{M}{\sum}} \, \eta_{\ell} ({\rm green})$ are introduced. [Remember also that the $H_{j \leq M}^b$ with their diagonal contacts are all in the passive part $Y_j$ of Figure 46.]

\bigskip

\noindent (147.2) \quad  The set of the parasitical $h_k$'s in (144.1) can only decrease, i.e. we have no new $R-B \ni h_i \subset {\rm LAVA} \, (t)$, touched by $\Sigma \, \eta_{\ell} ({\rm green})$, with respect to what we had already before.

\bigskip

\noindent (147.3) \quad  With our $\{$extended cocores$\}$ simply defined now by appealing to the PRODUCT PROPERTY OF $\{{\rm LAVA} , \delta \, {\rm LAVA}\}$, and see here the all-important (144.2) too, any possibe inclusion
$$
B^3 (r(j)) \subset \sum_{j=1}^M \{\mbox{extended cocore} \ H_j^b \} ,
$$
has been replaced by a harmless inclusion
$$
B^3 (b(j)) \subset \sum_{j=1}^M \{\mbox{extended cocore} \ H_j^b \} .
$$
The point here is that while with $B^3 (r(j)) \subset \{\mbox{extended cocore} \ H_j^b \} $ can come with its dangerous contacts $\eta_{\ell} ({\rm green}) \cdot B^3 (r(j)) \subset \eta_{\ell} ({\rm green}) \cdot \{\mbox{extended cocore} \ H_j^b \}$, the LITTLE BLUE diagonalization forbids the contacts
$$
\eta_{\ell} ({\rm green}) \cdot B^3 (b(j)) \ne \emptyset.
$$

\noindent (147.4) \quad  Since on the $S^3 \times [r,b](+)$-side of the COLOUR-CHANGING cylinder all the curves have to stay put, the undesirable (145.6) does not create any new unwanted terms
$$
\eta_{\ell} ({\rm green}) \cdot b(j) \subset \eta_{\ell} ({\rm green}) \cdot \{\mbox{extended cocore} \ H_j^b\}.
$$
So the (145.6) has been taken care of, more precisely the (145.6) becomes {\ibf irrelevant} once $r(j)$ has been replaced by $b(j)$.

\bigskip

\noindent (147.5) \quad  But there is also a piece to be paid for all this. All the $B \cap X_j^0$ (see (145.4-bis) and Figure 46) are prospective passive BLUE spectators $b$ for our move $\{B_a^3 ({\rm LAVA})$ SLIDES over $B_p^3 ({\rm LAVA})\}$, coming with unwanted contacts like (145.5). And these $b$'s could happily be among the $\{$past or future $b(j)$'s$\} \cap \underset{1}{\overset{M}{\sum}} \, \{$extended cocore $H_j^1\}$. So the (145.5) {\ibf has} to be demolished. This will be achieved by STEP III with which our $j-1 \Rightarrow j$ ends. This will also enforce the important condition, part of the LITTLE BLUE DIAGONALIZATION
$$
\eta_{\ell} ({\rm green}) \cdot \left( B - \sum_{j=1}^M H_j^b \right) = 0 .\leqno (147.6)
$$

At this point, we will ANTICIPATE a bit. The STEP III inside our big $j-1 \Rightarrow j$, soon to be described explicitly, achieves among other things the following items:

\medskip

$\bullet$) It destroys the (145.5) which the STEP I inside the same $j-1 \Rightarrow j$ has created.

\medskip

$\bullet$$\bullet$) It stays far from $b(j)$ and does not bring any $\eta_{\ell} ({\rm green})$ on top of the $b(j)$ in question.

\medskip

$\bullet$$\bullet$$\bullet$) It preserves the PRODUCT PROPERTY of LAVA and also the (144.2). Hand in hand with this, comes the next item.

\medskip

$\bullet$$\bullet$$\bullet$$\bullet$) It does not increase the set of parasitical $h_k \in R-B$ touched by the $\eta ({\rm green})$ (144.1).

\medskip

At the beginning of each step $j-1 \Rightarrow j$ we assume, inductively, that (147.6) is satisfied, then with the things said it is also verified at the end of $j-1 \Rightarrow j$.

\smallskip

With this, assuming the STEP III inside $j-1 \Rightarrow j$ has all the stated virtues, when we get to time $t=1$ the (147.6) is with us. Also, the totality of the STEPS II realizes condition (144.1) and hence, also, all the parasitical RED terms in the geometric intersection matrix
$$
\eta_i ({\rm green}) \cdot \{\mbox{extended cocore} \ H_j^b \} \, , \quad \mbox{when} \ 1 \leq i,j \leq M
$$
have disappeared. But to get from this to the desired grand BLUE diagonality
$$
\eta_i ({\rm green}) \cdot \{\mbox{extended cocore} \ H_j^b \} = \delta_{ij} ,
$$
we certainly need the (147.6) too. This will make sure that there are no $b(j-1) \subset \{$extended cocore $H_{i \leq M}^p\}$, coming with $\eta_{\ell} ({\rm green})$ on top.

\bigskip

\noindent {\bf Step III.} The aim of the present step is two-fold: we want to restore $N^4 (\Gamma (2\infty))$ as it was before Step I i.e. we want to find back our Figure 46, as it stood initially, but of course now with $r(j-1)$ a mere ``site'' without any other geometrical meaning. Our LAVA $(j-\varepsilon)$ contains now the $B^3 (b(j))$ in lieu of $B^3 (r(j))$. And then, we also want to repair the damage (145.5). For that purpose we perform now essentially, the Step~I {\ibf in reverse}, i.e. schematically
$$
\mbox{The} \ (t=j-\varepsilon) \underset{{\rm STEP \, III}}{\xRightarrow{ \qquad \qquad }} (t=j) \ \mbox{is, essentially, STEP III $=$ (STEP I)$^{-1}$}.
$$

All this is supposed to happen without undoing the colour change $r(j) \Rightarrow b(j)$. Since the Step I stayed far from $b(j)$, the same is true for STEP III.

\smallskip

Here are more details concerning out STEP III. Although, schematically speaking we perform $({\rm STEP})^{-1}$, the $r(j)$ is now no longer a $1$-handle, it is just a ``spot''. So, unlike what has happened in the STEP I, we no longer have now slidings of some interesting 1-handles over some other interesting 1-handles, like the kind of things which were displayed in the Figures 40 to 44. As far as the handles are concerned, namely the $B_a^3 ({\rm LAVA})$'s and the $G$'s, we have now a simple isotopic move without any incidence at the level of the geometric incidence matrices like $C \cdot h$. There are no moves of 1-handles sliding over other 1-handles now.

\smallskip

We proceed again in a factorized manner and we want both to destroy the (145.5) and to make sure that the PRODUCT PROPERTY of LAVA continues to be with us.

\medskip

i) When we are on the $\oplus$ side, then together with the moving 1-handles $B_a^3 ({\rm LAVA})$, we drag back to their initial (pure Step I) position the external curves $\eta_{\ell} ({\rm green})$. This way we dispose of the unwanted (145.5), which disappears, in this process. On the same side $\oplus$ we also have other contacts $\{$curves$\} \cap (\Sigma \, B_a^3 ({\rm LAVA}) + \{{\rm Gates} \ G\})$, with the $B_a^3 ({\rm LAVA}) + \{{\rm Gates} \ G\})$ attached now on the other side of our spot $r(j)$, outside of the COLOUR-CHANGING CYLINDER), and these contacts get dragged together with the $B_a^3 ({\rm LAVA}) + G$, back to their initial positions too.

\medskip

The Step III happens far from $b(j)$, on top of which it leaves intact, all the things, curves and LAVA bridges, brought by STEP II.

\smallskip

With all this, ${\rm LAVA}_b$ keeps its PRODUCT PROPERTY. We may change around glueing terms, but that is O.K. Of course, in this STEP III, the $\{G\} + B_a^3 ({\rm LAVA})$ slide  back over the {\ibf site} $r(j)$, to get back to their original positions.

\smallskip

On the $\ominus$ side we have both curves $(\{C ({\rm remaining})\}$'s) and LAVA bridges. STEP II has transformed the contacts
$$
(\{{\rm curves} \, ({\rm remaining})\} + \mbox{LAVA BRIDGES}) \cap r(j) \leqno (148)\mbox{-}r(j)
$$
to isomorphic contacts 
$$
(\{{\rm curves} \, ({\rm remaining})\} + \mbox{LAVA BRIDGES}) \cap b(j). \leqno (148)\mbox{-}b(j)
$$

ii) When we are on the $\ominus$ side and we perform our isotopic move of bringing the $B_a^3 ({\rm LAVA}) + G$ back to their original position, we solidarily displace, by dragging along OR pushing in front in a snake-like manner, like in Figure 51-(A) $\Rightarrow$ 51-(B), the curves and LAVA BRIDGES, so that the contacts $B_a^3 ({\rm LAVA}) \cap ({\rm curves} + \mbox{LAVA bridges})$ (and similarly for $G$, but we do not really care about that) should be isotopically preserved, without bumping (148)-$b(j)$. This way, the ${\rm LAVA}_r$ continues to have its STRONG PRODUCT property.

\smallskip

All this means that the global lava continues to have its product property.

\bigskip

\noindent {\bf Final Remark.} In all this story, the BLUE spectators stay passive, curves do not stick to them, even when touching them. [The ``sticking'' refers here to connections, inside ${\rm LAVA}_r$ or ${\rm LAVA}_b$, which when undone destroy the PRODUCT PROPERTY.] Also, and importantly, the curves $\eta_{\ell} ({\rm green})$ do not touch the blue spectators any longer, at the end of the STEP III. It can be checked that the STEP III, as described above, has the features $\bullet$), $\bullet$$\bullet$), $\bullet$$\bullet$$\bullet$), $\bullet$$\bullet$$\bullet$$\bullet$), listed immediately after the formula (147.6).

\bigskip

\noindent {\bf Lemma 20. (The time $t=1$ compactification.)} {\it Proceeding with steps $(t=j-1) \Rightarrow (t=j)$ for all the successive times in $(144.3)$ and with each of these $j-1 \Rightarrow j$ divided itself into its intermediary {\rm STEPS I, II, III} on the lines just described, we have inside $N_1^4 (2X_0^2)^{\wedge}$ a transformation
$$
\left(N^4 (2\Gamma (\infty)) \right) \left( t = \frac12 \right) , N^4 (\Gamma (3)) \Longrightarrow (N^4 (2\Gamma (\infty)(t=1) , N^4 (\Gamma (3)), \leqno (148) 
$$
which is the identity on $N^4 (\Gamma(3))$, with the following features.}

\medskip

1) {\it This transformation drags along the link from {\rm (101.1-bis)}, both the internal and the external curves, and also the LAVA bridges too; hence it comes with a transformation which conserves the PRODUCT PROPERTY}
$$
{LAVA} \left( t = \frac12 \right) \Longrightarrow {LAVA} \, (t=1) . \leqno (149)
$$

2) {\it The analogue of $(137)$ is valid at $t=1$ too, i.e. we have
$$
\left\{\underbrace{\left[ N^4 (2\Gamma (\infty))(t=1) - \sum_1^{\infty} h_n (t=1) \right] \cup \mbox{LAVA} (t=1)^{\wedge}
}_{\hat Z (t=1)}
\right\} \sum_1^M \{\mbox{extended cocore} \ H_i^b \}^{\wedge} =
\leqno (150)
$$
$$
\underset{\rm DIFF}{=} \ \underset{i=1}{\overset{M}{\#}} \ (S_i^1 \times B_i^3 , (*) \times B_i^3).
$$

The analogue of $(136)$ is also valid at $(t=1)$,
$$
\Delta^4_{\rm Schoenflies} \underset{\rm DIFF}{=} N^4 (2X_0^2)^{\wedge} (t=1) = \left[ N^4 (2\Gamma (\infty))(t=1) - \sum_1^{\infty} h_n (t=1) \cup \mbox{LAVA} (t=1)^{\wedge} \right] +
\leqno (151)
$$
$$
+ \sum_1^{\overset{=}{n}} D^2 (\Gamma_j) \subset N^4_1 (2X_0^2)^{\wedge} (t=1) = N^4 (2X_0^2)^{\wedge} (t=1) \cup \partial (N^4 (2X_0^2)^{\wedge} (t=1)) \times [0,1] \underset{\rm DIFF}{=} \Delta^4 \cup (\partial \Delta^4  \times [0,1]).
$$
If we disregard all the special subtelties of the COLOUR-CHANGING process, one goes from $(136)$ to $(151)$ by ambient isotopy and, moreover, under the transformation $\left( t = \frac12 \right) \Rightarrow (t=1)$, the $\underset{1}{\overset{\overset{=}{n}}{\sum}} \, D^2 (\Gamma_j)$ does not change otherwise than being dragged via the covering isotopy theorem. These things just said, should be enough for the proof of $(151)$.}

\medskip

3) {\it We have now, for $1 \leq i,j \leq M$ the BIG BLUE DIAGONALIZATION:
$$
\eta_i (\mbox{green}) \cdot \underbrace{\{\mbox{extended cocore} \ H_j^b \}^{\wedge}}_{\mbox{\footnotesize as defined by the LAVA $(t=1)^{\wedge}$}} = \delta_{ij} \leqno (152)
$$
and, of course, we continue to have also
$$
\eta_i (\mbox{green}) \cdot \left( B - \sum_1^M H_j^b \right) = 0 .
$$
}

Combining this time $(t=1)$-Lemma with the very first Lemma 3 we have our proof that $N^4 (\Delta^2)$, albeit now presented as
$$
N^4 (\Delta^2) \underset{\rm DIFF}{=} \left[ \left( N^4 (2\Gamma (\infty))(t=1) - \sum_1^h h_i (t=1)\right) \cup {\rm LAVA} \, (t=1)^{\wedge} \right] + \sum_1^{\overset{=}{n}} D^2 (\Gamma_j)
$$
is {\ibf geometrically simply-connected}. So our Theorem 2, and then Theorem 1 too, are by now proved.

\bigskip

\noindent {\bf Additional explanations concerning the change $\left( t = \frac12 \right) \Rightarrow (t=1)$.}

\medskip

A) In the context of the Figure 46, we certainly have $\Gamma (3) - R(j) \subset Y_j$, but this is not necessarily true for the $\underset{1}{\overset{\overset{=}{n}}{\sum}} \, \Gamma_j$ (see here the PROMOTION TABLE, after the (116)). On the other hand the $\{{\rm little} \, \Gamma_j\}$'s may happily occur like the $\{{\rm little} \, C\}$'s in Figure 49 and this is without consequence, but, in the context of the $\{\{$Figures 40, 44$\}$ for $B_a^3 = G (\mbox{non LAVA})\}$ the $\Sigma \, \{\Gamma_j \ {\rm remaining}\}$ may well occur too. We treat them like in the usual Figures 40 to 44, they live on the $\ominus$ side, and they do not stick on any of the 1-handles which are now involved. That means that they can be just ignored as far as LAVA is concerned. Otherwise we drag them along and then finally put them back as they were at $t = \frac12$, together with the other curves.

\medskip

B) Clearly, the $B_a^3 ({\rm LAVA})$ which has to stick over $r(j)$(LAVA) in the Figure 47 is in $R(j-1) - \{ r(j)\} - B$ and using $\{G\}$ makes that only a finite subfamily of $R(j-1)$ is concerned.

\medskip

C) When it comes to STEP II, in its translation move to the right along the $x$-axis (in Figure 49) the $r(j)$ which moves to the position $b(j)$ has, of course, to brush the sites $b^3 (\beta) , b^3 (Q)$, and the ghostly $B^3$'s. This brushing on the $\oplus$ side, leaves the curves $c(r)$, $\{$little $C\}$ in place and is without further consequence. Actually, to be precise, the $b^3 (Q)$ have, just like the $B^3$ (ghost) a $\oplus$ and a $\ominus$ half. Only the $\oplus$ half has been explicitly drawn. As one can see in Figure 45, the $b^3 (\beta)$'s live completely inside $N_+^4 (2\Gamma (\infty))$ and they are not factorized.

\medskip

D) We will explain now the factorized little detail
$$
\tilde S^2 = \frac12 \, \tilde B^2 (-) \cup \frac12 \, \tilde B^2 (+)
$$
which occurs in the Figure 47 surrounded by thick dotted lines $(= \bm{- - -})$ and which bounds the 3-ball:
$$
\tilde B^3 = \frac12 \, \tilde B^3 (-) \cup \frac12 \, \tilde B^3 (+) \subset \partial N^4 (2\Gamma (\infty)) ,
$$
which occurs, also, in the Figure 50 below. This is the place via which the cylinder which connects $B_p^3$ to $A_1^3$ at the level of the Figure 47, connects with the COLOUR-CHANGING CYLINDER (145.7), after the $B_a^3$ (LAVA) has slided over $B_p^3$ (LAVA). This is displayed, with some detail, albeit schematically, in Figure 50.

\smallskip

After that move, $B_p^3 = r(j)$ leaves the space of the Figure 47 and, via $\tilde S^2$, which is actually a section of (146) (lateral surfaces of the COLOUR-CHANGING cylinder), enters the COLOUR-CHANGING cylinder in question and makes its move $r(j) \xRightarrow{ \ \ {\rm STEP \ II} \ \ } b(j)$. [In terms of (146) the $\tilde S^2$ is located between $S^2 \times r$ and $S^2 \times b$.]

\smallskip

Coming back to STEP II, our COLOUR-CHANGE is NOT a $1$-handle slide, coming globally with some global isotopic move of an $N^4 (2\Gamma (\infty))$-like object (some $Z^4 (t)$), but an internal isotopic repositionning of a 1-handle cocore, in a fixed background.

\medskip

E) The condition $C \cdot h = {\rm id} + {\rm nilpotent}$ which we had at time $t=\frac12$ has gotten bumped as we saw, and the $\{$extended cocores $H_j^b \}$'s are defined now by appealing directly to the STRONG PRODUCT PROPERTY of LAVA. And, inside the big step $t=\frac12 \Longrightarrow t=1$, we have a lot of intermediary smaller steps $j-1 \Rightarrow j$ each of them divided into a
$$
\mbox{STEP I $+$ STEP II $+$ STEP III.}
$$
The STEP I and STEP III do not modify the $\{$extended cocore $H_j^b\}$ which they find (at the beginning time $t=j-1$ or just after STEP II ($t = j-\varepsilon$)). But STEP II achieves the following at $t=j$: The set $\underset{\ell = 1}{\overset{M}{\sum}} \ \{$extended cocore $H_{\ell}^b \} \cap \{$the parasitical 1-handles from (144.1), at time $t=j\}$, NO LONGER CONTAINS the $r(j)$ and this $r(j)$ will never reappear again inside any of the $\{$extended cocore $H_{\ell}^b\}$, $\ell \leq M$. So, when we get to $t=1$ and {\ibf all} the STEPS II will have been put into effect (and see here the (144.2) too), then we get that $\underset{\ell = 1}{\overset{M}{\sum}} \ \{$extended cocore $H_{\ell}^b \} \cap \{$all the parasitical $h_k \in R-B$, and see here (144.1)$\} = \emptyset$, actually the $h_k$'s in question have all turned BLUE. This, together with the (147.6) leads to the BIG BLUE DIAGONALIZATION.

\medskip

F) Here are some additional comments concerning the COLOUR-CHANGING STEP II, from Figure 49, inside $j-1 \Rightarrow j$. We certainly have
$$
\{ B^3 \times [r,b] \} \cap R(j-1) \equiv \{{\rm our} \ r(j)\}.
$$
At the beginning of STEP II, this $r(j)$ is LAVA $(t=j-1)$ and, at the end, the $t=j-\varepsilon$, it remains a mere ghostly spot which is non LAVA, while $b(j) \subset {\rm LAVA} \, (t=j-\varepsilon)$. This $b(j)$ becoming LAVA at $t=j-\varepsilon$, is actually a {\ibf decree.} Remember that $t = j-\varepsilon$ is the time when the STEP II has been completed, and the STEP III, the last piece of $j-1 \Rightarrow j$ is ready to start. The STEP III does not touch any longer the composition of the $\{$extended cocore $H_j^b \}$, which has changed during the STEP II, and so there is no worry  concerning the E) above.

\medskip

G) In all our story, inside the $Z(t)$'s $(\approx$ ``$N^4 (2\Gamma (\infty))$'') we have strands $\{C$ remaining$\}$, with $\{D^2 (C)$'s $\subset$ LAVA$\}$ sticking out of them, possibly glued on LAVA BRIDGES ($+$ DILATATIONS) and also strands $\{\Gamma_j$ remaining$\}$, with $D^2 (\Gamma_j)$'s, which are NON-LAVA sticking out. The $\Gamma_i$'s stick on the 1-handles of $\Delta^4$, but they are never glued to the LAVA bridges or dilatations. They run closely and parallel to their mates $C$ (and LAVA bridges) in the same sheaf.

\medskip

H) When one looks at Figure 46, at time $t=j-1$ both $\Gamma - R (j-1)$ and $\Gamma - B$ are trees. At $t=j$ the $R(j) = R(j-1) - \{ r(j)\} + \{ b(j)\}$ appears, and now $\Gamma - R(j)$ and $\Gamma - B$ are trees. This is the global view. Now, our COLOUR-CHANGING Step II is NOT a 1-handle sliding, like for instance the one in the Figures 47 $+$ 48 but an ``isotopy'' of 1-handle cocores, an internal and not an external process. In order to be  able to realize this, we have to go to $Z^4 (t=1-2\varepsilon)$, and to server the connections with the outer world at the site $b^3 (\beta)$, which are now to be brushed through, in the COLOUR-CHANGING cylinder which we have rather artificially created. Our ``isotopy'' of 1-handle cocores $r(j) \Rightarrow b(j)$, which clearly cannot take place, as such, in the real world, is a convenient recipee for changing in an appropriate way the connections between 1-handles and (2-handles) $+$ LAVA BRIDGES, so as to give consistency to the {\ibf colour-changing decree}.

\medskip

I) Figure 50, which is a companion to Figure 49 is a schematical representation of the creation of our COLOUR-CHANGING  cylinder.

$$
\includegraphics[width=175mm]{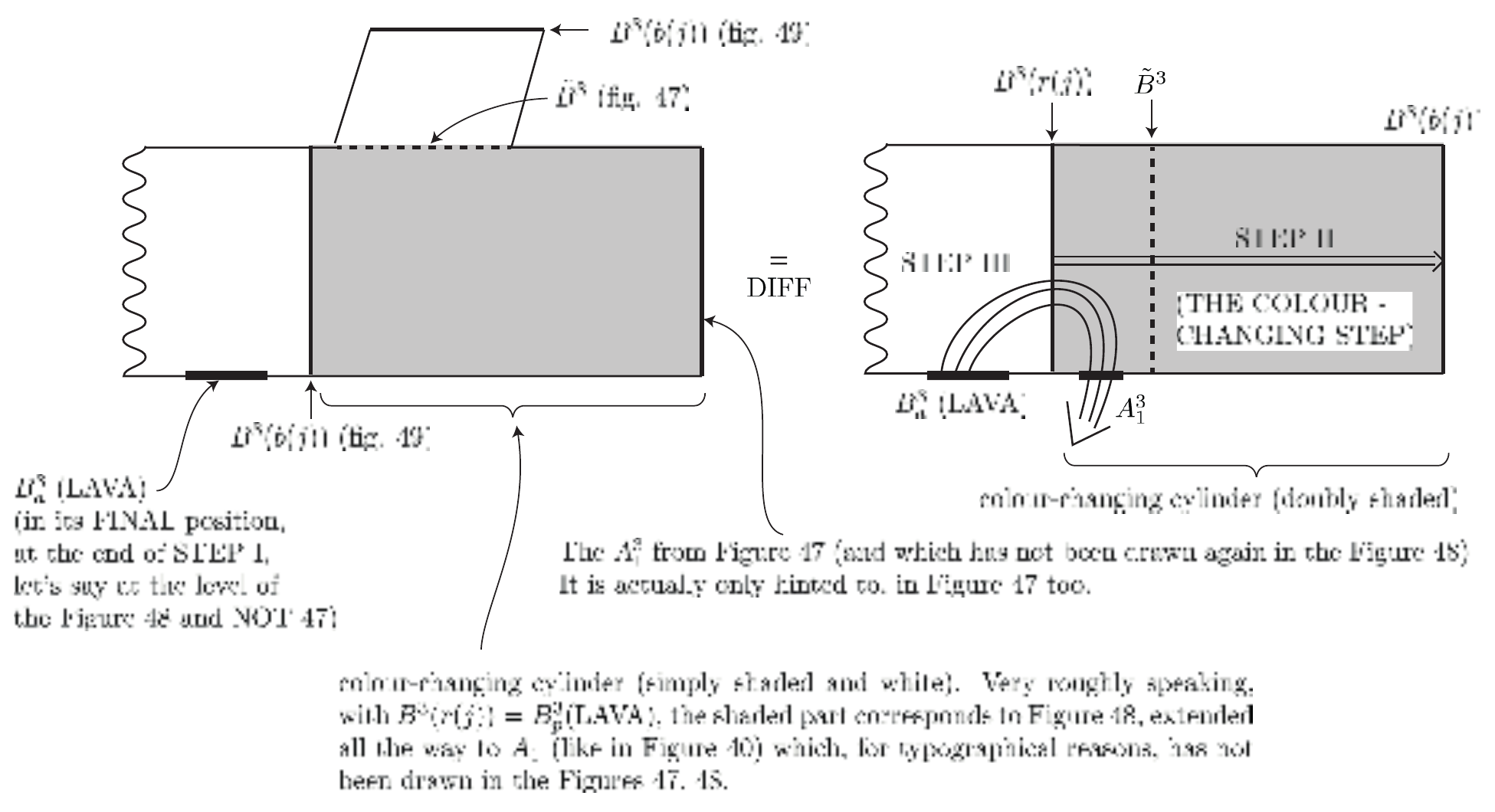}
$$
\label{fig50}
\centerline{\bf Figure 50.}
\begin{quote}
Geometry of the colour-changing cylinder (A schematical explanation). We can also explain now the $\tilde B^3$ which occurs as a fat dotted lines ({\bf - - -}) here, and also in Figure 47. Consider Figure 46 and let us say that the $B_a^3 ({\rm LAVA})$ in Figure 47 is $\{ (B_a^3)'$ OR $(B_a^3)''\}$, which has gotten now on the same edge as the 1-hanle cocore $r(j)$ and is ready to slide over it (a 1-handle sliding move). Then, on the road from $(B_a^3)'$ to $r(j)$ we find the $U$ (Figure 46), on the way, separating let us say, $r(j) + (B_a^3)'$ from $b(j)$. Then, for $B_a^3$ (Figure 47) $=$ $(B_a^3)'$, the $U$ is the $\tilde B^3$. When it comes to $(B_a^3)''$ or $(B_a^3)'''$ in Figure 46, it is $V$ which separates now $r(j) + (B_a^3)^{''(''')}$ from $b(j)$, and this is now the $\tilde B^3$ in the corresponding Figure 47.
\end{quote}

\bigskip

J) We will talk now about STEP III which restores things geometrically (the moving $Z^4 (t)$ is changed back into $N^4 (2\Gamma (\infty)))$, staying all the time far from $b(j)$. Figure 51 suggests what step III is supposed to do to the $B_a^3 ({\rm LAVA})$ as it stands after Step I, in its final position from the Figure 48, and it certainly does not budge from there, during the STEP II. The move of $B_a^3 ({\rm LAVA})$ sliding over $(B_p^3 = r(j))({\rm LAVA})$, which we consider here, is the one from the Figures 47 $+$ 48.

$$
\includegraphics[width=150mm]{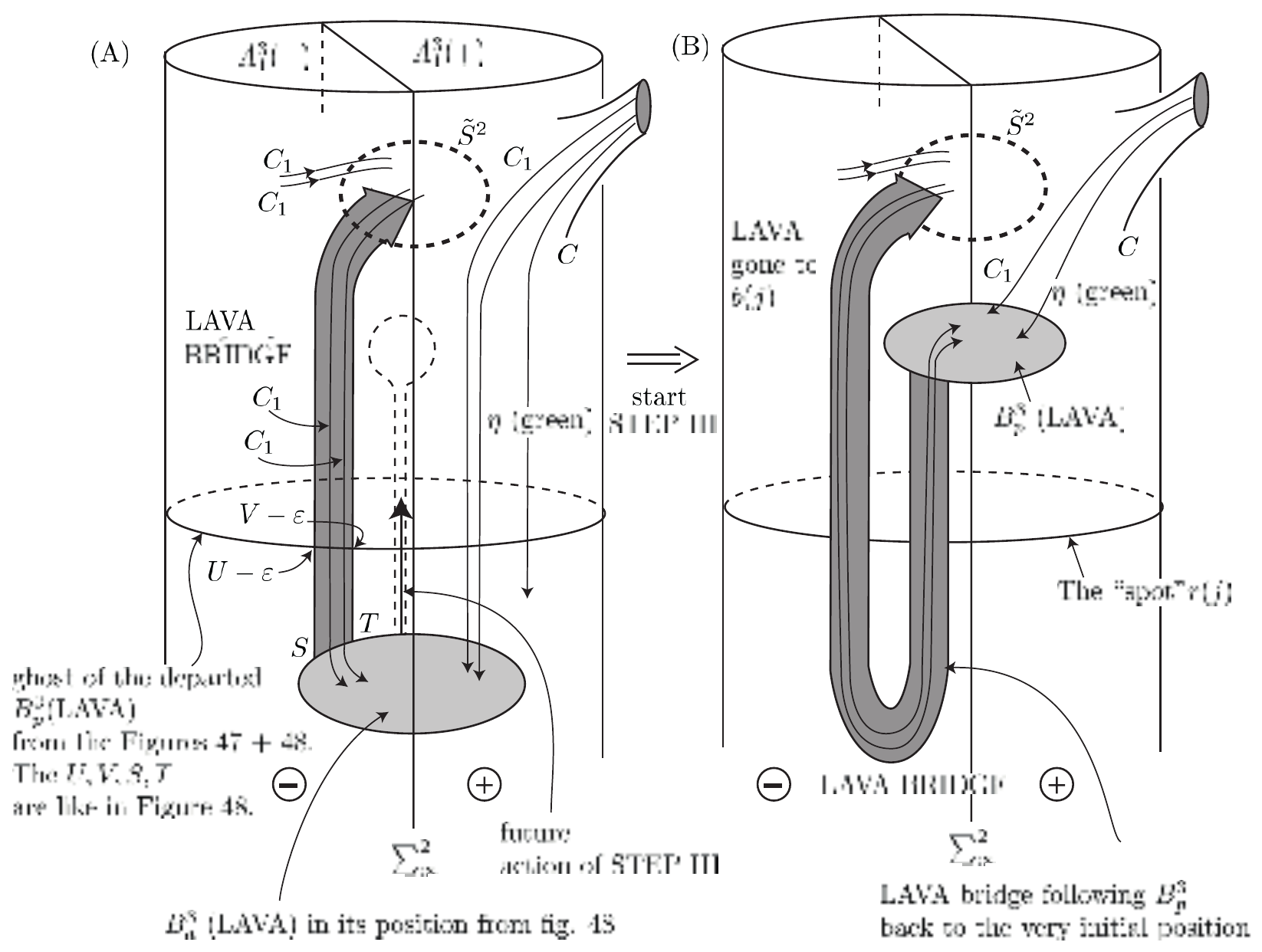}
$$
\label{fig51}
\centerline{\bf Figure 51.}
\begin{quote}
The STEP III of $j-1 \Rightarrow j$ at the level of Figure 48. At the level (A). The STEP II (following the Figure 47 $\underset{\rm STEP \, II}{\xRightarrow{ \qquad\qquad }}$ 48) has been performed, but not yet the STEP III. The $B^3 (p) = r(j)$ from Figure 48 has departed towards the position $b(j)$, going through $\tilde S^2$, with the LAVA BRIDGE dragged along after it. The $U,V$ from Figure 48 have gone with $B_p^3$ through $\tilde S^2$, leaving only the ghostly $U-\varepsilon$, $V-\varepsilon$ behind them. At the level of (B) the STEP III for $B_a^3$ has just started, taking it, for the time being just over the spot $r(j)$. We have tried to render graphically, the contorsion of the LAVA bridges. In (A) we see, dotted, the future move from (B).
\end{quote}

\newpage

\section*{References}

\begin{enumerate}
\item[{[C]}] A. {\sc Connes}, ``Non-Commutative Geometry'', Academic Press (1994).
\item[{[Ke-M]}] M. {\sc Kervaire} and J. {\sc Milnor}, ``Groups of homotopy spheres I'', {\sl Ann. of Math.} {\bf 77}, pp. 504-537 (1963).
\item[{[Ma]}] B. {\sc Mazur}, ``On embedding on spheres'', {\sl BAMS} {\bf 65} (1959), pp. 59-65.
\item[{[Po1]}] V. {\sc Po\'enaru}, ``Geometric simple connectivity in 4-dimensional Topology'', Pr\'epublications M/10/45, http://www.ihes.fr/PREPRINTS-M01/Resu/resu-M01-45.htm (2001).
\item[{[Po2]}] V. {\sc Po\'enaru}, ``Geometric simple connectivity and low-dimensional topology'', {\sl Proc. Steklov Inst. of Math.} {\bf 247} (2004), pp. 195-208.
\item[{[Po3]}] V. {\sc Po\'enaru}, ``A glimpse into the problems of the fourth dimension'' (to appear) preprint (2016).
\item[{[S]}] S. {\sc Smale}, ``On structure of manifolds'', {\sl Amer. J. of Math.} {\bf 84}, pp. 387-399 (1962).
\item[{[O-Po-Ta]}] D. {\sc Otera}, V. {\sc Po\'enaru}, C. {\sc Tanasi}, ``On Geometric Simple Connectivity'', {\sl Bull. Math. Soc. Math. Roum.} {\bf 53}, N$^{\rm o}$ 2, pp. 157-176 (2010).

\end{enumerate}

\end{document}